\documentclass[11pt,twoside]{article}

\usepackage{geometry}

\input{commands}

\begin{document}
\pagestyle{empty}
\begin{center}

  \textbf{\LARGE{Efficient reductions from a Gaussian source with applications to statistical-computational tradeoffs}}

\vspace*{.1in}

{\large{
\begin{tabular}{ccc}
Mengqi Lou$^{\star}$, Guy Bresler$^{\ddagger}$, Ashwin Pananjady$^{\star, \dagger}$
\end{tabular}
}}
\vspace*{.1in}

\begin{tabular}{c}
Schools of $^\star$Industrial and Systems Engineering and
$^\dagger$Electrical and Computer Engineering, \\
Georgia Tech \\
$^\ddagger$Department of Electrical Engineering and Computer Science, MIT
\end{tabular}

\vspace*{.2in}

\today

\vspace*{.2in}

\begin{abstract}
Given a single observation from a Gaussian distribution with unknown mean $\theta$, we design computationally efficient procedures that can approximately generate an observation from a different target distribution $\mathcal{Q}_{\theta}$ uniformly for all $\theta$ in a parameter set. We leverage our technique to establish reduction-based computational lower bounds for several canonical high-dimensional statistical models under widely-believed conjectures in average-case complexity. 
In particular, we cover cases in which:
\begin{enumerate}
	\item $\mathcal{Q}_{\theta}$ is a general location model with non-Gaussian distribution, including both light-tailed examples (e.g., generalized normal distributions) and heavy-tailed ones (e.g., Student's $t$-distributions). As a consequence, we show that computational lower bounds proved for spiked tensor PCA with Gaussian noise are universal, in that they extend to other non-Gaussian noise distributions within our class.
	\item $\mathcal{Q}_{\theta}$ is a normal distribution with mean $f(\theta)$ for a general, smooth, and nonlinear link function $f:\mathbb{R} \rightarrow \mathbb{R}$. Using this reduction, we construct a reduction from symmetric mixtures of linear regressions to generalized linear models with link function $f$, and establish computational lower bounds for solving the $k$-sparse generalized linear model when $f$ is an even function. This result constitutes the first reduction-based confirmation of a $k$-to-$k^2$ statistical-to-computational gap in $k$-sparse phase retrieval, resolving a conjecture posed by~\citet{cai2016optimal}. 
  As a second application, we construct a reduction from the sparse rank-1 submatrix model to the planted submatrix model, establishing a pointwise correspondence between the phase diagrams of the two models that faithfully preserves regions of computational hardness and tractability.
\end{enumerate}
\end{abstract}
\end{center}

\tableofcontents

\newpage
\pagestyle{plain}

\setcounter{page}{1}

\section{Introduction and background}\label{sec:introduction}

For some unknown real value $\theta$ from an uncountably large parameter set $\Theta$, suppose we observe a random draw $X_{\theta} \sim \mathcal{N}(\theta, 1)$. Consider two simple-to-state puzzles, both of a flavor that dates back to classical questions in mathematical statistics~\citep[e.g.,][]{blackwell1951comparison,blackwell1953equivalent,karlin1956theory,lindley1956measure,stone1961non,hansen1974comparison}:
\begin{enumerate}
  \item[(P1)] Can we generate a sample that is close in total variation distance to a random variable drawn from a prespecified, non-Gaussian location model with location parameter $\theta$?
  \item[(P2)] For a prespecified nonlinear function $f$, can we generate a sample from a distribution that is close in total variation distance to $\mathcal{N}(f(\theta), \sigma^2)$?
\end{enumerate}
The difficulty of these puzzles lies in the fact that $\theta$ is unknown a priori and impossible to estimate from a single observation---nevertheless, our procedure must work for infinitely many values of $\theta$.

A naive approach in both cases is the plug-in estimator. For puzzle (P2) above, this would amount to outputting $f(X_{\theta}) + Z$ where $Z \sim \mathcal{N}(0, \sigma^2)$ is an independent random variable. Indeed, this is a reasonable strategy for large $\sigma^2$: We have
\[
  \| \mathcal{N}(f(X_{\theta}), \sigma^2) - \mathcal{N}(f(\theta), \sigma^2) \|_{\mathsf{TV}} \overset{\1}{\leq} \sqrt{ \frac{1}{2} \mathsf{KL}( \mathcal{N}(f(X_{\theta}), \sigma^2) \; \| \; \mathcal{N}( f(\theta), \sigma^2) )} \leq \frac{\EE| f(X_{\theta}) - f(\theta) |}{2\sigma} =: \frac{C_{f, \Theta}}{\sigma},
\]
where step $\1$ holds by Pinsker's inequality, and $C_{f, \Theta}$ is uniformly bounded provided $f$ is well-behaved on the parameter set $\Theta$ and $\Theta$ is bounded.
Said a different way, for a natural class of functions $f$, if the variance $\sigma^2$ scales polynomially in $1/\epsilon$, then the total variation distance between $f(X_{\theta}) + Z$ and the desired $Y_{\theta} \sim \mathcal{N}(f(\theta), \sigma^2)$ is bounded above by $\epsilon$. We will show later in the paper that even for simple functions $f$, such a variance blow-up is necessary for the plug-in approach. 

One of the main results of this paper is a more sophisticated procedure that achieves this total variation guarantee---for a large class of functions $f$---with $\sigma^2$ scaling only polylogarithmically in $1/\epsilon$. This presents an exponential improvement over the plug-in approach, and together with our other main result, which is a procedure answering puzzle (P1) for large class of non-Gaussian location models, yields several consequences for establishing new statistical-computational tradeoffs in canonical high-dimensional statistical models.

To set up the problem formally in high dimensions, suppose we have two parametric statistical models. The first is the model $\mathcal{U} = (\mathbb{R}^D, \{ \mathcal{P}_{\xi} \}_{\xi \in \Xi})$, also called the source, and the second is $\mathcal{V} = (\mathbb{R}^D, \{ \mathcal{Q}_{\xi} \}_{\xi \in \Xi})$, also called the target. Suppose both the source model $\mathcal{U}$ and the target model $\mathcal{V}$ share the same sample space $\mathbb{R}^D$ and (possibly structured) parameter space $\Xi \subseteq \mathbb{R}^D$, but differ in their respective distributions $\mathcal{P}_{\xi}$ and $\mathcal{Q}_{\xi}$. For the moment, this is abstract notation for a high-dimensional model and accommodates isomorphisms of $\real^D$---for instance, the parameter set and observations could be matrices or tensors.
Our focus in this paper is specifically on the setting where the source statistical model is the \emph{Gaussian location model}, with $\mathcal{P}_{\xi} = \mathcal{N}(\xi, \sigma^{2} I_D)$. Specifically, for a variance parameter $\sigma^{2}$ that should be thought of as fixed and known, we have
\begin{align}\label{source-gaussian-location-vector}
\mathcal{U} = (\mathbb{R}^D, \{ \mathcal{N}(\xi, \sigma^{2} I_D) \}_{\xi \in \Xi}).
\end{align}

Now suppose there is some \textit{unknown} $\xi \in \Xi$ fixed by nature. We draw a \textit{single} observation (or sample) from the source distribution, and denote this random variable by $X_{\xi} \sim \NORMAL(\xi,\sigma^{2} I_D)$. We would like to transform this high-dimensional source random vector to a $D$-dimensional random vector over the target sample space that is close to $Y_{\xi} \sim \mathcal{Q}_{\xi}$. In particular, our goal is to design a randomized algorithm $\mathsf{K}:\real^D \rightarrow \real^D$ that additionally takes some $\epsilon \in (0, 1)$ as input and achieves the following twofold objective:
\begin{itemize}
\item[(A)] \textbf{Statistical:} Uniformly over $\xi \in \Xi$, the random variable $\mathsf{K}(X_{\xi})$ must be $\epsilon$-close in total variation (TV) distance to an output random variable $Y_{\xi} \sim \mathcal{Q}_{\xi}$, i.e., 
\[
\sup_{\xi \in \Xi} \; \mathsf{d_{TV}}(\textsc{K}(X_{\xi}), Y_{\xi}) \leq \epsilon, \qquad \text{and}
\] 
\item[(B)] \textbf{Computational:} The transformation algorithm $\mathsf{K}$ is computationally efficient, in that it must have (pointwise) runtime polynomial in the dimension $D$, and in $1/\epsilon$. 
\end{itemize}
We refer to such a transformation algorithm as a polynomial-time \emph{reduction}. It is worth noting that for some target models, there may not exist any reduction satisfying the statistical desideratum~(A) above due to data processing considerations~\citep[see, e.g.,][Chapter 7]{polyanskiy2023book}, and also that it may be possible that there are reductions that satisfy the statistical desideratum (A), but that a computationally efficient one satisfying desideratum~(B) is unlikely to exist under widely believed conjectures from average-case complexity (see \citet[Appendix A]{lou2025computationally} for examples).

One of the main applications of such a polynomial-time reduction is that it can be used to transfer computational hardness from the source model to the target model.
We briefly explain the underlying idea at a high level\footnote{Rigorous formulations of these computational implications can be found in~\citep[e.g.,][]{berthet2013sPCA,ma2015computational,brennan2018reducibility}.}. Suppose that---owing to evidence from average-case complexity---we believe that computationally efficient and accurate estimators do not exist for the source model. Now say we have a computationally efficient reduction as above, which transforms a sample from the source model into one that resembles a sample from the target model while keeping the underlying estimand the same. Then computationally efficient and accurate estimation for the target model must be impossible. To see this, suppose (for the sake of contradiction) that there exists a computationally efficient and accurate estimator for the target model. Then we can use this to produce an accurate estimator for the source model as follows: First draw a source sample, then use the reduction to transform it into a target sample and apply the estimator for the target model. Since the reduction and the estimator for the target model run in polynomial time, so does this entire procedure for the source model. But this is impossible by assumption and we have a  contradiction. Thus, any accurate estimator for the target model cannot run in polynomial time.

\subsection{Motivating examples}\label{sec:motivating_examples}

We now present concrete high-dimensional examples that motivate our paper.
These examples serve as a preview; they are formally presented again in Section~\ref{sec:applications} along with explicit guarantees. In all these examples,  reductions from the high-dimensional source model~\eqref{source-gaussian-location-vector} are key. As we will see shortly, however, it suffices in these cases to design reductions from the \emph{scalar} source
\begin{align} \label{source-gaussian-location}
\mathcal{U} = (\mathbb{R}, \{ \mathcal{N}(\theta, \sigma^{2}) \}_{\theta \in \Theta}),
\end{align}
where $\Theta \subseteq \real$ is some parameter set consisting of real values, and $\sigma^2$ is a known variance parameter. Indeed, at the core of these applications are puzzles of the form (P1) and (P2) mentioned above.

\paragraph{Motivation 1.} In spiked tensor Principal Component Analysis (PCA)~\citep{montanari2014statistical}, we observe a single order-$s$ tensor $T \in \mathbb{R}^{n^{\otimes s}} = \mathbb{R}^{n \times \cdots \times n}$ given by
\begin{align}\label{tensor-PCA-intro}
T = \beta (\sqrt{n} v)^{\otimes s} + \zeta,
\end{align}
where $\beta \in \mathbb{R}$ denotes the signal strength, $v \in \mathbb{S}^{n-1}$ is the unknown spike, and $\zeta \in \mathbb{R}^{n^{\otimes s}}$ is a symmetric Gaussian noise tensor (see Eq.~\eqref{def-tensor-gaussian-noise}). Based on a reduction from the secret leakage planted clique conjecture~\citep{brennan2020reducibility} or hypergraph planted clique conjecture~\citep{luo2022tensor}, it is believed that no polynomial-time algorithm can succeed at either detection or recovery when $\beta = \widetilde{o}(n^{-s/4})$. In contrast, recovery is information-theoretically possible for a much weaker signal strength $\beta = \widetilde{\omega}(n^{(1-s)/2})$ via exhaustive search~\citep{montanari2014statistical,lesieur2017statistical,chen2019phase,jagannath2020statistical,perry2020statistical}. 

A key question concerns the universality of the computational barrier: Does the same computational lower bound persist when entry-wise Gaussian noise $\zeta$ is replaced by a zero-mean, non-Gaussian distribution, for instance one having lighter or heavier tails? As argued above, one way to affirmatively answer this question is to reduce the original tensor $T$ to one with non-Gaussian noise. In addition, observe that for each index of the tensor $(i_1, \dots, i_s)$, the entry satisfies $T_{i_1, \dots, i_s} \sim \mathcal{N}(\theta, 1)$, where $\theta = \beta n^{s/2} v_{i_1} \cdots v_{i_s}$ is an unknown mean. Since the noise is independent across entries, the core problem is to design a one-dimensional reduction from the source model~\eqref{source-gaussian-location} (with $\mathcal{P}_{\theta} = \mathcal{N}(\theta, \sigma^2)$) to a general non-Gaussian target location model $\mathcal{Q}_{\theta}$, as posed in puzzle (P1). Executing these one-dimensional reductions entry-by-entry is one way to arrive at a polynomial-time, high-dimensional reduction.

\begin{figure}[ht!]
  \centering
  \includegraphics[width=0.8\textwidth]{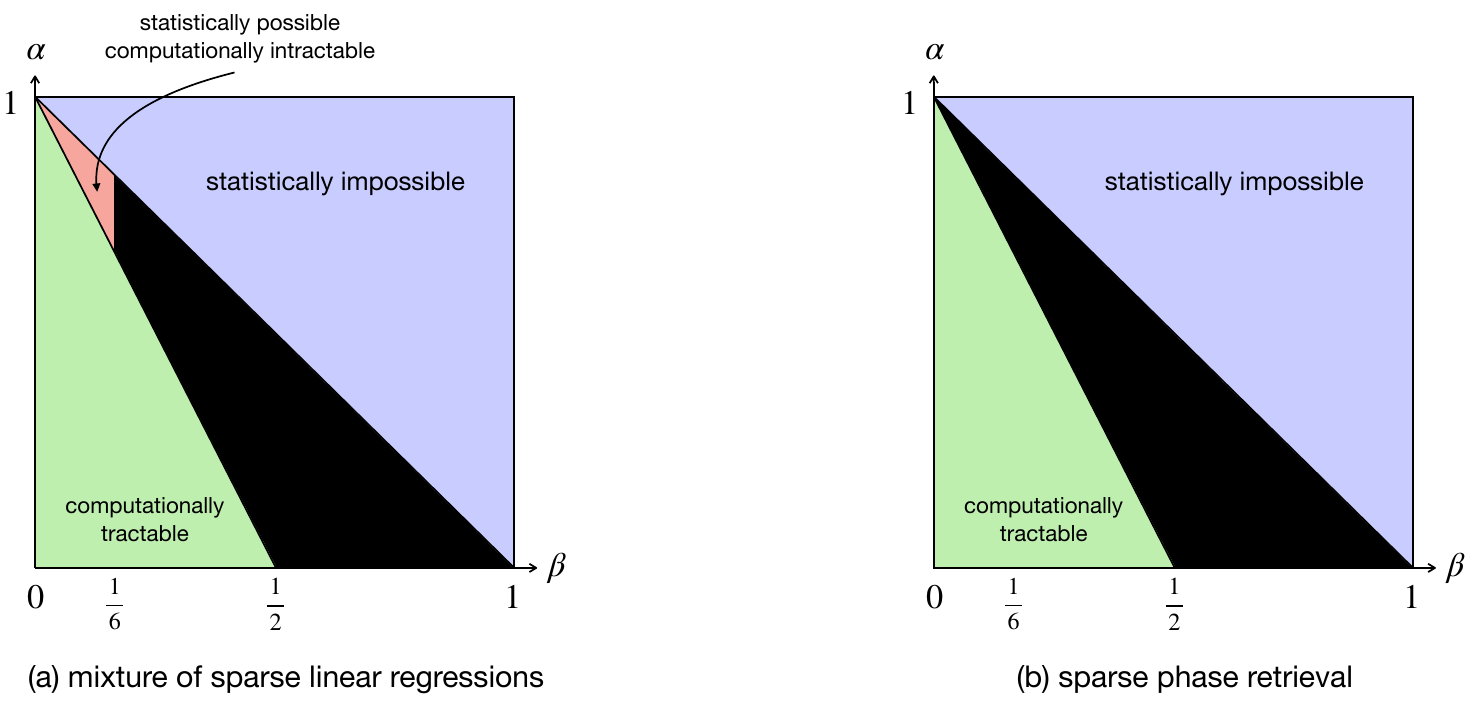}
  \caption{Phase diagrams for the sparse mixture of linear regressions ($\MixSLR$; panel (a)) and sparse phase retrieval ($\SPR$; panel (b)), under the parameterization $k=\widetilde{\Theta}(n^{\beta})$ and $\eta^{4}=\widetilde{\Theta}(n^{-\alpha})$, with constants $\alpha,\beta\in(0,1)$. 
  The purple region denotes the statistically impossible regime; the red region denotes the statistically possible but computationally intractable regime; the green region denotes the computationally tractable regime; and the black region denotes a statistically possible yet conjectured computationally hard regime, whose proof remains open. Our goal is to transfer all computationally hard regimes confirmed in panel (a) to panel (b).
  }
  \label{fig:SMLR-SPR-intro}
\end{figure}

\paragraph{Motivation 2.} In symmetric mixtures of sparse linear regressions ($\MixSLR$)~\citep{quandt1978estimating,stadler2010l1}, we observe i.i.d. samples $\{(x_{i},y_{i})\}_{i=1}^{n}$ generated by
\[
	y_{i} = R_{i} \cdot \langle x_{i}, v \rangle + \zeta_{i}, \quad \text{where} \quad  \zeta_{i} \sim \NORMAL(0,1).
\]
Here, $R_{i} \in \{-1,1\}$ is a hidden label and $v \in \mathbb{R}^{d}$ is the unknown $k$-sparse parameter. 
There is a so-called $k$-to-$k^{2}$ statistical-to-computational gap in $\MixSLR$. Information theoretically, we require $n = \widetilde{\Theta}\big( k \log(d) /\|v\|_{2}^{4} \big)$ for estimating $v$ to within $\ell_2$ distance $\| v \|_2/2$~\citep{fan2018curse}, yet there is reduction-based evidence from the secret leakage planted clique conjecture that no polynomial-time algorithm can succeed when $n = \widetilde{o}\big( k^{2}/\|v\|_{2}^{4} \big)$~\citep{brennan2020reducibility}. 

A related but distinct model is the sparse generalized linear model ($\SpGLM$), in which samples follow
\[
	\widetilde{y}_{i} = f\big( \langle x_{i}, v \rangle \big) + \zeta_{i}, \quad \text{where} \quad \zeta_{i} \sim \NORMAL(0,1),
\]
where $f:\mathbb{R} \to \mathbb{R}$ is the link function. Note that the computational gap in $\MixSLR$ arises due to the latent signs---in light of this, it is natural to ask: Does the same $k^{2}$ computational barrier persist in $\SpGLM$ with even link functions, which also destroy sign information? In particular, when $f(t) = t^{2}$ we obtain the familiar sparse phase retrieval model~\citep{li2013sparse}, and it was conjectured by~\citet{cai2016optimal} that there is a $k^2$ computational barrier, i.e., that no polynomial time algorithm can succeed when $n = \widetilde{o}(k^{2})$. However, no reduction-based evidence for this conjecture has been established. 
A natural way by which one might establish such computational hardness is by transforming samples $\{x_{i},y_{i}\}_{i=1}^{n}$ from $\MixSLR$ to samples $\{x_{i},\widetilde{y}_{i}\}_{i=1}^{n}$ from $\SpGLM$, thereby allowing us to transfer any computationally hard regimes that have been confirmed in Figure~\ref{fig:SMLR-SPR-intro}(a) to Figure~\ref{fig:SMLR-SPR-intro}(b). Note that $y_{i} \sim \mathcal{N}(\theta,1)$ and $\widetilde{y}_{i} \sim \mathcal{N}(f(\theta),1)$ for some unknown mean $\theta = R_{i} \cdot \langle x_{i}, v\rangle$, and the noise is independent and Gaussian across samples.  
Thus, a one-dimensional reduction from the Gaussian source~\eqref{source-gaussian-location} suffices: We must design a reduction from $\mathcal{P}_{\theta} = \mathcal{N}(\theta,1)$ to target $\mathcal{Q}_{\theta} = \mathcal{N}(f(\theta),1)$, as alluded to in puzzle (P2).

\begin{figure}[ht!]
  \centering
  \includegraphics[width=0.85\textwidth]{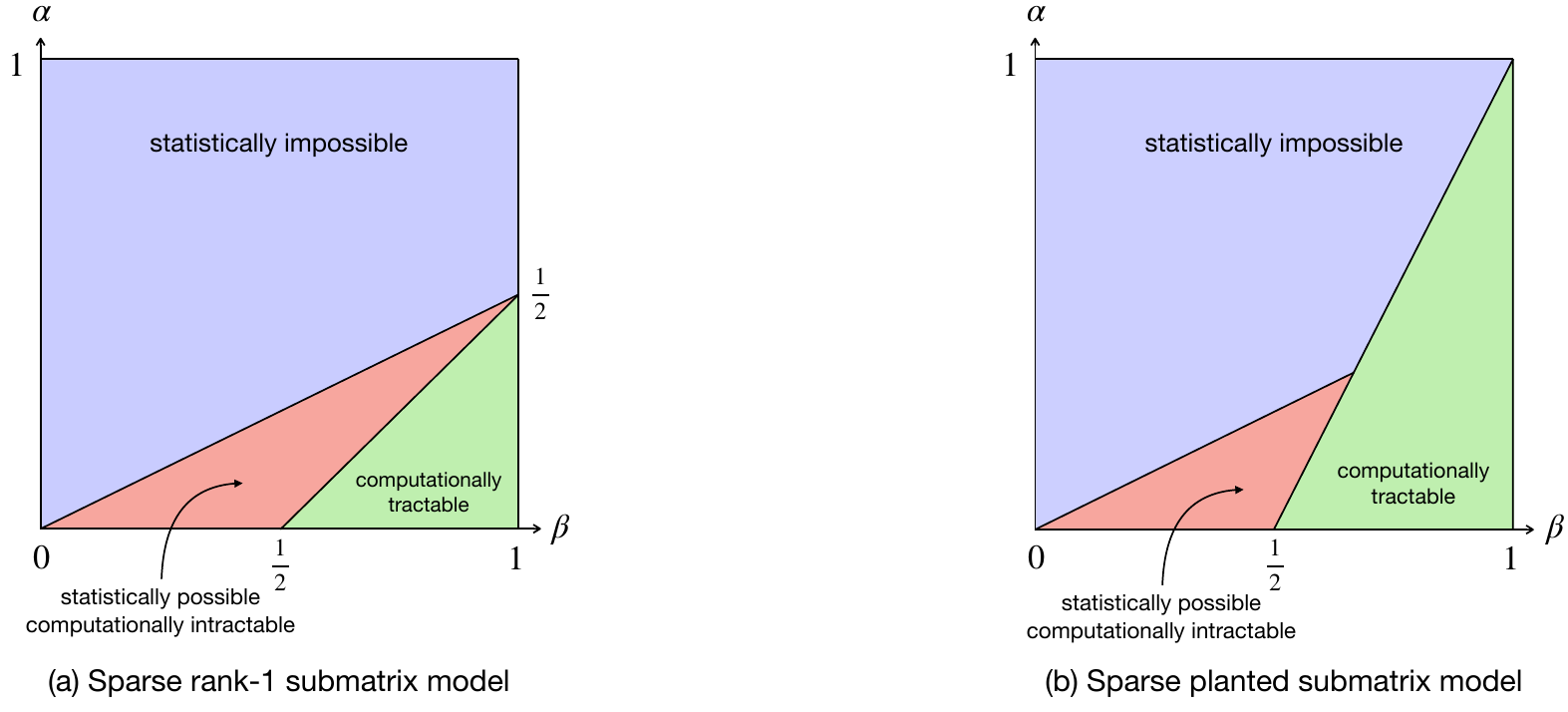}
  \caption{Phase diagrams for the sparse rank-1 submatrix model ($\SubMat$; panel (a)) and the planted submatrix model ($\PlantSubMat$; panel (b)), under the parameterization $\mu/k = \widetilde{\Theta}(n^{-\alpha})$ and $k = \widetilde{\Theta}(n^{\beta})$ for $\alpha, \beta \in (0,1)$. The purple region indicates the statistically impossible regime. The green region corresponds to the computationally tractable regime where polynomial-time algorithms are known to succeed. The red region indicates the statistically possible but computationally intractable regime. In panel (a), the computational threshold is given by the line $\alpha = \beta - 1/2$, while in panel (b), the threshold is the line $\alpha = 2\beta - 1$.}
  \label{fig:phase-diagram-submatrix}
\end{figure}

\paragraph{Motivation 3.} In the sparse rank-one submatrix model ($\SubMat$), we observe
\[
	M = \mu\, r c^{\top} + G, \quad \text{where} \quad G_{ij} \overset{\mathsf{i.i.d.}}{\sim} \NORMAL(0,1).
\]
Here, $\mu \in \mathbb{R}$ is the signal strength and $r,c \in \mathbb{R}^n$ are $k$-sparse unit-norm vectors with both positive and negative entries. In the planted submatrix model ($\PlantSubMat$), $M$ has the same form but $r,c$ are restricted to be $k$-sparse unit-norm vectors with only nonnegative entries. Existing reduction-based computational hardness results for $\PlantSubMat$ and $\SubMat$ are based on the planted clique conjecture~\citep{brennan2020reducibility,ma2015computational}; see Figure~\ref{fig:phase-diagram-submatrix} for their phase diagrams. However, computational hardness for each of these problems is established \emph{individually} through reductions from the planted clique model.  A natural question is whether, without relying on the planted clique conjecture, one can construct a \emph{direct} reduction from $\SubMat$ to $\PlantSubMat$. Note that once again, this would involve an entry-by-entry reduction from a Gaussian source. In principle, such a reduction would allow one to transfer computational hardness directly from instances of $\SubMat$ to instances of $\PlantSubMat$, allowing us to move beyond computational hardness that can be established through reductions from planted clique. Concretely, while estimation in the planted clique model can be performed using a quasi-polynomial time algorithm, it is conjectured that $\SubMat$ requires exponential time for certain regimes of signal-to-noise ratio~\citep{ding2024subexponential}. Stronger computational hardness of this flavor cannot be transferred to $\PlantSubMat$ via a direct reduction from planted clique.
Also note that unlike in the previous examples, Figure~\ref{fig:phase-diagram-submatrix} suggests that the two models exhibit different computational thresholds. Can a reduction between them be designed to preserve both computational tractability and intractability (i.e., map tractable regimes to tractable regimes and computationally hard regimes to hard regimes)? 

\subsection{Setup and contributions}

Having motivated the study of reductions from a Gaussian source, we now state our contributions. We begin with some setup. Given a real-valued parameter $\theta \in \Theta$, we let $u(\cdot;\theta):\real \rightarrow \real$ denote the density of the source distribution $\mathcal{P}_{\theta} = \NORMAL(\theta,\sigma^{2})$. Specifically, we let
\begin{align}\label{source-gaussian-density}
  u(x; \theta) = \frac{1}{\sqrt{2\pi} \sigma} \exp\left( - \frac{(x - \theta)^2}{2\sigma^{2}} \right), \quad \text{for all } x, \theta \in \mathbb{R}.
\end{align}
Similarly, we let $v(\cdot;\theta):\real \rightarrow \real$ denote the density of the target distribution $\mathcal{Q}_{\theta}$ with respect to Lebesgue measure. A mapping $\mathcal{T}(\cdot \mid \cdot):\real \times \real \rightarrow \real$ is called a Markov kernel if $\mathcal{T}(\cdot\mid x)$ is a density on $\real$ for all $x \in \real$. 
For a Markov kernel $\mathcal{T}$, we let 
\begin{align}\label{TV-deficiency}
  \delta(\mathcal{U},\mathcal{V};\mathcal{T}) : = \sup_{\theta \in \Theta}  \Big\| \int_{\real} \mathcal{T}(\cdot \mid x) \cdot u(x;\theta) \mathrm{d}x - v(\cdot ; \theta) \Big\|_{\mathsf{TV}}
\end{align}
denote the TV-deficiency attained by $\mathcal{T}$ when mapping the source statistical model $\mathcal{U}$ to the target statistical model $\mathcal{V}$.

\paragraph{Main contributions.} We now provide an overview of our main contributions. 
\begin{itemize}
  \item[(i)] \textbf{Design and analysis of two reductions.} 
  We construct a Markov kernel $\mathcal{T}$ from the Gaussian location model $\mathcal{U}$ in Eq.~\eqref{source-gaussian-location} to a general class of target models $\mathcal{V}$ and establish an upper bound on the TV-deficiency $\delta(\mathcal{U},\mathcal{V};\mathcal{T})$ in terms of properties of $\mathcal{V}$; see Section~\ref{sec:general_markov_kernel}. We use this Markov kernel to construct our single-sample reduction algorithm $\mathsf{K}$. 

  We then apply this framework to two concrete classes of target models and show that the reduction $\mathsf{K}$ is both statistically and computationally efficient. In particular:
  \begin{itemize}
    \item When $\mathcal{V}$ corresponds to a class of non-Gaussian location models with a suitably proxy for signal-to-noise ratio decaying polylogarithmically in $1/\epsilon$, Theorem~\ref{thm:non-gaussian} guarantees that our reduction achieves $\epsilon$-TV deficiency and runs in time polynomial in $1/\epsilon$ for arbitrarily small $\epsilon \in (0,1)$. We allow for the parameter to be any real value, and cover light-tailed target models such as the generalized normal as well as heavy-tailed models such as Student's $t$ distributions, among others. See Corollaries~\ref{example:bounded-psi} and~\ref{example:generalized-normal} for details.

    \item When $\mathcal{V}$ is a Gaussian model with mean $f(\theta)$ and a suitably proxy for signal-to-noise ratio decays polylogarithmically in $1/\epsilon$, Theorem~\ref{thm:gaussian-mean-nonlinear} guarantees that our reduction achieves $\epsilon$-TV deficiency and runs in time polynomial in $1/\epsilon$ for arbitrarily small $\epsilon \in (0,1)$. Under only the assumption that the parameter set is suitably bounded, we cover functions that are constant-degree polynomials or globally analytic, among others. See Corollaries~\ref{example:monomial} and~\ref{example:analytical} for details. 
  \end{itemize}

  \item[(ii)] \textbf{Applications to computational hardness.} 
  Using these reductions, we obtain affirmative answers to several open questions on computational hardness in high-dimensional statistics:
  \begin{itemize}
    \item We prove a reduction-based universality result for Tensor PCA, showing that the same computational barrier as in the Gaussian case persists for several non-Gaussian noise distributions (including both lighter-tailed and heavier-tailed cases); see Corollary~\ref{thm:test-hardness-pca-nongaussian}. 
    \item We construct a reduction from $\MixSLR$ to $\SpGLM$ for a general class of even link function, thereby confirming the existence of a $k^{2}$ computational barrier for $k$-sparse generalized linear models with an even link; see Corollary~\ref{thm:hardness-glm}. This provides the first reduction-based resolution of the conjectured $k$-to-$k^{2}$ statistical-to-computational gap in sparse phase retrieval posed by \citet{cai2016optimal}. See Corollary~\ref{corollary:spr} for details.
    \item We construct a reduction from $\SubMat$ to $\PlantSubMat$, yielding a pointwise correspondence between their phase diagrams (Figure~\ref{fig:phase-diagram}) and showing that the reduction faithfully preserves both computational tractability and intractability; see Theorem~\ref{thm:reduction-ROS-BC}.
  \end{itemize}
\end{itemize}

\subsection{Related work}

We review related work under two verticals. Given that our focus is on reductions, the first vertical will discuss other reduction-based approaches to proving computational hardness in high-dimensional statistical models. In our second vertical, we provide details on other approaches to proving hardness that are not reduction-based, wherein computational hardness is proved within some (restricted) classes of algorithms. Given that this latter literature is vast, we focus on papers that provide lower bounds for the high-dimensional problems outlined in Section~\ref{sec:motivating_examples}. 

\paragraph{Reduction-based computational hardness.} The program of transferring computational hardness between problems in high-dimensional statistics was initiated roughly a decade ago~\citep{berthet2013complexity,berthet2013sPCA,wainwright2014constrained,ma2015computational,hajek2015computational}. Since then, the framework of reductions has been applied to an extensive list of high-dimensional statistical problems~\citep{brennan2018reducibility,brennan2020reducibility}. Specific examples include principal component analysis with structure~\citep{brennan2019optimal,wang2023algorithms,bresler2025computational}, mixture models with adversarial corruptions~\citep{brennan2020reducibility}, block models in matrix and tensor settings~\citep{brennan2020reducibility,luo2022tensor,han2022exact}, tensor PCA~\citep{brennan2020reducibility}, and adaptation lower bounds in
shape-constrained regression problems~\citep{shah2019feeling,pananjady2022isotonic}. Recently, reductions have been constructed from hard problems that have seen a long line of work in cryptography~\citep{gupte2022continuous,bruna2021continuous,gupte2024sparse,bangachev2025near}.

This body of work introduces many powerful tools and techniques to transfer structure between problems, but it cannot handle continuous-valued inputs or parameter spaces like the Gaussian location model---indeed, all the distributional reductions in this literature are from source distributions taking either two or three values~\citep{hajek2015computational,ma2015computational,brennan2018reducibility,brennan2020reducibility}. There are some classical reductions that handle Gaussian sources, but they are only effective when the target is also Gaussian (see the books by~\cite{shiryaev2000statistical} and~\cite{torgersen1991comparison} for examples). To our knowledge, the only other single-sample reductions from continuous valued source models---apart from a few specialized source-target pairs~\citep{hansen1974comparison,lehmann2011comparing}---are from our own prior work~\citep{lou2025computationally}, which addresses reductions from Laplace, Erlang, and uniform source models to general target models. However, the tools introduced in that paper do not apply to the Gaussian source model~\eqref{source-gaussian-location}.

Other literature that is related in spirit to the reduction framework includes the problems of simultaneous optimal transport~\citep{wang2022simultaneous} and sample amplification~\citep{axelrod2024statistical}. However, both these problems differ in their focus.

\paragraph{Computational hardness for restricted algorithm classes.} There are several techniques that have been developed to prove computational hardness over restricted classes of algorithms, including those based on low-degree polynomials~\citep{kunisky2019notes,wein2025computational}, sum-of-squares hierarchies~\citep{raghavendra2018high}, and the statistical query model~\citep{kearns1998efficient}. As mentioned before, each is a vast field, so we focus on particular papers that address our problems of interest in Section~\ref{sec:motivating_examples}.

\smallskip

\noindent \underline{Tensor PCA.} For tensor PCA with Gaussian noise~\eqref{tensor-PCA-intro}, a large body of work has demonstrated that many natural classes of algorithms fail when the signal strength falls below the computational threshold $\beta = \widetilde{o}(n^{-s/4})$. This includes low-degree polynomial algorithms~\citep{kunisky2019notes}, sum-of-squares relaxations~\citep{hopkins2015tensor,hopkins2016fast,hopkins2017power,hopkins2018statistical,kim2017community}, approximate message passing (AMP) algorithms~\citep{montanari2014statistical,lesieur2017statistical}, algorithms based on the Kikuchi hierarchy~\citep{wein2019kikuchi}, statistical query algorithms~\citep{dudeja2021statistical,brennan2021statistical}, and Langevin dynamics applied to the maximum likelihood objective~\citep{arous2020algorithmic}. Additional evidence for computational hardness has been shown via landscape analysis~\citep{arous2019landscape,ros2019complex}, predictions from statistical physics~\citep{bandeira2018notes}, and communication-complexity lower bounds~\citep{dudeja2024statistical}. Universality of hardness for other noise distributions was recently shown by~\citet{kunisky2025low} for the class of low-coordinate degree algorithms.

\smallskip

\noindent \underline{Sparse generalized linear models.}
For sparse generalized linear models ($\SpGLM$), several works have established computational barriers. In particular, low-degree polynomial methods require at least $k^{2}$ samples for exact support recovery in $k$-sparse phase retrieval with Gaussian covariates in the noiseless setting~\citep{arpino2023statistical}. Note that this problem is identical to the noiseless $\MixSLR$ problem with Gaussian covariates. Using the statistical query framework, it has also been shown that no computationally efficient algorithm can attain the information-theoretic minimax risk for $\SpGLM$ with even link functions~\citep{wang2019statistical}. This work builds on a general framework based on showing statistical-query lower bounds for \emph{non-Gaussian component analysis} problems~\citep{diakonikolas2017statistical,diakonikolas2023algorithmic}. Beyond sparsity, statistical-computational gaps in generalized linear models arise can depending on the choice of link function, and lower bounds have been established in this setting via both statistical query and low-degree polynomial methods~\citep[see, e.g.,][]{damian2024computational}.

\smallskip

\noindent \underline{Sparse submatrix models.}
For sparse submatrix detection and recovery problems ($\SubMat$ and $\PlantSubMat$), evidence of computational hardness comes from multiple fronts: low-degree polynomial lower bounds~\citep{ding2024subexponential,schramm2022computational,wein2025computational,kunisky2019notes}, failure of AMP~\citep{lesieur2015phase,lesieur2017constrained,hajek2018submatrix,macris2020all}, sum-of-squares lower bounds~\citep{ma2015sum,hopkins2017power}, and the overlap gap property~\citep{gamarnik2021overlap,arous2023free}. Statistical query lower bounds have also been connected to the low-degree approach in such problems~\citep{brennan2021statistical}. Notably, in $\PlantSubMat$, there exists a gap between detection and recovery thresholds~\citep{chen2016statistical,schramm2022computational,arous2023free}.

\begin{remark}
In spite of the many different methods now in use to prove computational hardness in statistical problems, reductions remain an ironclad route to establishing fundamental relationships between the average-case computational complexity of different problems. 
 As discussed above, the reduction approach differs from proving computational hardness for restricted classes of algorithms---instead, it directly addresses the question of polynomial-time solvability \emph{by any algorithm} assuming hardness of some base problem. While showing hardness for a restricted class of algorithms does not depend on assumed hardness for any other problem, the ramification of such a result depends crucially on the belief that that class of algorithms is actually optimal for the problem at hand. It is by now well-known that each of these restricted classes of algorithms can be suboptimal for simple problems, and it is still unclear what type of problem each algorithm class solves optimally. For instance, the Lenstra--Lenstra--Lovasz (LLL) lattice basis reduction algorithm can be used to solve some of these problems in polynomial time with state-of-the-art guarantees~\citep[e.g.,][]{andoni2017correspondence,zadik2022lattice,diakonikolas2023algorithmic}, but it is not a low-degree polynomial algorithm or a statistical query algorithm. For another example, it has recently been observed that algorithmically simple problems such as shortest path on random graphs can exhibit the overlap gap property~\citep{li2024some}. Reductions, however, are immune to such specialized counterexamples. More generally, a key property of reductions is that they can work synergistically with analyses for restricted algorithmic approaches by clarifying what transformations preserve polynomial-time solutions to problems.
\end{remark}

\subsection{Notation and organization} 
For $a,b\in \real$, we let $a\vee b = \max\{a,b\}$ and $a\wedge b = \min\{a,b\}$. We use $\mathcal{L}[X]$ to denote the law of a random variable $X$. For two (possibly signed) measures $\mu$ and $\nu$, we use $\frac{\mathrm{d} \nu}{\mathrm{d} \mu}$ to denote the Radon--Nikodym derivative of $\nu$ with respect to $\mu$. We use $\NORMAL(\theta, \sigma^2)$ to denote a normal distribution with mean $\theta$ and variance $\sigma^2$. We use $\mathsf{GN}(\theta,\varT,\beta)$ to denote a Generalized normal distribution with mean $\theta$, scale $\varT$, and shape $\beta$. Let $\Gamma(z) = \int_{0}^{\infty} t^{z-1}e^{-t} \mathrm{d}t$ for all $z>0$ be the Gamma function. Let $\mathbb{N} = \{1,2,\dots\}$ denote the set of natural numbers and let $\PosInt = \{0,1,2,\dots\}$ denote the set of non-negative integers. Let $[n]$ denote the set of natural numbers less than or equal to $n$. 
For $n\in \mathbb{N}$, we define double factorial $n!! = \prod_{k=1}^{n/2} (2k)$ if $n$ is even and $n!! = \prod_{k=1}^{(n+1)/2} (2k-1)$ if $n$ is odd. For a function of two variables $v(y; \theta)$ and $k\in \PosInt$, we write $\nabla_{\theta}^{(k)} v(y; \theta) \big \vert_{\theta = \widetilde{\theta}}$ to denote the $k$-th derivative of $v$ with respect to its second argument $\theta$, evaluated at $\theta = \widetilde{\theta}$. Similarly, for a function $f:\real \rightarrow \real$, we let $f^{(k)}(x)$ to denote the $k$-th derivative of $f$ evaluated at $x$, i.e., $f^{(k)}(x) = \frac{\mathrm{d}^{k}\,f(t)}{\mathrm{d}t^{k}} \big \vert_{t = x}$. For a positive integer $s$, we let $\mathfrak{S}_{s}$ denote the set of all permutations of $[s]$, i.e., $\mathfrak{S}_{s} = \{ \pi :[s] \rightarrow [s]\;|\; \pi \text{ is a bijection} \}$. Following the notation in~\cite{raginsky2011shannon}, we use $\| \mu - \nu \|_{\mathsf{TV}}$ to denote the total variation distance between two measures $(\mu, \nu)$ defined on the same space. If $X \sim \mu$ and $Y \sim \nu$, then we use the notation $\mathsf{d_{TV}}(X, Y) := \| \mu - \nu \|_{\mathsf{TV}}$.  We define the (physicist's) Hermite polynomials $H_{n}:\real \rightarrow \real$ as
\begin{align}\label{def:hermite}
  H_{n}(x) = (-1)^{n} e^{x^{2}} \cdot \frac{\mathrm{d}^{n}\,e^{-x^{2}} }{\mathrm{d}x^{n}}, \quad \text{for all} \quad x\in\real \text{ and } n = 0,1,2,\dots.
\end{align}
It will be convenient later to work with the normalized Hermite polynomials
\begin{align}\label{def:hermite-normalized}
  \Nherm_{n}(x) = \frac{H_{n}(x)}{\big( \sqrt{\pi}2^{n}n! \big)^{1/2} }, \quad \text{for all} \quad x\in \real \text{ and } n = 0,1,2,\dots.
\end{align}
For two functions $f,g:\real \rightarrow \real$, define their inner product as 
\begin{align} \label{eq:in-prod}
  \langle f, g \rangle : = \int_{\real} f(x) g(x) \exp(-x^{2}) \mathrm{d}x.
\end{align}
Note that the normalized Hermite polynomials are orthonormal in that $\big \langle \Nherm_{n}(\cdot), \Nherm_{m}(\cdot)  \big\rangle = \delta_{nm}$, where $\delta_{nm}$ is the Kronecker delta with $\delta_{nm} = \ind{n = m}$.
We use $(c, c', C, C')$ to denote universal positive constants that may take different values in each instantiation. For two sequences of non-negative reals $\{f_n\}_{n\geq 1}$ and $\{g_n \}_{n \geq 1}$, we use $f_n \lesssim g_n$ to indicate that there is a universal positive constant $C$ such that $f_n \leq C g_n$ for all $n \geq 1$. The relation $f_n \gtrsim g_n$ indicates that $g_n \lesssim f_n$,
and we say that $f_n \asymp g_n$ if both $f_n \lesssim g_n$ and $f_n \gtrsim g_n$ hold simultaneously. We also use standard order notation $f_n = \order
(g_n)$ to indicate that $f_n \lesssim g_n$ and $f_n = \ordertil(g_n)$ to indicate that $f_n \lesssim g_n \log^c n$, for a universal constant $c>0$. We say that $f_n = \Omega(g_n)$ (resp. $f_n = \widetilde{\Omega}(g_n)$) if $g_n = \order(f_n)$ (resp. $g_n = \ordertil(f_n)$). We say
$f_n = \Theta(g_n)$ (resp. $f_n = \widetilde{\Theta}(g_n)$) if $f_n = \Omega(g_n)$  and $f_n = \mathcal{O}(g_n)$ hold simultaneously (resp. $f_n = \widetilde{\Omega}(g_n)$ and $f_n = \widetilde{\mathcal{O}}(g_n)$). 
The notation $f_n = o(g_n)$ is
used when $\lim_{n \to \infty} f_n / g_n = 0$, and $f_n = \omega(g_n)$ when $g_n = o(f_n)$. 
We write $\Omega_{n}(1)$ to denote a sequence that is bounded below by some universal and positive constant for all $n$ large enough. Similarly, $o_{n}(1)$ denotes a non-negative sequence that converges to $0$ as $n$ goes to infinity. For $x>0$, we write $\mathsf{poly}(x)$ to denote a polynomial function of $x$ and write $\mathsf{polylog}(x)$ to denote a polynomial function of $\log(x)$. 
All logarithms are to the natural base $e$ unless otherwise stated. We write $\sgn(z) = 1$ if $z \geq 0$ and $-1$ otherwise. We use $\ind{\cdot}$ to denote the indicator function.

\paragraph{Organization.} The remainder of this paper is organized as follows. 
In Section~\ref{sec:concrete_reductions}, we present concrete reductions from the Gaussian source model to general target models, focusing on the cases where the target is either a general location model or a Gaussian model with a nonlinear mean transformation. 
In Section~\ref{sec:applications}, we apply these reductions to establish computational lower bounds for several high-dimensional models, thereby addressing the questions raised in Motivations~1--3. 
We conclude in Section~\ref{sec:discussion} with a discussion of open problems. 
Proofs of all results are provided in Sections~\ref{sub:pf_general_markov_kernel},~\ref{sec:pf_thm_non_gaussian},~\ref{sub:pf_thm_gaussian_mean_nonlinear}, and~\ref{sec:proof_applications}.


\section{General families of reductions from the Gaussian source}\label{sec:concrete_reductions}
In this section, we present a general family of reductions from the Gaussian location model. Specifically, in Section~\ref{sec:general_markov_kernel}, we develop methodology for constructing a Markov kernel $\mathcal{T}$ that maps the Gaussian source $\mathcal{U}$~\eqref{source-gaussian-location} to a general target $\mathcal{V}$, and supply a bound on its TV deficiency---in Lemma~\ref{lemma:TV-deficiency-decompose}---as a function of various properties of the target model. In Section~\ref{sec:general_location_target}, we apply our Markov kernel to cases where the target $\mathcal{V}$ is a non-Gaussian location model and prove explicit bounds on the TV deficiency $\delta(\mathcal{U}, \mathcal{V}; \mathcal{T})$. In Section~\ref{sec:gaussian_target_nonlinear}, we consider the case where the target is Gaussian with mean given by a nonlinear transformation $f:\real \to \real$, i.e., $\mathcal{V} = (\real, \{ \NORMAL(f(\theta),1)\}_{\theta \in \Theta})$, and again provide explicit bounds on $\delta(\mathcal{U}, \mathcal{V}; \mathcal{T})$.

\subsection{Constructing a general Markov kernel}\label{sec:general_markov_kernel}
In this section, we describe the main idea behind constructing a reduction from the Gaussian location model~\eqref{source-gaussian-location} to a general target statistical model
\[
\mathcal{V} = (\mathbb{R}, \{ \mathcal{Q}_{\theta} \}_{\theta \in \Theta}).
\]
Since our approach is based on taking higher-order derivatives, we require some mild regularity conditions on the target density $v(\cdot;\theta)$.
\begin{assumption}\label{assump_1:target}
The following conditions hold:
\begin{itemize}
    \item[(a)] The higher-order derivatives
    $
    \nabla_{x}^{(k)} v(y; x)
    $
    exist for all $x, y \in \mathbb{R}$ and $k \in \PosInt$.

    \item[(b)] The higher-order derivatives vanish at infinity. In particular, for all $k \in \PosInt$, $y \in \real$, and $\theta \in \Theta$, 
    \[
      \lim_{x \uparrow +\infty} \big[ \nabla_{x}^{(k)} v(y;x) \big] \cdot e^{-\frac{(x-\theta)^{2}}{2\sigma^{2}}} = \lim_{x \downarrow -\infty} \big[ \nabla_{x}^{(k)} v(y;x) \big] \cdot e^{-\frac{(x-\theta)^{2}}{2\sigma^{2}} } = 0.
    \]
    \item[(c)] The function $v\big(y;\sqrt{2}\sigma x+\theta\big)$ satisfies, for all $y\in \real$ and $\theta \in \Theta$, 
    \[
      \int_{\real} \Big( v\big(y;\sqrt{2}\sigma x+\theta \big) \Big)^{2} \exp(-x^{2}) \mathrm{d}x < +\infty. 
    \]
\end{itemize}
\end{assumption}
In words, Assumption~\ref{assump_1:target}(a) requires $v(y;x)$ to be infinitely differentiable with respect to $x$; Assumption~\ref{assump_1:target}(b) requires that the partial derivatives $\nabla_{x}^{(k)} v(y;x)$ do not grow faster than exponentially, which is usually satisfied in most applications of interest; and Assumption~\ref{assump_1:target}(c) places a boundedness condition that typically holds since $v(y;x)$ is uniformly bounded, which in turn ensures that the integral is finite.

Let us now describe the core observations behind our construction. Clearly, any reduction $\textsc{K}$ corresponds to a Markov kernel $\mathcal{T}(\cdot \mid \cdot): \mathbb{R} \times \mathbb{R} \to \mathbb{R}$---given a sample $X_{\theta} \sim \mathcal{N}(\theta, \sigma^{2})$, the reduction generates the random variable $\textsc{K}(X_\theta)$ by sampling from the conditional density $\mathcal{T}(\cdot \mid X_\theta)$. As a result, the marginal density of $\textsc{K}(X_\theta)$ is given by
$
\int_{\mathbb{R}} \mathcal{T}(\cdot \mid x) \cdot u(x; \theta) \, \mathrm{d}x.
$
Minimizing the total variation (TV) distance between the distribution of $\textsc{K}(X_\theta)$ and the target distribution $\mathcal{Q}_{\theta}$ thus amounts to solving the following optimization problem:
\begin{align}\label{optimal-markov-kernel}
\begin{split}
\delta(\mathcal{U},\mathcal{V}): =  \inf_{\mathcal{T}(\cdot \mid \cdot): \mathbb{R} \times \mathbb{R} \to \mathbb{R}} \quad  \sup_{\theta \in \Theta}\quad & \frac{1}{2} \left\| \int_{\mathbb{R}} \mathcal{T}(\cdot \mid x) \cdot u(x; \theta) \, \mathrm{d}x - v(\cdot; \theta) \right\|_{1}, \\
\text{subject to} \quad & \mathcal{T}(y \mid x) \geq 0 \quad \text{and} \quad \int_{\mathbb{R}} \mathcal{T}(y' \mid x) \, \mathrm{d}y' = 1 \quad \text{for all } x, y \in \mathbb{R}. 
\end{split}
\end{align}
The optimization problem~\eqref{optimal-markov-kernel}, while linear, is infinite-dimensional. It can thus be computationally inefficient to solve, and analyzing its optimal value is challenging due to the complexity of the constraint set. 
Our approach is thus to pursue a feasible (instead of optimal)
Markov kernel $\mathcal{T}$ and to bound its TV deficiency $\delta(\mathcal{U}, \mathcal{V}; \mathcal{T})$. We proceed in a few steps:

\subsubsection*{Step 1: Construct a good signed kernel.}
Let $\mathcal{S}(\cdot \mid \cdot): \real \times \real \rightarrow \real$ be a bivariate function such that $S(\cdot \mid x)$ is a measurable function for all $x \in \real$. For any such \emph{signed kernel} $\mathcal{S}$, define
\begin{align}\label{deficiency-signed-kernel}
  \SignedDef( \mathcal{U}, \mathcal{V}; \mathcal{S}) := \sup_{\theta \in \Theta} \; \frac{1}{2} \left\| 
  \int_{\real}  \mathcal{S}( \cdot\mid x) \cdot u(x;\theta) \mathrm{d}x - v(\cdot;\theta)
  \right\|_{1},
\end{align}
which is analogous to $\delta(\mathcal{U},\mathcal{V}; \mathcal{T})$ defined in Eq.~\eqref{TV-deficiency}. The analogous optimization problem to Eq.~\eqref{optimal-markov-kernel} is then given by the unconstrained problem
\begin{align} \label{eq:signed-def}
\SignedDef( \mathcal{U}, \mathcal{V}) : =  \inf_{\mathcal{S}(\cdot \mid \cdot): \mathbb{R} \times \mathbb{R} \to \mathbb{R}} \quad \sup_{\theta \in \Theta} \quad \frac{1}{2} \left\| \int_{\mathbb{R}} \mathcal{S}(\cdot \mid x) \cdot u(x; \theta) \, \mathrm{d}x - v(\cdot; \theta) \right\|_{1}.
\end{align}
Our first observation is that the optimization problem~\eqref{eq:signed-def} admits a closed-form solution. We state this result as the following proposition.
\begin{proposition}\label{prop:optimal-signed-kernel}
Consider the source model $\mathcal{U} = (\mathbb{R}, \{ \mathcal{N}(\theta, \sigma^{2}) \}_{\theta \in \Theta})$ and general target model $\mathcal{V} = (\real,\{\mathcal{Q}_{\theta}\}_{\theta \in \Theta})$. Let $u(\cdot;\theta)$ be defined in Eq.~\eqref{source-gaussian-density} and let $v(\cdot;\theta):\real \rightarrow \real$ be the density of distribution $\mathcal{Q}_{\theta}$. Suppose Assumption~\ref{assump_1:target} holds for $\Theta \subseteq \real$ and for all $\theta \in \Theta$,
\begin{align}\label{condition-switch-int-sum}
  \sum_{k=0}^{+\infty} \int_{\mathbb{R}} \frac{\sigma^{2k}}{(2k)!!} \big| \nabla_{x}^{(2k)} v(y;x) \big| \cdot u(x;\theta) \, \mathrm{d}x < \infty.
\end{align} 
Consider the signed kernel defined as
\begin{align}\label{optimal-signed-kernel}
  \mathcal{S}^{\star}(y \mid x) := \sum_{k=0}^{+\infty} \frac{(-1)^{k} \sigma^{2k}}{(2k)!!} \nabla_{\theta}^{(2k)}v(y;\theta) \big \vert_{\theta = x}   \quad \text{for all } x,y\in \real.
\end{align}
Then we have $\SignedDef( \mathcal{U}, \mathcal{V}; \mathcal{S}^{\star}) = 0$.
\end{proposition}
We provide the proof of Proposition~\ref{prop:optimal-signed-kernel} in Section~\ref{sec:pf-prop1}. Note that the regularity condition~\eqref{condition-switch-int-sum} is a technical assumption imposed to justify the interchange of the integral and the infinite sum via Fubini's theorem. In later sections, we will make use of a modified signed kernel that does not require this condition. 

While Proposition~\ref{prop:optimal-signed-kernel} provides a closed-form expression for the optimal signed kernel, computing the signed kernel $\mathcal{S}^{\star}$ defined in Eq.~\eqref{optimal-signed-kernel} is computationally intractable since it requires evaluating $2k$-th order derivatives $\nabla_{\theta}^{(2k)} v(y; \theta)$ for infinitely many values of $k$. To address this, we consider a truncated version of the kernel. Specifically, for $N \in \mathbb{N}$, we define the truncated signed kernel as
\begin{align}\label{truncated-signed-kernel}
  \mathcal{S}^{\star}_{N}(y \mid x) = \sum_{k=0}^{N} \frac{(-1)^{k} \sigma^{2k}}{(2k)!!}  \nabla_{\theta}^{(2k)} v(y; \theta) \big \vert_{\theta = x} \quad \text{for all } x,y\in \real.
\end{align}
In light of Proposition~\ref{prop:optimal-signed-kernel}, we should expect that the signed deficiency $\SignedDef(\mathcal{U}, \mathcal{V}; \mathcal{S}^{\star}_{N})$ is small for large $N$. We quantify the computational cost of evaluating the signed kernel $\mathcal{S}^{\star}_{N}$~\eqref{truncated-signed-kernel} using the following definition.
\begin{definition}\label{assump:Teval}
For each $N \in \mathbb{N}$, define $T_{\mathsf{eval}}(N)$ as the maximum time taken to evaluate $\mathcal{S}^{\star}_{N}(y\mid x)$ for any $x, y \in \real$.
\end{definition}
In the specific applications that follow, we will provide explicit bounds on $T_{\mathsf{eval}}(N)$.

\subsubsection*{Step 2: Jordan decomposition to find a valid Markov kernel.}
Having motivated (through Proposition~\ref{prop:optimal-signed-kernel}) why the signed kernel $\mathcal{S}^{\star}_{N}$ defined in Eq.~\eqref{truncated-signed-kernel} may serve as a good choice to solve the optimization problem~\eqref{eq:signed-def}, we now construct a valid Markov kernel by retaining only the positive part of $\mathcal{S}^{\star}_{N}$ and normalizing. Specifically, since we have the Jordan decomposition of the signed measure $\mathcal{S}^{\star}_{N} (\cdot \mid x ) = [\mathcal{S}^{\star}_{N}(y \mid x) \vee 0] + [\mathcal{S}^{\star}_{N}(\cdot \mid x) \land 0]$, we define
\begin{align}\label{close-markov-kernel}
\mathcal{T}^{\star}_{N} (y \mid x) = \frac{ \mathcal{S}^{\star}_{N}(y \mid x) \vee 0 }{ \int_{\mathbb{R}} \left[ \mathcal{S}^{\star}_{N}(y' \mid x) \vee 0 \right] \, \mathrm{d}y' }, \quad \text{for all } x, y \in \mathbb{R}.
\end{align}
Note that by construction, for each $x \in \real$, $\mathcal{T}^{\star}_{N} (\cdot \mid x)$ is a density.

Lemma 1 of \citet{lou2025computationally} yields the Markov kernel deficiency bound
\begin{align}\label{ineq1:lemma-TV}
  \delta\big(\mathcal{U}, \mathcal{V};\mathcal{T}^{\star}_{N}\big)
  &\leq \SignedDef\big( \mathcal{U}, \mathcal{V}; \mathcal{S}_{N}^{\star} \big)
    + \frac{1}{2}\sup_{\theta\in \Theta}\int_{\mathbb{R}} \big(|p(x) - 1| + q(x)\big)\,u(x;\theta)\,\mathrm{d}x,
\end{align}
where the functions $p,q: \real \rightarrow [0,+\infty)$ are given by
\begin{align}\label{px-qx-def}
p(x) : = \int_{\real} \big[ \mathcal{S}^{\star}_{N}(y' \mid x) \vee 0 \big] \, \mathrm{d}y' \quad \text{and} \quad q(x): = - \int_{\real} \big[ \mathcal{S}^{\star}_{N}(y' \mid x) \wedge 0 \big] \, \mathrm{d}y'.
\end{align}
Specifically, the second term on the right-hand side of Eq.~\eqref{ineq1:lemma-TV} quantifies the discrepancy between the signed kernel $\mathcal{S}_{N}^{\star}$ and the Markov kernel $\mathcal{T}_{N}^{\star}$ introduced by the truncation and normalization in Eq.~\eqref{close-markov-kernel}. 
The following lemma provides an upper bound on the first term $\SignedDef\big( \mathcal{U}, \mathcal{V}; \mathcal{S}_{N}^{\star} \big)$.

\begin{lemma}\label{lemma:TV-deficiency-decompose}
  Consider the source model $\mathcal{U} = (\mathbb{R}, \{ \mathcal{N}(\theta, \sigma^{2}) \}_{\theta \in \Theta})$ and general target model $\mathcal{V} = (\real,\{\mathcal{Q}_{\theta}\}_{\theta \in \Theta})$. Let $u(\cdot;\theta)$ be defined in Eq.~\eqref{source-gaussian-density} and let $v(\cdot;\theta):\real \rightarrow \real$ denote the density of distribution $\mathcal{Q}_{\theta}$. For each $N \in \mathbb{N}$, consider the signed kernel $\mathcal{S}^{\star}_{N}$ defined in Eq.~\eqref{truncated-signed-kernel}.
  Suppose Assumption~\ref{assump_1:target} holds. Then
  \begin{align}\label{ineq2:lemma-TV}
    \SignedDef\big( \mathcal{U}, \mathcal{V}; \mathcal{S}_{N}^{\star} \big) &\leq \frac{1}{2}\sup_{\theta \in \Theta} \; \int_{\real} \sum_{k=N+1}^{+\infty} \bigg| \int_{\real} v\big(y;\theta + \sqrt{2}\sigma x\big) \cdot \Nherm_{2k}(x) e^{-x^2} \mathrm{d}x \bigg| \mathrm{d}y,
  \end{align} 
  where $\Nherm_{2k}(x)$ is the normalized Hermite polynomial defined in Eq.~\eqref{def:hermite-normalized}.
\end{lemma}
We provide the proof of Lemma~\ref{lemma:TV-deficiency-decompose} in Section~\ref{sec:pf-lem1}. 
Recall that the Hermite polynomials $\big\{ \Nherm_{n}(x) \big\}_{n \in \PosInt}$ form an orthonormal basis for the function space 
\[
  \bigg\{ f : \mathbb{R} \to \mathbb{R} \,\Big|\, \int_{\mathbb{R}} |f(x)|^{2} e^{-x^{2}} \, \mathrm{d}x < \infty \bigg \}.
\]
According to Ineq.~\eqref{ineq2:lemma-TV}, the signed deficiency $\SignedDef\big( \mathcal{U}, \mathcal{V}; \mathcal{S}_{N}^{\star} \big)$ can be bounded by the Hermite tail coefficients of the function $x \mapsto v\big(y;\theta + \sqrt{2}\sigma x\big)$, onto the subspace spanned by $\left\{ \Nherm_{2n}(x) \right\}_{n = N+1}^{\infty}$. These coefficients typically decay exponentially in $\sqrt{N}$ when $v\big(y;\theta + \sqrt{2}\sigma x\big)$ is sufficiently smooth in $x$; see Theorem~\ref{thm:non-gaussian} and Theorem~\ref{thm:gaussian-mean-nonlinear} for details.

\subsubsection*{Step 3: Approximately sample from the Markov kernel.}
Our final step is to use a rejection sampling algorithm to approximately sample from the family of distributions $\{\mathcal{T}_{N}^{\star}(\cdot \mid x)\}_{x \in \mathbb{R}}$. We use the particular algorithm from~\citet[Section~3]{lou2025computationally}, and denote it by $\textsc{RK}$. We do not reproduce a description of the entire algorithm here, but note that the algorithm $\textsc{RK}$ requires:
\begin{enumerate} 
\item A family of base measures $\{\mathcal{D}(\cdot \mid x)\}_{x \in \mathbb{R}}$ such that the ratio between $\mathcal{S}_{N}^{\star}(y \mid x)$ and $\mathcal{D}(y \mid x)$ is uniformly bounded. That is, for some constant $M > 0$,
\[
    \frac{\mathcal{S}_{N}^{\star}(y \mid x) \vee 0}{\mathcal{D}(y \mid x)} \leq M,
    \quad \text{for all } x,y \in \mathbb{R}.
\]
\item Moreover, the base measures $\{\mathcal{D}(\cdot \mid x)\}_{x \in \mathbb{R}}$ should be chosen so that it is computationally efficient to sample from them. We let $T_{\mathsf{samp}}$ denote the maximum time taken to sample from any base measure $\mathcal{D}(\cdot \mid x)$. 
\end{enumerate}
Given these, the output of the algorithm is a sample that is close in total variation to one having density proportional to $\mathcal{S}_{N}^{\star}(y \mid x) \vee 0$.
Concretely, denote the output of the algorithm is by $\textsc{RK}(x,T_{\mathsf{iter}},M,y_{0})$, where $T_{\mathsf{iter}} \in \mathbb{N}$ is the number of iterations and $y_{0} \in \mathbb{R}$ is an arbitrary initialization. As shown in~\citet{lou2025computationally}, the law of this random variable satisfies 
\begin{align}\label{error-rk-sampling}
    \big\| \mathcal{L}\big[ \textsc{RK}(x,T_{\mathsf{iter}},M,y_{0}) \big] - \mathcal{T}_{N}(\cdot \mid x) \big\|_{\mathsf{TV}} 
    \leq 2 \exp\!\left( -\frac{T_{\mathsf{iter}}}{M} \, p(x) \right), 
    \quad \text{for all } x \in \mathbb{R},
\end{align}
where $p(x)$ is defined in Eq.~\eqref{px-qx-def}. Thus, provided that $p(x) \geq c$ for some universal constant $c>0$, the error from sampling procedure decreases geometrically as the number of iterations $T_{\mathsf{iter}}$ increases. Computationally, the implementation of $\textsc{RK}(x,T_{\mathsf{iter}},M,y_{0})$ requires, in each iteration, generating a sample from $\mathcal{D}(\cdot \mid x)$ and evaluating the signed kernel $\mathcal{S}_{N}^{\star}$ once. Hence, $\textsc{RK}$ can be implemented in time 
$\mathcal{O} \big( T_{\mathsf{iter}} \cdot T_{\mathsf{samp}} \cdot T_{\mathsf{eval}}(N) \big)$.

\subsubsection*{Putting together the steps.}
Our reduction algorithm is defined as
\[
    \mathsf{K}(X_{\theta}) := \textsc{RK}(X_{\theta},T_{\mathsf{iter}},M,y_{0}),
\]
where $X_{\theta} \sim \mathcal{N}(\theta,\sigma^{2})$ is the input sample from the source. The initialization $y_0$ can be arbitrary\footnote{As discussed in~\cite{lou2025computationally}, an arbitrary initialization suffices for a TV guarantee. Specific initializations may be desired if we want to prove bounds in other metrics besides TV.}.
Recall that our goal was to control the total variation distance between $\mathsf{K}(X_{\theta})$ and the target sample $Y_{\theta} \sim v(\cdot;\theta)$, which we denoted by $\mathsf{d}_{\mathsf{TV}} \big( \mathsf{K}(X_{\theta}), Y_{\theta} \big)$. Error arises from two sources: (i) the sampling error of the rejection sampling procedure $\textsc{RK}$, which can be bounded as in Eq.~\eqref{error-rk-sampling}; and (ii) the deficiency of the Markov kernel $\mathcal{T}_{N}$, which is quantified by $\delta(\mathcal{U},\mathcal{V}; \mathcal{T}_{N})$ and can be bounded by using Ineq.~\eqref{ineq1:lemma-TV} and Ineq.~\eqref{ineq2:lemma-TV}. 

Note that while the reduction is general, its TV deficiency may not always be small.
We next turn to two specific class of target statistical models and evaluate both the sampling error and the deficiency $\delta(\mathcal{U},\mathcal{V}; \mathcal{T}_{N})$ in each case.

\subsection{General location target models}\label{sec:general_location_target} \label{sec:location-general}
In this section, we consider non-Gaussian target models that share the same location parameter as the Gaussian source model~\eqref{source-gaussian-location}. Specifically, we consider a target model $\mathcal{V} = (\mathbb{R}, \{\mathcal{Q}_{\theta}\}_{\theta \in \mathbb{R}})$, where each distribution $\mathcal{Q}_{\theta}$ has a density of the form
\begin{align}\label{density-non-gaussian-target}
  v(y; \theta) = \frac{1}{\varT}\exp\left(- \psi\left( \frac{y - \theta}{\varT} \right) \right), \quad \text{for all } y, \theta \in \mathbb{R}.
\end{align}
Here, $\psi : \mathbb{R} \to \mathbb{R}$ denotes the negative log-density function, and $\varT > 0$ is a scale parameter that controls (among other things) the variance of the distribution. For our reduction to be effective, we require some mild regularity conditions on the negative log-density function $\psi$.
\begin{assumption}\label{assump:log-density}
Suppose $\psi: \real \rightarrow \real$ satisfies the following conditions:
\begin{enumerate}
  \item[(a)] $\int_{\real} \exp\big(-\psi(t) \big) \mathrm{d}t = 1$ and $\int_{\real} \exp\big(-2\psi(t) \big) \mathrm{d}t<+\infty$.
  \item[(b)] $\psi$ is infinitely differentiable on $\real$ and there exist constants $\Cdensity \in \PosInt$ and $\CdensityT>0$ depending only on $\psi$ such that for all $n \in \mathbb{N}$ and $k \in [n]$, we have
  \begin{align*}
          \int_{\real} \exp\big(-\psi(t)\big) \bigg( \sum_{j=1}^{n} \frac{\big| \psi^{(j)}(t) \big|}{ \big( \CdensityT \big)^{j} j! } \bigg)^{k} \mathrm{d}t \leq (\Cdensity n)!.
  \end{align*}
\end{enumerate} 
\end{assumption}
Note that Assumption~\ref{assump:log-density} is very mild and encompasses a broad class of smooth densities. In particular, Assumption~\ref{assump:log-density}(a) ensures that the function $v(y;\theta)$ defined in Eq.~\eqref{density-non-gaussian-target} is indeed a valid density function. 
To interpret Assumption~\ref{assump:log-density}(b), recall the definition of analytic functions: a function $\psi: \mathbb{R} \rightarrow \mathbb{R}$ is called \emph{analytic} if for every compact set $D \subset \mathbb{R}$, there exists a constant $C > 0$ (depending only on $D$) such that
\begin{align}\label{condition-analytical}
    |\psi^{(j)}(x)| \leq C^{j+1} j! \quad \text{for all } x \in D \text{ and } j \in \PosInt.
\end{align}
Assumption~\ref{assump:log-density}(b) requires that $\psi$ is analytic with respect to the measure $e^{-\psi(x)} \, \mathrm{d}x$, which includes a rich class of functions. In particular, it contains all \emph{globally analytic} functions, which are functions that satisfy this condition uniformly over $\mathbb{R}$:
\begin{definition}\label{def:analytical-functions}
For a parameter $R > 0$, we define a class of globally analytic functions:
\begin{align}\label{function-class-bound-derivative}
   \mathcal{F}(R) := \bigg\{ \psi: \mathbb{R} \rightarrow \mathbb{R} \;\bigg|\, \sum_{j=1}^{\infty} \frac{|\psi^{(j)}(x)|}{R^{j} j!} \leq 1 \quad \text{for all } x \in \mathbb{R} \bigg\}.
\end{align}
\end{definition}
See Corollary~\ref{example:generalized-normal} and Example~\ref{examples:nongaussian} to follow, which illustrate that many canonical distributions satisfy Assumption~\ref{assump:log-density} and belong to the class $\mathcal{F}(R)$.

Before stating our theorem showcasing the reduction guarantee, it is useful to define a few quantities that will appear in the theorem statement. For the constants $(\Cdensity, \CdensityT)$ appearing in Assumption~\ref{assump:log-density}(b), we define the constant $\CdensityB$ for convenience as
\begin{align}\label{def-C-B}
    \CdensityB := 16 \pi^{2}(2e)^{4(1+\Cdensity)} e^{4(1+\Cdensity)^{2}}.
\end{align}
For $\lambda>0$ and $n \in \mathbb{N}$, we define a set $\Set(\psi,n,\lambda) \subseteq \real$ and also its complement as
\begin{align}\label{def-set-complement}
    \Set(\psi, n, \lambda): = \bigg\{ t \in \real \;\Big|\; \sum_{j=1}^{n}\frac{ |\psi^{(j)}(t)| }{\lambda^{j} j!} \leq 1 \bigg\} \quad \text{and} \quad \Set^{\complement}(\psi, n, \lambda) : = \real \setminus  \Set(\psi, n, \lambda).
\end{align}
Using these objects, we now define two quantities that will be used to measure the deficiency and computational complexity of the reduction, respectively. First, for $\lambda>0$ and $n \in \mathbb{N}$, define
\begin{align}\label{def-error-nongaussian}
    \err(\psi,n,\lambda) : =  \int_{\Set^{\complement}(\psi, n, \lambda)} e^{-\psi(t)} \mathrm{d}t + \max_{ k \in [n]} \int_{\Set^{\complement}(\psi, n, \lambda)} e^{-\psi(t)} \bigg( \sum_{j=1}^{n} \frac{|\psi^{(j)}(t)|}{\lambda^{j} j!} \bigg)^{k} \mathrm{d}t.
\end{align}
Second, for $\lambda>0$ and $n \in \mathbb{N}$, define
\begin{align}\label{def-M-nongaussian}
  M(\psi,n,\lambda): = \max_{k \in [n]}\; \sup_{t \in \real}\;\Bigg\{ e^{-\psi(t) + \psi(t/2)} \cdot \Bigg[ 1 + \bigg( \sum_{j=1}^{n} \frac{|\psi^{(j)}(t)|}{\lambda^{j} j!} \bigg)^{k}\Bigg] \Bigg\}.
\end{align}
Note that both $\err(\psi,n,\lambda)$ and $M(\psi,n,\lambda)$ decrease as $\lambda$ increases. In applications, one can choose a sufficiently large $\lambda$ to ensure that these two quantities are small. We are now ready to state the main result.
\begin{theorem}\label{thm:non-gaussian}
  Consider source model $\mathcal{U} = (\mathbb{R}, \{ \mathcal{N}(\theta, \sigma^{2}) \}_{\theta \in \real})$ and target model $\mathcal{V} = (\real, \{\mathcal{Q}_{\theta}\}_{\theta \in \real})$, where $\mathcal{Q}_{\theta}$ has density function $v(\cdot;\theta)$ defined in Eq.~\eqref{density-non-gaussian-target}. For each $N \in \mathbb{N}$, consider the signed kernel $\mathcal{S}^{\star}_{N}$ defined in Eq.~\eqref{truncated-signed-kernel} and the corresponding Markov kernel $\mathcal{T}^{\star}_{N}$ defined in Eq.~\eqref{close-markov-kernel}. For $\lambda>0$ and $n \in \mathbb{N}$, let $\err(\psi,n,\lambda)$ and $M(\psi,n,\lambda)$ be defined as in Eq.~\eqref{def-error-nongaussian} and Eq.~\eqref{def-M-nongaussian}, respectively.
  Suppose the negative log-density function $\psi:\real \rightarrow \real$ satisfies Assumption~\ref{assump:log-density} with a pair of constants $(\Cdensity,\CdensityT)$, and that the density $e^{-\psi}$ can be sampled in time $T_{\mathsf{samp}}$. Suppose the tuple of parameters $(N,\varT, \lambda,\sigma)$ satisfies
  \begin{align}\label{assump:thm-nongaussian}
     2N+1 \geq \CdensityB, \quad \varT \geq 8e\sigma \max\big\{ \CdensityT, \lambda N \big \}, \quad \text{and} \quad  \sigma,\lambda>0, 
  \end{align}
  where $\CdensityB$ is defined in Eq.~\eqref{def-C-B}. Then the following statements are true:
  \begin{enumerate}
  \item[(a)] We have
  \begin{align}\label{ineq-nongaussian-term1}
        \delta \big(\mathcal{U}, \mathcal{V} ; \mathcal{T}_{N}^{\star}\big) \leq \CdensityB \exp\bigg( - (2e)^{-1} (2N+1)^{\frac{1}{2(1+\Cdensity)}} \bigg)   + 4\err(\psi,2N,\lambda).
  \end{align}
  \item[(b)]
Consequently, for each $\epsilon \in (0,1)$, if the tuple $(N,\varT,\lambda)$ satisfies Ineq.~\eqref{assump:thm-nongaussian} and also 
\begin{align}\label{parameter-N-varT-nongaussian}
\begin{split}
    2N+1 \geq  \Big(2e\log\big(3\CdensityB/\epsilon\big) \Big)^{2(1+\Cdensity)},
\end{split}
\end{align}
then there is a reduction algorithm $\mathsf{K}:\real \rightarrow \real$ that runs in time 
\[
  \mathcal{O}\Big( T_{\mathsf{iter}}\cdot T_{\mathsf{samp}} \cdot T_{\mathsf{eval}} (N) \Big) \quad \text{for any} \quad T_{\mathsf{iter}} \in \mathbb{N},
\]
and satisfies that, for $X_{\theta} \sim \mathcal{N}(\theta, \sigma^{2})$ and $Y_{\theta} \sim Q_{\theta}$,
\begin{align}\label{ineq-nongaussian-term2}
\sup_{\theta \in \real}\; \mathsf{d}_{\mathsf{TV}} \big( \mathsf{K}(X_{\theta}), Y_{\theta} \big) \leq \exp\Bigg( - \frac{ T_{\mathsf{iter}} }{ 4M(\psi,2N,\lambda) } \Bigg) +  \frac{\epsilon}{3} + 4\err(\psi,2N,\lambda).
\end{align}
\end{enumerate}
\end{theorem}

We provide the proof of Theorem~\ref{thm:non-gaussian} in Section~\ref{sec:pf_thm_non_gaussian}. 
A few remarks are in order. To achieve $\epsilon$-TV deficiency, we must choose $N$, $\lambda$, and $T_{\mathsf{iter}}$ such that the right-hand sides of Ineqs.~\eqref{ineq-nongaussian-term1} and~\eqref{ineq-nongaussian-term2} are both smaller than~$\epsilon$. 
As a first step, observe that since the constants $(\Cdensity, \CdensityT)$ are independent of $\epsilon$, it suffices to set
\[
    N = \mathcal{O}\big(\mathsf{polylog}(1/\epsilon)\big).
\]

Given such a choice of $N$ specified by Ineq.~\eqref{parameter-N-varT-nongaussian}, we next select a sufficiently large $\lambda$ so that the quantities $\err(\psi,2N,\lambda)$ and $M(\psi,2N,\lambda)$ become small. Ideally, if we are able to choose
\[
\lambda = \mathcal{O}\big(\mathsf{polylog}(1/\epsilon)\big) \quad \text{such that} \quad \err(\psi, 2N, \lambda) \leq \frac{\epsilon}{12} \quad \text{and} \quad M(\psi, 2N, \lambda) = \mathcal{O}(1),
\]
then, by taking $T_{\mathsf{iter}} = \mathcal{O}(\log(3/\epsilon))$, Ineq.~\eqref{ineq-nongaussian-term2} ensures
\[
\sup_{\theta \in \mathbb{R}}\, \mathsf{d}_{\mathsf{TV}} \big( \mathsf{K}(X_{\theta}), Y_{\theta} \big) \leq \epsilon.
\]

Ultimately, the choices of $N$ and $\lambda$ necessitate an increase in the scale parameter $\varT^2$ of the target distribution, as stipulated by condition~\eqref{assump:thm-nongaussian}. In other words, the reduction incurs a cost in terms of the scale parameter (or equivalently, a decrease in signal-to-noise ratio). However, this cost is mild: we only require $\varT \gtrsim \sigma \, \mathsf{polylog}(1/\epsilon)$, which is merely \emph{polylogarithmic} in~$1/\epsilon$. This is a significant improvement over the plug-in approach, which typically requires $\varT$ to grow \emph{polynomially} in $1/\epsilon$; see~\citet[Proposition~1]{lou2025computationally} for a discussion of the inefficiency of the plug-in approach in a related context.

Furthermore, provided that the evaluation time $T_{\mathsf{eval}}(N)$ of the signed kernel (see Definition~\ref{assump:Teval}) is polynomial in $1/\epsilon$, the overall reduction can be carried out in time $\mathcal{O}(\mathsf{poly}(1/\epsilon))$, thereby ensuring that the reduction is computationally efficient.

Operationally, Theorem~\ref{thm:non-gaussian} can be viewed as a universality result in the setting of a \emph{single} observation, and affirmatively answers puzzle (P1) from the introduction. Specifically, given only \emph{one} observation $X_{\theta} \sim \NORMAL(\theta,\sigma^{2})$—--without knowledge of the true mean $\theta$, which can take any real value--—one can generate a new observation with the same mean, but whose distribution is transformed into a general non-Gaussian form. 
Such a reduction technique is particularly valuable for establishing the universality of computational hardness: lower bounds derived under Gaussian observations can be transferred to settings with general non-Gaussian observations. We demonstrate the use of this reduction in Section~\ref{sec:PCA-nongaussian} through the example of tensor PCA with non-Gaussian noise.

While Theorem~\ref{thm:non-gaussian} is abstractly stated, it applies to many concrete families of target densities. We prove two broad examples of such families, beginning with log-densities that are globally analytic according to Definition~\ref{def:analytical-functions}.
\begin{corollary}\label{example:bounded-psi}
  For $R>0$, let the class of functions $\mathcal{F}(R)$ be defined as in Definition~\ref{def:analytical-functions}. Suppose $\psi \in \mathcal{F}(R)$ and define $M := \sup_{t\in \real} \big\{ e^{-\psi(t) + \psi(t/2)} \big\}$. Then Assumption~\ref{assump:log-density} holds with $\Cdensity = 0$ and $\CdensityT = R$, and 
  \[
      \err(\psi,n,R) = 0 \quad \text{and} \quad M(\psi,n,R) = 2M \quad \text{for all } n \in \mathbb{N}.
  \]
  Moreover, if $\psi^{(n)}(x)$ can be evaluated in computational time $\mathcal{O}\big(e^{\mathcal{O}(\sqrt{n})} \big)$ for all $x\in \real$ and $n \in \mathbb{N}$, then $T_{\mathsf{eval}}(n) = \mathcal{O}\big(e^{\mathcal{O}(\sqrt{n})} \big)$.
  Consequently, for each $\epsilon \in (0,1)$, there exists a universal and positive constant $\widetilde{C}$ such that for 
  \[
      \varT \geq \widetilde{C} R \sigma \log^{2}\big(\widetilde{C}/ \epsilon\big),
  \]
  there is a reduction algorithm $\mathsf{K}:\real \rightarrow \real$ that runs in time 
  \[
      \mathcal{O}\Big( M \cdot T_{\mathsf{samp}} \cdot \mathsf{poly}(1/\epsilon) \Big) \quad \text{and satisfies} \quad \sup_{\theta \in \real}\; \mathsf{d}_{\mathsf{TV}} \big( \mathsf{K}(X_{\theta}), Y_{\theta} \big) \leq \epsilon,
  \]
  where $X_{\theta} \sim \mathcal{N}(\theta, \sigma^{2})$ and the random variable $Y_{\theta}$ has density $v(\cdot;\theta)$ defined in Eq.~\eqref{density-non-gaussian-target}.
\end{corollary}
We provide the proof of Corollary in Section~\ref{sec:pf_corollary_bounded_psi}. Note that it is common for the higher-order derivatives $\psi^{(n)}(x)$ to be computable in time $\mathcal{O}(e^{\mathcal{O}(\sqrt{n})})$. For a canonical example, consider the case where $\psi(x) = g \circ h(x)$ is the composition of two basic functions. By Faà di Bruno's formula, we have
\begin{align}\label{fadi-bruno-psi}
    \psi^{(n)}(x) = \sum_{m_{1},\dots,m_{n}} \frac{n!}{m_{1}! \,m_{2}! \cdots m_{n}! } \cdot g^{(m_{1}+\cdots+m_{n})}(h(x)) \cdot \prod_{j=1}^{n} \left( \frac{h^{(j)}(x)}{j!} \right)^{m_{j}},
\end{align}
where the summation is over all $n$-tuples of nonnegative integers $(m_{1},\dots,m_{n})$ satisfying $\sum_{j=1}^{n} j m_{j} = n$. The total number of such tuples is at most $\mathcal{O}(e^{\pi \sqrt{n}})$. If each term in the summation can be computed in time $\mathcal{O}(n)$—as $g$ and $h$ are assumed to be simple functions—then the total time to evaluate $\psi^{(n)}(x)$ is $\mathcal{O}(n \cdot e^{\pi \sqrt{n}})$, which satisfies our computational assumption. 

In the following example, we further specialize Corollary~\ref{example:bounded-psi} to illustrate that it covers natural target distributions.

\begin{example}\label{examples:nongaussian}
The following distributions have globally analytic log-densities and so Corollary~\ref{example:bounded-psi} applies. In particular, we obtain the following results.
\begin{enumerate}
  \item Student't t-distribution: let $\nu \geq 1$ be the degree of freedom and consider
  \[
      \psi(t) = \frac{\nu+1}{2} \log\bigg(1+\frac{t^{2}}{\nu}\bigg) - \log\bigg(\frac{\Gamma(\frac{v+1}{2})}{\sqrt{\pi \nu} \Gamma(\frac{\nu}{2})} \bigg)\quad \text{and} \quad v(y;\theta) = \frac{1}{\varT} \frac{\Gamma(\frac{v+1}{2})}{\sqrt{\pi \nu} \Gamma(\frac{\nu}{2})} \bigg( 1+\frac{(\frac{y-\theta}{\varT})^{2}}{\nu} \bigg)^{-\frac{\nu +1}{2}}.
  \]
  Then $\psi \in \mathcal{F}\big(6e(\nu+1) \big)$ and $M \leq 1$. Computing $\psi^{(n)}(x)$ takes time $\mathcal{O}(\mathsf{poly}(n))$ for all $x\in \real$, $n\in \mathbb{N}$.

  \item Hyperbolic secant distribution: consider 
  \[
      \psi(t) = \log\Big( e^{\frac{\pi t}{2}} + e^{-\frac{\pi t}{2}} \Big) \quad \text{and} \quad v(y;\theta) = \frac{1}{\varT} \cdot \frac{1}{\exp\big(\frac{\pi}{2} \frac{y-\theta}{\varT} \big) +  \exp\big(-\frac{\pi}{2} \frac{y-\theta}{\varT} \big) }.
  \]
  Then $\psi \in \mathcal{F}\big(2\pi e\big)$ and $M \leq 2$. Computing $\psi^{(n)}(x)$ takes time $\mathcal{O}\big(e^{5\sqrt{n}}\big)$ for all $x\in \real$, $n\in \mathbb{N}$.

  \item Logistic distribution: consider 
  \[
      \psi(t) = 2\log\Big( e^{t/2} + e^{-t/2} \Big) \quad \text{and} \quad v(y;\theta) = \frac{1}{\varT} \cdot \frac{1}{\Big[ \exp\big(\frac{y-\theta}{2\varT} \big) +  \exp\big(- \frac{y-\theta}{2\varT} \big) \Big]^{2}}.
  \]
  Then $\psi \in \mathcal{F}\big(2\pi e\big)$ and $M \leq 4$. Computing $\psi^{(n)}(x)$ takes time $\mathcal{O}\big(e^{5\sqrt{n}}\big)$ for all $x\in \real$, $n\in \mathbb{N}$.
\end{enumerate}
\end{example}

\begin{remark}\label{remark:eval}
Note that evaluating the signed kernel $\mathcal{S}_{N}^{\star}(y\mid x)$~\eqref{truncated-signed-kernel} for the target density $v(y;\theta)$~\eqref{density-non-gaussian-target} requires computing the higher-order derivative $\psi^{(2N)}(\cdot)$. Condition~\eqref{parameter-N-varT-nongaussian} dictates that $N \gtrsim \big(\log(1/\epsilon)\big)^{2(1+\Cdensity)}$. In the setting of Corollary~\ref{example:bounded-psi}, we have $\Cdensity = 0$, so we take $N \asymp \log^{2}(1/\epsilon)$. In this case, evaluating $\psi^{(2N)}(\cdot)$ using formula~\eqref{fadi-bruno-psi} requires computational time
\[
    \mathcal{O}\big(e^{\mathcal{O}(\sqrt{N})}\big) \;=\; \mathsf{poly}(1/\epsilon).
\]
However, in the more general case $\Cdensity \geq 1$, a naive application of formula~\eqref{fadi-bruno-psi} leads to computational time of order
\[
    e^{\sqrt{N}} \;=\; \exp \Big( \big( \log(1/\epsilon) \big)^{1+\Cdensity} \Big),
\]
which grows quasi-polynomially in $1/\epsilon$. Therefore, when $\Cdensity \geq 1$, a more computationally efficient algorithm for evaluating higher-order derivatives is required. 
The following Corollary~\ref{example:generalized-normal} serves as an illustration.
\end{remark}

Note that Theorem~\ref{thm:non-gaussian} also applies when the negative log density $\psi$ is not globally analytic. We next illustrate this with the case where $\psi$ is a monomial of fixed degree, which corresponds to the family of generalized normal distributions. 

\begin{corollary}[Generalized normal target]\label{example:generalized-normal}
Let $\beta \in \mathbb{N}$ be a positive and even integer. Consider 
\[
  \psi(t) = t^{\beta} - \log \bigg( \frac{\beta}{2\Gamma(1/\beta)} \bigg) \quad \text{and} \quad v(y;\theta) = \frac{\beta}{2\varT \Gamma(1/\beta)} \exp\bigg( - \frac{|y-\theta|^{\beta}}{\varT^{\beta}} \bigg),
\]
where $\Gamma(z) = \int_{0}^{\infty} t^{z-1} e^{-t} \mathrm{d}t$ for all $z>0$. Assumption~\ref{assump:log-density} holds with $\Cdensity = 1$ and $\CdensityT = 2\beta^{2}$, and 
\[
  \err(\psi,n,\lambda) \leq  3e^{-\sqrt{\lambda}}, \quad  M(\psi,n,\lambda) \leq 3, \quad T_{\mathsf{eval}}(n) =\mathcal{O} \big(n^{\beta+2}\big) \quad \text{for all } n\in \mathbb{N}, \lambda \geq 4n^{2} \beta^{2}.
\]
Consequently, for each $\epsilon \in (0,1)$, there exists a pair of universal and positive constants $(C_{1},C_{2})$ such that if
\[
      \varT \geq C_{1} \sigma \beta^{2} \log^{12}\bigg( \frac{C_{2}}{\epsilon}\bigg),
\]
then there is a reduction algorithm $\mathsf{K}:\real \rightarrow \real$ that runs in time 
\[
   \mathcal{O}\Big(\log^{4\beta+5}\big(C_{2}/\epsilon\big)\Big) \quad \text{and satisfies} \quad  \sup_{\theta \in \real}\; \mathsf{d}_{\mathsf{TV}} \big( \mathsf{K}(X_{\theta}), Y_{\theta} \big) \leq \epsilon,
\]
where $X_{\theta} \sim \mathcal{N}(\theta, \sigma^{2})$ and $Y_{\theta} \sim \mathsf{GN}(\theta,\varT,\beta)$.
\end{corollary}
We provide the proof of Corollary~\ref{example:generalized-normal} in Section~\ref{sec:pf_corollary_generalized_normal}.
Note that in Corollary~\ref{example:generalized-normal}, the target distribution has much lighter tails than the source distribution when $\beta \geq 4$. Also note that fast computation is still possible by using specific formulas for the derivatives (cf. Remark~\ref{remark:eval}).

\begin{remark}  
First, Corollary~\ref{example:generalized-normal} does not cover the case when $\beta$ is odd, e.g., the Laplace distribution with $\beta = 1$. When $\beta$ is odd, the target density is non-differentiable at one point, which violates Assumption~\ref{assump:log-density}(b), and our technique relies crucially on the differentiability of the target density. One could apply a mollification technique to smooth the target density; however, calculations show that to achieve $\epsilon$-TV deficiency, this would require inflating the scale parameter $\varT$ \emph{polynomially} in $1/\epsilon$. Second, when $\beta$ is chosen on the order of $1/\epsilon$, the generalized normal distribution is $\epsilon$-close to the uniform distribution in TV distance. In this case, Corollary~\ref{examples:nongaussian} also requires inflating $\varT$ polynomially in $1/\epsilon$. Thus, Corollary~\ref{examples:nongaussian} does not cover targets with extremely light tails such as the uniform distribution. Interestingly, a reduction from the uniform to the Gaussian location model is known~\citep{lou2025computationally}.
\end{remark}

\subsection{Gaussian target models with mean nonlinearly transformed}\label{sec:gaussian_target_nonlinear}
In this section, we consider a Gaussian target model where the mean parameter $\theta$ is transformed nonlinearly. Specifically, we study
\begin{align}\label{target:Gaussian}
  \mathcal{V} = \left( \Theta, \left\{ \NORMAL\left(\link\left(\theta/\varT\right), 1 \right) \right\}_{\theta \in \Theta} \right),
\end{align}
where $\link : \mathbb{R} \rightarrow \mathbb{R}$ is a link function that nonlinearly transforms the mean parameter $\theta$, and $\varT > 0$ is a fixed parameter that adjusts the signal-to-noise ratio. Since $\varT$ already serves this purpose, we fix the variance of the target distribution to be $1$ for convenience.\footnote{Both the reduction and the analysis can be readily extended to the general case where the target model~\eqref{target:Gaussian} has arbitrary variance. However, fixing the variance to $1$ simplifies exposition and suffices for our applications.} Note that the target distribution $\mathcal{Q}_{\theta}$ has the density
\begin{align}\label{target:density}
  v(y;\theta) = \frac{1}{\sqrt{2\pi}} \exp \left( -\frac{\left(y - \link\left(\frac{\theta}{\varT}\right)\right)^2}{2} \right), \quad \text{for all } y \in \mathbb{R},\ \theta \in \Theta.
\end{align}

We next state some mild regularity conditions on the link function $\link:\real \rightarrow \real$. 
\begin{assumption}\label{assump:link}
Suppose $f:\real \rightarrow \real$ is infinitely differentiable on $\real$ and there exists a pair of universal constants $(\Clink,\ClinkT)$, where $\Clink \in \PosInt$ and $\ClinkT>0$, that only depends on the link function $\link:\real \rightarrow \real$ such that 
\begin{align}\label{ineq:link}
    \int_{\real} \Bigg( \sum_{j=1}^{n} \frac{\big| \link^{(j)}\big( \big(\theta + \sqrt{2} x \big)/\varT \big) \big| }{\big(\ClinkT \big)^{j} j!} \Bigg)^{n} \frac{\exp(-x^{2}/2)}{\sqrt{2\pi}} \mathrm{d}x \leq \big( \Clink n \big)!,
\end{align}
for all $n \in \mathbb{N}$ and $\varT \geq 2\sqrt{2} (1 \vee |\theta|)$.
\end{assumption}
Note that Assumption~\ref{assump:link} imposes the analyticity of the function $f$ with respect to the standard Gaussian measure, which is analogous to Assumption~\ref{assump:log-density}(b). This condition encompasses a broad class of functions, including polynomials of constant degree (see Corollary~\ref{example:monomial}) as well as any globally analytic function (see Corollary~\ref{example:analytical}).

To state our theorem, we require formal definitions of several quantities derived from the link function $\link$. Let $(\Clink,\ClinkT)$ be the same pair of constants as in Assumption~\ref{assump:link} and for all $n \in \mathbb{N}$ and $\varT>0$, define
\begin{align}\label{def-quantity-link}
\begin{split}
  \ClinkB: = 32\pi^{2}(2e)^{4(1+\Clink)} &e^{4(1+\Clink)^{2}}, \qquad 
  \Llink(n) : = \sup_{|x| \leq 1}\; \sum_{j=1}^{n} \frac{ |\link^{(j)}(x)|}{\big(\ClinkT\big)^{j}j!},\\
  \text{and} \quad \LlinkT(n,\varT) &:= \sup_{|x|\geq 1}\; \sum_{j=1}^{n} \frac{ \big|\link^{(j)}(x)\big|}{\big(\ClinkT\big)^{j}j!} \cdot \exp\bigg(-\frac{\varT^{2}x^{2}}{8n} \bigg).
\end{split}
\end{align}
Note that the quantity $\Llink(n)$ ought to be bounded for any reasonable function $f$, given that the supremum is calculated over the range $|x| \leq 1$. Moreover, by increasing $\varT$, we can effectively ensure that $\LlinkT(n, \varT)$ also remains bounded. We are now ready to state the theorem.
\begin{theorem}\label{thm:gaussian-mean-nonlinear}
  Consider Gaussian source model $\mathcal{U} = (\mathbb{R}, \{ \mathcal{N}(\theta, 1) \}_{\theta \in \Theta})$ and Gaussian target model $\mathcal{V} = \big( \real, \big \{ \NORMAL\big(\link(\theta/\varT),1 \big) \big\}_{\theta \in \Theta}  \big)$. For each $N \in \mathbb{N}$, consider the signed kernel $\mathcal{S}^{\star}_{N}$ defined in Eq.~\eqref{truncated-signed-kernel} and the corresponding Markov kernel $\mathcal{T}^{\star}_{N}$ defined in Eq.~\eqref{close-markov-kernel}. Suppose the link function $f:\real \rightarrow \real$ satisfies Assumption~\ref{assump:link} with a pair of constants $(\Clink,\ClinkT)$. For $\varT>0$ and $n\in \mathbb{N}$, let $\ClinkB$, $\Llink(n)$ and $\LlinkT(n,\varT)$ be defined in Eq.~\eqref{def-quantity-link}. There exists a universal and positive constant $C$ such that if
  \begin{align}\label{assump:thm-nonlinear}
     2N+1 \geq \ClinkB \quad \text{and} \quad \varT^{2} \geq  C \bigg[ \big(  \ClinkT \vee 1 \big)^{4}  \Big( \Llink(2N) \vee \LlinkT(2N,\varT) \vee 1\Big)^{2}N^{2} \; \vee \; \sup_{\theta \in \Theta} \theta^{2} \bigg],
  \end{align}
  then the following statements hold.
  \begin{enumerate}
  \item[(a)] We have
  \begin{align}\label{ineq-gaussian-nonlinear-term0}
        \delta \big(\mathcal{U}, \mathcal{V} ; \mathcal{T}_{N}^{\star}\big) \leq \ClinkB \exp\bigg(-\frac{(2N+1)^{\frac{1}{2(1+\Clink)}}}{2e} \bigg) + 
        12e^{-\varT/8} + 15\sup_{\theta \in \Theta} \; \exp\bigg( \frac{\theta^{2}}{2} -\frac{\varT^{2}}{8} \bigg).
  \end{align}
  \item[(b)] Consequently, for each $\epsilon \in (0,1)$, if $N$ and $\varT$ are setting to satisfy Eq.~\eqref{assump:thm-nonlinear} and also 
  \begin{align}\label{parameter-N-varT-nonlinear}
  \begin{split}
         2N+1 \geq \Big( 2e \log \big(4\ClinkB/\epsilon\big) \Big)^{2(1+\Clink)}, \quad  
         \varT^{2} \geq 64 \log^{2}\big(96/\epsilon \big) \vee \; \Big[ 8 \log\big(120/\epsilon\big) + 4\sup_{\theta \in \Theta} \theta^{2} \Big],
  \end{split}
  \end{align}
  then there exists a reduction algorithm $\mathsf{K}:\real \rightarrow \real$ that runs in time 
  \begin{align}\label{ineq-gaussian-nonlinear-term3}
    \mathcal{O}\big(\log(4/\epsilon) T_{\mathsf{eval}} (N) \big) \quad \text{and satisfies} \quad \sup_{\theta \in \Theta}\; \mathsf{d}_{\mathsf{TV}} &\big( \mathsf{K}(X_{\theta}), Y_{\theta} \big) \leq \epsilon,
  \end{align}
  where $X_{\theta} \sim \mathcal{N}(\theta, 1)$ and $Y_{\theta} \sim \NORMAL\big(\link(\theta/\varT),1 \big)$.
\end{enumerate}
\end{theorem}
We provide the proof of Theorem~\ref{thm:gaussian-mean-nonlinear} in Section~\ref{sub:pf_thm_gaussian_mean_nonlinear}. A few remarks about Theorem~\ref{thm:gaussian-mean-nonlinear} are in order. To achieve $\epsilon$-TV deficiency, it suffices to set $N = \Theta(\mathsf{polylog}(1/\epsilon))$, since the constants $(\Clink, \ClinkT)$ are independent of $\epsilon$. Consequently, if the link function $f$ satisfies 
\begin{align} \label{eq:link-condn}
\Llink(2N) = \mathcal{O}(\mathsf{polylog}(1/\epsilon)), \quad \LlinkT(2N,\varT) = \mathcal{O}(\mathsf{polylog}(1/\epsilon))\quad  \text{for some } \varT = \mathcal{O}(\mathsf{polylog}(1/\epsilon) ),
\end{align}
and if the parameter set $\Theta \subset \mathbb{R}$ is bounded, then the reduction can attain $\epsilon$-TV deficiency provided that
\[
\varT^{2} \gtrsim \mathsf{polylog}(1/\epsilon).
\]
We will see shortly that condition~\eqref{eq:link-condn} is satisfied for many canonical link functions that obey natural regularity conditions.

As in Theorem~\ref{thm:non-gaussian}, this reduction comes at the cost of requiring a lower signal-to-noise ratio (SNR) in the target model~\eqref{target:Gaussian}---via inflating $\tau$---compared to the source model~\eqref{source-gaussian-location}. However, the degradation remains mild: the SNR is reduced only by a $\mathsf{polylog}(1/\epsilon)$ multiplicative factor. This represents a significant improvement over the plug-in approach, which typically requires reducing the SNR by a polynomial factor in $1/\epsilon$. For concreteness, we show in Appendix~\ref{sec:plug_in_approach} that the plug-in approach must incur this blow-up even when $f(\theta) = \theta^2$.

Operationally, Theorem~\ref{thm:gaussian-mean-nonlinear} provides an answer to puzzle (P2) posed in the introduction for a large class of functions. 
In Section~\ref{sec:hardness_glm}, we demonstrate a concrete application by establishing the computational hardness of the $k$-sparse generalized linear model with even link functions via reduction from sparse mixtures of linear regressions. In Section~\ref{sec:sparse_spiked_wigner_BCering}, we employ Theorem~\ref{thm:gaussian-mean-nonlinear} to construct a reduction between the two high-dimensional models $\SubMat$ and $\PlantSubMat$.

Let us now instantiate the abstract Theorem~\ref{thm:gaussian-mean-nonlinear} for several canonical function families. 

\begin{corollary}[Monomial]\label{example:monomial}
Consider link function $\link(\theta) = \theta^{B}$ for some $B \in \mathbb{N}$. Then Assumption~\ref{assump:link} holds with $\Clink = 2B+1$ and $\ClinkT = 4eB$. Moreover, 
we have for all $n \in \mathbb{N}$ and $\varT^{2} \geq 8n B$ that
\[
   \Llink(n) \leq 1, \quad \LlinkT(n,\varT) \leq 1, \quad \text{and} \quad T_{\mathsf{eval}}(n) = \mathcal{O}(n^{B+2}). 
\]
Consequently, for each $\epsilon \in (0,1)$, there exists a universal constant $C_{B} \geq 1$ that only depends on $B$ such that if
\[
      \varT^{2} \geq C_{B} \cdot \log^{8(B+1)}(C_{B}/\epsilon) + 4 \sup_{\theta \in \Theta} \theta^{2},
\]
then there is a reduction algorithm $\mathsf{K}:\real \rightarrow \real$ that runs in time 
\[
    \mathcal{O}\Big(C_{B} \log^{4(B+2)^{2}}(C_{B}/\epsilon) \Big) \quad \text{and satisfies} \quad \sup_{\theta \in \Theta}\; \mathsf{d}_{\mathsf{TV}} \big( \mathsf{K}(X_{\theta}), Y_{\theta} \big) \leq \epsilon,
\]
where $X_{\theta} \sim \mathcal{N}(\theta, 1)$ and $Y_{\theta} \sim \NORMAL\big(\link(\theta/\varT),1 \big)$.
\end{corollary}
We provide the proof of Corollary~\ref{example:monomial} in Section~\ref{sec:pf_corollary_monomial}. 

To produce other examples, we use the following lemma, which can be used to verify Assumption~\ref{assump:link} and bound the quantities defined in Eq.~\eqref{def-quantity-link} when the link function is a linear combination of a set of basic functions. 
\begin{lemma}\label{lemma:linear-combination}
Let $f_{k}:\real \rightarrow \real$ be a function satisfying Assumption~\ref{assump:link} with a pair of constants $\big(\widetilde{C}_{f_{k}},C_{f_{k}}\big)$ for all $k \in [K]$. Consider $f(x) = \sum_{k=1}^{K} c_{k} f_{k}(x)$ for all $x\in \real$ with linear coefficients $c_{k} \in \real$ for all $k\in [K]$. Then $f$ satisfies Assumption~\ref{assump:link} with the pair of constants $(\Clink,\ClinkT)$ defined as
\[
    \Clink = \max_{k \in  [K]} \{ C_{f_{k}}\} \quad \text{and} \quad \ClinkT = \sum_{\ell = 1}^{K} |c_{\ell}| \cdot \max_{k \in [K]} \{ \widetilde{C}_{f_{k}}\}.
\]
Moreover, we have for all $n \in \mathbb{N}$ and $\varT>0$, 
\[
    \Llink(n) \leq \frac{ \sum_{k=1}^{K} |c_{k}| \cdot L_{f_{k}}(n) }{\sum_{\ell = 1}^{K} |c_{\ell}|}  \quad \text{and} \quad \LlinkT(n,\varT) \leq 
    \frac{ \sum_{k=1}^{K} |c_{k}| \cdot \widetilde{L}_{f_{k}}(n,\varT) }{\sum_{\ell = 1}^{K} |c_{\ell}|}.
\]
\end{lemma}
The proof of Lemma~\ref{lemma:linear-combination} is provided in Appendix~\ref{sec:pf_lemma_linear_combi}.

Combining Corollary~\ref{example:monomial} with Lemma~\ref{lemma:linear-combination}, we obtain a reduction when the function $f$ is a polynomial of constant degree. This already encompasses several natural nonlinear transformations. To cover others that do not have a natural polynomial representation, we consider the class of link functions that are globally analytic in the sense of Definition~\ref{def:analytical-functions}.

\begin{corollary}\label{example:analytical}
For $R>0$, let the class of functions $\mathcal{F}(R)$ be defined as in Definition~\ref{def:analytical-functions}. Let $f: \real \rightarrow \real$ be any element in $\mathcal{F}(R)$ and suppose $f^{(n)}(x)$ can be evaluated in computational time $\mathcal{O}\big(e^{\mathcal{O}(\sqrt{n}) } \big)$ for all $x\in \real$ and $n \in \mathbb{N}$. Then Assumption~\ref{assump:link} holds with $\Clink = 0$ and $\ClinkT = 2R$. We have for all $n \in \mathbb{N}$ and $\varT>0$,
\[
    \Llink(n) \leq 1, \quad \LlinkT(n,\varT) \leq 1, \quad T_{\mathsf{eval}}(n) = \mathcal{O}\big(e^{\mathcal{O}(\sqrt{n}) } \big).
\]
Consequently, there exists a universal and positive constant $C$ such that if 
\[
    \varT^{2} \geq C R^{2} \log^{2}(C/\epsilon) + 4 \sup_{\theta \in \Theta} \theta^{2},
\]
there is a reduction algorithm $\mathsf{K}:\real \rightarrow \real$ that runs in time 
\[
    \mathcal{O}\Big( \mathsf{poly}(1/\epsilon) \Big) \quad \text{and satisfies} \quad \sup_{\theta \in \Theta}\; \mathsf{d}_{\mathsf{TV}} \big( \mathsf{K}(X_{\theta}), Y_{\theta} \big) \leq \epsilon,
\]
where $X_{\theta} \sim \mathcal{N}(\theta, 1)$ and $Y_{\theta} \sim \NORMAL\big(\link(\theta/\varT),1 \big)$.
\end{corollary}
There are many interesting link functions satisfying the conditions in Corollary~\ref{example:analytical}. We specify a few concrete examples below.
\begin{example}\label{examples}
The following functions are all globally analytic in the sense of Definition~\ref{def:analytical-functions}.
\begin{enumerate}
  \item \underline{Trigonometric function: $f(\theta) = \cos(a\theta+b)$ or $f(\theta) = \sin(a\theta + b)$ for any $a,b\in \real$.} We have $f \in \mathcal{F}(2|a|)$, and one can evaluate $f^{(n)}(x)$ in time $\mathcal{O}(1)$ for all $x\in \real$ and $n \in \mathbb{N}$.

  \item \underline{Hyperbolic function: $f(\theta) = \mathsf{sech}(\theta) = \frac{2}{e^{\theta}+e^{-\theta}}$.} We have $f \in \mathcal{F}(8e)$, and one can evaluate $f^{(n)}(x)$ in time $\mathcal{O}(e^{\pi \sqrt{n}})$ for all $n \in \mathbb{N}$ and $x\in \real$.

  \item \underline{Other examples: $f(\theta) = (1+\theta^{2})^{-1}$ or $f(\theta) = \exp(-\theta^{2})$ or $f(\theta) = (1+e^{-\theta})^{-1}$.} We have $f \in \mathcal{F}(8e)$, and one can evaluate $f^{(n)}(x)$ in time $\mathcal{O}\big( e^{\pi \sqrt{n}} \big)$ for all $n \in \mathbb{N}$ and $x\in \real$.
\end{enumerate}
\end{example}
We provide proof of the above claims in Section~\ref{sec:pf_concrete_examples}.


\section{Applications to statistical-computational tradeoffs}\label{sec:applications}
We now present three applications of our reductions in Section~\ref{sec:concrete_reductions} to the study of statistical-computational tradeoffs in some canonical high-dimensional statistics and machine learning models. Our computational hardness results are based on the $k$-partite hypergraph planted clique ($k$-$\mathrm{HPC}^{s}$) conjecture and the $k$-part bipartite planted clique conjecture ($k$-$\mathrm{BPC}$); see Appendix~\ref{sec:hardness_conjectures} for a detailed description. The two conjectures are also used in~\citet{brennan2020reducibility}. We will take these conjectures as given and state consequences of these hardness results as appropriate in the sequel. Our goal in this section is to cover the motivating examples mentioned in Section~\ref{sec:introduction} through the lens of reductions.

\subsection{Tensor principal component analysis with non-Gaussian noise}\label{sec:PCA-nongaussian}
In Tensor PCA~\citep{montanari2014statistical}, we observe a single order $s$ tensor $T(\beta,v)$ with dimensions $n^{\otimes s} = \underbrace{n \times \cdots \times n}_{s \text{ times}}$ given by
\begin{align}\label{tensor-pca-gaussian}
  T(\beta,v) = \beta (\sqrt{n} v)^{\otimes s} + \zeta.
\end{align}
The positive scalar $\beta$ denotes the signal-to-noise ratio, the unknown spike $v \in \mathbb{R}^n$ is typically sampled from a prior distribution, and the observations are corrupted by a symmetric Gaussian noise tensor $\zeta$. Specifically, let $ G \sim \NORMAL(0,1)^{\otimes n^{\otimes s} }$ be an asymmetric tensor with entries that are $\mathsf{i.i.d.}$ $\mathcal{N}(0,1)$, and define the symmetrized noise tensor $\zeta$ by
\begin{align}\label{def-tensor-gaussian-noise}
	\zeta_{i_{1},\dots,i_{s}} = \frac{1}{s!} \sum_{\pi \in \mathfrak{S}_{s}} G_{i_{\pi(1)},\dots,i_{\pi(s)} } \quad \text{for all } (i_{1},\dots,i_{s}) \in [n]^s.
\end{align}
For notational convenience, let $\sigma_{i_{1},\dots,i_{s}}^2$ denote the variance of $ \zeta_{i_{1},\dots,i_{s}} $. Since $s$ is a constant, all entries of $\zeta$ have constant-order variance. In particular, the minimum variance is $\sigma_{i_{1},\dots,i_{s}}^2 = 1/s!$ when all indices $(i_{1},\dots,i_{s})$ are distinct, while the maximum variance is $\sigma_{i_{1},\dots,i_{s}}^2 = 1$ when all indices are identical.

The task of recovering $v$ to within nontrivial $\ell_2$-error $o(1)$ exhibits a statistical-computational gap. On the one hand, it is information-theoretically possible to recover $v$ when $\beta = \widetilde{\omega}(n^{(1-s)/2})$ via exhaustive search~\citep{montanari2014statistical,lesieur2017statistical,chen2019phase,jagannath2020statistical,perry2020statistical}. On the other hand, all known polynomial-time algorithms require a stronger signal, namely $\beta = \widetilde{\Omega}(n^{-s/4})$, including spectral methods~\citep{montanari2014statistical}, the sum-of-squares hierarchy~\citep{hopkins2015tensor,hopkins2017power}, and spectral algorithms based on the Kikuchi hierarchy~\citep{wein2019kikuchi}. This statistical-computational gap in tensor PCA is believed to be fundamental, and the following reduction-based computational hardness result from~\citet[Theorem~10 and Lemma~102]{brennan2020reducibility} makes this precise.

\begin{proposition}\label{thm:test-hardness-pca-gaussian}
Let $n$ be a parameter and $s\geq 3$ be a constant. Let $\zeta \in \real^{n^{\otimes s}}$ be the symmetric Gaussian tensor defined in Eq.~\eqref{def-tensor-gaussian-noise} and define
\begin{align}\label{tensor-continuous-spike}
   T = \beta (\sqrt{n} v)^{\otimes s} + \zeta  \quad \text{where} \quad v \sim \mathsf{Unif}(\mathbb{S}^{n-1}). 
\end{align}
Suppose the $k$-$\mathrm{HPC}^{s}$ conjecture holds. If $\beta = \widetilde{o}\big(n^{-s/4}\big)$ and $\beta = \omega\big(n^{-s/2} \sqrt{s\log(n)}\big)$, then there \emph{cannot} exist an algorithm $\mathcal{A}: \real^{n^{\otimes s}} \rightarrow \mathbb{S}^{d-1}$ that runs in time $\mathsf{poly}(n)$ and satisfies
\begin{align}\label{estimation-requirement-TPCA-Gaussian}
    \big\langle \mathcal{A} (T), v \big \rangle = \Omega(1) \quad \text{with probability at least} \quad \Omega_{n}(1).
\end{align}
\end{proposition}
A few comments on Proposition~\ref{thm:test-hardness-pca-gaussian} are in order. First, the condition $\beta = \omega\big(n^{-s/2} \sqrt{s\log(n)}\big)$ is not an essential feature of the result and can be ignored, since it is weaker than the information-theoretic threshold $\beta = \widetilde{\omega}(n^{(1-s)/2})$ necessary for recovering $v$ with constant probability~\citep{montanari2014statistical,chen2019phase}.
Second, we remark that the requirement 
\[
  \big\langle \mathcal{A}(T), v \big\rangle = \Omega(1)
\] 
is weaker than the condition 
\[
  \big\| \mathcal{A}(T) - v \big\|_{2} = o(1).
\] 
Thus, any estimation algorithm $\mathcal{A}$ that achieves $\ell_{2}$ error $o(1)$ automatically satisfies requirement~\eqref{estimation-requirement-TPCA-Gaussian}. Third, Proposition~\ref{thm:test-hardness-pca-gaussian} as stated differs from the results of \citet[Theorem~10 and Lemma~102]{brennan2020reducibility} in two ways: In their setting, (i) the spike $\sqrt{n}v$ is a discrete vector in $\{-1,1\}^{n}$, and (ii) the Gaussian noise $\zeta$ is asymmetric. Nevertheless, the model we consider exhibits the same gap as the model in~\citet{brennan2020reducibility} because there is an exact reduction from the Tensor PCA model in \citet[Theorem~10]{brennan2020reducibility} to our model~\eqref{tensor-pca-gaussian}. This reduction proceeds by multiplying a uniformly random rotation matrix to each mode of the tensor, and achieves zero total variation deficiency. For completeness, we provide the proof of this claim in Appendix~\ref{sec:pf_claims_TPCA} and a 
formal proof of Proposition~\ref{thm:test-hardness-pca-gaussian} in Section~\ref{sec:pf_thm_hardness_pca_gaussian}.

The key takeaway from Proposition~\ref{thm:test-hardness-pca-gaussian} is that computational hardness result holds under the Gaussian noise model. This restriction is somewhat unsatisfactory, as real-world observations can be affected by various types of non-Gaussian noise---indeed, this has also been noticed in recent work~\citep{kunisky2025low}. We now use our results from Section~\ref{sec:concrete_reductions} to extend reduction-based hardness results for tensor PCA to a larger class of non-Gaussian distributions.

To be concrete, consider the same tensor PCA model as in Eq.~\eqref{tensor-pca-gaussian}, but with observations corrupted by general non-Gaussian noise. Specifically, we study the following model:
\begin{align}\label{tensor-pca-nongaussian}
  \widetilde{T}(\beta,v) = \beta (\sqrt{n} v)^{\otimes s} + \widetilde{\zeta},
\end{align}
Here the noise tensor $\widetilde{\zeta} \in \real^{n^{\otimes s}}$ is still symmetric---specifically, let $\mathcal{Q}_{0}$ denote the distribution of the noise, then
the entries $\widetilde{\zeta}_{i_{1},\dots,i_{s}} \sim \mathcal{Q}_{0}$ and satisfy $\widetilde{\zeta}_{i_{1},\dots,i_{s}} = \widetilde{\zeta}_{i_{\pi(1)},\dots,i_{\pi(s)}}$ for any permutation $\pi$ on $[s]$. We write the density of the noise distribution $\mathcal{Q}_{0}$ as
\begin{align}\label{density-Q0}
  f_{\mathcal{Q}_{0}}(y) = \frac{1}{\varT} \exp\big(-\psi\big(t/\varT \big) \big) \quad  \text{for all } y\in \real,
\end{align}
where $\varT>0$ is the scale parameter of distribution $\mathcal{Q}_{0}$ and $\psi:\real \rightarrow \real$ is a general function. Recall the definitions of $\err(\psi,2N,\lambda)$ in Eq.~\eqref{def-error-nongaussian}, $M(\psi,2N,\lambda)$ in Eq.~\eqref{def-M-nongaussian}, as well as the evaluation time $T_{\mathsf{eval}}(N)$ from Definition~\ref{assump:Teval}. We will make use of the signed kernel $\mathcal{S}_{N}^{\star}$~\eqref{truncated-signed-kernel}, with target density defined in Eq.~\eqref{density-non-gaussian-target}.
We make the following assumption on the function $\psi$. 
\begin{assumption}\label{assump:noise-tensor-pca}
The function $\psi$ satisfies Assumption~\ref{assump:log-density} with a pair of constants $(\Cdensity,\CdensityT)$. For each $\epsilon \in (0,1)$, if the pair of parameters $(N,\lambda)$ are set to
\[
    N = \Theta\Big( \CdensityB \log^{(2(1+\Cdensity))}\big(3\CdensityB / \epsilon\big) \Big) \quad \text{and} \quad \lambda = \Theta\Big(\mathsf{polylog}\big(1/\epsilon\big) \Big),
\] 
\[
   \text{then} \quad \err(\psi,2N,\lambda) \leq \epsilon/6, \quad \quad M(\psi,2N,\lambda) = \mathcal{O}(1), \quad \text{and} \quad T_{\mathsf{eval}}(N) = \mathcal{O}\big(\mathsf{poly}(1/\epsilon) \big).
\]
\end{assumption}
Note that a broad class of noise distributions satisfies Assumption~\ref{assump:noise-tensor-pca}. As illustrated in Corollary~\ref{example:generalized-normal}, Corollary~\ref{example:bounded-psi}, and Example~\ref{examples:nongaussian}, this class includes the generalized normal distribution $\mathsf{GN}(0, \varT, \beta)$, which has lighter tails than the Gaussian distribution, as well as the Student’s t-distribution (including the Cauchy distribution), which has heavier tails. It also includes any distribution whose negative log-density $\psi$ is a globally analytic function—such as the logistic distribution, the hyperbolic secant distribution, etc. We are now ready to state the result that relates the Gaussian noise model~\eqref{tensor-pca-gaussian} to the general noise model~\eqref{tensor-pca-nongaussian}.
\begin{theorem}\label{thm:tensor-pca-gaussian-nongaussian}
Let $n$ be a parameter and $s\geq 3$ be a constant. Let $\beta \in \real$ denote the signal strength and $v \in \real^{n}$ denote the spike vector. Let tensor $T(\beta,v) \in \real^{n^{\otimes s}}$ be distributed according to model~\eqref{tensor-pca-gaussian}. Let tensor $\widetilde{T}(\beta,v) \in \real^{n^{\otimes s}}$ be distributed according to model~\eqref{tensor-pca-nongaussian} with noise distribution $\mathcal{Q}_{0}$. Suppose Assumption~\ref{assump:noise-tensor-pca} holds. Then there is a reduction algorithm $\mathcal{R}:\real^{n^{\otimes s}} \rightarrow \real^{n^{\otimes s}}$ satisfying the following:
For each $\delta \in (0,1)$, if the scale of $\mathcal{Q}_{0}$ is set to $\varT^{2} = \Theta\big(\mathsf{polylog}(n^{s}/\delta) \big)$, then $\mathcal{R}$ runs in time
\[
   \mathcal{O}\Big( \mathsf{poly}\big( n^{s}/\delta \big)  \Big) \quad \text{and} \quad \sup_{\beta \in \real, v \in \real^{n}} \; \mathsf{d}_{\mathsf{TV}}\Big( \mathcal{R}\big(T(\beta,v)\big), \widetilde{T}(\beta,v)  \Big) \leq \delta.
\]
\end{theorem}
We provide the proof of Theorem~\ref{thm:tensor-pca-gaussian-nongaussian} in Section~\ref{sec:pf-thm-tensor-pca}. A few comments on Theorem~\ref{thm:tensor-pca-gaussian-nongaussian} are in order. First, the reduction guarantee does not impose any assumptions on the signal strength $\beta$ or the spike vector $v$. In particular, there is no need to assume any prior distribution on $\beta$ or $v$. Second, the reduction algorithm $\mathcal{R}$ is constructed by applying the entrywise transformation $\mathsf{K}:\mathbb{R} \rightarrow \mathbb{R}$ from Theorem~\ref{thm:non-gaussian} to each entry of the tensor $T(\beta, v)$. Specifically, recall that the distribution $\mathcal{Q}_{\theta}$ defined in Eq.~\eqref{density-non-gaussian-target} is equivalent to a location model centered at $\theta$, but with noise distribution $\mathcal{Q}_{0}$. Then, conditional on $\beta$ and $v$, for each index tuple $(i_{1}, \dots, i_{s})$, we have
\[
    T(\beta, v)_{i_{1}, \dots, i_{s}} \sim \mathcal{N}(\theta, \sigma^{2}_{i_{1}, \dots, i_{s}})
    \quad \text{and} \quad
    \widetilde{T}(\beta, v)_{i_{1}, \dots, i_{s}} \sim \mathcal{Q}_{\theta}, 
    \quad \text{where } \theta = \beta n^{s/2} v_{i_{1}} \times v_{i_{2}} \times \cdots \times v_{i_{s}} .
\]
This setup exactly matches the source and target distributions considered in Theorem~\ref{thm:non-gaussian}. 

The following hardness result immediately follows by combing the hardness of the source problem (Proposition~\ref{thm:test-hardness-pca-gaussian}) with the reduction (Theorem~\ref{thm:tensor-pca-gaussian-nongaussian}).
\begin{corollary}\label{thm:test-hardness-pca-nongaussian}
Suppose $s\geq 3$ is a constant. Consider the tensor $\widetilde{T} = \beta (\sqrt{n} v)^{\otimes n} + \widetilde{\zeta}$, where $v \sim \mathsf{Unif}(\mathbb{S}^{n-1})$ and $\widetilde{\zeta} \in \real^{n^{\otimes s}}$ is the symmetric noise tensor whose entries have density $f_{\mathcal{Q}_{0}}$ defined in Eq.~\eqref{density-Q0}. Suppose the noise distribution $\mathcal{Q}_{0}$ satisfies Assumption~\ref{assump:noise-tensor-pca} and the scale parameter satisfies $\varT^{2} = \Theta\big(\mathsf{polylog}(n)\big)$. Suppose the $k$-$\mathrm{HPC}^{s}$ conjecture holds. If $\beta = \widetilde{o}\big(n^{-s/4} \big)$ and $\beta = \omega\big(n^{-s/2} \sqrt{s\log(n)}\big)$, then there \emph{cannot} exist an algorithm $\mathcal{A}: \real^{n^{\otimes s}} \rightarrow \{0,1\}$ that runs in time $\mathsf{poly}(n)$ and satisfies
\[
    \big \langle \mathcal{A}\big( \widetilde{T} \big), v \big \rangle = \Omega(1) \quad \text{with probability at least} \quad \Omega_{n}(1).
\]
\end{corollary}
We prove Corollary~\ref{thm:test-hardness-pca-nongaussian} in Section~\ref{sec:pf-test-hardness-pca-nongaussian}. A few comments on Corollary~\ref{thm:test-hardness-pca-nongaussian} are in order. First, note that the scale parameter (which roughly corresponds to variance inflation) $\varT^{2} = \Theta\big(\mathsf{polylog}(n)\big)$ only affects the statistical-computational gap up to a polylogarithmic factor. This is because the signal-to-noise ratio in the non-Gaussian tensor PCA model~\eqref{tensor-pca-nongaussian} is given by $\beta/\varT$, and we only require $\varT$ to be \emph{polylogarithmic} in $n$. Note that the gap itself is polynomial, so this polylogarithmic change should not be thought of as substantial. Second, to the best of our knowledge, Corollary~\ref{thm:test-hardness-pca-nongaussian} is the first reduction-based hardness result for tensor PCA with non-Gaussian noise. It demonstrates a striking universality phenomenon: the same computational hardness persists for a wide class of noise distributions, including both light-tailed (e.g., generalized normal) and heavy-tailed (e.g., Student's $t$ or Cauchy) distributions. 

\begin{remark}
While~\citet{kunisky2025low} also establish a form of universality, their result applies only to a restricted class of algorithms, namely low coordinate degree algorithms. In contrast, our result applies to \textit{all} polynomial-time algorithms, a strictly broader class. Similarly,~\citet{brennan2019universality} demonstrate universality of computational lower bounds for the submatrix model; however, their techniques apply only when the matrix entries are drawn from at most two distributions. In contrast, our results hold under a continuous prior on the rank-one tensor, so that the entries of the tensor may be drawn from infinitely many possible distributions.
\end{remark}

\subsection{Symmetric mixtures of linear regressions and generalized linear models}\label{sec:hardness_glm}
In the symmetric version of the mixture of linear regressions model~\citep{quandt1978estimating,stadler2010l1}, there is some unknown vector $v \in \real^{d}$ and we observe $\mathsf{i.i.d.}$ data $\{(x_{i},y_{i})\}_{i=1}^{n}$ generated according to 
\begin{align}\label{model-MSLR}
  y_{i} = R_{i} \cdot \langle x_{i}, v \rangle + \zeta_{i}, \quad \text{where} \quad \zeta_{i} \sim \NORMAL(0,1).
\end{align}
Here, the random variable $R_{i} \in \{-1,1\}$ and the noise $\zeta_{i}$ are unobserved and independent of each other. For convenience, we denote $\eta = \|v\|_{2}$. In the setting that $v \in \real^{d}$ is $k$-sparse, the task is of estimating $v$ in $\ell_{2}$-error $o(\eta)$ exhibits a statistical-computational gap. In particular, the information-theoretic limit occurs at sample complexity $n = \widetilde{\Theta}(k\log(d)/\eta^{4})$; see~\citet[Section~E.3]{fan2018curse} for details\footnote{Although the information-theoretic lower bound in~\citet{fan2018curse} is stated for a hypothesis testing problem, it also applies to the estimation setting, since any estimator achieving $\ell_{2}$-error $o(\eta)$ can be used to solve the corresponding testing problem; see~\citet[Appendix P.2]{brennan2020reducibility} for a rigorous argument.}. On the other hand, there is a computational lower bound of $n = \widetilde{\Theta}(k^{2}/\eta^{4})$, below which no polynomial-time algorithm can succeed at the estimation task, under the $k$-$\mathrm{BPC}$ conjecture; see~\citet[Theorem~8]{brennan2020reducibility} 
for the precise statement. A refined computational lower bound---along with runtime tradeoffs---was obtained by~\citet{arpino2023statistical} for the class of low-degree polynomial algorithms. We restate the reduction-based hardness result of~\citet[Theorem~8]{brennan2020reducibility} as the following proposition.

\begin{proposition}\label{thm:hardness-MSLR}
Let $k,d,n \in \mathbb{N}$ and $\norm>0$ be polynomial in each other, $k = o(\sqrt{d})$, and $k = o(n^{1/6})$. Consider the setting $x_{i} \overset{\mathsf{i.i.d.}}{\sim} \NORMAL(0,I_{d})$, $\{R_{i}\}_{i=1}^{n} \overset{\mathsf{i.i.d.}}{\sim} \mathsf{Rad}$, $v = \norm \boldsymbol{1}_{S}/\sqrt{k}$ for some $k$-subset $S \subseteq [d]$, and $y_{i}$ follows model~\eqref{model-MSLR} for all $i \in [n]$. Suppose the $k$-$\mathrm{BPC}$ conjecture holds. If $n = \widetilde{o}\big(k^{2}/\norm^{4} \big)$, then there cannot exist an algorithm $\mathcal{A}:(\real^{d} \times \real)^{n} \rightarrow \real^{d}$ that runs in time $\mathsf{poly}(n)$ and satisfies 
\[ 
    \big\| \mathcal{A}\big( \{x_{i},y_{i}\}_{i=1}^{n} \big) - v  \big\|_{2} = o(\norm) \quad \text{with probability at least} \quad 1-o_{n}(1).
\]
\end{proposition}

\begin{figure}
  \centering
  \includegraphics[width=0.9\textwidth]{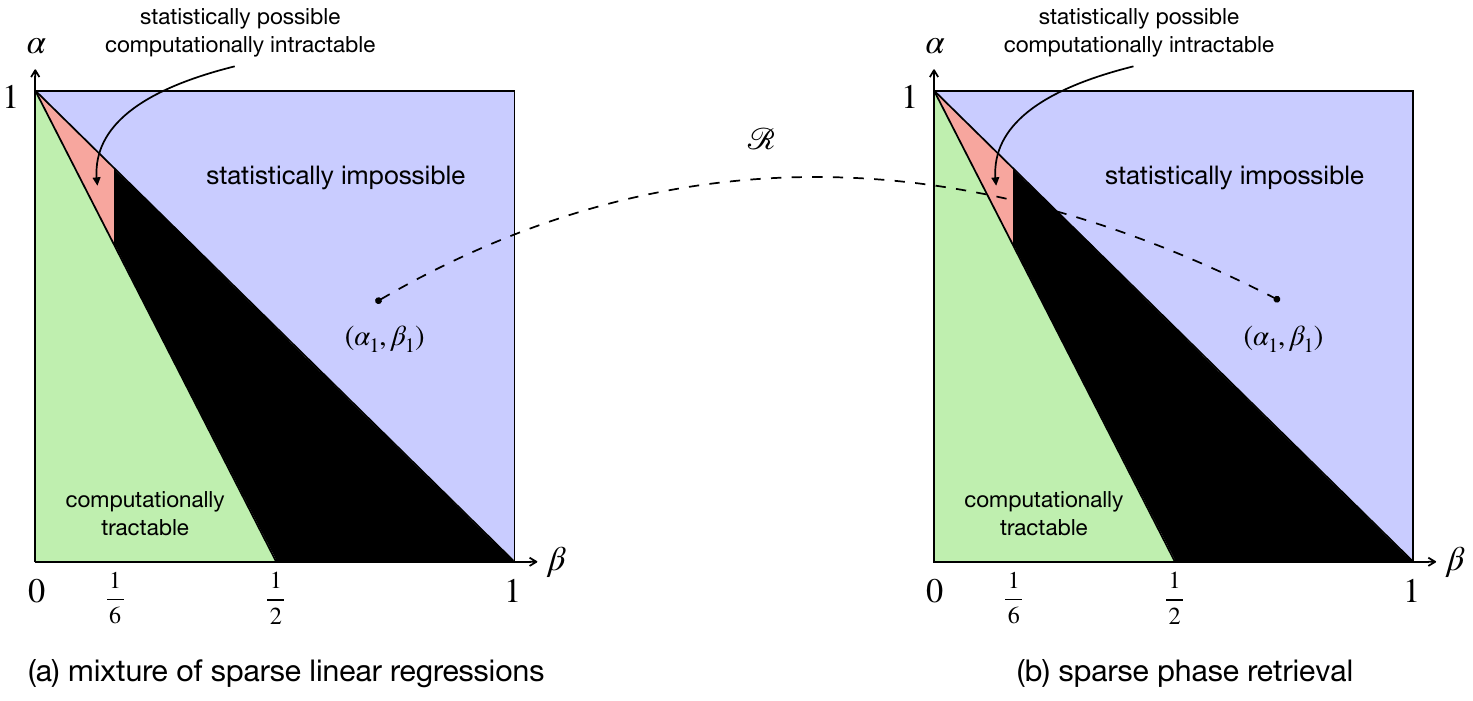}
  \caption{Phase diagrams for the mixture of sparse linear regressions ($\MixSLR$; panel (a)) and sparse phase retrieval ($\SPR$; panel (b)), under the parameterization $k=\widetilde{\Theta}(n^{\beta})$ and $\eta^{4}=\widetilde{\Theta}(n^{-\alpha})$, with constants $\alpha,\beta\in(0,1)$. 
  The purple region denotes the statistically impossible regime; the red region denotes the statistically possible but computationally intractable regime; the green region denotes the computationally tractable regime; and the black region denotes a statistically possible yet conjectured computationally hard regime, whose proof remains open. 
  The reduction $\mathcal{R}$ maps $\MixSLR$ instances with parameters $(\alpha,\beta)$ to $\SPR$ instances with the same $(\alpha,\beta)$ for all $(\alpha, \beta)$ in the phase diagram, so any computational hardness regime established in panel (a) transfers directly to panel (b).
  }
  \label{fig:SMLR-SPR-application}
\end{figure}

While Proposition~\ref{thm:hardness-MSLR} imposes the constraint $k=o(\sqrt{d})$ and $k = o(n^{1/6})$ on the sparsity level for technical reasons, it still establishes computational hardness in a nontrivial regime: By reparametrizing with
\begin{align}\label{parameterization-MSLR}
  k = \widetilde{\Theta}\big(n^{\beta}\big)  
  \quad \text{and} \quad  
  \eta^{4} = \widetilde{\Theta}\big(n^{-\alpha}\big),
\end{align}
for constants $\beta \in (0,1)$ and $\alpha > 0$, we obtain that every instance in the regime
\begin{align}\label{regime-hard-SMLR}
  \big\{(\alpha,\beta): 0 < \beta \leq 1/6, \; 1-2\beta \leq \alpha \leq 1-\beta \big\}
\end{align}
is statistically possible but remains computationally intractable for estimation; see Figure~\ref{fig:SMLR-SPR-application} for an illustration. 

At a heuristic level, the $k^{2}$ computational lower bound arises because in the model~\eqref{model-MSLR}, 
the Rademacher random variable $R_{i}$ destroys information about the sign of the regression term $\langle x_{i}, v \rangle$. Contrast this with the case of
sparse linear regression, where no such statistical-computational gap appears. This observation naturally raises the following question: Does the $\Omega(k^{2})$ computational lower bound arise for any signless observation model in variants of sparse linear regression?

Toward answering this question, we consider the class of generalized linear models~\citep{barbier2019optimal} with even link functions. Specifically, suppose there is an unknown vector $v \in \real^{d}$ and we observe $\mathsf{i.i.d.}$ data $\{(x_{i},\widetilde{y}_{i})\}_{i=1}^{n}$ generated according to 
\begin{align}\label{model-general-linear}
  \widetilde{y}_{i} = f\big( \langle x_{i}, v \rangle / \varT \big) + \widetilde{\zeta}_{i}, \quad \text{where} \quad \widetilde{\zeta}_{i} \sim \NORMAL(0,1).
\end{align}
In the above, $f:\real \rightarrow \real$ is a general even link function and the parameter $\varT>0$ is a proxy for the signal-to-noise ratio, which will be specified in concrete applications. Also suppose that the signal $v$ is $k$-sparse. 
We now use our results from Section~\ref{sec:concrete_reductions} to extend reduction-based hardness results for mixture of sparse linear regressions to a class of generalized sparse linear models.

We make the following assumption on the link function $f$. Recall the definitions of $\Llink(2N)$ and $\LlinkT(2N,\varT)$ in Eq.~\eqref{def-quantity-link} as well as the evaluation time $T_{\mathsf{eval}}(N)$ from Definition~\ref{assump:Teval}. We will make use of the signed kernel $\mathcal{S}_{N}^{\star}$~\eqref{truncated-signed-kernel}, with target density defined in Eq.~\eqref{target:density}.

\begin{assumption}\label{assump:link-linear-model}
 The link function $f:\real \rightarrow \real$ is even, i.e., $f(x) = f(-x)$ for all $x\in \real$, and it satisfies Assumption~\ref{assump:link} with a pair of constants $(\Clink,\ClinkT)$. For each $\epsilon \in (0,1)$, if $N = \Theta\Big( \ClinkB \log^{2(1+\Clink)}\big(4\ClinkB/\epsilon \big) \Big)$ where $\ClinkB$ is defined in Eq.~\eqref{def-quantity-link}, then we have for some $\varT = \mathcal{O}\big(\mathsf{polylog}(1/\epsilon)\big)$,
\[
    \Llink(2N) = \mathcal{O}\Big( \mathsf{polylog}(1/\epsilon)  \Big), \quad \LlinkT(2N,\varT) = \mathcal{O}\Big( \mathsf{polylog}(1/\epsilon)  \Big), \quad \text{and} \quad T_{\mathsf{eval}}(N) = \mathcal{O}\Big( \mathsf{poly}(1/\epsilon)  \Big).
\] 
\end{assumption}
Note that a broad class of functions satisfies Assumption~\ref{assump:link-linear-model}. As illustrated in Corollary~\ref{example:monomial}, Corollary~\ref{example:analytical}, and Example~\ref{examples}, this includes monomials of constant degree---for instance, the function $f(\theta) = \theta^{2}$ is covered, which corresponds to the phase retrieval model in the context of generalized linear models. Moreover, Assumption~\ref{assump:link-linear-model} also encompasses the class of globally analytic functions (see Definition~\ref{def:analytical-functions}), which includes many canonical and practically relevant link functions such as $f(\theta) = \cos(\theta)$, $f(\theta) = \mathsf{sech}(\theta)$, and $f(\theta) = (1+\theta^{2})^{-1}$, among others. Additionally, by Lemma~\ref{lemma:linear-combination}, any linear combination of these basic functions is also covered. We are now ready to relate mixtures of linear regressions~\eqref{model-MSLR} to the generalized linear model~\eqref{model-general-linear} with a link function satisfying Assumption~\ref{assump:link-linear-model}.
\begin{theorem}\label{thm:reduction-mslr-glm}
Let $v \in \real^{d}$ denote the unknown vector and let $\{R_{i}\}_{i=1}^{n} \subseteq \{-1,1\}^{n}$ be any vector of signs. Conditioned on $\{ R_{i}\}_{i=1}^{n}$ and $v$,
suppose samples $\{(x_{i},y_{i})\}_{i=1}^{n}$ are drawn $\mathsf{i.i.d.}$ according to model~\eqref{model-MSLR} and samples $\{(x_{i},\widetilde{y}_{i})\}_{i=1}^{n}$ are drawn $\mathsf{i.i.d.}$ according to model~\eqref{model-general-linear}. 
Suppose the link function $f:\real \rightarrow \real$ satisfies Assumption~\ref{assump:link-linear-model}. Then there is a reduction algorithm $\mathcal{R}: (\real^{d} \times \real)^{n} \rightarrow (\real^{d} \times \real)^{n}$ satisfying the following:
For each $\delta \in (0,1)$, if we set 
\begin{align}\label{tau-condition-prop-MSLR-GLM}
    \varT = \mathsf{polylog}(n/\delta) \quad \text{and large enough so that} \quad \tau \geq  2 \max_{ i \in [n]}\; \big|\langle x_{i},v \rangle \big|,
\end{align}
then the reduction $\mathcal{R}$ runs in time
\[ 
    \mathcal{O} \Big( \mathsf{poly}(n/\delta) \Big) \quad \text{and satisfies} \quad \mathsf{d}_{\mathsf{TV}} \Big(  \mathcal{R} \big( \{(x_{i},y_{i})\}_{i=1}^{n}  \big), \{(x_{i},\widetilde{y}_{i})\}_{i=1}^{n}  \Big)  \leq \delta.
\]
\end{theorem}
We provide the proof in Section~\ref{sec:pf_reduction_mslr_glm}. Note that the reduction guarantee does not impose any assumptions on the sign variables $\{R_{i}\}_{i=1}^{n}$ and the covariates $\{x_{i}\}_{i=1}^{n}$. It also does not place any structural assumptions on $v$ such as sparsity. In particular, there is no need to assume any prior distribution on $\{R_{i}\}_{i=1}^{n}$, $\{x_{i}\}_{i=1}^{n}$, and $v$. This is because the reduction algorithm $\mathcal{R}$ is constructed by applying the entrywise transformation $\mathsf{K}:\mathbb{R} \rightarrow \mathbb{R}$ from Theorem~\ref{thm:gaussian-mean-nonlinear} to each sample $y_{i}$. Specifically, conditional on the tuple $(R_{i}, x_{i}, v)$, we have
\[
    y_{i} \sim \NORMAL(\theta,1) \quad \text{and} \quad \widetilde{y}_{i} \sim \NORMAL\big(f(\theta/\varT),1\big), \quad \text{where } \theta = R_{i} \cdot \langle x_{i}, v \rangle.
\]
This setup exactly matches the source and target distributions considered in Theorem~\ref{thm:gaussian-mean-nonlinear}.

\begin{remark}\label{remark:MSLR-GLM}
Note that any instance of the mixture of sparse linear regressions ($\MixSLR$)~\eqref{model-MSLR} with signal $v$ is mapped under the reduction $\mathcal{R}$ to an instance of the sparse generalized linear model ($\SpGLM$)~\eqref{model-general-linear} with signal $v/\varT$. In most applications of Theorem~\ref{thm:reduction-mslr-glm}, it suffices to take $\delta$ to be an inverse polynomial of $n$ and in the Gaussian covariance setting $x_{i} \overset{\mathsf{i.i.d.}}{\sim} \NORMAL(0,I_{d})$, standard Gaussian tail bounds and a union bound yield
\[
    \mathbb{P}\Big\{ \max_{i \in [n]} |\langle x_{i}, v \rangle| \leq 2 \sqrt{\log(n)}\Big\} \geq 1 - n^{-1} \quad \text{for} \quad \|v\|_{2} \leq 1.
\]
Thus, condition~\eqref{tau-condition-prop-MSLR-GLM} is satisfied by letting $\varT = \Theta(\mathsf{polylog}(n))$. Therefore, the reduction $\mathcal{R}$ preserves the sparsity $k$ and only changes the signal strength $\eta$ by a polylogarithmic factor. 
\end{remark}

The hardness result of $\MixSLR$ (Proposition~\ref{thm:hardness-MSLR}) and the reduction (Theorem~\ref{thm:reduction-mslr-glm}) imply the following computational hardness of $\SpGLM$.
\begin{corollary}\label{thm:hardness-glm}
Let $k,d,n \in \mathbb{N}$ and $\norm>0$ be polynomial in each other, $k = o(\sqrt{d})$, and $k = o(n^{1/6})$. Consider the setting $x_{i} \overset{\mathsf{i.i.d.}}{\sim} \NORMAL(0,I_{d})$, $v = \norm \boldsymbol{1}_{S}/\sqrt{k}$ for some $k$-subset $S \subseteq [d]$, and $\widetilde{y}_{i}$ follows model~\eqref{model-general-linear} for all $i \in [n]$. Suppose in model~\eqref{model-general-linear}, the link function $f$ satisfies Assumption~\ref{assump:link-linear-model} and the scale parameter satisfies $\varT = \Theta\big(\mathsf{polylog}(n) \big)$. Suppose the $k$-$\mathrm{BPC}$ conjecture holds. If $n = \widetilde{o}\big(k^{2}/\norm^{4} \big)$, then there cannot exist an algorithm $\mathcal{A}:(\real^{d} \times \real)^{n} \rightarrow \real^{d}$ that runs in time $\mathsf{poly}(n)$ and satisfies 
\[ 
    \big\| \mathcal{A}\big( \{x_{i},\widetilde{y}_{i} \} \big) - v \big\|_{2} = o(\norm) \quad \text{with probability at least} \quad 1 - o_{n}(1).
\]
\end{corollary}
We provide the proof of Corollary~\ref{thm:hardness-glm} in Section~\ref{sec:pf_thm_hardness_glm}, and record a few remarks here. 
First, the constraints $k=o(\sqrt{d})$ and $k = o(n^{1/6})$ are inherited from Proposition~\ref{thm:hardness-MSLR}, which is the source of the computational hardness. Nevertheless, this still establishes computational hardness in the nontrivial regime defined in Eq.~\eqref{regime-hard-SMLR} under the parameterization~\eqref{parameterization-MSLR}; see the illustration of this computationally intractable regime in Figure~\ref{fig:SMLR-SPR-application}. Second, our reduction itself imposes no constraints on $k$, so any hardness shown for $\MixSLR$ automatically transfers to $\SpGLM$. In particular, to understand the computational complexity of $\SpGLM$ in the black region of Figure~\ref{fig:SMLR-SPR-application}(b), it suffices to understand the corresponding black region of Figure~\ref{fig:SMLR-SPR-application}(a).
Third, the scaling factor $\varT$ does not affect the signal strength beyond a polylogarithmic factor, since in model~\eqref{model-general-linear} the effective signal strength is $\|v\|_{2}/\varT = \norm/\varT$, and we set $\varT = \Theta(\mathsf{polylog}(n))$. 
We have thus demonstrated that the $\widetilde{\Omega}(k^{2})$ computational barrier is universal for models with signless observations arising from even link functions. This resolves an open question posed in~\citet[Appendix D]{brennan2020reducibility}. We remark that a $k$-to-$k^{2}$ gap has been established under the statistical query model~\citep{wang2019statistical} for a class of even link functions, but this only applies to a subset of algorithms and link functions.

We now highlight the computational hardness result for sparse phase retrieval ($\SPR$), where we have
the model
\begin{align}\label{model-phase-retrieval}
      \widetilde{y}_{i} = |\langle x_{i},v\rangle|^{2} + \widetilde{\zeta}_{i}, \quad \text{where} \quad \widetilde{\zeta}_{i} \sim \NORMAL(0,1).
\end{align}
The following computational hardness result immediately follows from Corollary~\ref{thm:hardness-glm}. 

\begin{corollary}[Sparse Phase Retrieval]\label{corollary:spr}
Let $k,d,n \in \mathbb{N}$ and $\norm>0$ be polynomial in each other, $k = o(\sqrt{d})$, and $k = o(n^{1/6})$. Consider the setting $x_{i} \overset{\mathsf{i.i.d.}}{\sim} \NORMAL(0,I_{d})$, $v = \norm \boldsymbol{1}_{S}/\sqrt{k}$ for some $k$-subset $S \subseteq [d]$, and $\widetilde{y}_{i}$ follows the phase retrieval model~\eqref{model-phase-retrieval} for all $i \in [n]$. Suppose the $k$-$\mathrm{BPC}$ conjecture holds. 
If $n = \widetilde{o}\big(k^{2}/\norm^{4} \big)$, then there cannot exist an algorithm $\mathcal{A}:(\real^{d} \times \real)^{n} \rightarrow \real^{d}$ that runs in time $\mathsf{poly}(n)$ and satisfies 
\[ 
    \big\| \mathcal{A}\big( \{x_{i},\widetilde{y}_{i} \} \big) - v \big\|_{2} = o(\norm) \quad \text{with probability at least} \quad 1 - o_{n}(1).
\]
\end{corollary}
Note that the square link function $f(t) = t^{2}$ satisfies Assumption~\ref{assump:link-linear-model} (see Corollary~\ref{example:monomial}). 
Hence, Corollary~\ref{corollary:spr} follows directly from Corollary~\ref{thm:hardness-glm} by setting $f(\theta) = \theta^{2}$ and absorbing $\varT$ into $\eta$. 
The information-theoretic limit for $k$-sparse phase retrieval is known to occur at sample complexity $n = \Theta\big(k \log(d)/\eta^{4}\big)$; see~\citet{cai2016optimal,lecue2015minimax}. In contrast, all polynomial-time algorithms require sample complexity in the order of $k^{2}$; see~\citet{li2013sparse,cai2016optimal,wang2017sparse,celentano2020estimation}. Motiveted by this gap,~\citet{cai2016optimal} conjectured that any computationally efficient method must require $n = \widetilde{\Omega}(k^{2}/\eta^{4})$ samples. Corollary~\ref{thm:reduction-mslr-glm} provides the first reduction-based confirmation of this conjecture by establishing a computational lower bound: no polynomial-time algorithm can succeed when $n = \widetilde{o}(k^{2}/\eta^{4})$. This demonstrates a statistical–computational gap between the information-theoretic rate $k$ and the computational barrier $k^{2}$. In particular, it establishes computational hardness in a nontrivial regime; see Figure~\ref{fig:SMLR-SPR-application}(b) for an illustration.

\subsection{Sparse rank-1 submatrix model and planted submatrix model} \label{sec:sparse_spiked_wigner_BCering}
In the sparse rank-1 submatrix model ($\SubMat$), we observe an $n\times n$ matrix
\begin{align} \label{eq:rank-1-model}
  M = \mu\, r c^{\top} + \zeta, \quad \text{where } r,c \in S_{n,k},\; \zeta \sim \mathcal{N}(0,1)^{\otimes n \times n}.
\end{align}
Here, $\mu \in \mathbb{R}$ is the signal strength, $\NORMAL(0,1)^{\otimes n \times n}$ denotes the distribution of an $n \times n$ matrix with $\mathsf{i.i.d.}$ $\mathcal{N}(0,1)$ entries, and $S_{n,k}$ denotes set of $k$-sparse unit vectors defined as
\begin{align}\label{set-sparse-vector}
  S_{n,k} := \left\{ x \in \mathbb{S}^{n-1}:\; \|x\|_{0} = k,\ x_{i} \in \left\{ 0, \pm \tfrac{1}{\sqrt{k}} \right\} \text{ for all } i \in [n] \right\}.
\end{align}
The corresponding detection problem can be formalized as distinguishing between the following two hypotheses upon observing the matrix $M$: 
\begin{subequations} \label{spike-wigner-model}
\begin{align}
  &\mathcal{H}_{0}: M \text{ follows model~\eqref{eq:rank-1-model} with } \mu = 0, \quad \text{and} \\ 
  &\mathcal{H}_{1}: M \text{ follows model~\eqref{eq:rank-1-model} for some } \mu > 0, r,c \in S_{n,k}.
\end{align}
\end{subequations}
The detection problem in Eq.~\eqref{spike-wigner-model} exhibits a computational barrier, which we briefly describe. Define the computational threshold
\begin{align}\label{threshold-spiWig}
  \mu_{\mathsf{Rank1}}^{\text{comp}} := \min \big\{k, \sqrt{n} \big\}.
\end{align}
When the signal strength satisfies $\mu = \widetilde{\omega}\big(\mu_{\mathsf{Rank1}}^{\text{comp}}\big)$, polynomial-time algorithms exist for detection, and these also succeed for recovering $r, c$ under $H_1$. Specifically, if $\mu = \widetilde{\omega}(k)$, then each entry of the signal matrix $\mu\, r c^{\top}$ dominates the noise $\zeta$. If $\mu = \widetilde{\omega}(\sqrt{n})$, then spectral methods succeed~\citep{montanari2015limitation}, since the spectral norm of the noise matrix is only $\mathcal{O}(\sqrt{n})$ and the top singular vectors of $M$ therefore reveal information about the signal $\mu\, r c^{\top}$.
Conversely, it has been shown via reductions from the planted clique conjecture that no polynomial-time algorithm can succeed when $\mu = \widetilde{o}\big(\mu_{\mathsf{Rank1}}^{\text{comp}}\big)$; see~\citet{brennan2018reducibility}.

\begin{figure}
  \centering
  \includegraphics[width=0.9\textwidth]{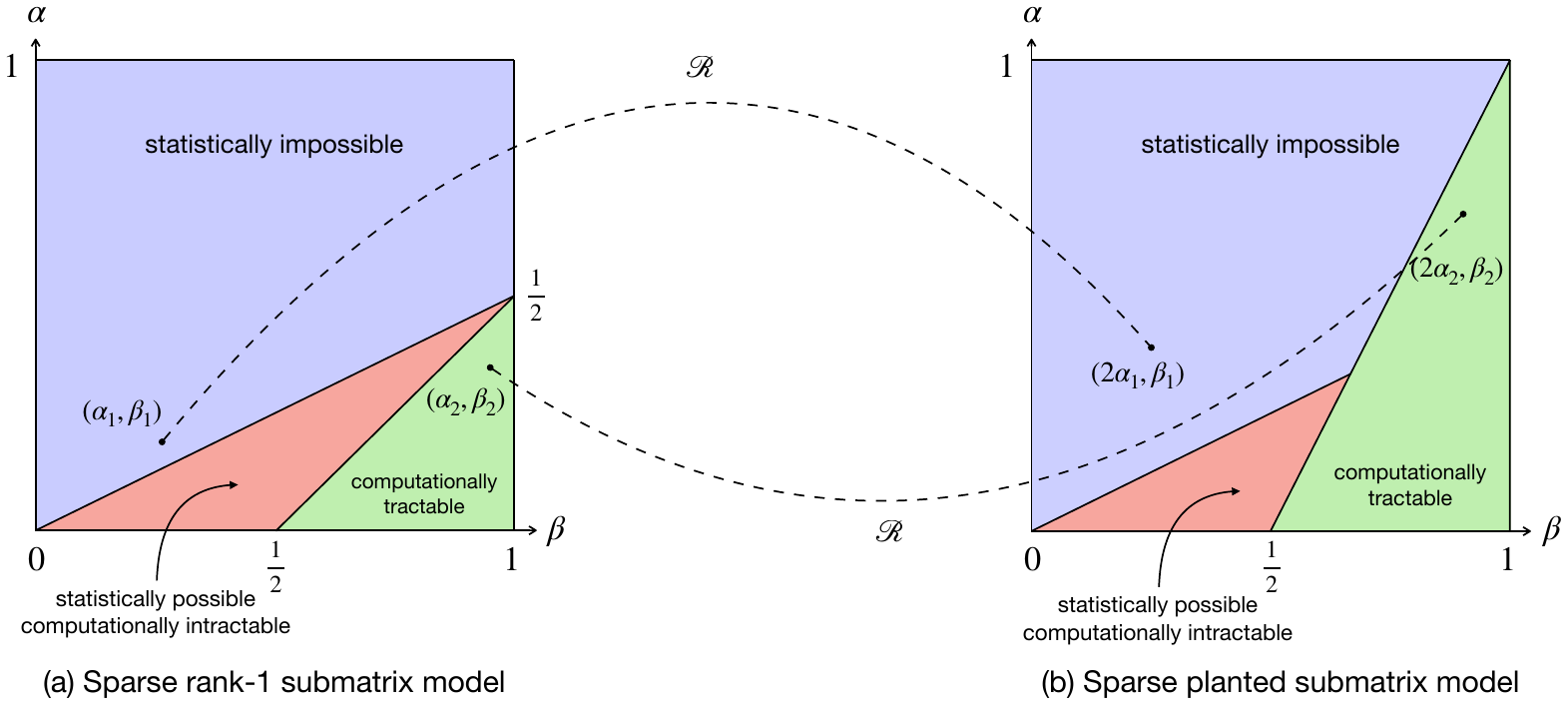}
  \caption{Phase diagrams for the sparse rank-1 submatrix model ($\SubMat$; panel (a)) and the planted submatrix model ($\PlantSubMat$; panel (b)), under the parameterization $\mu/k = \widetilde{\Theta}(n^{-\alpha})$ and $k = \widetilde{\Theta}(n^{\beta})$ for $\alpha, \beta \in (0,1)$. The purple region indicates the statistically impossible regime. The green region corresponds to the computationally tractable regime where polynomial-time algorithms are known to succeed. The red region indicates the statistically possible but computationally intractable regime. In panel (a), the computational threshold is given by the line $\alpha = \beta - 1/2$, while in panel (b), the threshold is the line $\alpha = 2\beta - 1$. The reduction $\mathcal{R}$ established in Theorem~\ref{thm:reduction-ROS-BC} maps instances with parameters $(\alpha, \beta)$ in $\SubMat$ to instances with parameters $(2\alpha, \beta)$ in $\PlantSubMat$. This transformation demonstrates a pointwise correspondence between the phase diagrams, preserving the computational tractability and intractability of problem instances.
  }
  \label{fig:phase-diagram}
\end{figure}

In the planted submatrix model $(\PlantSubMat)$, we observe an $n \times n$ matrix
\begin{align} \label{eq:plantedmat}
  M = \mu\, r c^{\top} + \zeta, 
  \quad \text{where } r,c \in S_{n,k}^{+},\; \zeta \sim \NORMAL(0,1)^{\otimes n \times n}.
\end{align}
Here, the spike vectors $r$ and $c$ are restricted to have nonnegative entries. Specifically, $S_{n,k}^{+}$ denotes the set of $k$-sparse unit vectors with nonnegative coordinates:
\begin{align}\label{set-sparse-vector-plus}
  S_{n,k}^{+} 
  \;:=\; \left\{ x \in \mathbb{S}^{n-1} : \|x\|_{0} = k,\; 
  x_{i} \in \left\{ 0, \tfrac{1}{\sqrt{k}} \right\} \ \forall i \in [n] \right\}.
\end{align}
The corresponding detection problem involves distinguishing between the following hypotheses upon observing $M$:
\begin{subequations} \label{BCering-model}
\begin{align}
  &\mathcal{H}_{0}: M \text{ follows model~\eqref{eq:plantedmat} with } \mu = 0, \quad \text{and} \\ 
  &\mathcal{H}_{1}: M \text{ follows model~\eqref{eq:plantedmat} for some } \mu > 0, r,c \in S_{n,k}^+.
\end{align}
\end{subequations}

Compared to the sparse rank-1 submatrix model, the planted submatrix model exhibits a different computational threshold, defined as
\begin{align}\label{threshold-BC}
  \mu_{\mathsf{Plant}}^{\text{comp}} := \min\left\{k, \tfrac{n}{k} \right\}.
\end{align}
When $\mu = \widetilde{\omega}\big(\mu_{\mathsf{Plant}}^{\text{comp}} \big)$, polynomial-time algorithms exist. For instance, when $\mu = \widetilde{\omega}(k)$, each entry of the signal matrix $\mu\, r c^{\top}$ again dominates the noise. Alternatively, when $\mu = \widetilde{\omega}(n/k)$, summing all entries of $\mu\, r c^{\top}$ yields a signal $\mu k = \widetilde{\omega}(n)$, while the total sum of entries in the noise matrix is only $\widetilde{O}(n)$, making the two hypotheses distinguishable. Conversely, when $\mu = \widetilde{o}\big(\mu_{\mathsf{Plant}}^{\text{comp}} \big)$, then no polynomial-time algorithm is known, and evidence of computational hardness has been provided through reductions from the planted clique conjecture~\citep{ma2015computational,brennan2018reducibility}.

As observed, the computational threshold $\mu_{\mathsf{Plant}}^{\mathrm{comp}}$ is significantly smaller than $\mu_{\mathsf{Rank1}}^{\mathrm{comp}}$ in the regime $k \geq \sqrt{n}$, indicating that $\PlantSubMat$ is computationally easier than $\SubMat$. To rigorously establish this relation, 
we construct a reduction from $\SubMat$ to $\PlantSubMat$, thereby showing that the computational threshold of $\SubMat$ directly implies the corresponding threshold of $\PlantSubMat$, without relying on the planted clique conjecture. We now formally state the result.

\begin{theorem}\label{thm:reduction-ROS-BC}
  Let $S_{n,k} \subseteq \real^{n}$ be defined as in Eq.~\eqref{set-sparse-vector}. There exists a universal and positive constant $C$ and a reduction algorithm $\mathcal{R}:\real^{n\times n} \rightarrow \real^{n\times n}$ such that the following holds: For each $\delta \in (0,1)$, if the tuple of parameters $(\mu,k,n,\varT)$ satisfies
  \begin{align}\label{condition-tau-ROS-BC}
      \varT^{2} \geq C\log^{24}(n^{2}/\delta) + \frac{4\mu^{2}}{k^{2}},
  \end{align}
  then the reduction $\mathcal{R}$ runs in time $\mathcal{O} \Big( n^{2}\log^{64}\big(Cn^{2}/ \delta \big) \Big)$ and
  \begin{align*}  
       \sup_{\mu \in \real,\; r,c \in S_{n,k}} \mathsf{d}_{\mathsf{TV}} \bigg( \mathcal{R}\Big(\mu rc^{\top} +  \zeta \Big) ,\;  \frac{\mu^{2}}{\varT^{2}k} |r| |c|^{\top} + \zeta \bigg) \leq \delta,
  \end{align*}
  where $\zeta \sim \NORMAL(0,1)^{\otimes n \times n}$ and $|r|,|c|\in \real^{n}$ denote the vectors obtained by taking the entrywise absolute value of $r, c$, respectively. 
\end{theorem}

We provide the proof of Theorem~\ref{thm:reduction-ROS-BC} in Section~\ref{sec:pf_thm_reduction_ROS_BC}. Note that it is typical to set $\delta = \mathcal{O}(n^{-C'})$ for some constant $C' > 0$ to achieve meaningful reductions, so that condition~\eqref{condition-tau-ROS-BC} is satisfied by choosing $\varT = \Theta(\mathrm{polylog}(n))$. 

Let us now discuss some consequences of the theorem, focusing on the regime where $\mu/k = \widetilde{\Theta}(n^{-\alpha})$ for some $\alpha \in (0,1)$. The reduction maps an instance of $\SubMat$ with parameters $(\mu, k, n)$ to an instance of $\PlantSubMat$ with parameters $\big(\mu^2/(\varT^2 k),\, k,\, n\big)$. Since $\varT = \Theta(\mathrm{polylog}(n))$, its effect on the signal strength can be neglected. This ensures a pointwise correspondence between the phase diagrams of the two models: if the $\SubMat$ instance with parameters $(\mu,k,n)$ is computationally tractable, then its image under the reduction, an instance of $\PlantSubMat$ with parameters $(\mu^{2}/(\varT^{2}k), k, n)$, is also tractable, and the same holds for intractability.  
This correspondence can be formalized by parameterizing $\mu/k = \widetilde{\Theta}(n^{-\alpha})$ and $k = \widetilde{\Theta}(n^{\beta})$ for $\alpha, \beta \in (0,1)$. Under this parameterization, a point $(\alpha, \beta)$ in the $\SubMat$ phase diagram is mapped to $(2\alpha, \beta)$ in the $\PlantSubMat$ phase diagram, confirming that the reduction preserves computational feasibility pointwise; see Figure~\ref{fig:phase-diagram} for an illustration.  

Finally, we emphasize that Theorem~\ref{thm:reduction-ROS-BC} establishes computational lower bounds without relying on the planted clique conjecture. Indeed, even if the planted clique conjecture were false, a computational lower bound for $\SubMat$ would still imply the corresponding lower bound for $\PlantSubMat$.  
Moreover, reductions from $\SubMat$ provide evidence of hardness beyond polynomial time. Since existing hardness results for $\PlantSubMat$ are conditioned on the planted clique conjecture~\citep{brennan2018reducibility,ma2015computational}, they rule out only polynomial-time algorithms. However, these reductions lose strength once we go beyond polynomial time, since the planted clique problem is known to be solvable in quasi-polynomial time, specifically $\mathcal{O}\!\left(n^{\mathcal{O}(\log n)}\right)$, whenever the clique size exceeds the information-theoretic threshold. In contrast, reductions from $\SubMat$ remain meaningful in this regime, since $\SubMat$ is believed to require subexponential time in the hard regime. More concretely, if we fix $\mu = \widetilde{\Theta}(\sqrt{n})$ and let the sparsity ratio be $\rho = k/n$, then in the statistically feasible but computationally intractable regime $n^{-1/2} \ll \rho \ll 1$,~\citet{ding2024subexponential} showed that recovery 
 of the spike requires time $\mathcal{O}\!\left(\exp(\rho^{2} n)\right)$, and evidence for the optimality of this recovery algorithm was also provided. Consequently, hardness based on $\SubMat$ can exclude a strictly larger class of algorithms than hardness based on planted clique.

\section{Discussion}\label{sec:discussion}

Motivated by the goal of characterizing fundamental statistical-computational tradeoffs in canonical high-dimensional models, we constructed a Markov kernel that maps the Gaussian location model (the source) to general target models within small total variation deficiency. The target models we covered include (i) general non-Gaussian location models and (ii) Gaussian models with a nonlinearly transformed location parameter. Building on this kernel, we designed a computationally efficient reduction that transforms a single source observation into a single target observation without knowledge of the underlying location parameter. This single-sample reduction allows us to establish new computational lower bounds for several high-dimensional problems including tensor PCA and sparse generalized linear models, as well as novel connections between the sparse rank-one submatrix model and the planted submatrix model.

Our work leaves several open questions, which we briefly outline here. First, our construction of the optimal signed kernel relies on the fact that, under the Gaussian measure, the Hermite polynomials form an orthonormal basis for a suitable function space. Many other distributions also admit orthonormal polynomial bases—for instance, the Laguerre polynomials with respect to the Gamma distribution, and the Gegenbauer polynomials with respect to the Beta distribution; see~\citet{goldstein2005distributional} for more examples. Can our design of the optimal signed kernel be generalized to these source distributions as well? Furthermore, for Laplace, Uniform, and Erlang sources, the optimal signed kernels have been derived using different techniques, such as Fourier inversion; see~\citet{lou2025computationally}. Can we unify these approaches under a single framework based on orthonormal polynomials? 

Second, our technique relies critically on the differentiability of the target density with respect to the location parameter. For instance, when the target is the Laplace location model, the reduction fails due to the density being non-differentiable at a single point. Nevertheless, Gaussian-to-Laplace reductions are of particular interest, for example in differential privacy. Can the optimal Markov kernel defined in Eq.~\eqref{optimal-markov-kernel} overcome this limitation, or does a lower bound—perhaps based on the data processing inequality—rule out such reductions entirely? Furthermore, our reduction requires the target density to satisfy smoothness assumptions such as Assumption~\ref{assump:log-density} and Assumption~\ref{assump:link}. While these are met by a broad class of distributions, an important open question is whether these conditions can be relaxed or even removed, or conversely, whether one can establish lower bounds showing that they are indeed necessary.

Finally, we believe the Gaussian-source reduction framework has applications well beyond the settings explored in this paper. For instance, we chose to focus our discussion of universality of computational hardness on tensor PCA. However, the reduction in Theorem~\ref{thm:non-gaussian} is quite general, and could be adapted to regression and denoising models where the observations are corrupted by Gaussian noise. By transferring these to analogous models with non-Gaussian noise satisfying Assumption~\ref{assump:noise-tensor-pca}, one could potentially extend universality results across a broad class of problems, complementing the directions studied in~\citet{brennan2018reducibility,brennan2020reducibility}.

\subsection*{Acknowledgments}
ML was supported in part by the ARC-ACO fellowship through the Algorithms and Randomness Center at Georgia Tech. ML and AP were supported in part by the NSF under grants CCF-2107455 and DMS-2210734, a Google Research Scholar award, and by research awards/gifts from Adobe, Amazon, and Mathworks. 
GB gratefully acknowledges support from NSF award CCF-2428619.


\section{Proofs of results from Section~\ref{sec:general_markov_kernel}}\label{sub:pf_general_markov_kernel}
We provide proofs of Proposition~\ref{prop:optimal-signed-kernel} and Lemma~\ref{lemma:TV-deficiency-decompose} below. We begin by reviewing some basic properties of the Hermite polynomials defined in Eqs.~\eqref{def:hermite} and~\eqref{def:hermite-normalized}, which will be used frequently in our proofs. 

First, note that the normalized Hermite polynomials are orthonormal in that $\big \langle \Nherm_{n}(\cdot), \Nherm_{m}(\cdot)  \big\rangle = \delta_{nm}$, where $\delta_{nm}$ is the Kronecker delta with $\delta_{nm} = \ind{n = m}$. The other properties we use are:
\begin{subequations}\label{Hermite-Eq-Ineqs}
\begin{align}
  \label{Hermite-boundness}
  &\big| H_{n}(x) \exp(-x^{2}/2) \big| \leq \big( \sqrt{\pi} 2^{n}n! \big)^{1/2}, \qquad \qquad \big| \Nherm_{n}(x) \exp(-x^{2}/2) \big| \leq 1, \\
  \label{Hermite-derivative}
  &\frac{\mathrm{d} H_{n+1}(x)}{ \mathrm{d} x} = 2(n+1)H_{n}(x), \qquad \qquad \frac{\mathrm{d} \Nherm_{n+1}(x)}{ \mathrm{d} x} = \sqrt{2(n+1)} \cdot \Nherm_{n}(x),\\
  \label{Hermite-closed-form}
  &\text{and} \qquad \qquad H_{n}(x) = \sum_{m=0}^{\left\lfloor \frac{n}{2} \right\rfloor } \frac{n!}{m! (n-2m)!} \cdot (-1)^{m} (2x)^{n-2m}.
\end{align}
\end{subequations} 
Inequality~\eqref{Hermite-boundness} can be found in~\citet[page~190]{szeg1939orthogonal}, and Eqs.~\eqref{Hermite-derivative} and~\eqref{Hermite-closed-form} correspond to~\citet[Eq.~(5.5.10) and Eq.~(5.5.4)]{szeg1939orthogonal}, respectively.

\subsection{Proof of Proposition~\ref{prop:optimal-signed-kernel}} \label{sec:pf-prop1}
Recall the signed $\mathcal{S}^{\star}$ from Eq.~\eqref{optimal-signed-kernel} and the source density $u(x;\theta)$ from Eq.~\eqref{source-gaussian-density}. It suffices to verify
\begin{align} \label{eq:deconv}
	\int_{\real} \mathcal{S}^{\star}(y \mid x) \cdot u(x;\theta) \mathrm{d}x = v(y;\theta) \quad \text{for all} \quad y,\theta \in \real.
\end{align}
To reduce the notational burden, we let $c_{k} = (-1)^{k} \sigma^{2k}/(2k)!!$ for $k \in \PosInt$ and $\phi(x) = e^{-x^2}/\sqrt{\pi}$ for $x\in\real$. In our proof, we separately evaluate the LHS and RHS of Eq.~\eqref{eq:deconv} in terms of Hermite projections, and then show that these are equal.

\paragraph{LHS of Eq.~\eqref{eq:deconv}.}
By applying integration by parts and using Assumption~\ref{assump_1:target}(b), we obtain
\begin{align*}
  \int_{\real} \nabla_{x}^{(2k)} v(y;x)  \cdot \phi\bigg(\frac{x-\theta }{\sqrt{2}\sigma}\bigg)  \mathrm{d}x &= \nabla_{x}^{(2k-1)} v(y;x) \cdot \phi\bigg(\frac{x-\theta }{\sqrt{2}\sigma}\bigg) \bigg \vert_{-\infty}^{+\infty} - \int_{\real} \nabla_{x}^{(2k-1)} v(y;x)  \cdot \phi'\bigg(\frac{x-\theta }{\sqrt{2}\sigma}\bigg) \cdot \frac{1}{\sqrt{2}\sigma}  \mathrm{d}x \nonumber \\
  & = -\frac{1}{\sqrt{2}\sigma} \int_{\real} \nabla_{x}^{(2k-1)} v(y;x)  \cdot \phi'\bigg(\frac{x-\theta }{\sqrt{2}\sigma}\bigg)  \mathrm{d}x \nonumber \\
  & \overset{\1}{=} \Big( \frac{1}{\sqrt{2}\sigma} \Big)^{2k} \cdot \int_{\real}  v(y;x)  \cdot \phi^{(2k)}\bigg(\frac{x-\theta }{\sqrt{2}\sigma}\bigg) \mathrm{d}x \\
  & \overset{\2}{=} \frac{\sqrt{2}\sigma}{2^{k} \sigma^{2k}} \int_{\real} v\big(y;\theta+\sqrt{2}\sigma x\big) \cdot \phi^{(2k)} (x) \mathrm{d}x \\
  & \overset{\3}{=} \frac{\sqrt{2}\sigma}{2^{k} \sigma^{2k}} \int_{\real} v\big(y;\theta+\sqrt{2}\sigma x\big) \cdot H_{2k}(x) \frac{e^{-x^{2}}}{\sqrt{\pi}} \mathrm{d}x,
\end{align*}
where in step $\1$ we perform the integration by parts from the previous step $2k$ times. In step $\2$ we use a change of variables, and in step $\3$ we use the definition of Hermite polynomials~\eqref{def:hermite}.
Consequently, we obtain
\begin{align*}
    \int_{\real} c_{k} \nabla_{x}^{(2k)} v(y;x)  \cdot u(x;\theta)  \mathrm{d}x &= \int_{\real} c_{k} \nabla_{x}^{(2k)} v(y;x)  \cdot \frac{1}{\sqrt{2}\sigma} \phi\bigg( \frac{x-\theta}{\sqrt{2}\sigma} \bigg)  \mathrm{d}x \\
    & = \frac{c_{k}}{2^{k} \sigma^{2k}} \int_{\real} v\big(y;\theta+\sqrt{2}\sigma x\big) \cdot H_{2k}(x) \frac{e^{-x^{2}}}{\sqrt{\pi}} \mathrm{d}x.
\end{align*}
Now define $f(x;\theta,y) = v\big(y;\theta + \sqrt{2}\sigma x\big)$ for convenience, and recall our definition of the inner product operation $\langle \cdot, \cdot \rangle$ from Eq.~\eqref{eq:in-prod}. Then
\begin{align}\label{eq:integral-to-hermite}
    \int_{\real} c_{k} \nabla_{x}^{(2k)} v(y;x)  \cdot u(x;\theta)  \mathrm{d}x 
    & = \frac{c_{k}}{\sqrt{\pi}2^{k} \sigma^{2k}} \big \langle f(\cdot;\theta,y), H_{2k}(\cdot)  \big \rangle \nonumber \\
    & = \big \langle f(\cdot;\theta,y), \Nherm_{2k}(\cdot) \big \rangle \cdot \frac{c_{k}}{\sqrt{\pi} 2^{k}\sigma^{2k}} \cdot ( \sqrt{\pi} 2^{2k} (2k)! )^{1/2}  \nonumber \\
   & =  \big \langle f(\cdot;\theta,y), \Nherm_{2k}(\cdot) \big \rangle \cdot \Nherm_{2k}(0).
\end{align}
To justify Eq.~\eqref{eq:integral-to-hermite}, recall that $c_{k} = (-1)^{k} \sigma^{2k}/(2k)!!$. Then by definition of $\Nherm_{2k}(x)$ in Eq.~\eqref{def:hermite}, we have
\[
   \Nherm_{2k}(0) =  \frac{H_{2k}(0)}{\big( \sqrt{\pi} 2^{2k} (2k)! \big)^{1/2}} \overset{\1}{=} \frac{(-1)^{k} (2k)! / k! }{\big( \sqrt{\pi} 2^{2k} (2k)! \big)^{1/2}} = \frac{c_{k}}{\sqrt{\pi} 2^{k} \sigma^{2k}} \cdot ( \sqrt{\pi} 2^{2k} (2k)! )^{1/2},
\]
where in step $\1$ we use~\citet[Eq.~(5.5.5)]{szeg1939orthogonal} so that $H_{2k}(0) = (-1)^{k} (2k)! / k!$. 

Putting together the pieces, we obtain
\begin{align*}
  \int_{\real} \mathcal{S}^{\star}(y \mid x) \cdot u(x;\theta) \mathrm{d}x &= \int_{\real} \sum_{k=0}^{+\infty} c_{k} \nabla_{x}^{(2k)} v(y;x) \cdot u(x;\theta)  \mathrm{d}x \\
  &\overset{\1}{=} \sum_{k=0}^{+\infty} \int_{\real} c_{k} \nabla_{x}^{(2k)} v(y;x)  \cdot u(x;\theta)  \mathrm{d}x \\
  &\overset{\2}{=}  \sum_{k=0}^{+\infty} \big \langle f(\cdot;\theta,y), \Nherm_{2k}(\cdot) \big \rangle \cdot \Nherm_{2k}(0),
\end{align*}
where in step $\1$ we use Assumption~\ref{assump_1:target}(d) and Fubini's theorem so that switching the summation and integral is valid, and in step $\2$ we use Eq.~\eqref{eq:integral-to-hermite}.
 
\paragraph{RHS of Eq.~\eqref{eq:deconv}.}
By Assumption~\ref{assump_1:target}(c), $f(x;\theta,y)$ admits a Hermite expansion, that is,
\[
	f(x;\theta,y) = \sum_{k=0}^{+\infty} \big \langle f(\cdot ;\theta,y)  , \Nherm_{k}(\cdot) \big \rangle \cdot \Nherm_{k}(x), \quad \text{for all} \quad x,y\in \real.
\]
Setting $x = 0$ yields for all $y\in \real$,
\begin{align}\label{eq:hermite-expand-target}
	v(y;\theta) = f(0;\theta,y) = \sum_{k=0}^{+\infty} \big \langle f(\cdot ;\theta,y)  , \Nherm_{k}(\cdot) \big \rangle \cdot \Nherm_{k}(0) = \sum_{k=0}^{+\infty} \big \langle f(\cdot ;\theta,y)  , \Nherm_{2k}(\cdot) \big \rangle \cdot \Nherm_{2k}(0),
\end{align}
where in the last step we use the fact that $\Nherm_{k}(0) = 0$ when $k$ is odd; see~\citet[Eq.~(5.5.5)]{szeg1939orthogonal}. 

\paragraph{Combing the pieces.} Our two steps yield
\[
	\int_{\real} \mathcal{S}^{\star}(y \mid x) \cdot u(x;\theta) \mathrm{d}x = f(0;\theta,y) = v(y;\theta) \quad \text{for all } y\in \real, \theta \in \Theta.
\]
This concludes the proof. \qed

\subsection{Proof of Lemma~\ref{lemma:TV-deficiency-decompose}} \label{sec:pf-lem1}
Before proving the lemma, let us establish the claimed inequality~\eqref{ineq1:lemma-TV} for completeness. Applying the triangle inequality, we have
\begin{align}\label{ineq0-pr-lemma1}
	\Big \| \int_{\real} \mathcal{T}_{N}^{\star}(\cdot \mid x) \cdot u(x;\theta) \mathrm{d}x - v(\cdot;\theta)\Big\|_{1}  &\leq 
	\Big \| \int_{\real} \mathcal{T}_{N}^{\star}(\cdot \mid x) \cdot u(x;\theta) \mathrm{d}x - \int_{\real} \mathcal{S}_{N}^{\star}(\cdot \mid x) \cdot u(x;\theta) \mathrm{d}x \Big\|_{1} \nonumber \\
	&\quad + \Big \| \int_{\real} \mathcal{S}_{N}^{\star}(\cdot \mid x) \cdot u(x;\theta) \mathrm{d}x - v(\cdot;\theta) \Big\|_{1}.
\end{align}
By following the same steps as in~\citet[Section 7.2]{lou2025computationally}, we obtain
\begin{align}\label{ineq1-pr-lemma1}
	\Big \| \int_{\real} \mathcal{T}_{N}^{\star}(\cdot \mid x) \cdot u(x;\theta) \mathrm{d}x - \int_{\real} \mathcal{S}_{N}^{\star}(\cdot \mid x) \cdot u(x;\theta) \mathrm{d}x \Big\|_{1} \leq \int_{\real} \big(|p(x)-1| + q(x) \big) \cdot u(x;\theta) \mathrm{d}x.
\end{align}
Consequently, substituting Ineq.~\eqref{ineq1-pr-lemma1} into the RHS of Ineq.~\eqref{ineq0-pr-lemma1} and taking a supremum over $\theta \in \Theta$ yields inequality~\eqref{ineq1:lemma-TV}. 

We now turn to proving Ineq.~\eqref{ineq2:lemma-TV}; the reader is advised to go through the proof of Proposition~\ref{prop:optimal-signed-kernel} for context. In that proof, we defined $f(x;\theta,y) = v\big(y;\theta + \sqrt{2}\sigma x\big)$ and $c_{k} = (-1)^{k} \sigma^{2k}/(2k)!!$ for convenience. By Assumption~\ref{assump_1:target}(c) and Eq.~\eqref{eq:hermite-expand-target}, we have
\[
  v(y;\theta)  = \sum_{k=0}^{+\infty} \big \langle f(\cdot ;\theta,y)  , \Nherm_{2k}(\cdot) \big \rangle \cdot \Nherm_{2k}(0) \quad \text{for all } y\in \real, \theta \in \Theta.
\]
Next,
\begin{align*}
  \int_{\real} \mathcal{S}_{N}^{\star}(y \mid x) \cdot u(x;\theta) \mathrm{d}x & = \int_{\real} \sum_{k=0}^{N} c_{k} \nabla_{x}^{(2k)} v(y;x) \cdot u(x;\theta)  \mathrm{d}x \\
  & =  \sum_{k=0}^{N} \int_{\real} c_{k} \nabla_{x}^{(2k)} v(y;x) \cdot u(x;\theta)  \mathrm{d}x \\
  & \overset{\1}{=}  \sum_{k=0}^{N} \big \langle f(\cdot;\theta,y), \Nherm_{2k}(\cdot) \big \rangle \cdot \Nherm_{2k}(0),
\end{align*}
where in step $\1$ we use Eq.~\eqref{eq:integral-to-hermite}. Combining the two equations above yields
\begin{align*}
  \Big \| \int_{\real} \mathcal{S}_{N}^{\star}(\cdot \mid x) \cdot u(x;\theta) \mathrm{d}x - v(\cdot;\theta) \Big\|_{1} &= \int_{\real} \bigg| \int_{\real}
  \mathcal{S}_{N}^{\star}(y \mid x) u(x;\theta) \mathrm{d}x - v(y;\theta)   \bigg| \mathrm{d}y \\
  & =  \int_{\real} \bigg| \sum_{k=0}^{N} \big \langle f(\cdot;\theta,y), \Nherm_{2k}(\cdot) \big \rangle \cdot \Nherm_{2k}(0) - \sum_{k=0}^{+\infty} \big \langle f(\cdot ;\theta,y)  , \Nherm_{2k}(\cdot) \big \rangle \cdot \Nherm_{2k}(0)  \bigg| \mathrm{d}y \\
  & = \int_{\real} \bigg| \sum_{k=N+1}^{+\infty} \big \langle f(\cdot;\theta,y), \Nherm_{2k}(\cdot) \big \rangle \cdot \Nherm_{2k}(0) \bigg| \mathrm{d}y \\
  & \overset{\1}{\leq}  \int_{\real}  \sum_{k=N+1}^{+\infty} \Big|  \big \langle f(\cdot;\theta,y), \Nherm_{2k}(\cdot) \big \rangle \Big| \mathrm{d}y \\
  & = \int_{\real}  \sum_{k=N+1}^{+\infty} \bigg|  \int_{\real} v(y;\theta+\sqrt{2}\sigma x) \cdot \Nherm_{2k}(x) \exp(-x^{2}) \mathrm{d}x \bigg| \mathrm{d}y,
\end{align*}
where in step $\1$ we apply the triangle inequality and use the sequence of bounds
\[
    | \Nherm_{2k}(0) | = \frac{| H_{2k}(0) |}{\big( \sqrt{\pi} 2^{2k} (2k)! \big)^{1/2}} =  \frac{ (2k)! / k! }{\big( \sqrt{\pi} 2^{2k} (2k)! \big)^{1/2}} = \frac{(2k-1)!!}{ \big( \sqrt{\pi} (2k)! \big)^{1/2} } \leq 1.
\]
Taking supremum over $\theta \in \Theta$ on both sides of the inequality in the display above yields the claim. \qed

\section{Proofs of results from Section~\ref{sec:location-general}} \label{sec:pf_thm_non_gaussian}

We prove Theorem~\ref{thm:non-gaussian}(a) in Section~\ref{sec:pf_thm_nongaussian_term1} and Theorem~\ref{thm:non-gaussian}(b) in Section~\ref{sec:construct_reduction_nongaussian}. 
We provide the proof of the corollaries and examples following Theorem~\ref{thm:non-gaussian} in Sections~\ref{sec:pf_corollary_bounded_psi}--\ref{sec:pf_corollary_generalized_normal}.

\subsection{Proof of Theorem~\ref{thm:non-gaussian}(a)}\label{sec:pf_thm_nongaussian_term1}
We will use Lemma~\ref{lemma:TV-deficiency-decompose} to bound $\delta(\mathcal{U},\mathcal{V};\mathcal{T}_{N}^{\star})$ by the signed deficiency and other terms, and then proceed to bounding those quantities.
The first part of the proof is dedicated to verifying Assumption~\ref{assump_1:target} for the target density $v(y;\theta)$ defined in Eq.~\eqref{density-non-gaussian-target}. This will then allow us to apply Lemma~\ref{lemma:TV-deficiency-decompose}.

Assumption~\ref{assump_1:target}(a) immediately holds by Assumption~\ref{assump:log-density}(b). To verify Assumption~\ref{assump_1:target}(b),
note that by definition, $v(y;x) = g(h(x))$ where $g(x) = \frac{1}{\varT}\exp(-x)$ and $h(x) = \psi((y-x)/\varT)$ for all $x,y\in \real$. Applying Faà di Bruno's formula yields
\begin{align}\label{eq-target-derivative-nongaussian-0}
    \nabla_{x}^{(n)}v(y;x) &= \sum_{m_{1},\dots,m_{n}} \frac{n!}{m_{1}!\,m_{2}!\cdots m_{n}!} \cdot g^{(m_{1}+\cdots+m_{n})}(h(x)) \cdot \prod_{j=1}^{n} \bigg( \frac{h^{(j)}(x)}{j!} \bigg)^{m_{j}} \nonumber \\
    &=  \frac{\exp(-h(x))}{\varT} \cdot (n!) \cdot \sum_{m_{1},\dots,m_{n}} \frac{(-1)^{m_{1}+\cdots+m_{n}}}{m_{1}!\,m_{2}!\cdots m_{n}!}  \prod_{j=1}^{n} \bigg( \psi^{(j)}\bigg( \frac{y-x}{\varT} \bigg) \frac{(-1)^{j}}{\varT^{j} j!} \bigg)^{m_{j}},
\end{align} 
where the summation is over all $n$-tuples of nonnegative integers $(m_{1},\cdots,m_{n})$ satisfying $\sum_{j=1}^{n}jm_{j} = n$. Note that for all $n \in \mathbb{N}$, the following bound holds and will be used repeatedly:
\begin{align}\label{ineq:size-n-tuples}
   \bigg| \bigg\{ (m_{1},\dots,m_{n}) : \sum_{j=1}^{n} j m_{j} = n, m_{j} \geq 0, m_{j} \in \mathbb{Z} \text{ for all } j\in[n] \bigg\} \bigg| \leq \frac{n}{1} \cdot \frac{n}{2} \cdots \frac{n}{n} \leq \frac{n^{n}}{n!} \leq e^{n}. 
\end{align}
By applying the triangle inequality, we obtain
\begin{align*}
  \big| \nabla_{x}^{(n)}v(y;x) \big| &\leq \frac{\exp(-h(x))}{\varT} \cdot (n!) \cdot \sum_{m_{1},\dots,m_{n}} \frac{\prod_{j=1}^{n} \big| \psi^{(j)}\big( \frac{y-x}{\varT} \big) \frac{1}{\varT^{j} j!} \big|^{m_{j}} }{m_{1}!\,m_{2}!\cdots m_{n}!}  \\
  & \leq \frac{\exp(-h(x))}{\varT} \cdot (n!) \cdot \sum_{m_{1},\dots,m_{n}} \bigg( \sum_{j=1}^{n} \bigg| \psi^{(j)}\bigg( \frac{y-x}{\varT} \bigg) \frac{1}{\varT^{j} j!} \bigg|  \bigg)^{m_{1}+\cdots+m_{n}},
\end{align*}
where in the last step we use the multinomial theorem so that 
\[
    \frac{\prod_{j=1}^{n} \big| \psi^{(j)}\big( \frac{y-x}{\varT} \big) \frac{1}{\varT^{j} j!} \big|^{m_{j}} }{m_{1}!\,m_{2}!\cdots m_{n}!}  \leq \bigg( \sum_{j=1}^{n} \bigg| \psi^{(j)}\bigg( \frac{y-x}{\varT} \bigg) \frac{1}{\varT^{j} j!} \bigg|  \bigg)^{m_{1}+\cdots+m_{n}}.
\]
Integrating over $x \in \real$ and using the change of variables $t = (y-x)/\varT$ yields
\begin{align*}
  \int_{\real}  \big| \nabla_{x}^{(n)}v(y;x) \big| \mathrm{d}x &\leq n! \cdot \sum_{m_{1},\dots,m_{n}} \int_{\real} \exp(-\psi(t)) \bigg( \sum_{j=1}^{n}  \frac{ |\psi^{(j)}(t)|}{\varT^{j} j!}   \bigg)^{m_{1}+\cdots+m_{n}} \mathrm{d}t \\
  & \overset{\1}{\leq} n! \cdot \sum_{m_{1},\dots,m_{n}} (\Cdensity n)! \overset{\2}{\leq} n! \cdot e^{n} \cdot (\Cdensity n)!<+\infty,
\end{align*}
where in step $\1$ we use Assumption~\ref{assump:log-density}(b), $\varT \geq \CdensityT$, and $1\leq m_{1}+\cdots+m_{n} \leq n$, and in step $\2$ we use the fact that there are at most $e^{n}$ terms in the summation; see Ineq.~\eqref{ineq:size-n-tuples}. Note that $\big| \nabla_{x}^{(n)}v(y;x) \big|$ is a continuous function in $x$ by definition~\eqref{density-non-gaussian-target} and Assumption~\ref{assump:log-density}. Using the bounded integral in the display above, we obtain
\[
    \lim_{x \uparrow +\infty} \big| \nabla_{x}^{(n)}v(y;x) \big| =  \lim_{x \downarrow -\infty} \big| \nabla_{x}^{(n)}v(y;x) \big| = 0.
\]
This verifies Assumption~\ref{assump_1:target}(b). We next turn to prove Assumption~\ref{assump_1:target}(c). By definition~\eqref{density-non-gaussian-target},
\[
    v\big(y;\sqrt{2}\sigma x + \theta\big) = \frac{1}{\varT} \exp\bigg( - \psi\bigg( \frac{y-\sqrt{2}\sigma x - \theta}{\varT} \bigg) \bigg).
\]
Thus,
\begin{align*}
  \int_{\real} \big(v\big(y;\sqrt{2}\sigma x + \theta\big)\big)^{2} e^{-x^{2}} \mathrm{d}x &\leq  \int_{\real} \big(v\big(y;\sqrt{2}\sigma x + \theta\big) \big)^{2} \mathrm{d}x \\
  &= \int_{\real} \frac{1}{\sqrt{2} \sigma \varT} \exp\big(-2\psi(t)\big) \mathrm{d}t<+\infty,
\end{align*}
where in the second step we use the change of variables $t = \frac{y-\sqrt{2}\sigma x - \theta}{\varT}$, and in the last step we use Assumption~\ref{assump:log-density}(a). This proves Assumption~\ref{assump_1:target}(c). 

Having verified Assumption~\ref{assump:log-density}, we can now apply Lemma~\ref{lemma:TV-deficiency-decompose} to obtain
\begin{align}\label{ineq-pf-nongaussian-term3}
  \delta\big(\mathcal{U}, \mathcal{V};\mathcal{T}^{\star}_{N}\big) &\leq \SignedDef\big( \mathcal{U}, \mathcal{V}; \mathcal{S}_{N}^{\star} \big)  + \frac{1}{2}\sup_{\theta\in \real}\; \int_{\real} \big(|p(x) - 1| + q(x)\big) \cdot u(x;\theta) \mathrm{d}x. 
\end{align}
To bound the terms on the RHS of the above display, we use the following lemma.
\begin{lemma}\label{lemma:aux-pf-thm1-part1}
Let the source and target statistical models $\mathcal{U}$ and $\mathcal{V}$ be given by Eq.~\eqref{source-gaussian-location} and Eq.~\eqref{density-non-gaussian-target}, respectively. Let the signed kernel $\mathcal{S}_{N}^{\star}$ be defined in Eq.~\eqref{truncated-signed-kernel} and let $\SignedDef\big( \mathcal{U}, \mathcal{V}; \mathcal{S}_{N}^{\star} \big)$ be defined in Eq.~\eqref{deficiency-signed-kernel}. Suppose the tuple of parameters $(N,\varT,\lambda,\sigma)$ satisfies condition~\eqref{assump:thm-nongaussian} and Assumption~\ref{assump:log-density} holds. Let $\CdensityB$ be defined in Eq.~\eqref{def-C-B}. Then the following two parts hold.
\begin{enumerate} 
\item[(a)] We have
\begin{align}\label{ineq-gaussian-nongaussian-term4}
  \SignedDef\big( \mathcal{U}, \mathcal{V}; \mathcal{S}_{N}^{\star} \big) \leq \CdensityB \exp\bigg( - (2e)^{-1} (2N+1)^{\frac{1}{2(1+\Cdensity)}} \bigg),
\end{align}
\item[(b)] Let $p(x)$ and $q(x)$ be defined in Eq.~\eqref{px-qx-def}. Let $\err(\psi,n,\lambda)$ be defined in Eq.~\eqref{def-error-nongaussian}. Then
\begin{align} 
\label{ineq-gaussian-nongaussian-term5}
    \frac{1}{2}\sup_{\theta\in \real}\; \int_{\real} \big(|p(x) - 1| + q(x)\big) \cdot u(x;\theta) \mathrm{d}x \leq 4 \err(\psi,2N,\lambda)
\end{align}
\end{enumerate}
\end{lemma}
Indeed, substituting Ineq.~\eqref{ineq-gaussian-nongaussian-term4} and Ineq.~\eqref{ineq-gaussian-nongaussian-term5} into the RHS of Ineq.~\eqref{ineq-pf-nongaussian-term3} proves Theorem~\ref{thm:non-gaussian}(a). \qed

It remains to prove Lemma~\ref{lemma:aux-pf-thm1-part1}; we do so in the following subsections.

\subsubsection{Proof of Lemma~\ref{lemma:aux-pf-thm1-part1}(a)}\label{sec:pf_gaussian_nongaussian_term4}
By Lemma~\ref{lemma:TV-deficiency-decompose} (Ineq.~\eqref{ineq2:lemma-TV}), we have
\begin{align}\label{ineq0:pf-thm-nongaussian}
    \SignedDef\big( \mathcal{U}, \mathcal{V}; \mathcal{S}_{N}^{\star} \big) \leq \frac{1}{2}\sup_{\theta \in \real} \; \int_{\real} \sum_{k=N+1}^{+\infty} \bigg| \int_{\real} v\big(y;\theta + \sqrt{2}\sigma x\big) \cdot \Nherm_{2k}(x) e^{-x^2} \mathrm{d}x \bigg| \mathrm{d}y.
\end{align}
For each $n \in \mathbb{N}$, let $a_{n} = \sqrt{2(n+1)}$ for convenience.
En route to bounding the RHS of the inequality in the display above, let us apply integration by parts and use Eq.~\eqref{Hermite-derivative}. We obtain that for all $n \in \mathbb{N}$,
\begin{align*}
    &\int_{\real} v\big(y;\theta + \sqrt{2}\sigma x\big) \cdot \Nherm_{n}(x) e^{-x^2} \mathrm{d}x  \\
    & = \frac{1}{a_{n}}\int_{\real} v\big(y;\theta + \sqrt{2}\sigma x\big) \cdot e^{-x^2} \mathrm{d} \Nherm_{n+1}(x) \\
    & =  \frac{1}{a_{n}} \Big[ v\big(y;\theta + \sqrt{2}\sigma x\big)  e^{-x^2} \Nherm_{n+1}(x) \Big] \Big\vert_{-\infty}^{+\infty} - 
    \frac{1}{a_{n}}\int_{\real} \frac{ \mathrm{d} \big[ v\big(y;\theta + \sqrt{2}\sigma x\big) \cdot e^{-x^2} \big] }{ \mathrm{d}x} \Nherm_{n+1}(x) \mathrm{d}x.
\end{align*}
Performing the same integration by parts $r$ times yields
\begin{align}\label{ineq1:pf-thm-nongaussian}
  \int_{\real} v\big(y;\theta + \sqrt{2}\sigma x\big) \cdot \Nherm_{n}(x) e^{-x^2} \mathrm{d}x &= \sum_{\ell=0}^{r-1} \frac{(-1)^{\ell}}{\prod_{j=0}^{\ell} a_{n+j} } \bigg[ \frac{ \mathrm{d}^{\ell} \big[ v\big(y;\theta + \sqrt{2}\sigma x\big)  e^{-x^2} \big] }{\mathrm{d}x^{\ell}} \Nherm_{n+1+\ell}(x) \bigg] \bigg \vert_{-\infty}^{+\infty} \nonumber \\ 
  &  \quad + \frac{(-1)^{r}}{\prod_{j=0}^{r-1} a_{n+j} } \int_{\real} \frac{ \mathrm{d}^{r} \big[ v\big(y;\theta + \sqrt{2}\sigma x\big) \cdot e^{-x^2} \big] }{ \mathrm{d}x^{r}} \Nherm_{n+r}(x) \mathrm{d}x.
\end{align}
We claim the following equation holds for all $n\in \mathbb{N}$, $\ell \in \PosInt$, $y \in \real$, and $\varT \geq \sqrt{2} \sigma \CdensityT, \sigma>0$:
\begin{align}\label{claim1-pf-nongaussian-part1}
    \bigg[ \frac{ \mathrm{d}^{\ell} \big[ v\big(y;\theta + \sqrt{2}\sigma x\big)  e^{-x^2} \big] }{\mathrm{d}x^{\ell}} \Nherm_{n+1+\ell}(x) \bigg] \bigg \vert_{-\infty}^{+\infty} = 0,
\end{align}
where the evaluation at infinity is with respect to variable $x$. We defer the proof of claim~\eqref{claim1-pf-nongaussian-part1} to the end of this section.
Applying the product rule for higher-order derivative yields for all $r \in \PosInt$,
\begin{align} \label{eq:pd-rule-nongaussian}
    \frac{ \mathrm{d}^{r} \big[ v\big(y;\theta + \sqrt{2}\sigma x\big) \cdot e^{-x^2} \big] }{ \mathrm{d}x^{r}} &= \sum_{k = 0}^{r} \binom{r}{k} \nabla_{x}^{(r-k)} v\big(y;\theta + \sqrt{2}\sigma x\big) \cdot \frac{\mathrm{d}^{k} e^{-x^2}}{\mathrm{d}x^{k}} \nonumber \\
    &= \sum_{k = 0}^{r} \binom{r}{k} \nabla_{x}^{(r-k)} v\big(y;\theta + \sqrt{2}\sigma x\big) \cdot (-1)^{k} H_{k}(x) e^{-x^{2}}.
\end{align}
Substituting Eq.~\eqref{eq:pd-rule-nongaussian} and Eq.~\eqref{claim1-pf-nongaussian-part1} into Eq.~\eqref{ineq1:pf-thm-nongaussian} and applying the triangle inequality yields
\begin{align}\label{ineq2:pf-thm-nongaussian}
    &\int_{\real} \bigg| \int_{\real} v\big(y;\theta + \sqrt{2}\sigma x\big) \cdot \Nherm_{n}(x) e^{-x^2} \mathrm{d}x \bigg| \mathrm{d}y \nonumber \\
    &\leq  \frac{1}{\prod_{j=0}^{r-1} a_{n+j}} \sum_{k = 0}^{r} \binom{r}{k} \int_{\real}\int_{\real} \bigg| \nabla_{x}^{(r-k)} v\big(y;\theta + \sqrt{2}\sigma x\big) \cdot  H_{k}(x) e^{-x^{2}} \cdot  \Nherm_{n+r}(x)  \bigg| \mathrm{d}x\mathrm{d}y \nonumber \\
    & \leq \frac{1}{\prod_{j=0}^{r-1} a_{n+j}} \sum_{k = 0}^{r} \binom{r}{k} \cdot \sqrt{\gamma_{k}} \int_{\real}\int_{\real} \bigg| \nabla_{x}^{(r-k)} v\big(y;\theta + \sqrt{2}\sigma x\big) \cdot  e^{-x^{2}/2}   \bigg| \mathrm{d}x\mathrm{d}y,
\end{align}
where in the last step we let $\gamma_{k} = \sqrt{\pi} 2^{k} k!$ and use $\big|H_{k}(x)e^{-x^{2}/2}\big| \leq \sqrt{\gamma_{k}}$ and also $\big|\Nherm_{n+r}(x)\big| \leq 1$ for all $x\in\real$, by Ineq.~\eqref{Hermite-boundness}. We further claim that for all $r,k\in \PosInt$ and $k \leq r$, and for all $\varT \geq 4\sqrt{2}e\sigma \CdensityT$, we have
\begin{align}\label{claim2-pf-nongaussian-part1}
  \int_{\real}\int_{\real} \bigg| \nabla_{x}^{(r-k)} v\big(y;\theta + \sqrt{2}\sigma x\big) \cdot  e^{-x^{2}/2}   \bigg| \mathrm{d}x\mathrm{d}y \leq 4\pi (r-k)! \cdot \frac{1}{4^{r-k}} \cdot \big(\Cdensity(r-k)\big)!
\end{align}
We defer the proof of Claim~\eqref{claim2-pf-nongaussian-part1} to the end of this section.
Substituting Ineq.~\eqref{claim2-pf-nongaussian-part1} into Ineq.~\eqref{ineq2:pf-thm-nongaussian} yields
\begin{align}\label{ineq3:pf-thm-nongaussian}
    \int_{\real} \bigg| \int_{\real} v\big(y;\theta + \sqrt{2}\sigma x\big) \cdot \Nherm_{n}(x) e^{-x^2} \mathrm{d}x \bigg| \mathrm{d}y &\leq 
    \frac{4\pi}{\prod_{j=0}^{r-1} a_{n+j}} \sum_{k = 0}^{r} \binom{r}{k} \cdot \sqrt{\gamma_{k}} \cdot \frac{(r-k)!}{4^{r-k}} \cdot \big(\Cdensity(r-k)\big)! \nonumber \\
    & \overset{\1}{\leq}  \frac{4\pi^{5/4} \cdot r! \cdot \big(\Cdensity \cdot r\big)! \cdot 2^{r/2} }{\prod_{j=0}^{r-1} a_{n+j}}  \sum_{k=0}^{r} \frac{1}{4^{r-k}} \nonumber \\
    & \overset{\2}{\leq} \frac{8\pi^{5/4} \cdot r! \cdot \big(\Cdensity \cdot r\big)! \cdot 2^{r/2} }{2^{r/2} \cdot n^{r/2} } = \frac{8\pi^{5/4} \cdot r! \cdot \big(\Cdensity \cdot r\big)!  }{ n^{r/2} },
\end{align}
where in the step $\1$ we use $\binom{r}{k} = r! / \big((r-k)!\cdot k! \big)$ and $\gamma_{k} = \sqrt{\pi} 2^{k} \cdot k!$, and in step $\2$ we recall that $a_{n+j} = \sqrt{2(n+j+1)} \geq \sqrt{2n}$ for all $j$ and $\sum_{k=0}^{r} 4^{-(r-k)} \leq 2$. Note that the inequality in the display above holds for all $r \in \PosInt$. We now set
\[
    r = \left \lfloor n^{\frac{1}{2(1+\Cdensity)}}/\big(e \big(\Cdensity \vee 1\big) \big) \right \rfloor \quad \text{so that} \quad \frac{ \big( \big(\Cdensity \vee 1 \big) r \big)^{1+\Cdensity} }{\sqrt{n}} \leq e^{-(1+\Cdensity)}.
\]
Consequently, we obtain
\begin{align}\label{ineq4:pf-thm-nongaussian}
  \frac{ r! \cdot \big(\Cdensity \cdot r\big)!  }{ n^{r/2} } \leq \frac{ r^{r} \cdot \big( \big( \Cdensity \vee 1 \big) r\big)^{\Cdensity  r}  }{ n^{r/2} } &= \bigg( \frac{\big( \big(\Cdensity \vee 1 \big)  r\big)^{1+\Cdensity}}{\sqrt{n}} \bigg)^{r} \nonumber \\
  & \leq e^{-(1+\Cdensity)r} \leq e^{1+\Cdensity} \exp\Big(-e^{-1} n^{\frac{1}{2(1+\Cdensity)}}  \Big),
\end{align}
where in the last step we use that by definition, we have $r \geq \frac{n^{\frac{1}{2(1+\Cdensity)}}}{e (\Cdensity \vee 1)} - 1$. Substituting Ineq.~\eqref{ineq4:pf-thm-nongaussian} into the RHS of Ineq.~\eqref{ineq3:pf-thm-nongaussian} yields that for all $n \in \PosInt$,
\begin{align*}
    \int_{\real} \bigg| \int_{\real} v\big(y;\theta + \sqrt{2}\sigma x\big) \cdot \Nherm_{n}(x) e^{-x^2} \mathrm{d}x \bigg| \mathrm{d}y \leq 8\pi^{5/4} e^{1+\Cdensity} \exp\Big(-e^{-1} n^{\frac{1}{2(1+\Cdensity)}}  \Big).
\end{align*}
Substituting the inequality in the display above into the RHS of Ineq.~\eqref{ineq0:pf-thm-nongaussian} yields
\begin{align*}
    \SignedDef\big( \mathcal{U}, \mathcal{V}; \mathcal{S}_{N}^{\star} \big) &\leq 4\pi^{5/4} e^{1+\Cdensity} \sum_{n = 2N+2}^{+\infty} \exp\Big(-e^{-1} n^{\frac{1}{2(1+\Cdensity)}}  \Big) \\
    &\leq 4\pi^{5/4} e^{1+\Cdensity} \int_{2N+1}^{+\infty} \exp\Big(-e^{-1} t^{\frac{1}{2(1+\Cdensity)}}  \Big) \mathrm{d}t.
\end{align*}
To conclude, we must bound the integral in the display above. To reduce the notational burden, let $a = 2N+1$, $b = 2(1+\Cdensity)$, and $C = e^{e^{-1}} >1$. Note that $\frac{ba^{-1/b}}{\log(C)} \leq 1$ by assumption. We obtain 
\begin{align*}
  \int_{2N+1}^{+\infty}  \exp\Big(-e^{-1} t^{\frac{1}{2(1+\Cdensity)}} \Big) \mathrm{d}t &= \int_{a}^{+\infty} C^{-t^{1/b}} \mathrm{d}t = b\int_{a^{1/b}}^{+\infty} z^{b-1}e^{-\log(C) z} \mathrm{d}z,
\end{align*}
where in the last step we use the change of variables $z = t^{1/b}$. Integrating by parts then yields
\begin{align*}
  \int_{a^{1/b}}^{+\infty} z^{b-1}e^{-\log(C) z} \mathrm{d}z &= e^{-\log(C) a^{1/b}} \cdot \bigg( \frac{a^{\frac{b-1}{b}}}{\log(C)} + \frac{(b-1)a^{\frac{b-2}{b}}}{\log(C)^{2}} 
  + \frac{(b-1)(b-2)a^{\frac{b-3}{b}}}{\log(C)^{3}} + \cdots + \frac{b!}{\log(C)^{b}} \bigg) \\
  &\leq e^{-\log(C) a^{1/b}} \cdot \frac{a^{\frac{b-1}{b}}b}{\log(C)}.
\end{align*}
Here the second step uses the fact that since $\frac{ba^{-1/b}}{\log(C)} \leq 1$, the series is decreasing, and thus the total sum can be bounded by the first term multiplied by the total number of terms. Putting all the  pieces together, we have obtained for $2N+1 \geq (2e(1+\Cdensity))^{2(1+\Cdensity)}$,
\begin{align*}
   \SignedDef\big( \mathcal{U}, \mathcal{V}; \mathcal{S}_{N}^{\star} \big) &\leq 16 \pi^{\frac{5}{4}} e^{2+\Cdensity} (1+\Cdensity)^{2} (2N+1) \exp\bigg(-e^{-1}(2N+1)^{\frac{1}{2(1+\Cdensity)}} \bigg) \\
   & = 16 \pi^{\frac{5}{4}} e^{2+\Cdensity} (1+\Cdensity)^{2} \exp\bigg(-(2e)^{-1}(2N+1)^{\frac{1}{2(1+\Cdensity)}} \bigg) \\ 
   & \quad \times (2e)^{2(1+\Cdensity)} \cdot \frac{2N+1}{(2e)^{2(1+\Cdensity)}} \cdot \exp\bigg(-(2e)^{-1}(2N+1)^{\frac{1}{2(1+\Cdensity)}} \bigg) \\
   & \overset{\1}{\leq} 16 \pi^{\frac{5}{4}} e^{2+\Cdensity} (1+\Cdensity)^{2} (2e)^{2(1+\Cdensity)} \exp\bigg(-(2e)^{-1}(2N+1)^{\frac{1}{2(1+\Cdensity)}} \bigg) ,
\end{align*}
where in step $\1$ use the numerical inequality
\begin{align}\label{ineq:numerical-nongaussian}
    \frac{t}{(2e)^{2(1+\Cdensity)}} \leq \exp\Big( (2e)^{-1} t^{\frac{1}{2(1+\Cdensity)}} \Big) \quad \text{for all } t \geq (2e)^{2(1+\Cdensity)} e^{4(1+\Cdensity)^{2}}.
\end{align}
In particular, we apply this inequality for $t = 2N+1$, which is in the desired range by assumption. To conclude, note that Ineq.~\eqref{ineq-gaussian-nongaussian-term4} holds by recalling the definition $\CdensityB$ in Eq.~\eqref{def-C-B}.

It remains to verify the claims~\eqref{claim1-pf-nongaussian-part1},~\eqref{claim2-pf-nongaussian-part1}, and~\eqref{ineq:numerical-nongaussian}.

\paragraph{Proof of Claim~\eqref{claim1-pf-nongaussian-part1}.} 
By definition~\eqref{density-non-gaussian-target}, we obtain $v\big(y;\theta + \sqrt{2}\sigma x\big) = g(h(x))$, where $g(x) = \exp(-x)/\varT$ and $h(x) = \psi\big(\big(y-\theta - \sqrt{2}\sigma x\big)/\varT \big)$. 
Applying the Faà di Bruno's formula yields
\begin{align}\label{eq-target-derivative-nongaussian}
    \nabla_{x}^{(n)}v\big(y;\theta + \sqrt{2}\sigma x\big) &= \sum_{m_{1},\dots,m_{n}} \frac{n!}{m_{1}!\,m_{2}!\cdots m_{n}!} \cdot g^{(m_{1}+\cdots+m_{n})}(h(x)) \cdot \prod_{j=1}^{n} \bigg( \frac{h^{(j)}(x)}{j!} \bigg)^{m_{j}} \nonumber \\
    &=  \frac{\exp(-h(x))}{\varT} \cdot (n!) \cdot \sum_{m_{1},\dots,m_{n}} \frac{(-1)^{m_{1}+\cdots+m_{n}}}{m_{1}!\,m_{2}!\cdots m_{n}!}  \prod_{j=1}^{n} \bigg(  \frac{\psi^{(j)}\big( \frac{y-\theta-\sqrt{2} \sigma x}{\varT} \big) \cdot \big(-\sqrt{2}\sigma\big)^{j}}{\varT^{j} j!} \bigg)^{m_{j}},
\end{align} 
where the summation is over all $n$-tuples of nonnegative integers $(m_{1},\cdots,m_{n})$ satisfying $\sum_{j=1}^{n}jm_{j} = n$. Applying the triangle inequality and using the multinomial theorem yield
\begin{align*}
   \big| \nabla_{x}^{(n)}v\big(y;\theta + \sqrt{2}\sigma x\big) \big| &\leq \frac{\exp(-h(x))}{\varT} \cdot (n!) \cdot \sum_{m_{1},\dots,m_{n}} \bigg( \sum_{j=1}^{n} \frac{\big| \psi^{(j)}\big( \frac{y-\theta-\sqrt{2} \sigma x}{\varT} \big) \big| \cdot \big(\sqrt{2}\sigma\big)^{j}}{\varT^{j} j!}  \bigg)^{m_{1}+\cdots+m_{n}} \\
   &\leq \frac{\exp(-h(x))}{\varT} \cdot (n!) \cdot \sum_{m_{1},\dots,m_{n}} \bigg( \sum_{j=1}^{n} \frac{\big| \psi^{(j)}\big( \frac{y-\theta-\sqrt{2} \sigma x}{\varT} \big) \big| }{\big( \CdensityT \big)^{j} j!}  \bigg)^{m_{1}+\cdots+m_{n}},
\end{align*}
where in the last step we use $\varT \geq \sqrt{2}\sigma \CdensityT$. Integrating over $x\in \real$ and using change of variables $t = \frac{y-\theta-\sqrt{2} \sigma x}{\varT}$ yield that for all $n \in \PosInt$ and $y \in \real$,
\begin{align*}
  \int_{\real} \big| \nabla_{x}^{(n)}v\big(y;\theta + \sqrt{2}\sigma x\big) \big| \mathrm{d}x &\leq \frac{n!}{\sqrt{2}\sigma} \sum_{m_{1},\dots,m_{n}} \int_{\real} \exp(-\psi(t)) \bigg( \sum_{j=1}^{n} \frac{| \psi^{(j)}( t ) | }{\big( \CdensityT \big)^{j} j!}  \bigg)^{m_{1}+\cdots+m_{n}} \mathrm{d}t
  \\ & \leq \frac{n!}{\sqrt{2}\sigma} \cdot e^{n} \cdot (\Cdensity n)! < +\infty,
\end{align*}
where we have used Assumption~\ref{assump:log-density} and the fact that there are at most $e^{n}$ terms in the summation. Combining the inequality in the display with Eq.~\eqref{eq:pd-rule-nongaussian} and applying the triangle inequality yield
\begin{align*}
    &\int_{\real}\bigg| \frac{ \mathrm{d}^{\ell} \big[ v\big(y;\theta + \sqrt{2}\sigma x\big) \cdot e^{-x^2} \big] }{ \mathrm{d}x^{\ell}} \cdot \Nherm_{n+1+\ell}(x) \bigg|  \mathrm{d}x \\
    &\leq \sum_{k = 0}^{\ell} \binom{\ell}{k} \int_{\real} \big| \nabla_{x}^{(\ell-k)} v\big(y;\theta + \sqrt{2}\sigma x\big) \big| \cdot \big| \Nherm_{n+1+\ell}(x) \cdot H_{k}(x) e^{-x^{2}} \big| \mathrm{d}x \\
    &\overset{\1}{\leq} \sum_{k = 0}^{\ell} \binom{\ell}{k} \big( \sqrt{\pi} 2^{k} k! \big)^{1/2} \int_{\real} \big| \nabla_{x}^{(\ell-k)} v\big(y;\theta + \sqrt{2}\sigma x\big) \big| \mathrm{d}x <+\infty
\end{align*}
where in step $\1$ we use the boundness property of Hermite polynomial~\eqref{Hermite-boundness}. Consequently, for a positive and continuous function, if the integral over $\real$ is finite, then this function must converge to zero at infinity. Thus, 
\begin{align*}
  \lim_{|x| \uparrow +\infty}\; \frac{ \mathrm{d}^{\ell} \big[ v\big(y;\theta + \sqrt{2}\sigma x\big) \cdot e^{-x^2} \big] }{ \mathrm{d}x^{\ell}} \cdot \Nherm_{n+1+\ell}(x) = 0.
\end{align*}
This proves Claim~\eqref{claim1-pf-nongaussian-part1}.

\paragraph{Proof of Claim~\eqref{claim2-pf-nongaussian-part1}.}
Using Eq.~\eqref{eq-target-derivative-nongaussian} with the triangle inequality yields for all $n \in \mathbb{N}$, 
\begin{align*}
   \big| \nabla_{x}^{(n)}v\big(y;\theta + \sqrt{2}\sigma x\big) \big| 
    & \leq  \frac{\exp(-h(x))}{\varT} \sum_{m_{1},\dots,m_{n}} \frac{n!}{m_{1}!\,m_{2}!\cdots m_{n}!}  \prod_{j=1}^{n} \bigg(  \frac{ \big| \psi^{(j)}\big( \frac{y-\theta-\sqrt{2} \sigma x}{\varT} \big) \big| \cdot \big(\sqrt{2}\sigma\big)^{j}}{\varT^{j} j!} \bigg)^{m_{j}} \\
    & =\bigg( \frac{\sqrt{2}\CdensityT \sigma}{\varT} \bigg)^{n} \cdot \frac{\exp(-h(x))}{\varT} \sum_{m_{1},\dots,m_{n}} \frac{n!}{m_{1}!\,m_{2}!\cdots m_{n}!}  \prod_{j=1}^{n} \bigg(  \frac{ \big| \psi^{(j)}\big( \frac{y-\theta-\sqrt{2} \sigma x}{\varT} \big) \big|}{\big(\CdensityT\big)^{j} j!} \bigg)^{m_{j}},
\end{align*}
where in the second step we use $\sum_{j=1}^{n}jm_{j} = n$ so that
\[
    \bigg( \frac{\sqrt{2}\CdensityT \sigma}{\varT} \bigg)^{n} \prod_{j=1}^{n} \bigg(\frac{1}{\big(\CdensityT\big)^{j} }\bigg)^{m_{j}} = \bigg( \frac{\sqrt{2} \sigma}{\varT} \bigg)^{n} = \prod_{j=1}^{n} \bigg( \frac{\big(\sqrt{2} \sigma\big)^{j}}{\varT^{j}} \bigg)^{m_{j}}.
\]
Continuing and applying the multinomial theorem yields
\begin{align*}
  \big| \nabla_{x}^{(n)}v\big(y;\theta + \sqrt{2}\sigma x\big) \big| \leq \bigg( \frac{\sqrt{2}\CdensityT \sigma}{\varT} \bigg)^{n} \cdot (n!) \cdot \frac{\exp(-h(x))}{\varT} \cdot
  \sum_{m_{1},\dots,m_{n}}   \bigg(  \sum_{j=1}^{n}  \frac{ \big| \psi^{(j)}\big( \frac{y-\theta-\sqrt{2} \sigma x}{\varT} \big) \big|}{\big(\CdensityT\big)^{j} j!} \bigg)^{m_{1}+\cdots+m_{n}}.
\end{align*}
Recall that we let $h(x) = \psi\big(\big(y-\theta-\sqrt{2}\sigma x\big)/\varT\big)$. Integrating the above expression over $x,y \in \real$ then yields
\begin{align*}
  &\int_{\real} \int_{\real}  \big| \nabla_{x}^{(n)}v\big(y;\theta + \sqrt{2}\sigma x\big) \big| e^{-x^{2}/2} \mathrm{d}x \mathrm{d}y \\
  &\leq \bigg( \frac{\sqrt{2}\CdensityT \sigma}{\varT} \bigg)^{n} \cdot (n!) \cdot \sum_{m_{1},\dots,m_{n}} \int_{\real} \int_{\real} \frac{1}{\varT} e^{-\psi\big(\frac{y-\theta-\sqrt{2}\sigma x}{\varT}\big)} \bigg(  \sum_{j=1}^{n}  \frac{ \big| \psi^{(j)}\big( \frac{y-\theta-\sqrt{2} \sigma x}{\varT} \big) \big|}{\big(\CdensityT\big)^{j} j!} \bigg)^{m_{1}+\cdots+m_{n}} \mathrm{d}y \cdot e^{-x^{2}/2} \mathrm{dx} \\
  &\overset{\1}{=} \bigg( \frac{\sqrt{2}\CdensityT \sigma}{\varT} \bigg)^{n} \cdot (n!) \cdot \sum_{m_{1},\dots,m_{n}} \int_{\real} \int_{\real} \exp(-\psi(t))
  \bigg(  \sum_{j=1}^{n}  \frac{ \big| \psi^{(j)}(t) \big|}{\big(\CdensityT\big)^{j} j!} \bigg)^{m_{1}+\cdots+m_{n}} \mathrm{d}t \cdot e^{-x^{2}/2} \mathrm{d}x \\
  & \overset{\2}{\leq} \bigg( \frac{\sqrt{2}\CdensityT \sigma}{\varT} \bigg)^{n} \cdot (n!) \cdot \sum_{m_{1},\dots,m_{n}} (\Cdensity n)! \cdot \int_{\real} e^{-x^{2}/2} \mathrm{d}x \\
  & \leq \bigg( \frac{\sqrt{2}\CdensityT \sigma}{\varT} \bigg)^{n} \cdot (n!) \cdot e^{n} \cdot (\Cdensity n)! \cdot \sqrt{2\pi},
\end{align*}
where in step $\1$ we first integrate over $y \in \real$ and use the change of variables $t = (y-\theta-\sqrt{2}\sigma x)/\varT$, in step $\2$ we use Assumption~\ref{assump:log-density}, and in the last step we use there are at most $e^{n}$ terms in the summation (see Ineq.~\eqref{ineq:size-n-tuples}) and also that $\int_{\real} e^{-x^2/2}\mathrm{d}x = \sqrt{2\pi}$. Thus, by letting $\varT \geq 4\sqrt{2}e\sigma \CdensityT$, we obtain
\[
    \int_{\real} \int_{\real}  \big| \nabla_{x}^{(n)}v\big(y;\theta + \sqrt{2}\sigma x\big) \big| e^{-x^{2}/2} \mathrm{d}x \mathrm{d}y \leq \sqrt{2\pi} (n!) \cdot (\Cdensity n)! \cdot \frac{1}{4^{n}} \quad \text{for all } n\in \PosInt.
\]
This proves Claim~\eqref{claim2-pf-nongaussian-part1}.

\paragraph{Proof of Claim~\eqref{ineq:numerical-nongaussian}.} By letting $x = (2e)^{-1} t^{\frac{1}{2(1+\Cdensity)}} $, it is equivalent to show
\[
    x^{2(1+\Cdensity)} \leq e^{x} \quad \text{for all } x \geq e^{2(1+\Cdensity)},
\]
which is equivalent to, by taking the log on both sides and re-arranging the terms, 
\[
   h(x) :=  2(1+\Cdensity) \frac{ \log(x) }{ x} \leq 1 \quad \text{for all } x \geq e^{2(1+\Cdensity)}. 
\]

To show this, simply note that the function $h(x)$ is monotone decreasing for all $x\geq 1$. Thus,
\[
   h(x) \leq h\big(e^{2(1+\Cdensity)}\big) = \frac{(2(1+\Cdensity))^{2}}{e^{2(1+\Cdensity)}} \leq \sup_{a>0} \frac{a^{2}}{e^{a}} = \frac{2^{2}}{e^{2}} \leq 1  \quad \text{for all } x \geq e^{2(1+\Cdensity)},
\]
as claimed.
\qed

\subsubsection{Proof of Lemma~\ref{lemma:aux-pf-thm1-part1}(b)}\label{sec:pf_gaussian_nongaussian_term5}
We first state an auxiliary lemma and defer its proof to Appendix~\ref{sec:pf_lemma_aux_px_qx_nongaussian}.
\begin{lemma}\label{lemma:aux-px-qx-nongaussian}
Consider the signed kernel $\mathcal{S}_{N}^{\star}(y \mid x)$ defined in Eq.~\eqref{truncated-signed-kernel} for the target density $v(y;\theta)$ defined in Eq.~\eqref{density-non-gaussian-target}. Let the set $\Set(\psi,n,\lambda)$, its complement $\Set^{\complement}(\psi,n,\lambda)$, and the scalar $\err(\psi,n,\lambda)$ be defined as in Eq.~\eqref{def-set-complement} and Eq.~\eqref{def-error-nongaussian}, respectively. Suppose Assumption~\ref{assump:log-density} holds with the pair of constants $(\Cdensity, \CdensityT)$. Then the following holds for all $N \in \mathbb{N}$ and $\lambda>0$.
\begin{enumerate}
    \item[(i)] If $\varT \geq 8e\lambda \sigma N$, then $\mathcal{S}_{N}^{\star}(y \mid x) \geq 0$ as long as $(y-x)/\varT \in \Set(\psi,2N,\lambda)$.
    \item[(ii)] Let $\rho(t) = \exp(-\psi(t))$ for all $t \in \real$. Then for all $1\leq n \leq 2N$ and $x \in \real$, we have
    \begin{align}
        \label{ineq1:integral-derivative-density}
        &\int_{\real} \rho^{(n)}(t) \mathrm{d}t = 0
        \quad \text{and} \quad \int_{ \frac{y-x}{\varT} \in \Set^{\complement}(\psi,2N,\lambda) } \big| \nabla_{x}^{(n)} v(y;x) \big| \mathrm{d}y \leq \frac{n! \cdot \lambda^{n} \cdot e^{n}}{\varT^{n}} \cdot \err(\psi,2N,\lambda).
    \end{align}
\end{enumerate}
\end{lemma}

We now use the lemma above to prove inequality~\eqref{ineq-gaussian-nongaussian-term5}. Note that $v(y;x) = \frac{1}{\varT} \rho(\frac{y-x}{\varT})$. Consequently, we obtain for all $n\in \PosInt$ and $x,y\in \real$ that
\begin{align}\label{derivative-v-rho}
  \nabla_{x}^{(n)}v(y;x) = \frac{1}{\varT} \rho^{(n)}\bigg(\frac{y-x}{\varT}\bigg) \bigg(\frac{-1}{\varT}\bigg)^{n}.
\end{align}

Let us now upper bound $p(x)$~\eqref{px-qx-def}. Using part (i) of Lemma~\ref{lemma:aux-px-qx-nongaussian} and definition~\eqref{px-qx-def}, we obtain
\begin{align}\label{ineq-px-up-nongaussian}
  p(x) &\leq \int_{\frac{y-x}{\varT} \in \Set(\psi,2N,\lambda)} \mathcal{S}_{N}^{\star}(y \mid x) \mathrm{d}y + \int_{\frac{y-x}{\varT} \in \Set^{\complement}(\psi,2N,\lambda)} \big|  \mathcal{S}_{N}^{\star}(y \mid x) \big| \mathrm{d}y \nonumber \\
  & \leq  \int_{\real} \mathcal{S}_{N}^{\star}(y \mid x) \mathrm{d}y + 2\int_{ \frac{y-x}{\varT} \in \Set^{\complement}(\psi,2N,\lambda)} \big|  \mathcal{S}_{N}^{\star}(y \mid x) \big| \mathrm{d}y.
\end{align}
Next, by definition, we have
\begin{align}\label{eq-px-up-term1}
  \int_{\real} \mathcal{S}_{N}^{\star}(y \mid x) \mathrm{d}y = \sum_{k=0}^{N} \frac{(-1)^{k} \sigma^{2k}}{(2k)!!} \int_{\real} \nabla_{x}^{(2k)}v(y;x) \mathrm{d} y \overset{\1}{=} 
  \sum_{k=0}^{N} \frac{(-1)^{k} \sigma^{2k}}{(2k)!!} \int_{\real} \frac{\rho^{(2k)}(t)}{\varT^{2k}} \mathrm{d} t = 1,
\end{align}
where in step $\1$ we use Eq.~\eqref{derivative-v-rho} and use the change of variables $t = (y-x)/\varT$, and in the last step we use Ineq.~\eqref{ineq1:integral-derivative-density} and note that $\int_{\real} \rho(t) \mathrm{d}t = 1$. Towards bounding the second term on the RHS of Ineq.~\eqref{ineq-px-up-nongaussian}, we obtain by using the triangle inequality that
\begin{align}\label{ineq-px-up-term2}
    \int_{ \frac{y-x}{\varT} \in \Set^{\complement}(\psi,2N,\lambda)} \big|  \mathcal{S}_{N}^{\star}(y \mid x) \big| \mathrm{d}y &\leq \sum_{k=0}^{N} \frac{\sigma^{2k}}{(2k)!!} \int_{\frac{y-x}{\varT} \in \Set^{\complement}(\psi,2N,\lambda)} \big| \nabla_{x}^{(2k)}v(y;x) \big| \mathrm{d}y \nonumber \\
    &\overset{\1}{\leq} \int_{\frac{y-x}{\varT} \in \Set^{\complement}(\psi,2N,\lambda)} v(y;x) \mathrm{d}y  + \sum_{k=1}^{N} \frac{\sigma^{2k}}{(2k)!!} \frac{(2k)! \cdot \lambda^{2k} \cdot e^{2k}}{\varT^{2k}} \cdot \err(\psi,2N,\lambda) \nonumber \\
    & \overset{\2}{\leq} 2 \err(\psi,2N,\lambda),
\end{align}
where in step $\1$ we use Ineq.~\eqref{ineq1:integral-derivative-density}. In step $\2$ we use
\[
    \int_{\frac{y-x}{\varT} \in \Set^{\complement}(\psi,2N,\lambda)} v(y;x) \mathrm{d}y = \int_{\frac{y-x}{\varT} \in \Set^{\complement}(\psi,2N,\lambda)} \frac{\rho(\frac{y-x}{\varT})}{\varT} \mathrm{d}y = \int_{\Set^{\complement}(\psi,2N,\lambda)}e^{-\psi(t)} \mathrm{d}t \leq \err(\psi,2N,\lambda),
\]
and also use $\varT \geq 8e \lambda \sigma N$ so that we have the chain of bounds
\[
  \sum_{k=1}^{N}  \frac{\sigma^{2k} \cdot (2k)! \cdot (e\lambda)^{2k}}{(2k)!! \cdot \varT^{2k}} \leq \sum_{k=1}^{N} \bigg( \frac{ 2e \lambda \sigma k }{\varT} \bigg)^{2k} \leq \sum_{k=1}^{N} \frac{1}{4^{k}} \leq 1.
\] 
Substituting Eq.~\eqref{eq-px-up-term1} and Ineq.~\eqref{ineq-px-up-term2} into the RHS of Ineq.~\eqref{ineq-px-up-nongaussian} yields
\begin{align*}
  p(x) \leq 1 + 4 \err(\psi,2N,\lambda) \quad \text{for all } x\in\real. 
\end{align*}
We now bound $q(x)$~\eqref{px-qx-def}. We obtain by definition that
\[
    p(x) - q(x) = \int_{\real} \big[ \mathcal{S}_{N}^{\star}(y\mid x) \vee 0 \big] \mathrm{d}y + \int_{\real} \big[ \mathcal{S}_{N}^{\star}(y\mid x) \wedge 0 \big] \mathrm{d}y  = \int_{\real} \mathcal{S}_{N}^{\star}(y\mid x) \mathrm{d}y = 1,
\]
where in the last step we use Eq.~\eqref{eq-px-up-term1}. Note that $q(x) \geq 0$ by its definition~\eqref{px-qx-def} and thus, 
\begin{align}\label{ineq:px-geq-1}
   p(x) = q(x) + 1 \geq 1 \quad \text{for all}\quad  x \in \real. 
\end{align}
Putting the pieces together yields
\[
    \big|p(x)-1\big| + q(x) \leq 2\big|p(x)-1\big| \leq 8 \err(\psi,2N,\lambda) \quad \text{for all } x \in \real.
\]
Note that the RHS of the inequality in the display is independent of $x$. Integrating it over the density $u(x;\theta)$ and taking a supremum over $\theta$, we complete the proof of Ineq.~\eqref{ineq-gaussian-nongaussian-term5}.
\qed

\subsection{Proof of Theorem~\ref{thm:non-gaussian}(b)}\label{sec:construct_reduction_nongaussian}
We split the proof into several parts.
\paragraph{Constructing the reduction based on a general rejection kernel.}
We define a class of base measures $\left\{ \mathcal{D}(\cdot \mid x) \right\}_{x \in \mathbb{R}}$ as
\begin{align}\label{base-measure-nongaussian}
    \mathcal{D}(y \mid x) = \frac{1}{2\varT} \exp\bigg( -\psi\bigg(\frac{y-x}{2\varT}\bigg) \bigg), \quad \text{for all } x, y \in \mathbb{R}.
\end{align}
We recall the general rejection kernel $x \mapsto \textsc{rk}(x, T, M, y_{0})$ from~\citet[page 10]{lou2025computationally}, and employ it using the base measures $\mathcal{D}(y \mid x)$ defined in~\eqref{base-measure-nongaussian} and the signed kernel $\mathcal{S}_{N}^{\star}(y \mid x)$ defined in~\eqref{truncated-signed-kernel}. The rejection kernel $\textsc{rk}$ requires as input a source sample $x$, the total number of iterations $T$, and a scalar $M > 0$ that is required to satisfy, for input $x$, the bound 
\[
    \frac{\mathcal{S}_{N}^{\star}(y \mid x) \vee 0}{ \mathcal{D}(y \mid x)} \leq M, \quad \text{for all } y \in \mathbb{R}.
\]
An output sample point $y_{0}$ is also required, and serves as an initialization for the algorithm.

Staying faithful to the theorem statement, we consider $T_{\mathsf{iter}}$ to be any nonnegative integer, set $y_{0} = 0$, and recall $M(\psi,2N,\lambda)$ in Eq.~\eqref{def-M-nongaussian}. We will now construct our reduction algorithm from the rejection kernel and show that it satisfies the conditions required for $\textsc{rk}$ to succeed. For $X_{\theta} \sim \NORMAL(\theta,\sigma^{2})$, define
\begin{align}\label{reduction-nongaussian}
  \mathsf{K}(X_{\theta}): =  \textsc{rk}\big(X_{\theta}, T_{\mathsf{iter}}, 4 M(\psi,2N,\lambda), y_{0} \big).
\end{align}
Note that for all $x,y \in \real$, we have by definition of $\mathcal{S}_{N}^{\star}(y \mid x)$~\eqref{truncated-signed-kernel},
\begin{align*}
  \frac{\mathcal{S}_{N}^{\star}(y \mid x) \vee 0}{ \mathcal{D}(y \mid x)} &\leq \sum_{k=0}^{N} \frac{\sigma^{2k}}{(2k)!!} \frac{\big|\nabla_{x}^{(2k)}v(y;x)\big|}{\mathcal{D}(y \mid x)} \\ 
  & \overset{\1}{\leq} 2 e^{ -\psi(\frac{y-x}{\varT}) + \psi(\frac{y-x}{2\varT})}\Bigg( 1 +  \sum_{k=1}^{N} \frac{\sigma^{2k} (2k)! \lambda^{2k} }{(2k)!! \varT^{2k}} \sum_{m_{1},\dots,m_{2k}} \bigg( \sum_{j=1}^{2k} \frac{ |\psi^{(j)}(\frac{y-x}{\varT})|}{\lambda^{j}j!} \bigg)^{m_{1}+\cdots+m_{2k}} \Bigg) \\
  & \overset{\2}{\leq} 2M(\psi,2N,\lambda) + 2\sum_{k=1}^{N} \frac{\sigma^{2k} (2k)^{2k} \lambda^{2k} }{ \varT^{2k}} \sum_{m_{1},\dots,m_{2k}} M(\psi,2N,\lambda) \leq 4M(\psi,2N,\lambda),
\end{align*}
where in step $\1$ we use Ineq.~\eqref{ineq-up-derivative-vabs} and Eq.~\eqref{base-measure-nongaussian}, in step $\2$ we use definition~\eqref{def-M-nongaussian} so that
\begin{align*}
    &e^{ -\psi(\frac{y-x}{\varT}) + \psi(\frac{y-x}{2\varT})} \leq M(\psi,2N,\lambda), \quad \text{and}\\
    &e^{ -\psi(\frac{y-x}{\varT}) + \psi(\frac{y-x}{2\varT})}\bigg( \sum_{j=1}^{2k} \frac{ |\psi^{(j)}(\frac{y-x}{\varT})|}{\lambda^{j}j!} \bigg)^{m_{1}+\cdots+m_{2k}} \leq M(\psi,2k,\lambda) \leq M(\psi,2N,\lambda).
\end{align*}
In the last step, we use the fact that $\varT \geq 8e\sigma\lambda N$, and also that there are at most $e^{2k}$ many $2k$-tuples of nonnegative integers $(m_{1},\cdots,m_{2k})$ satisfying $\sum_{j=1}^{2k}jm_{j} = 2k$. Therefore, the reduction $\mathsf{K}(X_{\theta})$ in Eq.~\eqref{reduction-nongaussian} satisfies the conditions required by $\textsc{rk}$.

\paragraph{Computational complexity of the reduction.}
Note that from~\citet[page 10]{lou2025computationally}, the reduction procedure defined in~\eqref{reduction-nongaussian} requires, at each iteration, generating a sample form $\mathsf{Unif}([0,1])$ and a sample from base measures~\eqref{base-measure-nongaussian}, which takes time $\mathcal{O}(T_{\mathsf{samp}})$ by Assumption~\ref{assump:log-density}(c). It also requires evaluating $\mathcal{S}_{N}^{\star}$ once per iteration, which takes time $\mathcal{O}(T_{\mathsf{eval}}(N))$ by Assumption~\ref{assump:Teval}. Moreover, the total number of iterations is $T_{\mathsf{iter}}$. Consequently, the reduction $\mathsf{K}$~\eqref{reduction-nongaussian} has computational complexity
\[
    \mathcal{O}\Big( T_{\mathsf{iter}} \cdot T_{\mathsf{samp}} \cdot T_{\mathsf{eval}}(N)  \Big).
\]

\paragraph{Proof of TV deficiency guarantee of the reduction.} We now turn to prove Eq.~\eqref{ineq-nongaussian-term2}. Applying the triangle inequality as in Ineq.~\eqref{ineq:triangle-TV-KtoV} yields
\begin{align}\label{ineq:triangle-TV-KtoV-nongaussian}
    \sup_{\theta \in \real}\; \mathsf{d}_{\mathsf{TV}} \big( \mathsf{K}(X_{\theta}), Y_{\theta} \big) \leq \frac{1}{2}\sup_{\theta \in \real} \bigg\| \mathsf{K}(X_{\theta}) - \int_{\real} \mathcal{T}_{N}^{\star}(\cdot \mid x) \cdot u(x;\theta) \mathrm{d}x  \bigg\|_{1} + \delta\big(\mathcal{U},\mathcal{V};\mathcal{T}_{N}^{\star} \big).
\end{align}
We next turn to bound the two terms in the RHS of the inequality in the display. Applying Ineq.~\eqref{ineq-nongaussian-term1} and using $2N+1 \geq \big(2e\log\big(3\CdensityB/\epsilon\big) \big)^{2(1+\Cdensity)}$ yield
\begin{align}\label{term2-TV-guarantee-nongaussian}
    \delta\big(\mathcal{U},\mathcal{V};\mathcal{T}_{N}^{\star} \big) \leq \frac{\epsilon}{3} + 4\err(\psi,2N,\lambda).
\end{align}
Towards bounding the first term of the RHS of Ineq.~\eqref{ineq:triangle-TV-KtoV-nongaussian}, we let $Y = \textsc{rk}\big(x, T_{\mathsf{iter}}, 4 M(\psi,2N,\lambda), y_{0} \big)$ for each $x\in\real$.
Then from~\citet[Lemma 2]{lou2025computationally}, we obtain that the random variable $Y$ is drawn from the mixture distribution
\begin{align}\label{condition-density-nongaussian}
  f_{Y}(y \mid x) = \frac{\mathcal{S}_{N}^{\star}(y\mid x) \vee 0}{p(x)} \cdot \big( 1 - \ell(x) \big) + \delta_{y_{0}}(y) \cdot \ell(x), \quad \text{for all } y \in \real,
\end{align}
where $p(x)$ is defined in Eq.~\eqref{px-qx-def}, $\delta_{y_{0}}(y)$ is the Dirac delta function centered at $y_{0}$, and 
\begin{align}\label{def-lx-nongaussian}
    \ell(x) =  \bigg(1 - \frac{p(x)}{4 M(\psi,2N,\lambda)} \bigg)^{T_{\mathsf{iter}}}.
\end{align}
By the law of total probability and definition of $\mathsf{K}(X_{\theta})$ in Eq.~\eqref{reduction-nongaussian}, we obtain that the marginal distribution of $\mathsf{K}(X_{\theta})$ is given by
\begin{align}\label{marginal-density-nongaussian}
    f_{\mathsf{K}(X_{\theta})}(y) &= \int_{\real} f_{Y}(y\mid x) \cdot u(x;\theta) \mathrm{d}x \nonumber \\
    & = \int_{\real} \bigg( \mathcal{T}_{N}^{\star}(y \mid x) \big( 1 - \ell(x) \big) + \delta_{y_{0}}(y) \cdot \ell(x) \bigg) \cdot u(x;\theta) \mathrm{d}x,
\end{align}
where we use Eq.~\eqref{condition-density-nongaussian} and also note by definition~\eqref{close-markov-kernel}, $\mathcal{T}_{N}^{\star}(y \mid x) = \frac{\mathcal{S}_{N}^{\star}(y\mid x) \vee 0}{p(x)}$. Continuing, we obtain
\begin{align}\label{term1-abs-guarantee-nongaussian}
  \bigg\| \mathsf{K}(X_{\theta}) - \int_{\real} \mathcal{T}_{N}^{\star}(\cdot \mid x) \cdot u(x;\theta) \mathrm{d}x  \bigg\|_{1} &= \int_{\real} \bigg| f_{\mathsf{K}(X_{\theta})}(y) - \int_{\real} \mathcal{T}_{N}^{\star}(y \mid x) \cdot u(x;\theta) \mathrm{d}x  \bigg| \mathrm{d}y \nonumber \\
  & \overset{\1}{=} \int_{\real} \bigg| \int_{\real} \Big(- \ell(x) \cdot \mathcal{T}_{N}^{\star}(y \mid x) +   \delta_{y_{0}}(y) \cdot \ell(x) \Big)u(x;\theta)  \mathrm{d}x \bigg| \mathrm{d}y \nonumber \\
  & \overset{\2}{\leq} 2\int_{\real} |\ell(x)| u(x;\theta) \mathrm{d}x \nonumber \\
  & \overset{\3}{\leq} 2\exp\bigg( - \frac{ T_{\mathsf{iter}}   }{ 4M(\psi,2N,\lambda) } \bigg),
\end{align}
where in step $\1$ we use Eq.~\eqref{marginal-density-nongaussian}, and in step $\2$ we apply triangle inequality and use
\[
    \int_{\real} \mathcal{T}_{N}^{\star}(y \mid x) \mathrm{d}y = 1 \quad \text{and} \quad \int_{\real} \delta_{y_{0}}(y) \mathrm{d}y = 1.
\]
In step $\3$, we use Ineq.~\eqref{ineq:px-geq-1}: $p(x) \geq 1$ for all $x \in \real$, so that from Eq.~\eqref{def-lx-nongaussian},
\[
  \ell(x) \leq \exp\bigg( - \frac{ T_{\mathsf{iter}} \cdot p(x)}{ 4M(\psi,2N,\lambda) } \bigg) \leq \exp\bigg( - \frac{ T_{\mathsf{iter}}  }{ 4M(\psi,2N,\lambda) } \bigg)
\]
Substituting Ineq.~\eqref{term1-abs-guarantee-nongaussian} and Ineq.~\eqref{term2-TV-guarantee-nongaussian} into the RHS of Ineq.~\eqref{ineq:triangle-TV-KtoV-nongaussian} completes the proof of guarantee~\eqref{ineq-nongaussian-term2}.
\qed

\subsection{Proof of Corollary~\ref{example:bounded-psi}}\label{sec:pf_corollary_bounded_psi}
Note that since $\psi \in \mathcal{F}(R)$~\eqref{def:analytical-functions}, we immediately obtain that Assumption~\ref{assump:log-density}(b) holds with $\CdensityT = R$ and $\Cdensity = 0$. Similarly, we obtain $\Set(\psi,n,R) = \real$ and $\Set(\psi,n,R) = \emptyset$ for all $n\in \mathbb{N}$. Consequently, we have
\[
    \err(\psi,n,R) = 0 \quad \text{and} \quad M(\psi,n,R) \leq 2 \sup_{t \in\real} \Big\{ e^{-\psi(t)+\psi(t/2)} \Big\} = \mathcal{O}(1) \quad \text{for all } n \in \mathbb{N}.
\]
Applying Ineq.~\eqref{ineq-nongaussian-term2}, by setting $T_{\mathsf{iter}} = \mathcal{O}\big( \log(3/\epsilon) \big)$, we obtain $\sup_{\theta \in \Theta}\; \mathsf{d}_{\mathsf{TV}} \big( \mathsf{K}(X_{\theta}), Y_{\theta} \big) \leq \epsilon$. Note that eventually we take 
\[
    \lambda = R, \quad N = \mathcal{O}\big( \log^{2}(C/\epsilon) \big), \quad \text{and} \quad \varT = \mathcal{O}\big(\sigma R \log^{2}(C/\epsilon) \big),
\]
so that conditions~\eqref{assump:thm-nongaussian} and~\eqref{parameter-N-varT-nongaussian} are satisfied. 

Finally, we bound the computational complexity of the reduction. Note that we can use Eq.~\eqref{eq-target-derivative-nongaussian-0} to evaluate the signed kernel $\mathcal{S}_{N}^{\star}$~\eqref{truncated-signed-kernel}. There are at most $\mathcal{O}\big(e^{\pi\sqrt{n}} \big)$ many terms in the summation of Eq.~\eqref{eq-target-derivative-nongaussian-0} and each term in the summation can be computed in time $\mathcal{O}\big(e^{C\sqrt{n}} \big)$. Also, evaluating $\mathcal{S}_{N}^{\star}$~\eqref{truncated-signed-kernel} only requires computing $\nabla_{x}^{(n)}v(y;x)$ for $1\leq n \leq \mathcal{O}\big(\log^{2}(C/\epsilon)\big)$. Thus, we obtain
\[
    T_{\mathsf{eval}}(N) = \mathcal{O}\Big( e^{C'\sqrt{N}} \Big) = \mathsf{poly}(1/\epsilon).
\]
Combining with $T_{\mathsf{iter}} = \mathcal{O}\big( \log(3/\epsilon) \big)$ yields that the overall computational complexity is of the order $T_{\mathsf{samp}} \cdot \mathsf{poly}(1/\epsilon)$. 
\qed

\subsection{Proof of claims in Example~\ref{examples:nongaussian}} 

We provide proofs for the different families of distributions separately.

\paragraph{Student's t-distribution.} Note that $\psi(x) = \frac{\nu+1}{2} h(g(x))$ for $h(x) = \log(x)$ and $g(x) = 1+x^{2}/\nu$. Applying Faà di Bruno's formula yields
\begin{align*}
    \psi^{(n)}(x) &= \frac{\nu+1}{2} \sum_{m_{1}+2m_{2} = n} \frac{n!}{m_{1}! m_{2}!} h^{(m_{1}+m_{2})}\big( g(x) \big) \prod_{j=1}^{2} \bigg(\frac{g^{(j)}(x)}{j!}\bigg)^{m_{j}} \\
    &= \frac{\nu+1}{2} \cdot n! \cdot \sum_{m_{1}+2m_{2} = n} \frac{(-1)^{m_{1}+m_{2}-1}(m_{1}+m_{2}-1)!}{m_{1}! m_{2}!} \frac{(2x/\nu)^{m_{1}} (1/\nu)^{m_{2}} }{(1+x^{2}/\nu)^{m_{1}+m_{2}-1}},
\end{align*}
where the summation is over all nonnegative integers $m_{1}$ and $m_{2}$ such that $m_{1}+2m_{2} = n$. Using the above closed form equation, we can evaluate $\psi^{(n)}(x)$ in time $\mathcal{O}\big(\mathsf{poly}(n)\big)$. Applying the triangle inequality yields
\begin{align*}
   \big| \psi^{(n)}(x) \big| &\leq \frac{\nu+1}{2} \cdot n! \cdot \sum_{m_{1}+2m_{2} = n} \frac{(m_{1}+m_{2})!}{m_{1}! m_{2}!} \frac{(2|x|/\nu)^{m_{1}} (1/\nu)^{m_{2}} }{(1+x^{2}/\nu)^{m_{1}+m_{2}}} \\
   &\overset{\1}{\leq} \frac{\nu+1}{2} \cdot n! \cdot \sum_{m_{1}+2m_{2} = n} \bigg( \frac{ \frac{2|x|+1}{\nu} }{1+x^{2}/\nu} \bigg)^{m_{1}+m_{2}} \\
   & \leq \frac{\nu+1}{2} \cdot n! \cdot (3e)^{n},
\end{align*}
where in step $\1$ we use the multinomial theorem. In step $\2$ we use 
\[
    \frac{ \frac{2|x|+1}{\nu} }{1+x^{2}/\nu} \leq \frac{ x^{2}+2 }{\nu+x^{2}} \leq 1+ \frac{2}{\nu} \leq 3,
\]
the fact that $m_{1}+m_{2} \leq n$ and that there are at most $e^n$ terms in the summation. Consequently, by letting $R = 6e(\nu+1)$, we obtain
\[
    \sum_{j=1}^{n} \frac{|\psi^{(j)}(x)|}{R^{j}j!} \leq \sum_{j=1}^{n}\frac{1}{2^{j}} \leq 1 \quad \text{for all} \quad x\in \real, n\in \mathbb{N}.
\]
Thus, $\psi \in \mathcal{F}(6e(\nu+1))$. Note that for all $t\in \real$, we have
\[
    \psi(t)- \psi(t/2) = \frac{\nu+1}{2} \log\bigg( \frac{ 1+t^{2}/v}{1+t^{2}/4v} \bigg) \geq \frac{1}{2} \log(5/4).
\]
Consequently, $\sup_{t\in \real} \big\{ e^{-\psi(t) + \psi(t/2)} \big\} \leq \sqrt{4/5} \leq 1$. This completes the proof for Student's t-distribution. \hfill $\clubsuit$

\paragraph{Hyperbolic secant distribution.} By letting $h(t) = \log(t)$ and $g(t) = e^{\frac{\pi t}{2}} + e^{-\frac{\pi t}{2}}$, we have $\psi(t) = h(g(t))$. Now, applying Faà di Bruno's formula yields
\begin{align*}
  \big| \psi^{(n)}(t)\big| &\leq \sum_{m_{1},\dots,m_{n}} \frac{n!}{(m_{1}!)\cdots (m_{n}!)} \cdot \big| h^{(m_{1}+\cdots +m_{n})} (g(t))\big| \cdot \prod_{j=1}^{n} \bigg( \frac{\big|g^{(j)}(t)\big|}{j!} \bigg)^{m_{j}} \\
  & \leq n! \sum_{m_{1},\dots,m_{n}} \frac{(m_{1}+\cdots +m_{n})!}{(m_{1}!)\cdots (m_{n}!)} \frac{1}{(g(t))^{m_{1} + \cdots + m_{n}}} \prod_{j=1}^{n} \bigg( \frac{ (\frac{\pi}{2})^{j} g(t) }{j!} \bigg)^{m_{j}} \\
  & \leq n! \cdot (\pi/2)^{n} \sum_{m_{1},\dots,m_{n}} \frac{(m_{1}+\cdots +m_{n})!}{(m_{1}!)\cdots (m_{n}!)}  \prod_{j=1}^{n} \bigg( \frac{ 1 }{j!} \bigg)^{m_{j}} \\
  & \leq n! \cdot (\pi/2)^{n} \sum_{m_{1},\dots,m_{n}} \bigg( \sum_{j=1}^{n} \frac{ 1 }{j!} \bigg)^{m_{1}+\cdots + m_{n}} \leq n! \cdot (\pi/2)^{n} \cdot (2e)^{n}.
\end{align*}
Consequently, for $R = 2\pi e$ and uniformly for all $n \in \mathbb{N}$ and $t\in \real$, we have
\[
  \sum_{j=1}^{n}\frac{|\psi^{(j)}(t)|}{(R)^{j}j!} \leq \sum_{j=1}^{n} \frac{1}{2^{j}} \leq 1.
\]
Thus, $\psi \in \mathcal{F}(2\pi e)$. Finally, note that for all $t \in \real$, we have
\[
    \psi(t) - \psi(t/2) = \log\bigg( \frac{e^{\frac{\pi t}{2}} + e^{-\frac{\pi t}{2}}}{ e^{\frac{\pi t}{4}} + e^{-\frac{\pi t}{4}}  } \bigg) \geq \log(1/2).
\]
Consequently, we obtain $\sup_{t\in \real} \big\{ e^{-\psi(t) + \psi(t/2)} \big\} \leq  2$. One can use Faà di Bruno's formula to evaluate $\psi^{(n)}(t)$, since there are at must $\mathcal{O}(e^{\pi \sqrt{n}})$ terms in the summation and each term can be computed in time $\mathcal{O}(\mathsf{poly}(n))$. Thus, the total computational time is $\mathcal{O}(e^{5\sqrt{n}})$. \hfill $\clubsuit$

The proof for the logistic distribution is identical so we omit it.

\subsection{Proof of Corollary~\ref{example:generalized-normal}}\label{sec:pf_corollary_generalized_normal}
We must verify the assumption, bound the error quantities, and prove the consequence; we do so separately.
\paragraph{Verifying Assumption~\ref{assump:log-density}.} By definition, we have
\begin{align*}
\psi^{(j)}(t) = 
\begin{cases} 
\beta! \cdot t^{\beta-j}/(\beta-j)! \;\; &\text{ if } 1 \leq j \leq \beta \\
0 \;\; &\text{ if } j \geq \beta + 1. 
\end{cases}
\end{align*}
We thus obtain for all $n \in \mathbb{N}$ and $k\in [n]$ that
\begin{align*}
      \bigg(\sum_{j=1}^{n} \frac{|\psi^{(j)}(t)|}{\big(\CdensityT \big)^{j} j! }\bigg)^{k} \leq \bigg(\sum_{j=1}^{\beta} \frac{ \binom{\beta}{j} |t|^{\beta-j} }{\big(\CdensityT \big)^{j}  }\bigg)^{k} \leq \beta^{k-1} \sum_{j=1}^{\beta} \bigg( \frac{\beta}{\CdensityT} \bigg)^{jk} |t|^{(\beta-j)k} \leq \frac{1}{2\beta} \sum_{j=1}^{\beta} |t|^{(\beta-j)k},
\end{align*}
where in the last step we set $\CdensityT = 2\beta^{2}$. Integrating both sides over $t\in\real$ yields
\begin{align*}
    \int_{\real} e^{-\psi(t)} \cdot \bigg(\sum_{j=1}^{n} \frac{|\psi^{(j)}(t)|}{\big(\CdensityT \big)^{j} j! }\bigg)^{k} \mathrm{d}t \leq \frac{1}{2\beta} \sum_{j=1}^{\beta} \int_{\real} e^{-\psi(t)} |t|^{(\beta-j)k} \mathrm{d}t \overset{\1}{\leq} k! \leq n!.
\end{align*}
In step $\1$ we use that for any $j \in [\beta]$,
\begin{align*}
  \int_{\real} e^{-\psi(t)} |t|^{(\beta-j)k} \mathrm{d}t &= \frac{\beta}{\Gamma(1/\beta)} \int_{0}^{+\infty} e^{-t^{\beta}} t^{(\beta-j)k} \mathrm{d}t \\
  & \leq  \frac{\beta}{\Gamma(1/\beta)} \int_{0}^{+\infty} e^{-t^{\beta}}  \mathrm{d}t +  \frac{\beta}{\Gamma(1/\beta)} \int_{0}^{+\infty} e^{-t^{\beta}} t^{\beta k}  \mathrm{d}t \\
  &= \frac{1}{\Gamma(1/\beta)} \int_{0}^{+\infty} e^{-x} x^{1/\beta-1} \mathrm{d}x + \frac{1}{\Gamma(1/\beta)} \int_{0}^{+\infty} e^{-x} x^{k+1/\beta-1} \mathrm{d}x  \\
  &= 1+ \frac{\Gamma(k + 1/\beta)}{\Gamma(1/\beta)} = 1+ k! \leq 2 \cdot (k!)
\end{align*}
Thus, Assumption~\ref{assump:log-density} holds with $\CdensityT = 2\beta^{2}$ and $\Cdensity = 1$.

\paragraph{Bounding $\err(\psi,n,\lambda)$ and $M(\psi,n,\lambda)$.} Note that for $|t| \leq \lambda^{\frac{1}{2\beta}}$ and $\lambda \geq 4\beta^{2}$, we have
\begin{align*}
    \sum_{j=1}^{n} \frac{|\psi^{(j)}(t)|}{\lambda^{j} j! } = \sum_{j=1}^{\beta} \frac{\binom{\beta}{j} |t|^{\beta-j}}{\lambda^{j}} \leq \sum_{j=1}^{\beta} \frac{\beta^{j} \lambda^{1/2}}{\lambda^{j}} \leq \sum_{j=1}^{\beta} \frac{1}{2^{j}} \leq 1.
\end{align*}
Consequently, we obtain that for $\lambda \geq 4\beta^{2}$ and for all $n \in \mathbb{N}$,
\[
  \Big\{ t \in \real: |t| \leq \lambda^{\frac{1}{2\beta}} \Big\} \subseteq \Set(\psi,n,\lambda) \quad \text{and} \quad  \Set^{\complement}(\psi,n,\lambda) \subseteq \Big\{ t \in \real: |t| \geq \lambda^{\frac{1}{2\beta}} \Big\}.
\]
We now turn to bound $\err(\psi,n,\lambda)$~\eqref{def-error-nongaussian}. Using the relation in the display above, we obtain
\begin{align}\label{ineq1-tail-general-normal}
  \int_{\Set^{\complement}(\psi,n,\lambda)} e^{-\psi(t)} \mathrm{d}t &\leq  \frac{\beta}{\Gamma(1/\beta)} \int_{\lambda^{\frac{1}{2\beta}}}^{+\infty} e^{-t^{\beta}} \mathrm{d}t
   \nonumber \\
  &= \frac{\beta}{\Gamma(1/\beta)} \int_{0}^{+\infty} e^{-\big(t+ \lambda^{\frac{1}{2\beta}} \big)^{\beta}} \mathrm{d}t \leq \frac{\beta}{\Gamma(1/\beta)} \int_{0}^{+\infty} e^{-t^{\beta}} \mathrm{d}t \cdot e^{-\sqrt{\lambda}} = e^{-\sqrt{\lambda}}.
\end{align}
Continuing, we have that for all $n \in \mathbb{N}$ and $k \in [n]$,
\begin{align}\label{ineq2-tail-general-normal}
  \int_{\Set^{\complement}(\psi,n,\lambda)} e^{-\psi(t)} \bigg( \sum_{j=1}^{\beta} \frac{|\psi^{(j)}(t)|}{\lambda^{j} j! }  \bigg)^{k} \mathrm{d}t &\leq 
  \int_{\Set^{\complement}(\psi,n,\lambda)} e^{-\psi(t)} \bigg( \sum_{j=1}^{\beta} \frac{\beta^{j}|t|^{\beta-j}}{\lambda^{j}}  \bigg)^{k} \mathrm{d}t \notag \\
  & \leq \frac{\beta}{\Gamma(1/\beta)} \int_{\lambda^{\frac{1}{2\beta}}}^{+\infty}   e^{-t^{\beta}} \bigg( \sum_{j=1}^{\beta} \frac{\beta^{j}|t|^{\beta-j}}{\lambda^{j}}  \bigg)^{k} \mathrm{d}t \notag \\
  & \overset{\1}{\leq} \frac{\beta}{\Gamma(1/\beta)} \int_{\lambda^{\frac{1}{2\beta}}}^{+\infty}   e^{-t^{\beta}} \frac{t^{\beta k}}{\lambda^{k/2}} \mathrm{d}t \notag \\
  & \overset{\2}{\leq} e^{-\sqrt{\lambda}} \int_{0}^{+\infty}  \frac{\beta}{\Gamma(1/\beta)} e^{-\frac{t^{\beta}}{2}} \mathrm{d}t \leq 2e^{-\sqrt{\lambda}}.
\end{align}
Step $\1$ above holds for all $\lambda \geq 4\beta^{2}$. In this case, for $|t| \geq \lambda^{\frac{1}{2\beta}} \geq 1$, we have
$
  \sum_{j=1}^{\beta} \frac{\beta^{j}|t|^{-j}}{\lambda^{j/2}} \leq \sum_{j=1}^{\beta} \frac{1}{2^{j}} \leq 1.
$
Step $\2$ above holds for all $\lambda \geq 4k^{2}$; in this case
\begin{align*}
    \int_{\lambda^{\frac{1}{2\beta}}}^{+\infty}   e^{-t^{\beta}} t^{\beta k} \mathrm{d}t &\leq \sup_{t\geq \lambda^{\frac{1}{2\beta}}}\Big\{ t^{\beta k} e^{-t^{\beta}/2}  \Big\} \cdot  \int_{\lambda^{\frac{1}{2\beta}}}^{+\infty}  e^{-t^{\beta}/2} \mathrm{d}t \\
    & \leq \lambda^{k/2} e^{-\sqrt{\lambda}/2} \cdot \int_{0}^{+\infty} e^{-t^{\beta}/2} \mathrm{d}t \cdot e^{-\sqrt{\lambda}/2}  = \lambda^{k/2} e^{-\sqrt{\lambda}} \cdot \int_{0}^{+\infty} e^{-t^{\beta}/2} \mathrm{d}t.
\end{align*}
Combining Ineq.~\eqref{ineq1-tail-general-normal} and Ineq.~\eqref{ineq2-tail-general-normal} yields that 
\[
    \err(\psi,n,\lambda) \leq  3e^{-\sqrt{\lambda}} \quad \text{for all } \lambda \geq 4 n^{2} \vee 4\beta^{2} \text{ and } n \in \mathbb{N}.
\]    
We next bound $M(\psi,n,\lambda)$~\eqref{def-M-nongaussian}. Since $-\psi(t)+\psi(t/2) \leq -t^{\beta}/2$, we obtain
\begin{align*}
    \sup_{t \in \real} \bigg\{ e^{-\psi(t)+\psi(t/2)} \bigg( \sum_{j=1}^{n} \frac{|\psi^{(j)}(t)|}{\lambda^{j}j! }\bigg)^{k} \bigg\} &\leq \beta^{k} \sum_{j=1}^{\beta}  \frac{\beta^{jk}}{\lambda^{jk}} \cdot \sup_{t \in \real} \bigg\{ e^{-\frac{t^{\beta}}{2}} \cdot |t|^{(\beta-j)k} \bigg\} \\
    & \overset{\1}{\leq} \beta^{k} k^{k} \sum_{j=1}^{\beta} \bigg(\frac{\beta}{\lambda}\bigg)^{jk} \overset{\2}{\leq} \frac{\beta^{2k}k^{k} }{\lambda^{k}} \sum_{j=1}^{\beta}\frac{1}{2^{j-1}} \overset{\3}{\leq} 2.
\end{align*} 
In step $\1$ above, we use the chain of bounds
\[
    \sup_{t \in \real} \bigg\{ e^{-\frac{t^{\beta}}{2}} \cdot |t|^{(\beta-j)k} \bigg\} \leq 1 \vee \sup_{t>1} \Big\{ e^{-\frac{t^{\beta}}{2}} \cdot |t|^{\beta k} \Big\} \leq 1 \vee k^{k} \leq k^{k}.
\]     
Step $\2$ above holds for all $\lambda \geq \beta^{2}$, since $\beta\geq 2$. Finally, step $\3$ holds for all $\lambda \geq \beta^2 n$. Putting together the pieces, we obtain by definition~\eqref{def-M-nongaussian} that
\[
    M(\psi,n,\lambda) \leq 3, \quad \text{for all } n\in \mathbb{N} \text{ and } \lambda \geq \beta^{2}n.
\]

\paragraph{Proof of the consequence.} Note that since $\Cdensity = 1$, the requirement on $N$ from conditions~\eqref{assump:thm-nongaussian} and~\eqref{parameter-N-varT-nongaussian} is $N \geq C \log^{4}(C/\epsilon)$ for a universal and positive constant $C$. It suffices to set $\lambda \geq 4 \beta^2 N$ to satisfy the conditions from the previous parts. We thus have 
\[
    \err(\psi,2N,\lambda) \leq  \frac{\epsilon}{6}\quad \text{and} \quad M(\psi,2N,\lambda) \leq 3, \quad \text{if} \quad \lambda \gtrsim \beta^{2} \log^{8}(C/\epsilon).
\]  
Consequently, by letting $T_{\mathsf{iter}} \gtrsim \log(3/\epsilon)$, we obtain from guarantee~\eqref{ineq-nongaussian-term2} that
\[
    \sup_{\theta \in \real} \mathsf{d}_{\mathsf{TV}}\; (\mathsf{K}(X_{\theta}),Y_{\theta}) \leq \epsilon.
\]
Note that the eventual requirement on $\varT$ from condition~\eqref{assump:thm-nongaussian} reads as $\varT \gtrsim \sigma \beta^{2} \log^{12}(C/\epsilon)$. 

We now turn to the computational complexity of the reduction. Note that we can use Eq.~\eqref{eq-target-derivative-nongaussian-0} to evaluate the signed kernel $\mathcal{S}_{N}^{\star}$~\eqref{truncated-signed-kernel}. Since $\psi^{(j)}(t) = 0$ for $j\geq \beta+1$, then the summation in Eq.~\eqref{eq-target-derivative-nongaussian-0} is over all tuples of nonnegative integers $(m_{1},\cdots,m_{\beta})$ satisfying $\sum_{j=1}^{\beta} jm_{j} = n$. There are at most $n^{\beta}$ such tuples. Thus, evaluating $\nabla_{x}^{(n)}v(y;x)$ takes time $\mathcal{O}(n^{\beta+1})$. Since in addition we have $N =\mathcal{O}( \log^{4}(C/\epsilon))$, we can bound the evaluation time  $T_{\mathsf{eval}}(N) = \mathcal{O}(\log^{4\beta+4}(C/\epsilon))$. We view the sampling time as a constant ($T_{\mathsf{samp}} = \mathcal{O}(1)$), since the generalized normal distribution is straightforward to sample from. Putting this together with number of iterations $T_{\mathsf{iter}} \asymp \log(3/\epsilon)$ yields a total computational complexity of $\mathcal{O}\Big(\log^{4\beta+5}(C/\epsilon)\Big)$.
\qed


\section{Proofs of results from Section~\ref{sec:gaussian_target_nonlinear}} \label{sub:pf_thm_gaussian_mean_nonlinear}
We prove Theorem~\ref{thm:gaussian-mean-nonlinear}(a) in Section~\ref{sec:pf_gaussian_nonlinear_term0} and Theorem~\ref{thm:gaussian-mean-nonlinear}(b) in Section~\ref{sec:pf_gaussian_nonlinear_term3}.
We provide the proof of the corollaries and concrete examples of Theorem~\ref{thm:gaussian-mean-nonlinear} in Sections~\ref{sec:pf_corollary_monomial}--\ref{sec:pf_concrete_examples}.

\subsection{Proof of Theorem~\ref{thm:gaussian-mean-nonlinear}(a)}\label{sec:pf_gaussian_nonlinear_term0}
Throughout this section, to reduce the notational burden, we use the shorthand:
\begin{align}\label{shorthand-h-x}
  h_{j}(x) = \frac{  \link^{(j)}(x/\varT) }{ \big(\ClinkT \big)^{j} j!},\quad \phi(x) = \frac{e^{-x^{2}}}{\sqrt{2\pi}}, \quad g(x) = \frac{f(x/\varT) - y}{2}, \quad\text{for all } x,y\in \real, j\in \PosInt.
\end{align} 
We first state an auxiliary lemma, whose proof we defer to Appendix~\ref{sec:pf_lemma_aux_derivative_bound_nonlinear}.
\begin{lemma}\label{lemma:aux-derivative-bound-nonlinear}
Let $v(y;\theta)$ be defined as in Eq.~\eqref{target:density}. Let $h_{j}(x)$, $\phi(x)$, and $g(x)$ be defined in Eq.~\eqref{shorthand-h-x}.
Then we have for all $n \in \PosInt$ and $x,y\in \real$,
\small
\begin{align}\label{eq:derivative-gaussian-nonlinear}
    \nabla_{\theta}^{(n)} v(y;\theta) \big|_{\theta = x} = \frac{(\ClinkT)^{n}}{\varT^{n}} \cdot n! \cdot \sum_{m_{1}, \dots,m_{n}} \frac{(-1)^{m_{1} + \cdots+m_{n}} \cdot \phi(g(x)) \cdot H_{m_{1}+\cdots+m_{n}}(g(x)) }{m_{1}!\,m_{2}! \cdots m_{n}!}   \cdot \prod_{j=1}^{n} \bigg( \frac{h_{j}(x) }{2} \bigg)^{m_{j}},
\end{align}
\normalsize
where the summation is over all $n$-tuples of nonnegative integers $(m_{1},\cdots,m_{n})$ satisfying \mbox{$\sum_{j=1}^{n}jm_{j} = n$}. Moreover, we have
\begin{align}\label{ineq:derivative-gaussian-nonlinear}
    \Big|  \nabla_{\theta}^{(n)} v(y;\theta) \big|_{\theta = x} \Big| \leq  e^{-\frac{g(x)^{2}}{2}} \cdot \frac{(2\ClinkT)^{n}}{\varT^{n}} \cdot n! \cdot e^{n} \cdot \bigg[ \bigg( \sum_{j=1}^{n} |h_{j}(x)| \bigg)^{n} \vee 1 \bigg].
\end{align}
\end{lemma}

We now use the Lemma~\ref{lemma:aux-derivative-bound-nonlinear} to verify that the target density $v(y;\theta)$ defined in Eq.~\eqref{target:density} satisfies Assumption~\ref{assump_1:target}, ensuring that Lemma~\ref{lemma:TV-deficiency-decompose} is applicable. 

Part (a) of the assumption immediately holds by definition of $v(y;\theta)$~\eqref{target:density}. 
We now verify Assumption~\ref{assump_1:target}(b). Ineq.~\eqref{ineq:derivative-gaussian-nonlinear} implies that
\begin{align*}
  \int_{\real} \big| \nabla_{x}^{(n)} v(y;x) \big| e^{-\frac{(x-\theta)^{2}}{2}} \mathrm{d}x &\leq \frac{(2e\ClinkT)^{n}}{\varT^{n}} \cdot n!  \cdot \bigg[ \int_{\real} \bigg( \sum_{j=1}^{n} |h_{j}(x)| \bigg)^{n} e^{-\frac{(x-\theta)^{2}}{2}} \mathrm{d}x + \int_{\real}  e^{-\frac{(x-\theta)^{2}}{2}} \mathrm{d}x \bigg] \\
  & \overset{\1}{=}  \frac{(2e\ClinkT)^{n}}{\varT^{n}} \cdot n!  \cdot \bigg[ \sqrt{2}\int_{\real} \bigg( \sum_{j=1}^{n} \frac{|f^{(j)}((\theta+\sqrt{2}t )/\varT)|}{(\ClinkT)^{j}j!} \bigg)^{n} e^{-t^{2}} \mathrm{d}t + \sqrt{2\pi} \bigg] \\
  & \overset{\2}{\leq} \frac{(2e\ClinkT)^{n}}{\varT^{n}} \cdot n!  \cdot \Big( 2\sqrt{\pi} \cdot  (\Clink n)! + \sqrt{2\pi} \Big) < +\infty,
\end{align*}
where in step $\1$ we use the change of variable $t = (x-\theta)/\sqrt{2}$, and in step $\2$ we use Assumption~\ref{assump:link}. Note that $\big| \nabla_{x}^{(n)} v(y;x) \big| e^{-\frac{(x-\theta)^{2}}{2}}$ is a continuous and nonnegative function on $\real$, and its integral is bounded. This means it converges to zero at infinity, that is, 
\[
    \lim_{x \uparrow +\infty}  \big[ \nabla_{x}^{(n)} v(y;x) \big] e^{-\frac{(x-\theta)^{2}}{2}} = 
    \lim_{x \downarrow -\infty}  \big[ \nabla_{x}^{(n)} v(y;x) \big] e^{-\frac{(x-\theta)^{2}}{2}} = 0.
\]
This proves Assumption~\ref{assump_1:target}(b).
We now verify Assumption~\ref{assump_1:target}(c). Note that by definition~\eqref{target:density}, $v(y;\sqrt{2} x+ \theta) \leq 1/\sqrt{\pi}$ for all $x,y,\theta \in \real$ , and consequently
\[
    \int_{\real} \Big( v\big(y;\sqrt{2} x+\theta \big) \Big)^{2} \exp(-x^{2}) \mathrm{d}x  \leq \int_{\real} e^{-x^{2}}/\sqrt{\pi} \mathrm{d}x = 1 <+\infty.
\]
This proves that Assumption~\ref{assump_1:target}(c) holds.

Consequently, Lemma~\ref{lemma:TV-deficiency-decompose} holds and we obtain
\begin{align}\label{ineq-pf-nonlinear-term0}
  \delta\big(\mathcal{U}, \mathcal{V};\mathcal{T}^{\star}_{N}\big) &\leq \SignedDef\big( \mathcal{U}, \mathcal{V}; \mathcal{S}_{N}^{\star} \big)  + \frac{1}{2}\sup_{\theta\in \Theta}\; \int_{\real} \big(|p(x) - 1| + q(x)\big) \cdot u(x;\theta) \mathrm{d}x. 
\end{align}
To bound the terms appearing on the RHS above, we state an auxiliary lemma, whose proof we defer to the following subsections.
\begin{lemma}\label{lemma:aux-pf-thm2-part1}
Let the source and target statistical models $\mathcal{U}$ and $\mathcal{V}$ be given by Eq.~\eqref{source-gaussian-location} (with $\sigma = 1$) and Eq.~\eqref{target:Gaussian}, respectively. Let the signed kernel $\mathcal{S}_{N}^{\star}$ be defined in Eq.~\eqref{truncated-signed-kernel} and let $\SignedDef\big( \mathcal{U}, \mathcal{V}; \mathcal{S}_{N}^{\star} \big)$ be defined in Eq.~\eqref{deficiency-signed-kernel}. Suppose the pair of parameters $(N,\varT)$ satisfies condition~\eqref{assump:thm-nonlinear} and Assumption~\ref{assump:link} holds. Let $\ClinkB$ be defined in Eq.~\eqref{def-quantity-link}. Then the following statements hold:
\begin{enumerate} 
\item[(a)] We have
\begin{align}\label{ineq-gaussian-nonlinear-term1}
  \SignedDef\big( \mathcal{U}, \mathcal{V}; \mathcal{S}_{N}^{\star} \big) \leq \ClinkB \exp\bigg( - (2e)^{-1} (2N+1)^{\frac{1}{2(1+\Clink)}} \bigg).
\end{align}
\item[(b)] Let $p(x)$ and $q(x)$ be defined in Eq.~\eqref{px-qx-def}. Then
\begin{align} 
\label{ineq-gaussian-nonlinear-term2}
    \frac{1}{2}\sup_{\theta\in \real}\; \int_{\real} \big(|p(x) - 1| + q(x)\big) \cdot u(x;\theta) \mathrm{d}x \leq 12\exp(-\varT/8) + 15\sup_{\theta \in \Theta} \; \exp\bigg( \frac{\theta^{2}}{2} -\frac{\varT^{2}}{8} \bigg).
\end{align}
\end{enumerate}
\end{lemma}

Substituting Ineq.~\eqref{ineq-gaussian-nonlinear-term1} and Ineq.~\eqref{ineq-gaussian-nonlinear-term2} into the RHS of Ineq.~\eqref{ineq-pf-nonlinear-term0} proves part (a) of the theorem.
\qed

We next prove Ineq.~\eqref{ineq-gaussian-nonlinear-term1} in Section~\ref{sec:pf_gaussian_nonlinear_term1}, and prove Ineq.~\eqref{ineq-gaussian-nonlinear-term2} in Section~\ref{sec:pf_gaussian_nonlinear_term2}.

\subsubsection{Proof of Lemma~\ref{lemma:aux-pf-thm2-part1}(a)}\label{sec:pf_gaussian_nonlinear_term1}
The proof is identical to that of Lemma~\ref{lemma:aux-pf-thm1-part1}(a) in Section~\ref{sec:pf_gaussian_nongaussian_term4}. 
We state two key claims. First, for all $\varT \geq 2 \vee 2|\theta| \vee 2 \ClinkT$ and $n,\ell \in \PosInt$ and $y\in \real$, we have
\begin{align}\label{claim1-pf-nonlinear-part1}
    \bigg[ \frac{ \mathrm{d}^{\ell} \big[ v\big(y;\theta + \sqrt{2} x\big)  e^{-x^2} \big] }{\mathrm{d}x^{\ell}} \Nherm_{n+1+\ell}(x) \bigg] \bigg \vert_{-\infty}^{+\infty} = 0,
\end{align}
where the evaluation at infinity is with respect to variable $x$. Second, for all $\varT \geq 8\sqrt{2} e  \ClinkT$ and for all $r,k\in \PosInt$ and $k \leq r$, we have
\begin{align}\label{claim2-pf-nonlinear-part1}
  \int_{\real}\int_{\real} \bigg| \nabla_{x}^{(r-k)} v\big(y;\theta + \sqrt{2} x\big) \cdot  e^{-x^{2}/2}   \bigg| \mathrm{d}x\mathrm{d}y \leq 4 \sqrt{2}\pi (r-k)! \cdot \frac{1}{4^{r-k}} \cdot \big(\Clink(r-k)\big)!.
\end{align}
We prove the two claims above, which are analogous to Claims~\eqref{claim1-pf-nongaussian-part1} and~\eqref{claim2-pf-nongaussian-part1} in Section~\ref{sec:pf_gaussian_nongaussian_term4}, and we omit the remaining steps to prove the lemma, since these are exactly parallel to before.

\paragraph{Proof of Claim~\eqref{claim1-pf-nonlinear-part1}.}
By definition, we have for all $n \in \PosInt$ that
\begin{align*}
  \nabla_{x}^{(n)} v\big(y;\theta+\sqrt{2}x \big) = \nabla_{z}^{(n)} v\big(y;z \big) \big|_{z = \theta+\sqrt{2}x} \cdot (\sqrt{2})^{n}.
\end{align*}
Applying Ineq.~\eqref{ineq:derivative-gaussian-nonlinear}, we obtain
\begin{align}
  \big| \nabla_{x}^{(n)} v\big(y;\theta+\sqrt{2}x \big) \big| &\leq 2^{n/2} \cdot \Big| \nabla_{z}^{(n)} v\big(y;z \big) \big|_{z = \theta+\sqrt{2}x} \Big| \nonumber \\
  &\leq e^{-\frac{g(\theta+\sqrt{2}x)^{2}}{2}} \cdot \frac{(2\sqrt{2}e\ClinkT)^{n}}{\varT^{n}} \cdot n!  \cdot \bigg[ \bigg( \sum_{j=1}^{n} \big|h_{j}\big(\theta+\sqrt{2}x\big) \big| \bigg)^{n} \vee 1 \bigg]  \label{ineq0:derivative-bound-claim-nonlinear} \\
  & \leq  \frac{(2\sqrt{2}e\ClinkT)^{n}}{\varT^{n}} \cdot n!  \cdot \bigg[ \bigg( \sum_{j=1}^{n} \big|h_{j}\big(\theta+\sqrt{2}x\big) \big| \bigg)^{n} \vee 1 \bigg]. \label{ineq:derivative-bound-claim-nonlinear}
\end{align}
Using the product rule~\eqref{eq:pd-rule-nongaussian} and Ineq.~\eqref{Hermite-boundness}, we obtain
\begin{align*}
\bigg| \frac{ \mathrm{d}^{\ell} \big[ v\big(y;\theta + \sqrt{2} x\big)  e^{-x^2} \big] }{\mathrm{d}x^{\ell}} &\Nherm_{n+1+\ell}(x) \bigg|  \leq 
\sum_{k = 0}^{\ell} \binom{\ell}{k} \Big| \nabla_{x}^{(\ell-k)} v\big(y;\theta + \sqrt{2} x\big) \cdot H_{k}(x) e^{-x^{2}} \Big| \\
& \leq \sum_{k = 0}^{\ell} \binom{\ell}{k} \frac{(2\sqrt{2}e\ClinkT)^{\ell-k}}{\varT^{\ell-k}} \cdot (\ell-k)!  \cdot \bigg[ \bigg( \sum_{j=1}^{\ell-k} \big|h_{j}\big(\theta+\sqrt{2}x\big) \big| \bigg)^{\ell-k} \vee 1 \bigg] \cdot \sqrt{\gamma_{k}} \cdot e^{-x^{2}/2},
\end{align*}
where in the last step we have used Ineq.~\eqref{ineq:derivative-bound-claim-nonlinear} and $| H_{k}(x) e^{-x^{2}/2}| \leq \sqrt{\gamma_{k}}$ with $\gamma_{k} = \sqrt{\pi}2^{k} k!$ from Ineq.~\eqref{Hermite-boundness}. Note that from Assumption~\ref{assump:link} and Eq.~\eqref{shorthand-h-x}, we obtain
\begin{align*}
      \int_{\real} \bigg( \sum_{j=1}^{\ell-k} \big|h_{j}\big(\theta+\sqrt{2}x\big) \big| \bigg)^{\ell-k} e^{-x^{2}/2} \mathrm{d}x \leq \sqrt{2\pi}(\Clink(\ell-k))! <+\infty.
\end{align*}
If the integral of a continuous and nonnegative function over $\real$ is finite, then it must converge to zero at infinity. Thus, for all $\ell-k \geq 0$
\begin{align*}
  \lim_{x \rightarrow \infty} \bigg( \sum_{j=1}^{\ell-k} \big|h_{j}\big(\theta+\sqrt{2}x\big) \big| \bigg)^{\ell-k} e^{-x^{2}/2} = 0.
\end{align*}
Combining the pieces yields
\begin{align*}
  \lim_{x \rightarrow \infty} \frac{ \mathrm{d}^{\ell} \big[ v\big(y;\theta + \sqrt{2} x\big)  e^{-x^2} \big] }{\mathrm{d}x^{\ell}} \cdot \Nherm_{n+1+\ell}(x) = 0,
\end{align*}
which proves Claim~\eqref{claim1-pf-nonlinear-part1}.

\paragraph{Proof of Claim~\eqref{claim2-pf-nonlinear-part1}.}
Applying Ineq.~\eqref{ineq0:derivative-bound-claim-nonlinear} yields
\begin{align*}
  &\int_{\real} \int_{\real} \bigg| \nabla_{x}^{(r-k)} v\big(y;\theta + \sqrt{2} x\big) \cdot  e^{-x^{2}/2}  \bigg| \mathrm{d}y \mathrm{d}x \\
  &\leq \frac{\big(2\sqrt{2}e\ClinkT\big)^{r-k}}{\varT^{r-k}} \cdot (r-k)! \cdot \int_{\real} \int_{\real} 
  e^{-\frac{g(\theta+\sqrt{2}x)^{2}}{2}} \mathrm{d}y \cdot \bigg[ \bigg( \sum_{j=1}^{r-k} \big|h_{j}\big(\theta+\sqrt{2}x\big) \big| \bigg)^{r-k} + 1 \bigg]   e^{-x^{2}/2} \mathrm{d}x.
\end{align*}
Recall the function $g(x) = (f(x/\varT) - y)/2$ defined in Eq.~\eqref{shorthand-h-x}. We obtain
\[
    \int_{\real} e^{-\frac{g(\theta+\sqrt{2}x)^{2}}{2}} \mathrm{d}y = \int_{\real} e^{-y^{2}/4} \mathrm{d}y = 2\sqrt{\pi}.
\]
By Assumption~\ref{assump:link} and by definition of $h_{j}(x)$ in Eq.~\eqref{shorthand-h-x}, we obtain
\[
    \int_{\real} \bigg[ \bigg( \sum_{j=1}^{r-k} \big|h_{j}\big(\theta+\sqrt{2}x\big) \big| \bigg)^{r-k} + 1 \bigg]   e^{-x^{2}/2}  \mathrm{d}x \leq \sqrt{2\pi} \big(1 + (\Clink(r-k))! \big) \leq 2\sqrt{2\pi} \cdot (\Clink(r-k))!.
\]
Combining all the pieces and letting $\varT \geq 8\sqrt{2} e\ClinkT$, we have
\[
    \int_{\real} \int_{\real} \bigg| \nabla_{x}^{(r-k)} v\big(y;\theta + \sqrt{2} x\big) \cdot  e^{-x^{2}/2}  \bigg| \mathrm{d}y \mathrm{d}x \leq  \frac{4\sqrt{2}\pi}{4^{r-k}} \cdot  (\Clink(r-k))! \cdot (r-k)!.
\]
This concludes the proof of Claim.~\eqref{claim2-pf-nonlinear-part1}.
\qed

\subsubsection{Proof of Lemma~\ref{lemma:aux-pf-thm2-part1}(b)}\label{sec:pf_gaussian_nonlinear_term2}

We first state an auxiliary lemma and defer its proof to Appendix~\ref{sec:pf_lemma_positivity_nonlinear}.
\begin{lemma}\label{lemma:positivity-nonlinear}
Consider the signed kernel $\mathcal{S}^{\star}_{N}(y\mid x)$ defined in Eq.~\eqref{truncated-signed-kernel} for target density $v(y;\theta)$ defined in Eq.~\eqref{target:density}. For all $n \in \mathbb{N}$, let $\Llink(n)$ be defined as in Eq.~\eqref{def-quantity-link}. Recall the shorthand in Eq.~\eqref{shorthand-h-x}. Suppose Assumption~\ref{assump:link} holds with the pair of constants $(\Clink,\ClinkT)$. Then we have:
\begin{enumerate}
\item[(a)] If $\varT \geq 4 e^{3}  \big( \ClinkT \big)^{2} \big( \Llink(2N) \vee 1 \big) N$, then
\begin{align}\label{ineq2:aux-lemma-nonlinear}
  \mathcal{S}_{N}^{\star}(y\mid x) \geq 0 \quad \text{for all} \quad |x| \leq \varT, \quad |y - \link(x/\varT)| \leq \sqrt{\varT}, \quad N \in \PosInt.
\end{align}
\item[(b)] For all $|x| \leq \tau$ and $k\geq 1$ and $k \in \mathbb{N}$, 
\begin{align}\label{ineq3:px-nonlinear}
    &\int_{f(\frac{x}{\varT}) + \sqrt{\varT}}^{+\infty} \big| \nabla_{x}^{(2k)}v(y;x) \big| \mathrm{d}y + \int_{-\infty}^{f(\frac{x}{\varT}) - \sqrt{\varT}} \big| \nabla_{x}^{(2k)}v(y;x) \big| \mathrm{d}y \nonumber \\
    &\leq 2\sqrt{2\pi} \exp(-\varT/8) \cdot (2k)! \cdot \bigg( \frac{2e\ClinkT}{\varT}\bigg)^{2k} \cdot \big( \Llink(2k) \vee 1 \big)^{2k}.
\end{align}
\item[(c)] For all $x \in \real$ and $n\geq 1$ and $n\in \mathbb{N}$,
\begin{align}\label{ineq4:px-nonlinear}
      \int_{\real}  \nabla_{x}^{(n)}v(y;x) \mathrm{d}y  = 0.
\end{align}
\end{enumerate}
\end{lemma} 
Taking the above lemma as given, let us now prove Lemma~\ref{lemma:aux-pf-thm2-part1}(b).
We split the proof into two cases.

\paragraph{Case 1: bound $|p(x)-1|+q(x)$ for $|x| \leq \varT$.}
By definition~\eqref{px-qx-def}, we obtain for all $x\in \real$
\begin{align*}
  p(x) - q(x) &= \int_{\real} \big[\mathcal{S}_{N}^{\star}(y\mid x) \vee 0 \big] +  \big[\mathcal{S}_{N}^{\star}(y\mid x) \wedge 0 \big] \mathrm{d}y \\
  & = \int_{\real} \mathcal{S}_{N}^{\star}(y\mid x) \mathrm{d}y \\ 
  & = \int_{\real} v(y;x) \mathrm{d}y + \sum_{k=1}^{N} \frac{(-1)^{k}}{(2k)!!} \int_{\real} \nabla_{x}^{(2k)} v(y;x) \mathrm{d}y  \overset{\1}{=} 1 + 0 = 1, 
\end{align*}
where in step $\1$ we use Eq.~\eqref{ineq4:px-nonlinear} to control the second term and definition~\eqref{target:density} to control the first, whereby $\int_{\real} v(y;x) \mathrm{d}y = \int_{\real} e^{-y^2/2}/\sqrt{2\pi} \mathrm{d}y = 1$. We consequently obtain $q(x) = p(x) - 1$ for all $x\in\real$. Moreover, note that since $q(x) \geq 0$ by definition~\eqref{px-qx-def}, we obtain $p(x) \geq 1$ for all $x\in \real$. We thus obtain
\begin{align}\label{px-qx-nonlinear-relation}
    |p(x)-1| + q(x) \leq 2(p(x) -1) \quad \text{for all } x \in \real. 
\end{align}
We next turn to upper bound $p(x)$ when $|x| \leq \varT$. Using Ineq.~\eqref{ineq2:aux-lemma-nonlinear} and definition of $p(x)$~\eqref{px-qx-def}, we have
\begin{align}\label{ineq:px-upper-nonlinear-case1}
  p(x) &\leq \int_{f(\frac{x}{\varT}) - \sqrt{\varT}}^{f(\frac{x}{\varT}) + \sqrt{\varT}} \mathcal{S}_{N}^{\star}(y \mid x) \mathrm{d}y + 
  \int_{f(\frac{x}{\varT}) + \sqrt{\varT}}^{+\infty} \big| \mathcal{S}_{N}^{\star}(y \mid x) \big| \mathrm{d}y + 
  \int_{-\infty}^{f(\frac{x}{\varT}) - \sqrt{\varT}} \big| \mathcal{S}_{N}^{\star}(y \mid x) \big| \mathrm{d}y \nonumber \\
  & \leq \int_{\real} v(y;x) \mathrm{d}y + \sum_{k=1}^{N} \frac{1}{(2k)!!} \int_{f(\frac{x}{\varT}) - \sqrt{\varT}}^{f(\frac{x}{\varT}) + \sqrt{\varT}} \nabla_{x}^{(2k)}v(y;x)  \mathrm{d}y \nonumber \\
  & \quad + \sum_{k=1}^{N} \frac{1}{(2k)!!} \bigg(  \int_{f(\frac{x}{\varT}) + \sqrt{\varT}}^{+\infty} \big| \nabla_{x}^{(2k)}v(y;x) \big| \mathrm{d}y + \int_{-\infty}^{f(\frac{x}{\varT}) - \sqrt{\varT}} \big| \nabla_{x}^{(2k)}v(y;x) \big| \mathrm{d}y \bigg) \nonumber \\
  & \leq 1 + 2 \sum_{k=1}^{N} \frac{1}{(2k)!!} \bigg(  \int_{f(\frac{x}{\varT}) + \sqrt{\varT}}^{+\infty} \big| \nabla_{x}^{(2k)}v(y;x) \big| \mathrm{d}y + \int_{-\infty}^{f(\frac{x}{\varT}) - \sqrt{\varT}} \big| \nabla_{x}^{(2k)}v(y;x) \big| \mathrm{d}y \bigg)/
\end{align}
Here in the last step we use $\int_{\real} v(y;x) \mathrm{d}y = 1$ by definition~\eqref{target:density}. We also and use Eq.~\eqref{ineq4:px-nonlinear} so that for $k \geq 1$, we have
\begin{align*}
   \int_{f(\frac{x}{\varT}) - \sqrt{\varT}}^{f(\frac{x}{\varT}) + \sqrt{\varT}} \nabla_{x}^{(2k)}v(y;x)  \mathrm{d}y &= -\int_{f(\frac{x}{\varT}) + \sqrt{\varT}}^{+\infty}  \nabla_{x}^{(2k)}v(y;x)  \mathrm{d}y - \int_{-\infty}^{f(\frac{x}{\varT}) - \sqrt{\varT}}  \nabla_{x}^{(2k)}v(y;x)  \mathrm{d}y \\
   &\leq  \int_{f(\frac{x}{\varT}) + \sqrt{\varT}}^{+\infty} \big| \nabla_{x}^{(2k)}v(y;x) \big| \mathrm{d}y + \int_{-\infty}^{f(\frac{x}{\varT}) - \sqrt{\varT}} \big| \nabla_{x}^{(2k)}v(y;x) \big| \mathrm{d}y.
\end{align*}
Substituting Ineq.~\eqref{ineq3:px-nonlinear} into the RHS of Ineq.~\eqref{ineq:px-upper-nonlinear-case1} yields
\begin{align}\label{ineq:px-sharp-case1}
   p(x) &\leq  1 + 4\sqrt{2\pi} e^{-\varT/8} \sum_{k=1}^{N} \frac{(2k)!}{(2k)!!}  \cdot (4e^{2})^{k} \cdot \bigg( \frac{\ClinkT}{\varT} \bigg)^{2k} \cdot \big (\Llink(2k) \vee 1\big)^{2k} \nonumber \\
   & \overset{\1}{\leq} 1 + 4\sqrt{2\pi} e^{-\varT/8}  \sum_{k=1}^{N} \bigg( \frac{ 8k  e^{2} \big( \ClinkT \big)^{2} \big (\Llink(2k) \vee 1\big)^{2} }{\varT^{2}} \bigg)^{k} \nonumber \\
   &\overset{\2}{\leq} 1+4\sqrt{2\pi} e^{-\varT/8} \sum_{k=1}^{N} 2^{-k} \leq  1+12e^{-\varT/8}.
\end{align}
Above, in step $\1$ we use $(2k)!/(2k)!! = (2k-1)!! \leq (2k)^{k}$ and in step $\2$ we use
\begin{align*}
    \varT^{2} \geq 16N  e^{2} \big( \ClinkT \big)^{2} \big (\Llink(2N) \vee 1\big)^{2} \geq 16k  e^{2} \big( \ClinkT \big)^{2} \big (\Llink(2k) \vee 1\big)^{2}, \quad \text{for all } k \in [N].
\end{align*}
Combining Ineq.~\eqref{ineq:px-sharp-case1} with Ineq.~\eqref{px-qx-nonlinear-relation} yields
\begin{align}\label{px-qx-bound-case1}
  \big| p(x) -1 \big| + q(x) \leq 2 \big( p(x) -1 \big) \leq 24\exp(-\varT/8), \quad \text{for all } |x| \leq \varT.
\end{align} 
This concludes the proof for the first case. Now we consider the second case.

\paragraph{Case 2: bound $|p(x)-1|+q(x)$ for $|x| \geq \varT$.}
By definition, we obtain
\begin{align*}
    p(x) \vee q(x) &\leq \int_{\real} \big| \mathcal{S}_{N}^{\star}(y \mid x) \big| \mathrm{d}y \\
    &\leq 1 + \sum_{k=1}^{N} \frac{1}{(2k)!!} \int_{\real} \big| \nabla_{x}^{(2k)} v(y;x) \big| \mathrm{d}y \\
    & \overset{\1}{\leq} 1 + \sum_{k=1}^{N} \frac{(2k)!}{(2k)!!} \cdot (4e^{2})^{k}  \cdot \bigg( \frac{\ClinkT}{\varT} \bigg)^{2k} \cdot \bigg[ \bigg( \sum_{j=1}^{2k} \big| h_{j}(x) \big| \bigg) \vee 1 \bigg]^{2k} \int_{\real} e^{-t^{2}/8} \mathrm{d}t \\
    & \overset{\2}{\leq} 1 + 2\sqrt{2\pi} \sum_{k=1}^{N} \bigg( \frac{ 8e^{2}  k \big( \ClinkT \big)^{2} }{\varT^{2}} \bigg)^{k} \cdot \bigg[ \bigg( \sum_{j=1}^{2k} \big| h_{j}(x) \big| \bigg) \vee 1 \bigg]^{2k}.
\end{align*}
In step $\1$ we use Ineq.~\eqref{ineq:derivative-gaussian-nonlinear} and note that since $g(x) = \frac{f(x/\varT) - y}{2}$, we also use
\[
    \int_{\real} e^{-\frac{g(x)^{2}}{2}} \mathrm{d}y = \int_{\real} e^{-\frac{(y-f(x/\varT))^{2}}{8}} \mathrm{d}y = \int_{\real} e^{-t^{2}/8} \mathrm{d}t,
\]
where in the last step we use change of variables $t = y - f(x/\varT)$. In step $\2$, we use the inequalities $(2k)!/(2k)!! = (2k-1)!! \leq (2k)^{k}$ and $\int_{\real} e^{-t^{2}/8} \mathrm{d}t \leq 2\sqrt{2\pi}$. Note that $p(x),q(x)\geq 0$ for all $x\in\real$ by definition~\eqref{px-qx-def}. Consequently, we obtain 
\begin{align}\label{px-qx-bound-case2}
  |p(x)-1| + q(x) &\leq 1 + p(x) + q(x) \nonumber \\
  &\leq 3 + 4\sqrt{2\pi} \sum_{k=1}^{N} \bigg( \frac{ 8e^{2} k \big( \ClinkT \big)^{2} }{\varT^{2}} \bigg)^{k} \cdot \bigg[ \bigg( \sum_{j=1}^{2k} \big| h_{j}(x) \big| \bigg) \vee 1 \bigg]^{2k}, \quad \text{for all } x \in \real.
\end{align}

\paragraph{Putting the two cases together.} We have
\begin{align*}
  \int_{\real} \big( |p(x)-1| + q(x) \big) \cdot u(x;\theta) \mathrm{d}x &\leq \int_{|x| \leq \varT} \big( |p(x)-1| + q(x) \big) \cdot u(x;\theta) \mathrm{d}x \\
  & \quad + \int_{|x| \geq \varT} \big( |p(x)-1| + q(x) \big) \cdot u(x;\theta) \mathrm{d}x.
\end{align*}
We then bound the two integrals in the display above. Applying Ineq.~\eqref{px-qx-bound-case1} yields
\begin{align}\label{ineq8:px-qx-nonlinear}
  \int_{|x| \leq \varT} \big( |p(x)-1| + q(x) \big) \cdot u(x;\theta) \mathrm{d}x \leq 24\exp(-\varT/8) \cdot \int_{|x| \leq \varT} u(x;\theta) \mathrm{d}x \leq 24\exp(-\varT/8). 
\end{align}
Applying Ineq.~\eqref{px-qx-bound-case2} yields
\begin{align}\label{ineq4:px-qx-nonlinear}
    & \int_{|x| \geq \varT} \big( |p(x)-1| + q(x) \big) \cdot u(x;\theta) \mathrm{d}x  \nonumber \\
    & \leq 3\int_{|x| \geq \varT} u(x;\theta) \mathrm{d}x + 4\sqrt{2\pi} \sum_{k=1}^{N} \bigg( \frac{ 8e^{2} k \big( \ClinkT \big)^{2} }{\varT^{2}} \bigg)^{k} \cdot \int_{|x|\geq \varT}\bigg[ \bigg( \sum_{j=1}^{2k} \big| h_{j}(x) \big| \bigg) \vee 1 \bigg]^{2k} \cdot u(x;\theta) \mathrm{d}x.
\end{align}
By definition of $u(x;\theta)$~\eqref{source-gaussian-density} for $\sigma=1$, we obtain
\[
    u(x;\theta) = \frac{1}{\sqrt{2\pi}} e^{-\frac{(x-\theta)^{2}}{2} } \leq e^{ \frac{\theta^{2}}{2} } \cdot \frac{1}{\sqrt{2\pi} } e^{-\frac{x^{2}}{4}},
\]
where in the last step we use $(x-\theta)^{2} \geq x^{2}+\theta^{2} - 2|\theta x| \geq x^{2} + \theta^{2} - (x^{2}/2 + 2\theta^{2}) \geq x^{2}/2 - \theta^{2}$. Consequently, we obtain 
\begin{align}\label{ineq5:px-qx-nonlinear}
    \int_{|x| \geq \varT} u(x;\theta) \mathrm{d}x \leq 2 e^{ \frac{\theta^{2}}{2} }  \int_{\varT}^{+\infty} \frac{1}{\sqrt{2\pi} } e^{-\frac{x^{2}}{4}} \mathrm{d}x \leq \sqrt{2} e^{ \frac{\theta^{2}}{2} } e^{-\frac{\varT^{2}}{4}},
\end{align}
and similarly for $k\geq 1$, we have
\begin{align}\label{ineq6:px-qx-nonlinear}
    \int_{|x|\geq \varT}\bigg[ \bigg( \sum_{j=1}^{2k} \big| h_{j}(x) \big| \bigg) \vee 1 \bigg]^{2k} \cdot u(x;\theta) \mathrm{d}x &\leq e^{ \frac{\theta^{2}}{2} } \int_{|x|\geq \varT} \bigg[ \bigg( \sum_{j=1}^{2k} \big| h_{j}(x) \big| \bigg) \vee 1 \bigg]^{2k} \cdot \frac{ e^{-\frac{x^{2}}{4}} }{\sqrt{2\pi} } \mathrm{d}x \nonumber \\
    &\overset{\1}{\leq} e^{ \frac{\theta^{2}}{2} } \int_{|x|\geq \varT}\frac{ e^{-\frac{x^{2}}{8}} }{\sqrt{2\pi} } \mathrm{d}x \cdot \Big( \LlinkT(2k,\varT) \vee 1\Big)^{2k} \nonumber \\
    & \overset{\2}{\leq}  2 e^{ \frac{\theta^{2}}{2} } e^{-\frac{\varT^{2}}{8}} \cdot \Big( \LlinkT(2k,\varT) \vee 1\Big)^{2k}.
\end{align}
Above, in step $\1$ we use definition of $h_{j}(x)$~\eqref{shorthand-h-x} and $\LlinkT(n,\varT,\sigma)$~\eqref{def-quantity-link} so that for all $|x| \geq \varT$,
\begin{align*}
    \bigg[ \bigg( \sum_{j=1}^{2k} \big| h_{j}(x) \big| \bigg) \vee 1 \bigg] e^{-\frac{x^{2}}{8 (2k) }} &=  \bigg[ \bigg( \sum_{j=1}^{2k} \frac{|f(x/\varT)|}{\big( \ClinkT \big)^{j}j!} \bigg) \vee 1 \bigg] e^{-\frac{x^{2}}{8 (2k) }} \\
    & \leq \sup_{|t|\geq 1} \bigg[ \bigg( \sum_{j=1}^{2k} \frac{|f(t)|}{\big( \ClinkT \big)^{j}j!} \bigg) \vee 1 \bigg] e^{-\frac{\varT^{2}t^{2}}{8 (2k) }}
    \\ &\leq \LlinkT(2k,\varT) \vee 1.
\end{align*}
In step $\2$, we use
\[
    \int_{|x|\geq \varT}\frac{1}{\sqrt{2\pi} } e^{-\frac{x^{2}}{8}} \mathrm{d}x \leq 2\exp\big( - \varT^{2} / 8 \big).
\]
Substituting Ineq.~\eqref{ineq5:px-qx-nonlinear} and Ineq.~\eqref{ineq6:px-qx-nonlinear} into Ineq.~\eqref{ineq4:px-qx-nonlinear} yields
\begin{align}\label{ineq7:px-qx-nonlinear}
    &\int_{|x| \geq \varT} \big( |p(x)-1| + q(x) \big) \cdot u(x;\theta) \mathrm{d}x \nonumber \\
    &\leq 2 e^{ \frac{\theta^{2}}{2} } e^{-\frac{\varT^{2}}{8}} \cdot \bigg( 3 + 4\sqrt{2\pi} \sum_{k=1}^{N} \Bigg( \frac{ 8e^{2}  k \big( \ClinkT \big)^{2} }{\varT^{2}} \bigg)^{k} \cdot \Big( \LlinkT(2k,\varT) \vee 1\Big)^{2k}  \Bigg) \nonumber \\
    & \overset{\1}{\leq} 2 e^{ \frac{\theta^{2}}{2} } e^{-\frac{\varT^{2}}{8}} \cdot \bigg( 3 + 4\sqrt{2\pi} \sum_{k=1}^{N} \frac{1}{2^{k}}  \Bigg) \leq 30 e^{ \frac{\theta^{2}}{2} } e^{-\frac{\varT^{2}}{8}},
\end{align}
where in step $\1$ we use the condition
\[
    \varT^{2} \geq 16 e^{2}  N \big( \ClinkT \big)^{2}  \big( \LlinkT(2N,\varT) \vee 1\big)^{2} \geq 8 e^{2}  k \big( \ClinkT \big)^{2}  \big( \LlinkT(2k,\varT) \vee 1\big)^{2} \quad \text{for all } k \in [N].
\]
Finally, combining Ineq.~\eqref{ineq7:px-qx-nonlinear} with Ineq.~\eqref{ineq8:px-qx-nonlinear} yields
\begin{align*}
    \int_{\real} \big( |p(x)-1| + q(x) \big) \cdot u(x;\theta) \mathrm{d}x \leq 24\exp(-\varT/8) + 30 e^{ \frac{\theta^{2}}{2} } e^{-\frac{\varT^{2}}{8}}.
\end{align*}
Taking a supremum over $\theta \in \Theta$ proves the desired result.  \qed

\subsection{Proof of Theorem~\ref{thm:gaussian-mean-nonlinear}(b)}\label{sec:pf_gaussian_nonlinear_term3}
We split the proof into three parts.
\paragraph{Construction of the reduction based on rejection kernel.}
We define a class of base measures $\left\{ \mathcal{D}(\cdot \mid x) \right\}_{x \in \mathbb{R}}$ as
\begin{align}\label{base-measure-nonlinear}
    \mathcal{D}(y \mid x) = \frac{1}{2\sqrt{2\pi}} \exp\left( -\frac{(y - f(x/\varT))^2}{8} \right), \quad \text{for all } x, y \in \mathbb{R}.
\end{align}
As before, we use the general rejection kernel $x \mapsto \textsc{rk}(x, T, M, y_{0})$ from~\citet[page 10]{lou2025computationally} with the base measures $\mathcal{D}(y \mid x)$ defined in~\eqref{base-measure-nonlinear} and the signed kernel $\mathcal{S}_{N}^{\star}(y \mid x)$ defined in~\eqref{truncated-signed-kernel}. 
As before, $M > 0$ is a constant that is required to satisfy: for input $x$, 
\[
    \frac{\mathcal{S}_{N}^{\star}(y \mid x) \vee 0}{ \mathcal{D}(y \mid x)} \leq M, \quad \text{for all } y \in \mathbb{R},
\]
and $y_{0}$ is the initialization. Throughout this section, we set
\begin{align}\label{parameters-rk-nonlinear}
    T = \left \lceil 4\sqrt{2\pi} \log(4/\epsilon) \right \rceil , \quad M = 2\sqrt{2\pi}, \quad \text{and} \quad y_{0} = 0.
\end{align}
We are now ready to describe the algorithm $\mathsf{K}$. For $X_{\theta} \sim \NORMAL(\theta,\sigma^{2})$, define
\begin{align}\label{reduction-nonlinear}
  \mathsf{K}(X_{\theta}): = \begin{cases} \textsc{rk}(X_{\theta}, T, M, y_{0}), &\text{if } |X_{\theta}| \leq \varT,\\
  y_{0}, &\text{otherwise.}   \end{cases}
\end{align}
Note that for all $|x| \leq \varT$ and $y \in \real$, we have by definition of $\mathcal{S}_{N}^{\star}(y \mid x)$~\eqref{truncated-signed-kernel},
\begin{align*}
  \frac{\mathcal{S}_{N}^{\star}(y \mid x) \vee 0}{ \mathcal{D}(y \mid x)} &\leq \sum_{k=0}^{N} \frac{1}{(2k)!!} \frac{\big|\nabla_{x}^{(2k)}v(y;x)\big|}{\mathcal{D}(y \mid x)} \\ 
  & \overset{\1}{\leq} 2\sqrt{2\pi} \sum_{k=0}^{N}  \frac{(2k)!}{(2k)!!} \cdot \frac{(4e^{2})^{k} \big( \ClinkT \big)^{2k} }{\varT^{2k}} \cdot \bigg[ \bigg( \sum_{j=1}^{2k} \big|h_{j}(x)\big| \bigg) \vee 1 \bigg]^{2k} \\
  & \overset{\2}{\leq} 2\sqrt{2\pi} \sum_{k=0}^{N} \frac{ (2k)^{k} \cdot (4e^{2})^{k} \big( \ClinkT \big)^{2k} }{\varT^{2k}} \cdot \big( \Llink(2k) \vee 1 \big)^{2k} \\
  & \overset{\3}{\leq} 2\sqrt{2\pi} \sum_{k=0}^{N} \frac{1}{2^{k}} = M.
\end{align*}
In step $\1$ we use Ineq.~\eqref{ineq:upbound-derivative-abs-nonlinear} and definition of $\mathcal{D}(y \mid x)$~\eqref{base-measure-nonlinear}. In step $\2$ we use $(2k)!/(2k)!! = (2k-1)!! \leq (2k)^{k}$ and recall $h_{j}(x)$~\eqref{shorthand-h-x} and $\Llink(2k)$~\eqref{def-quantity-link}, so that
\[
    \sum_{j=1}^{2k} \big|h_{j}(x)\big| \leq \Llink(2k) \quad \text{for all } |x| \leq \varT.
\]
In step $\3$ we use condition~\eqref{assump:thm-nonlinear} that $\varT$ is large enough to guarantee that the summation is bounded. Therefore, the reduction $\mathsf{K}(X_{\theta}) = \textsc{rk}(X_{\theta}, T, M, y_{0})$ when $|X_{\theta}| \leq \varT$~\eqref{reduction-nonlinear} is well defined.

\paragraph{Computational complexity of the reduction.} The rejection kernel $x \mapsto \textsc{rk}(x, T, M, y_{0})$ requires, at each iteration, generating a sample from the uniform distribution $\mathsf{Unif}([0,1])$ and a sample from the Gaussian base measure~\eqref{base-measure-nonlinear}, both of which take $\mathcal{O}(1)$ time. Additionally, it requires evaluating the signed kernel $\mathcal{S}_{N}^{\star}$ once per iteration, which takes time $T_{\mathsf{eval}}(N)$. The procedure is repeated for at most $T = \mathcal{O}(\log(4/\epsilon))$ iterations. Therefore, the total computational complexity of the reduction $\mathsf{K}(X_{\theta})$ is
\[
\mathcal{O}\big( \log(4/\epsilon) \cdot T_{\mathsf{eval}}(N) \big).
\]

\paragraph{Proof of TV deficiency guarantee of the reduction.}
Having defined the reduction, we turn to prove the claimed guarantee~\eqref{ineq-gaussian-nonlinear-term3}. Applying the triangle inequality, we obtain
\begin{align*}
    \mathsf{d}_{\mathsf{TV}}\big( \mathsf{K}(X_{\theta}), Y_{\theta}\big) & = \frac{1}{2} \big\| \mathsf{K}(X_{\theta}) - v(\cdot;\theta) \big \|_{1}  \\
    & \leq \frac{1}{2} \left\| \mathsf{K}(X_{\theta}) - \int_{\real} \mathcal{T}_{N}^{\star}(\cdot \mid x) u(x;\theta) \mathrm{d}x  \right \|_{1} + 
    \frac{1}{2} \left\|  \int_{\real} \mathcal{T}_{N}^{\star}(\cdot \mid x) u(x;\theta) \mathrm{d}x - v(\cdot;\theta)  \right \|_{1}.
\end{align*}
Taking the supremum over $\theta \in \Theta$ on both sides yields
\begin{align}\label{ineq:triangle-TV-KtoV}
  \sup_{\theta \in \Theta}\; \mathsf{d}_{\mathsf{TV}}\big( \mathsf{K}(X_{\theta}), Y_{\theta}\big) \leq 
  \frac{1}{2} \sup_{\theta \in \Theta}\; \left\| \mathsf{K}(X_{\theta}) - \int_{\real} \mathcal{T}_{N}^{\star}(\cdot \mid x) u(x;\theta) \mathrm{d}x  \right \|_{1} + \delta\big(\mathcal{U}, \mathcal{V}; \mathcal{T}_{N}^{\star} \big).
\end{align}
We claim that under our setting of parameters~\eqref{parameter-N-varT-nonlinear}~\eqref{parameters-rk-nonlinear}, the following holds.
\begin{subequations}
\begin{align}
\label{part1-nonlinear-reduction}
    &\frac{1}{2} \sup_{\theta \in \Theta}\; \left\| \mathsf{K}(X_{\theta}) - \int_{\real} \mathcal{T}_{N}^{\star}(\cdot \mid x) u(x;\theta) \mathrm{d}x  \right \|_{1} \leq \frac{\epsilon}{2}, \\
\label{part2-nonlinear-reduction}
    & \text{and} \qquad \delta\big(\mathcal{U}, \mathcal{V}; \mathcal{T}_{N}^{\star} \big) \leq \frac{\epsilon}{2}.
\end{align}
\end{subequations}
The desired guarantee Ineq.~\eqref{ineq-gaussian-nonlinear-term3} follows immediately. We now prove the above two inequalities separately.

\noindent \underline{Proof of Ineq.~\eqref{part1-nonlinear-reduction}.}
For each $x$ satisfying $|x| \leq \varT$, we let $Y = \textsc{rk}(x, T, M, y_{0})$. Then from~\citet[Lemma 2]{lou2025computationally}, we obtain that the random variable $Y$ has the conditional distribution
\begin{align}\label{condition-density-nonlinear}
  f_{Y}(y \mid x) = \frac{\mathcal{S}_{N}^{\star}(y\mid x) \vee 0}{p(x)} \cdot \big( 1 - (1-p(x)/M)^{T} \big) + \delta_{y_{0}}(y) \cdot (1-p(x)/M)^{T}, \quad \text{for all } y \in \real,
\end{align}
where $p(x)$ is defined in Eq.~\eqref{px-qx-def} and $\delta_{y_{0}}(y)$ is the Dirac delta function centered at $y_{0}$. By the law of total probability and definition of $\mathsf{K}(X_{\theta})$ in Eq.~\eqref{reduction-nonlinear}, we obtain that $\mathsf{K}(X_{\theta})$ has the marginal distribution
\begin{align}\label{marginal-density-nonlinear}
    f_{\mathsf{K}(X_{\theta})}(y) & = \int_{|x| \leq \varT} f_{Y}(y\mid x) \cdot u(x;\theta) \mathrm{d}x  + \delta_{y_{0}}(y) \cdot \int_{|x| \geq \varT} u(x;\theta) \mathrm{d}x \nonumber \\
    & = \int_{|x| \leq \varT} \mathcal{T}_{N}^{\star}(y \mid x) \cdot (1-\ell(x)) \cdot u(x;\theta) \mathrm{d}x \nonumber \\
    & \hspace{3cm} + \delta_{y_{0}}(y) \cdot \bigg( \int_{|x| \geq \varT} u(x;\theta) \mathrm{d}x + \int_{|x| \leq \varT} \ell(x) u(x;\theta) \mathrm{d}x\bigg),
\end{align}
where in the last step we have used Eq.~\eqref{condition-density-nonlinear} along with the shorthand notation $\ell(x) = (1-p(x)/M)^{T}$ and  the previously defined notation for the Markov kernel~\eqref{close-markov-kernel}, $\mathcal{T}_{N}^{\star}(y \mid x) = \frac{\mathcal{S}_{N}^{\star}(y\mid x) \vee 0}{p(x)}$. Further recall from Ineq.~\eqref{ineq:px-sharp-case1} that $|p(x)-1| \leq 12\exp(-\varT/8) \leq \epsilon/2 \leq 1/2$ for all $|x| \leq \varT$. Consequently, we obtain that for all $|x| \leq \varT$,
\begin{align}\label{ineq:bound-lx}
    1/2 \leq p(x) \leq 3/2 \quad \text{and} \quad 0 \leq \ell(x) \leq e^{- p(x) T/M} \leq e^{- T/(2M)}.
\end{align}
Using Eq.~\eqref{marginal-density-nonlinear}, we obtain
\begin{align*}
    &\bigg\| \mathsf{K}(X_{\theta}) - \int_{\real} \mathcal{T}_{N}^{\star}(\cdot \mid x) u(x;\theta) \mathrm{d}x  \bigg \|_{1} \\
    & = \int_{\real} \bigg| \int_{|x| \leq \varT} \mathcal{T}_{N}^{\star}(y \mid x) \cdot (1-\ell(x)) \cdot u(x;\theta) \mathrm{d}x   - \int_{\real} \mathcal{T}_{N}^{\star}(y \mid x) u(x;\theta) \mathrm{d}x \bigg| \mathrm{d}y  \\
    &  \hspace{3cm} + \int_{|x| \geq \varT} u(x;\theta) \mathrm{d}x + \int_{|x| \leq \varT} \ell(x) u(x;\theta) \mathrm{d}x \\
    & \leq \int_{\real} \int_{|x| \leq \varT} \mathcal{T}_{N}^{\star}(y \mid x)\cdot \ell(x) u(x;\theta) \mathrm{d}x \mathrm{d}y + \int_{\real} \int_{|x| \geq \varT} \mathcal{T}_{N}^{\star}(y \mid x) \cdot u(x;\theta) \mathrm{d}x \mathrm{d}y  \\ 
    &  \hspace{3cm} + \int_{|x| \geq \varT} u(x;\theta) \mathrm{d}x + \int_{|x| \leq \varT} \ell(x) u(x;\theta) \mathrm{d}x,
\end{align*}
where in the last step we apply the triangle inequality. Continuing, we bound each of the terms in the display above
\begin{align*}
  \int_{\real} \int_{|x| \leq \varT} \mathcal{T}_{N}^{\star}(y \mid x)\cdot \ell(x) u(x;\theta) \mathrm{d}x \mathrm{d}y &=  \int_{|x| \leq \varT} \int_{\real} \mathcal{T}_{N}^{\star}(y \mid x) \mathrm{d}y \cdot \ell(x) u(x;\theta) \mathrm{d}x \mathrm{d}y \\
  &=  \int_{|x| \leq \varT} \ell(x) u(x;\theta) \mathrm{d}x \leq e^{-\frac{T}{2M}},
\end{align*}
where in the last step we use Ineq.~\eqref{ineq:bound-lx}. Also, we have
\[
  \int_{\real} \int_{|x| \geq \varT} \mathcal{T}_{N}^{\star}(y \mid x) \cdot u(x;\theta) \mathrm{d}x \mathrm{d}y = \int_{|x| \geq \varT } u(x;\theta) \mathrm{d}x \leq e^{\frac{\theta^{2}}{2}} \int_{|x| \geq \varT} \frac{e^{-x^{2}/(4)}}{\sqrt{2\pi}} \mathrm{d}x \leq \sqrt{2}e^{\frac{\theta^{2}}{2} - \frac{\varT^{2}}{4}}.
\]
Combining the three pieces together yields
\begin{align*}
    \bigg\| \mathsf{K}(X_{\theta}) - \int_{\real} \mathcal{T}_{N}^{\star}(\cdot \mid x) u(x;\theta) \mathrm{d}x  \bigg \|_{1} \leq 2e^{-\frac{T}{2M}} + 2\sqrt{2}e^{\frac{\theta^{2}}{2} - \frac{\varT^{2}}{4}} \leq \epsilon,
\end{align*}
where in the last step we use $T \geq 2M \log(4/\epsilon)$ and $\varT^{2} \geq 4 \log(4\sqrt{2}/\epsilon) + 2\sup_{\theta \in \Theta} \theta^{2}$. This concludes the proof of Ineq.~\eqref{part1-nonlinear-reduction}.

\noindent \underline{Proof of Ineq.~\eqref{part2-nonlinear-reduction}.}
Since we set 
\[
    2N+1 \geq  \Big( 2e \log \big(4\ClinkB/\epsilon\big) \Big)^{2(1+\Clink)}, \quad \varT^{2} \geq 64 \log^{2}\big(96/\epsilon \big) \vee \; \Big[ 8 \log\big(120/\epsilon\big) + 4\sup_{\theta \in \Theta} \theta^{2} \Big],
\]
we obtain
\begin{align*}
    &\ClinkB \exp\bigg(-(2e)^{-1}(2N+1)^{\frac{1}{2(1+\Clink)}} \bigg) \leq \frac{\epsilon}{4}, \quad \text{and} \\
    &12\exp(-\varT/8) + 15\sup_{\theta \in \Theta} \; \exp\bigg( \frac{\theta^{2}}{2} -\frac{\varT^{2}}{8} \bigg) \leq \frac{\epsilon}{4}.
\end{align*}
Substituting the two inequalities in the display above into Ineq.~\eqref{ineq-gaussian-nonlinear-term0} completes the proof.
\qed

\subsection{Proof of Corollary~\ref{example:monomial}}\label{sec:pf_corollary_monomial}
We verify each part separately.
\paragraph{Verifying Assumption~\ref{assump:link}.}
By definition, we obtain $f^{(j)}(x) = B! \cdot x^{B-j} / (B-j)!$ for all $x\in \real$ and $0 \leq j \leq B$, and $f^{(j)}(x) = 0$ for all $x\in \real$ and $j \geq B+1$. We now verify Assumption~\ref{assump:link} in two cases.

\noindent \underline{Case 1: $n \leq B$.}
For $n \leq B$ and $x \in \real$, we have
\begin{align*}
 \sum_{j=1}^{n} \frac{ \big| f^{(j)} \big( \big(\theta + \sqrt{2} x \big)/\varT \big) \big| }{\big(\ClinkT \big)^{j} j! } 
 &= \sum_{j=1}^{n} \frac{B! }{\big(\ClinkT \big)^{j} j! \cdot (B-j)!} \cdot \frac{ \big| \theta + \sqrt{2} x \big|^{B-j}}{\varT^{B-j}} \\
 & \leq \sum_{j=1}^{n} \frac{ \binom{B}{j} }{\big( \ClinkT \big)^{j}} \cdot \frac{ (2|\theta|)^{B-j} + (2\sqrt{2}|x|)^{B-j}  }{\varT^{B-j}} \\
 & \overset{\1}{\leq} \sum_{j=1}^{n} \frac{ 1 + |x|^{B-j}  }{(4e)^{j}},
\end{align*}
where in step $\1$ we use $\binom{B}{j} \leq B^{j}$ and let $\ClinkT = 4eB$, and we use $\varT \geq 2\sup_{\theta \in \Theta} |\theta|$ and $\varT \geq 2\sqrt{2}$. Consequently, we obtain for all $x \in \real$,
\begin{align}\label{ineq1-example-monomial}
     \bigg( \sum_{j=1}^{n} \frac{ \big| f^{(j)} \big( \big(\theta + \sqrt{2} x \big)/\varT \big) \big| }{\big(\ClinkT \big)^{j} j! }  \bigg)^{n} \leq \bigg( \sum_{j=1}^{n} \frac{ 1 + |x|^{B-j}  }{(4e)^{j}}  \bigg)^{n} &\leq n^{n} \sum_{j=1}^{n} \bigg( \frac{ 1 + |x|^{B-j}  }{(4e)^{j}} \bigg)^{n} \nonumber \\
     & \leq 2^{n}n^{n} \sum_{j=1}^{n} \frac{1 + |x|^{(B-j)n} }{(4e)^{jn}}.
\end{align}
Integrating over $x$ yields
\begin{align}\label{ineq2-example-monomial}
    \int_{\real} \bigg( \sum_{j=1}^{n} \frac{ \big| f^{(j)} \big( \big(\theta + \sqrt{2} x \big)/\varT \big) \big| }{\big(\ClinkT \big)^{j} j! }  \bigg)^{n} \frac{e^{-x^{2}/2}}{\sqrt{2\pi}} \mathrm{d}x & \leq \frac{2^{n} n^{n}}{(2e)^{n}} \sum_{j=1}^{n} \Big[ 2^{-j} + 2^{-j} \big( (B-j)n !! \big) \Big] \nonumber \\
    & \overset{\1}{\leq} n! \cdot (2B n)! \leq \big( (2B+1)n \big)!,
\end{align} 
where in step $\1$ we use $n^{n} \leq e^{n} n!$. Thus, Assumption~\ref{assump:link} holds with $\Clink = 2B+1$ and $\Clink = 4eB$ for all $n \leq B$. 

\noindent \underline{Case 2: $n \geq B+1$.} Note that $f^{(j)}(x) = 0$ for all $x \in \real$ when $j \geq B+1$. Consequently, following the same steps as in Ineq.~\eqref{ineq1-example-monomial} yields 
\begin{align*}
  \bigg( \sum_{j=1}^{n} \frac{ \big| f^{(j)} \big( \big(\theta + \sqrt{2} x \big)/\varT \big) \big| }{\big(\ClinkT \big)^{j} j! }  \bigg)^{n} &= 
  \bigg( \sum_{j=1}^{B} \frac{ \big| f^{(j)} \big( \big(\theta + \sqrt{2} x \big)/\varT \big) \big| }{\big(\ClinkT \big)^{j} j! }  \bigg)^{n} \\
  & \leq 2^{n}B^{n} \sum_{j=1}^{B} \frac{1 + |x|^{(B-j)n} }{(4e)^{jn}}.
\end{align*}
Integrating over $x$ and using $B \leq n$ yields the same bound as Ineq.~\eqref{ineq2-example-monomial}. 

\paragraph{Bounding the quantities $\Llink(n)$ and $\LlinkT(n,\varT)$.} We next turn to bound $\Llink(n)$ and $\LlinkT(n,\varT)$ defined in Eq.~\eqref{def-quantity-link}. We obtain by definition that for all $n \in \mathbb{N}$ and $n\geq 1$,
\begin{align*}
    \Llink(n) \leq \sup_{|x| \leq 1} \sum_{j=1}^{B} \frac{B! \cdot |x|^{B-j}}{(B-j)! \cdot \big( \ClinkT \big)^{j}j!} \leq \sum_{j=1}^{B} (2e)^{-j} \leq 1,
\end{align*}
where we use $\frac{B!}{(B-j)! \cdot j!} = \binom{B}{j} \leq B^{j}$ and $\ClinkT = 2eB$. Reasoning similarly, we obtain
\begin{align*}
  \LlinkT(n,\varT) \leq  \sup_{|x| \geq 1} \sum_{j=1}^{B} \frac{|x|^{B}}{(2e)^{j}} \cdot \exp\bigg( -\frac{\varT^{2}x^{2}}{8n} \bigg) \leq \sup_{|x|\geq 1} |x|^{B} \exp\bigg( -\frac{\varT^{2}x^{2}}{8n} \bigg) \leq 1,
\end{align*}
where in the last step we let $\varT^{2} \geq 8n B$, whereby the supremum is achieved at $|x| = 1$.

\paragraph{Bounding the evaluation time $T_{\mathsf{eval}}(n)$.} We next bound the evaluation time $T_{\mathsf{eval}}(n)$ of the signed kernel $\mathsf{S}_{n}^{\star}(y\mid x)$~\eqref{truncated-signed-kernel},
noting the closed-form expression~\eqref{eq:derivative-gaussian-nonlinear}. Note that $h_{j}(x) = 0$ for all $j\geq B+1$ and $x\in\real$. Thus, we obtain
\begin{align*}
  \nabla_{x}^{(n)}v(y;x) &= n! \sum_{m_{1},\dots,m_{B}} \frac{(-1)^{m_{1}+\cdots+m_{B}}}{m_{1}!\,m_{2}! \cdots m_{B}!} \cdot H_{m_{1}+\cdots+m_{B}}(g(x)) \phi(g(x)) \cdot  \prod_{j=1}^{B} \bigg( \frac{h_{j}(x)}{2} \bigg)^{m_{j}}.
\end{align*}
Here the summation is over all $B$-tuples of nonnegative integers $(m_{1},\cdots,m_{B})$ satisfying $\sum_{j=1}^{B} j m_{j} = n$; thus there is at most $n^{B}$ such $B$-tuples. Also, each term in the summation can be computed in time $\mathcal{O}(n^{2})$ using Eq.~\eqref{Hermite-closed-form}.
Consequently, evaluating $\nabla_{x}^{(n)}v(y;x)$ for any pair of $(x,y)$ takes time $\mathcal{O}(n^{B+2})$. Using the expression~\eqref{truncated-signed-kernel}, we conclude that $T_{\mathsf{eval}}(n) = \mathcal{O}(n^{B+2})$.

\paragraph{Proof of the consequence.} Finally, note that conditions~\eqref{assump:thm-nonlinear} and~\eqref{parameter-N-varT-nonlinear} are equivalent to
\[
  N \geq C_{B} \log^{4(B+1)}(C_{B}/\epsilon) \quad \text{and} \quad \varT^{2} \geq C_{B}  \cdot \log^{8(B+1)}(C_{B}/\epsilon) + 4 \sup_{\theta \in \Theta} \theta^{2},
\]
where $C_{B}$ is a constant that depends only on $B$. Applying Theorem~\ref{thm:gaussian-mean-nonlinear} with the guarantees above yields the desired result.
\qed

\subsection{Proof of Corollary~\ref{example:analytical}}

We verify each part separately.
\paragraph{Verifying Assumption~\ref{assump:link}.}
Let $\ClinkT = 2R$. Since $f \in \mathcal{F}(R)$, then we obtain for all $x\in \real$,
\begin{align*}
  \sum_{j=1}^{n} \frac{ \big| f^{(j)} \big( \big(\theta + \sqrt{2} x \big)/\varT \big) \big| }{\big(\ClinkT \big)^{j} j! }  \leq \sum_{j=1}^{n} \frac{1}{2^{j}} \leq 1.
\end{align*}
Consequently, Assumption~\ref{assump:link} holds with $\ClinkT = 2R$ and $\Clink = 0$.

\paragraph{Bounding the quantities $\Llink(n)$ and $\LlinkT(n,\varT)$.}
We obtain by definition~\eqref{def-quantity-link},
\begin{align*}
    \Llink(n) \leq \sum_{j=1}^{n} \frac{1}{2^{j}} \leq 1 \quad \text{and} \quad \LlinkT(n,\varT) \leq 1 \quad \text{for all } n \in \mathbb{N}, \varT >0.
\end{align*}

\paragraph{Bounding the evaluation time $T_{\mathsf{eval}}(N)$.} We first derive the computational time of evaluating $\nabla_{x}^{(n)}v(y;x)$. Using Faà di Bruno's formula yields
\begin{align*}
  \nabla_{x}^{(n)}v(y;x) &= n! \sum_{m_{1},\dots,m_{n}} \frac{(-1)^{m_{1}+\cdots+m_{n}}}{m_{1}!\,m_{2}! \cdots m_{n}!} \cdot H_{m_{1}+\cdots+m_{n}}(g(x)) \phi(g(x)) \cdot  \prod_{j=1}^{n} \bigg( \frac{ f^{(j)}(x/\varT)}{\varT^{j}j!} \bigg)^{m_{j}},
\end{align*}
where the summation is over all nonnegative integers $(m_{1},\dots,m_{n})$ satisfying $\sum_{j=1}^{n} jm_{j} = n$. Let $P(n)$ denote the number of $n$-tuples of nonnegative integers $(m_{1},\cdots,m_{n})$ satisfying $\sum_{j=1}^{n} j m_{j} = n$. Note that $P(n)$ corresponds to the partition function and asymptotically satisfies
\[
    P(n) \sim \frac{1}{4n\sqrt{3}} \exp\Big(\pi \sqrt{2n/3}\Big) \quad \text{as} \quad n \uparrow +\infty.
\]
Thus, there are at most $\mathcal{O}\big(\exp\big(\pi \sqrt{n} \big) \big)$ terms in the above summation, and each term in the summation can be computed in time $\mathcal{O}(\mathsf{poly}(n) \exp\big( C \sqrt{n} \big) \big)$ provided we can evaluate $f^{(n)}(x)$ in time $\mathcal{O}\big(\exp\big(C \sqrt{n} \big) \big)$, which is in turn true by assumption. Thus we can evaluate $\mathcal{S}_{N}^{\star}(y \mid x)$~\eqref{truncated-signed-kernel} in time
\[
  \mathcal{O} \Big( \mathsf{poly}(N) \exp\Big( (C+\pi) \sqrt{N} \Big)   \Big) = \mathcal{O}\big(\mathsf{poly}(1/\epsilon)\big),
\]
where in the last step we use the following logic: Since $\Clink = 0$, it suffices to take $N = \mathcal{O}(\log^{2}(1/\epsilon))$. Thus $N$ satisfies conditions~\eqref{assump:thm-nonlinear} and~\eqref{parameter-N-varT-nonlinear}. 

\paragraph{Proof of the consequence.} Finally, note that conditions~\eqref{assump:thm-nonlinear} and~\eqref{parameter-N-varT-nonlinear} are equivalent to
\[
  N \geq \widetilde{C} \log^{2}(\widetilde{C}/\epsilon) \quad \text{and} \quad \varT^{2} \geq \widetilde{C} R^{2}  \cdot \log^{2}(\widetilde{C}/\epsilon) + 4 \sup_{\theta \in \Theta} \theta^{2},
\]
where $\widetilde{C}$ is a universal constant. Applying Theorem~\ref{thm:gaussian-mean-nonlinear} with the guarantees above yields the desired result.
\qed

\subsection{Proofs of claims in Example~\ref{examples}} \label{sec:pf_concrete_examples}

We perform derivations for the different functions separately.

\paragraph{Trigonometric links.} Consider $f(\theta) = \cos(a\theta+b)$. Clearly, we have $|f^{(n)}(x)| \leq |a|^{n}$ for all $x \in \real$ and $n \in \mathbb{N}$. Thus by letting $R = 2|a|$, we have
\[
    \sum_{j=1}^{n} \frac{|f^{(j)}(x)|}{R^{j} j!} \leq \sum_{j=1}^{n} \frac{1}{2^{j}} \leq 1.
\]
Consequently, $f \in \mathcal{F}(2|a|)$. The same holds when $f(\theta) = \sin(a\theta+b)$. This concludes the proof.

\paragraph{Function $f(\theta) = \mathsf{sech}(\theta) = \frac{2}{e^{\theta} + e^{-\theta}}$.} We let $g(x) = \frac{2}{x}$ and $h(x) = e^{x} + e^{-x}$. Note that $f(x) = g(h(x))$. Applying Faà di Bruno's formula yields
\begin{align}\label{eq:derivative-sech}
    f^{(n)}(x) &= \sum_{m_{1},\dots,m_{n}} \frac{n!}{m_{1}!\,m_{2}! \cdots m_{n}!} g^{(m_{1}+\cdots+m_{n})}(h(x)) \prod_{j=1}^{n} \bigg(\frac{h^{(j)}(x)}{j!}  \bigg)^{m_{j}} \nonumber \\
    &= n!\sum_{m_{1},\dots,m_{n}} \frac{(-1)^{m_{1}+\cdots+m_{n}}}{m_{1}!\,m_{2}! \cdots m_{n}!} \cdot \frac{2(m_{1}+\cdots+m_{n})!  }{h(x)^{m_{1}+\cdots+m_{n}}} \prod_{j=1}^{n} \bigg(\frac{h^{(j)}(x)}{j!}  \bigg)^{m_{j}}.
\end{align}
Note that by definition $|h^{(j)}(x)| \leq h(x)$ for all $x\in \real$. Consequently, we obtain
\begin{align*}
    &\frac{(m_{1}+\cdots+m_{n})!}{m_{1}!\,m_{2}! \cdots m_{n}!} \cdot \frac{ 1 }{h(x)^{m_{1}+\cdots+m_{n}}} \prod_{j=1}^{n} \bigg(\frac{| h^{(j)}(x)| }{j!}  \bigg)^{m_{j}} \\
    & \leq \frac{(m_{1}+\cdots+m_{n})!}{m_{1}!\,m_{2}! \cdots m_{n}!} \prod_{j=1}^{n} \frac{1}{j!} \overset{\1}{\leq} \bigg( \sum_{j=1}^{n} \frac{1}{j!} \bigg)^{m_{1}+\cdots+m_{n}} \leq 2^{n},
\end{align*}
where in step $\1$ we use the multinomial theorem and in the last step we use $m_{1}+\cdots+m_{n} \leq n$. Putting the pieces together and applying the triangle inequality yields for all $n \in \mathbb{N}$ and $x\in \real$,
\begin{align*}
  |f^{(n)}(x)| \leq 2n! \sum_{m_{1},\dots,m_{n}} 2^{n} \leq 2 \cdot n! \cdot (2e)^{n},
\end{align*}
where in the last step we use that there are at most $e^{n}$ terms in the summation. Consequently, by letting $R = 8e$, we obtain
\[
   \sum_{j=1}^{n} \frac{|f^{(j)}(x)|}{R^{j}j!} \leq \sum_{j=1}^{n} 2^{-j} \leq 1.
\]
Note that one can evaluate $f^{(n)}(x)$ using Eq.~\eqref{eq:derivative-sech}, where one need to sum over at most $\mathcal{O}\Big(\exp\Big(\pi\sqrt{2n/3}\Big)\Big)$ terms and each term in the summation can be computed in time $\mathcal{O}(\mathsf{poly}(n))$. Thus, the evaluation time of $f^{(n)}(x)$ can be bounded by 
\[
\mathcal{O}\Big( \mathsf{poly}(n) \cdot \exp \Big(\pi\sqrt{2n/3} \Big)  \Big) = \mathcal{O}\Big( \exp \Big(\pi\sqrt{n} \Big)  \Big).
\]

\paragraph{Other functions.} For $f(\theta) = (1+\theta^{2})^{-1}$ and $f(\theta) = (1+e^{-\theta})^{-1}$. we can follow the identical steps as above to show $f \in \mathcal{F}(8e)$. For $f(\theta) = \exp(-\theta^{2})$, we obtain 
\[
    f^{(n)}(x) = (-1)^{n} H_{n}(x)\exp(-x^{2}).
\]
By using the properties of Hermite polynomials~\eqref{Hermite-boundness}, it is straightforward to show $f \in \mathcal{F}(8e)$. We omit the details for brevity. \qed


\section{Proofs of results from Section~\ref{sec:applications}}\label{sec:proof_applications}
In this section, we prove all the results related to the applications presented in this paper.

\subsection{Proof of Proposition~\ref{thm:test-hardness-pca-gaussian}}\label{sec:pf_thm_hardness_pca_gaussian}
We first state a result which directly follows by combining Theorem~10 and Lemma~102 of~\citet{brennan2020reducibility}.
\begin{proposition}\label{prop:TPCA-hardness-BB20}
Let $s\geq 3$ be a constant. Let 
\begin{align}\label{tensor-BB20}
	\widetilde{T} = \beta (\sqrt{n} \widetilde{v})^{\otimes s} + \widetilde{\zeta}, \quad \text{where } \widetilde{v} \sim \mathsf{Unif}\big(\big\{\tfrac{-1}{\sqrt{n}}, \tfrac{1}{\sqrt{n}} \big\}^{n} \big) \quad \text{and} \quad \widetilde{\zeta} \sim \NORMAL(0,1)^{\otimes n^{\otimes s}}.
\end{align}
Suppose the $k$-$\mathrm{HPC}^{s}$ conjecture holds. If $\beta = \tilde{o}\big(n^{-s/4}\big)$ and $\beta = \omega\big(n^{-s/2} \sqrt{s\log(n)}\big)$, then there exists \emph{no} algorithm $\mathcal{A}': \real^{n^{\otimes s}} \rightarrow \mathbb{S}^{d-1}$ that runs in time $\mathsf{poly}(n)$ and satisfies
\begin{align*}
    \big\langle \mathcal{A}' \big(\widetilde{T}\big), \widetilde{v} \big \rangle = \Omega(1) \quad \text{with probability at least} \quad \Omega_{n}(1)
\end{align*}
\end{proposition}
We next present an auxiliary lemma establishing that there exists a polynomial-time reduction that maps $\widetilde{T}$ defined in Eq.~\eqref{tensor-BB20} to the tensor $T$ defined in Eq.~\eqref{tensor-continuous-spike} and achieves zero TV deficiency.

\begin{lemma}\label{lemma:reduction-tensor-BB20}
Let $\widetilde{T}$ denote the tensor defined in Eq.~\eqref{tensor-BB20} with signal strength $\beta \in \real$ and spike $\sqrt{n} \widetilde{v}$. Then there exists a reduction $\mathcal{R}:\real^{n^{\otimes s}} \rightarrow \real^{n^{\otimes s}}$ that runs in time $\mathsf{poly}(n)$ and satisfies
\[
	 \mathcal{R}(\widetilde{T}) = \beta\big(\sqrt{n} U \widetilde{v}\big)^{\otimes s} + \zeta,
\]
where $U$ is an $n \times n$ unitary matrix chosen uniformly at random and $\zeta$ is as defined in Eq.~\eqref{def-tensor-gaussian-noise}.
\end{lemma}
We provide the proof Lemma~\ref{lemma:reduction-tensor-BB20} in Appendix~\ref{sec:pf_claims_TPCA}.
Note that $U \widetilde{v} \sim \mathsf{Unif}\big( \mathbb{S}^{n-1} \big)$ and thus $\mathcal{R}\big(\widetilde{T}\big)$ is equal in distribution to the tensor $T$ defined in Eq.~\eqref{tensor-continuous-spike}.

We now prove Proposition~\ref{thm:test-hardness-pca-gaussian} by combining Proposition~\ref{prop:TPCA-hardness-BB20} with Lemma~\ref{lemma:reduction-tensor-BB20}. We prove it by contradiction. Let $T$ be defined as in Eq.~\eqref{tensor-continuous-spike} and assume that there is an algorithm $\mathcal{A}: \real^{n^{\otimes s}} \rightarrow \mathbb{S}^{d-1}$ that runs in time $\mathsf{poly}(n)$ and satisfies
\begin{align*}
    \big\langle \mathcal{A} \big(T\big), v \big \rangle = \Omega(1) \quad \text{with probability at least} \quad \Omega_{n}(1).
\end{align*}
Then the algorithm $\mathcal{A} \circ \mathcal{R}: \real^{n^{\otimes s}} \rightarrow \mathbb{S}^{d-1}$ runs in time $\mathsf{poly}(n)$ and satisfies
\begin{align*}
    \big\langle \mathcal{A} \circ \mathcal{R} \big( \widetilde{T}\big), U \widetilde{v}  \big \rangle = \Omega(1) \quad \text{with probability at least} \quad \Omega_{n}(1).
\end{align*}
This implies that 
\begin{align*}
    \big\langle U^{\top} \mathcal{A} \circ \mathcal{R} \big( \widetilde{T}\big),  \widetilde{v}  \big \rangle = \Omega(1) \quad \text{with probability at least} \quad \Omega_{n}(1).
\end{align*}
In words, the algorithm $U^{\top}\mathcal{A} \circ \mathcal{R}: \real^{n^{\otimes s}} \rightarrow \mathbb{S}^{d-1}$ runs in time $\mathsf{poly}(n)$ and produces accurate recovery. However, this violates Proposition~\ref{prop:TPCA-hardness-BB20}. Consequently, our assumption that algorithm $\mathcal{A}$ runs in polynomial time must be false, thereby proving Proposition~\ref{thm:test-hardness-pca-gaussian}.
\qed

\subsection{Proof of Theorem~\ref{thm:tensor-pca-gaussian-nongaussian}}\label{sec:pf-thm-tensor-pca}
We partition the set of all index tuples
\[
\big\{(i_{1}, \dots, i_{s}) \; \big| \; i_{k} \in [n], \; \forall\, k \in [s] \big\}
\]
into disjoint subsets $S_1, S_2, \dots, S_L$, where each subset $S_\ell$ consists of index tuples that are equivalent up to permutation—that is, any two tuples in $S_k$ are permutations of each other. Let $\mathsf{K}:\real \rightarrow \real$ be the reduction in Theorem~\ref{thm:non-gaussian}. We design the tensor-to-tensor reduction algorithm $\mathcal{R}: \real^{n ^{\otimes s}} \rightarrow \real^{n ^{\otimes s}}$ as follows.\\
Input: $T \in \real^{n ^{\otimes s}}$. For each $\ell = 1,2,\dots,L$:  
\begin{enumerate}
\item Pick one index tuple $\big(i_{1}^{\ell},\dots,i_{s}^{\ell} \big) \in S_{\ell}$. Let $\mathcal{R}(T)_{i_{1}^{\ell},\dots,i_{s}^{\ell}} = \mathsf{K}\big( T_{i_{1}^{\ell},\dots,i_{s}^{\ell}} \big)$.
\item For all other index tuples $(i'_{1},\dots,i'_{s}) \in S_{\ell} \setminus \big\{ \big(i_{1}^{\ell},\dots,i_{s}^{\ell} \big) \big\}$, set $\mathcal{R}(T)_{(i'_{1},\dots,i'_{s})} = \mathcal{R}(T)_{i_{1}^{\ell},\dots,i_{s}^{\ell}}$ to ensure the output is a symmetric tensor.
\end{enumerate}
We next bound the TV distance between $\mathcal{R}(T(\beta,v))$ and $\widetilde{T}(\beta,v)$, where we use $T(\beta,v)$ to denote the tensor from model~\eqref{tensor-pca-gaussian} when $\beta$ and $v$ are the fixed parameters, and likewise $\widetilde{T}(\beta,v)$ is the tensor from model~\eqref{tensor-pca-nongaussian}. Recall the density of $\mathcal{Q}_{0}$ in Eq.~\eqref{density-Q0} and the density of $\mathcal{Q}_{\theta}$ in Eq.~\eqref{density-non-gaussian-target}. Conditional on $\beta$ and $v$, we have by definition for all $\ell \in [L]$ that
\[
    T(\beta, v)_{i_{1}^{\ell}, \dots, i_{s}^{\ell}} \sim \mathcal{N}(\theta, \sigma_{i_{1},\dots,i_{s}}^{2})
    \quad \text{and} \quad
    \widetilde{T}(\beta, v)_{i_{1}^{\ell}, \dots, i_{s}^{\ell}} \sim \mathcal{Q}_{\theta}, 
    \quad \text{where } \theta = \beta n^{s/2} \times v_{i_{1}^{\ell}} \times v_{i_{2}^{\ell}} \times \cdots \times v_{i_{s}^{\ell}}.
\]
Note that $\sigma_{i_{1},\dots,i_{s}}$ is a constant, since $s$ is a fixed parameter and $1/s!\leq \sigma_{i_{1},\dots,i_{s}}^{2} \leq 1$.
Using Assumption~\ref{assump:noise-tensor-pca} and Theorem~\ref{thm:non-gaussian}, we then obtain our reduction. In particular, for all $\epsilon \in (0,1)$, if $\varT =  \Theta \big(\mathsf{polylog}(1/\epsilon) \big)$ then setting $T_{\mathsf{iter}} =  \Theta \big( \log(1/\epsilon) \big)$ yields, for all $\ell \in [L]$, that
\[
	\sup_{\beta \in \real, v\in \real^{n}}\mathsf{d}_{\mathsf{TV}}\Big( \mathsf{K} \big( T(\beta, v)_{i_{1}^{\ell}, \dots, i_{s}^{\ell} } \big),  \widetilde{T}(\beta, v)_{i_{1}^{\ell}, \dots, i_{s}^{\ell}}  \Big) \leq \epsilon.
\]
Applying the triangle inequality for the TV distance (or equivalently, a union bound)~\citep[Lemma 6]{brennan2018reducibility}, we obtain
\begin{align*}
	\mathsf{d}_{\mathsf{TV}}\Big( \mathcal{R}\big(T(\beta,v)\big), \widetilde{T}(\beta,v)  \Big) &\leq  \sum_{\ell=1}^{L} \mathsf{d}_{\mathsf{TV}}\Big( \mathcal{R}\big(T(\beta,v)\big)_{i_{1}^{\ell},\dots,i_{s}^{\ell}}, \widetilde{T}(\beta,v)_{i_{1}^{\ell},\dots,i_{s}^{\ell}} \Big) \\
	&= \sum_{\ell=1}^{L} \mathsf{d}_{\mathsf{TV}}\Big( \mathsf{K} \big( T(\beta, v)_{i_{1}^{\ell}, \dots, i_{s}^{\ell} } \big) , \widetilde{T}(\beta,v)_{i_{1}^{\ell},\dots,i_{s}^{\ell}} \Big) \\
	& \leq L \cdot \epsilon \leq n^{s} \cdot \epsilon.
\end{align*}
For each $\delta \in (0,1)$, let $\epsilon = \delta/n^{s}$. Then taking the supremum over $\beta\in \real$ and $v \in \real^{n}$ of the inequality in the display above completes the proof of Theorem~\ref{thm:tensor-pca-gaussian-nongaussian}. Note that from Assumption~\ref{assump:noise-tensor-pca} and Theorem~\ref{thm:non-gaussian}, the reduction algorithm $\mathcal{R}$ runs in time 
\[
	\mathcal{O}\Big( n^{s} \cdot T_{\mathsf{samp}} \cdot T_{\mathsf{iter}} \cdot T_{\mathsf{eval}}(N) \Big) = \mathcal{O}\Big( \mathsf{poly}\big(n^{s}/\delta\big) \Big),
\]
where $T_{\mathsf{samp}}$ is independent of $\delta$ and $n$.
\qed

\subsection{Proof of Corollary~\ref{thm:test-hardness-pca-nongaussian}}\label{sec:pf-test-hardness-pca-nongaussian}
Let $T$ be a tensor defined as in Eq.~\eqref{tensor-continuous-spike} and let $\widetilde{T}$ be the tensor defined in Theorem~\ref{thm:test-hardness-pca-nongaussian}. Note that $T$ and $\widetilde{T}$ both have the same signal strength $\beta \in \real$ and spike $v \in \mathbb{S}^{n-1}$. Let $\mathcal{R}$ be the reduction defined in Proposition~\ref{thm:tensor-pca-gaussian-nongaussian} with TV deficiency $\delta = n^{-1}$. By Proposition~\ref{thm:tensor-pca-gaussian-nongaussian}, we have that $\mathcal{R}$ runs in time $\mathsf{poly}(n)$ and satisfies $\mathsf{d}_{\mathsf{TV}}\big( \mathcal{R}(T), \widetilde{T}\big) \leq n^{-1}$. We now prove Corollary~\ref{thm:test-hardness-pca-nongaussian} by contradiction. Assume that there exists an algorithm $\mathcal{A}:\real^{n^{\otimes s}} \rightarrow \mathbb{S}^{n-1}$ that runs in time $\mathsf{poly}(n)$ and satisfies
\[
		\big \langle \mathcal{A}\big( \widetilde{T} \big) , v\big \rangle = \Omega(1) \quad \text{with probability at least} \quad \Omega_{n}(1).
\]
Then $\mathcal{A} \circ \mathcal{R}: \real^{n^{\otimes s}} \rightarrow \mathbb{S}^{n-1}$ also runs in time $\mathsf{poly}(n)$ and satisfies
\[
		\big \langle \mathcal{A} \circ \mathcal{R}\big( T \big) , v\big \rangle = \Omega(1) \quad \text{with probability at least} \quad \Omega_{n}(1) - n^{-1}.
\]
However, this contradicts Proposition~\ref{thm:test-hardness-pca-gaussian}. Thus, our assumption that the algorithm $\mathcal{A}$ runs in polynomial time cannot hold, thereby proving Corollary~\ref{thm:test-hardness-pca-nongaussian}.
\qed

\subsection{Proof of Theorem~\ref{thm:reduction-mslr-glm}}\label{sec:pf_reduction_mslr_glm}
Let $\mathsf{K}: \real \rightarrow \real$ be the reduction in Theorem~\ref{thm:gaussian-mean-nonlinear}.
We design the reduction algorithm $\mathcal{R}: (\real^{d} \times \real)^{n} \rightarrow (\real^{d} \times \real)^{n}$ as follows: For inputs $\{(x_{i},y_{i})\}_{i=1}^{n}$ from mixtures of linear regressions, return
$\mathcal{R}\big( \{(x_{i},y_{i})\}_{i=1}^{n} \big) := \{(x_{i}, \mathsf{K}(y_{i})) \}_{i=1}^{n}$. Note that by definition, conditioned on $x_{i}$ and $R_{i}$,
\[
	y_{i} \sim \NORMAL(\theta,1) \quad \text{and} \quad \widetilde{y}_{i} \sim \NORMAL(f(\theta/\varT),1) \quad \text{where } \theta = R_{i} \cdot \langle x_{i}, v \rangle.
\]
Consequently, applying Assumption~\ref{assump:link-linear-model} and Theorem~\ref{thm:gaussian-mean-nonlinear}, we obtain that for each $\epsilon \in (0,1)$,
\[
	\text{if} \quad \varT = \mathsf{polylog}(1/\epsilon) + 2 \max_{i \in [n]}|\langle x_{i}, v \rangle|, \quad \text{then} \quad \mathsf{d}_{\mathsf{TV}}\big(\mathsf{K}(y_{i}),\widetilde{y}_{i} \big) \leq \epsilon.
\]
Applying the triangle inequality for the TV distance~\citep[Lemma 6]{brennan2018reducibility}, we obtain
\begin{align*}
	\mathsf{d}_{\mathsf{TV}} \Big( \mathcal{R}\big( \{(x_{i},y_{i})\}_{i=1}^{n} \big), \{(x_{i},\widetilde{y}_{i})\}_{i=1}^{n}  \Big) &= \mathsf{d}_{\mathsf{TV}} \Big(  \big\{\big(x_{i}, \mathsf{K}(y_{i} ) \big) \big\}_{i=1}^{n} , \{(x_{i},\widetilde{y}_{i})\}_{i=1}^{n}  \Big) \\
	& \leq \sum_{i=1}^{n}  \mathsf{d}_{\mathsf{TV}}\big(\mathsf{K}(y_{i}),\widetilde{y}_{i} \big) \leq n \epsilon.
\end{align*}
The desired TV guarantee follows by letting $\epsilon = \delta/n$. From Assumption~\ref{assump:link-linear-model} and Theorem~\ref{thm:gaussian-mean-nonlinear}, the total time taken to compute $\mathcal{R}$ is 
\begin{align*}	
	\mathcal{O}\Big( n \cdot \log(4n/\delta) \cdot T_{\mathsf{eval}}(N) \Big) = \mathcal{O}\Big( \mathsf{poly}(n/\delta) \Big),
\end{align*}
as claimed.
\qed

\subsection{Proof of Corollary~\ref{thm:hardness-glm}}\label{sec:pf_thm_hardness_glm}
We define the following event 
\[
	\mathcal{E} = \Big\{ \max_{ i \in [n]} \big|\langle x_{i}, v \rangle \big| \leq 10 \log^{2}(n) \Big\}.
\]
Note that $\langle x_{i},v \rangle \overset{\mathsf{i.i.d.}}{\sim} \NORMAL(0,\eta^{2})$ since $\|v\|_{2} = \eta$ and $x_{i} \overset{\mathsf{i.i.d.}}{\sim} \NORMAL(0,I_{d})$. Moreover, by the condition that $k = o(n^{1/6})$ and $n = \tilde{o}(k^{2}/ \eta^{4})$, we obtain $\eta = \widetilde{o}(k^{-1}) \leq 1$. Consequently, applying a standard bound on Gaussian maxima yields 
\begin{align}\label{high-prob-event}
	\mathbb{P}\big[\mathcal{E} \big] \geq 1 - n^{-3}.
\end{align}
By setting $\varT = \mathsf{polylog}(n/\delta)$ for any $\delta \in (0,1)$, the reduction $\mathcal{R}$ in Proposition~\ref{thm:reduction-mslr-glm} runs in time $\mathsf{poly}(n/\delta)$. The TV-deficiency guarantee in Proposition~\ref{thm:reduction-mslr-glm} requires condition~\eqref{tau-condition-prop-MSLR-GLM}, which holds with high probability according to Ineq.~\eqref{high-prob-event}. We now combine the high probability guarantee~\eqref{high-prob-event} with Proposition~\ref{thm:reduction-mslr-glm} to derive the desired TV-deficiency guarantee.

Let $\mathsf{K}:\real \rightarrow \real$ be the reduction in Theorem~\ref{thm:gaussian-mean-nonlinear} and note that
\[
		y_{i} \mid x_{i} \sim \NORMAL\big(R_{i} \cdot \langle x_{i}, v\rangle, 1\big) \quad \text{and} \quad \widetilde{y}_{i} \mid x_{i} \sim \NORMAL\big(f\big(\langle x_{i}, v\rangle/\varT\big), 1 \big).
\]
To reduce the notational burden, we let $\rho:\real^{d} \rightarrow \real$ denote the density function of $x_{i} \sim \NORMAL(0,I_{d})$, let $\gamma(\cdot;x), \widetilde{\gamma}(\cdot ;x):\real \rightarrow \real$ denote the conditional densities of $\mathsf{K}(y_{i}) \mid x_{i} = x$ and $\widetilde{y}_{i} \mid x_{i} = x$, respectively. Then we have
\begin{align*}
	&\mathsf{d}_{\mathsf{TV}} \left( \mathcal{R}\big( \{(x_{i},y_{i})\}_{i=1}^{n} \big), \{(x_{i},\widetilde{y}_{i})\}_{i=1}^{n}  \right) \\
	& \quad \overset{\1}{\leq} n \cdot \mathsf{d}_{\mathsf{TV}} \big( (x_{i}, \mathsf{K}(y_{i}) ), (x_{i}, \widetilde{y_{i}}) \big) \\
	& \quad = \frac{n}{2} \int_{\real^{d}} \rho(x) \cdot \int_{\real}\Big| \gamma(t;x) - \widetilde{\gamma}(t;x)\Big| \mathrm{d}t \, \mathrm{d}x \\
	& \quad = \frac{n}{2} \int_{\mathcal{E}} \rho(x) \cdot \int_{\real}\Big| \gamma(t;x) - \widetilde{\gamma}(t;x)\Big| \mathrm{d}t \, \mathrm{d}x + \frac{n}{2} \int_{\mathcal{E}^{\complement}} \rho(x) \cdot \int_{\real}\Big| \gamma(t;x) - \widetilde{\gamma}(t;x)\Big| \mathrm{d}t\, \mathrm{d}x,
\end{align*}
where step $\1$ follows since the samples are $\mathsf{i.i.d.}$ Towards bounding the two terms in the display above, we obtain
\begin{align*}
	\frac{n}{2}\int_{\mathcal{E}} \rho(x) \cdot \int_{\real}\Big| \gamma(t;x) - \widetilde{\gamma}(t;x)\Big| \mathrm{d}t \, \mathrm{d}x \leq n \cdot \int_{\mathcal{E}} \rho(x) \mathrm{d}x \cdot \frac{\delta}{n} \leq \delta,
\end{align*}
where we use the guarantee that on event $\mathcal{E}$, $\mathrm{d}_{\mathsf{TV}} ( \mathsf{K}(y_{i}), \widetilde{y_{i}}) \leq \delta/n$. Moreover, we obtain
\begin{align*}
\frac{n}{2} \int_{\mathcal{E}^{\complement}} \rho(x) \cdot \int_{\real}\Big| \gamma(t;x) - \widetilde{\gamma}(t;x)\Big| \mathrm{d}t\, \mathrm{d}x \leq n \cdot \int_{\mathcal{E}^{\complement}} \rho(x) \mathrm{d}x \leq n \cdot \mathbb{P} \Big[ \mathcal{E}^{\complement} \Big] \leq n^{-2}.
\end{align*}
Putting the three pieces together and setting $\delta = n^{-1}$ yields
\begin{align}\label{TV-H1-glm}
	\mathsf{d}_{\mathsf{TV}} \Big( \mathcal{R}\big( \{(x_{i},y_{i})\}_{i=1}^{n} \big), \{(x_{i},\widetilde{y}_{i})\}_{i=1}^{n}  \Big) \leq n^{-1}+n^{-2} \leq 2n^{-1}.
\end{align}

For the sake of contradiction, we assume that there exists an algorithm $\mathcal{A}: (\real^{d} \times \real)^{n} \rightarrow \real^{d}$ that runs in time $\mathsf{poly}(n)$ and for sample size $n = \tilde{o}\big(k^{2}/\norm^{4} \big)$, satisfies
\begin{align}\label{assump-guarantee-GLM}
	\big\| \mathcal{A}\big( \{x_{i},\widetilde{y}_{i}\}_{i=1}^{n} \big) - v  \big\|_{2} = o(\norm) \quad \text{with probability at least} \quad 1 - o_{n}(1).
\end{align}
Then we note that $\mathcal{A}\circ \mathcal{R}: (\real^{d} \times \real)^{n} \rightarrow \real^{d}$ also runs in polynomial time and it is an estimation algorithm for data $\{x_{i},y_{i}\}_{i=1}^{n}$ from $\MixSLR$~\eqref{model-MSLR}. Moreover, the TV-deficiency guarantee~\eqref{TV-H1-glm} and the guarantee of $\mathcal{A}$~\eqref{assump-guarantee-GLM} imply
\[
	\big\| \mathcal{A} \circ \mathcal{R}\big( \{x_{i},y_{i}\}_{i=1}^{n} \big) - v  \big\|_{2} = o(\norm) \quad \text{with probability at least} \quad 1 - o_{n}(1) - 2n^{-1} = 1 - o_{n}(1).
\]
This contradicts the computational hardness claimed in Proposition~\ref{thm:hardness-MSLR}. Thus, our assumption that algorithm $\mathcal{A}$ runs in polynomial time must be false, thereby proving Corollary~\ref{thm:hardness-glm}.
\qed

\subsection{Proof of Theorem~\ref{thm:reduction-ROS-BC}} \label{sec:pf_thm_reduction_ROS_BC}
We construct the reduction $\mathcal{R}: \mathbb{R}^{n \times n} \rightarrow \mathbb{R}^{n \times n}$ using the scalar reduction $\mathsf{K}: \mathbb{R} \rightarrow \mathbb{R}$ developed in Corollary~\ref{example:monomial} with the parameter setting $B = 2$.

\medskip
\noindent
\textbf{Input:} A matrix $M = \mu r c^{\top} + \NORMAL(0,1)^{\otimes n \times n}$ for $r,c\in S_{n,k}$.\\
\noindent
\textbf{Procedure:} Construct a new matrix $\widetilde{M} \in \mathbb{R}^{n \times n}$:
For all $i,j \in [n]$, set $\widetilde{M}_{ij} = \mathsf{K}(M_{ij})$.\\ 
\noindent
\textbf{Output:} The reduced matrix $\mathcal{R}(M) := \widetilde{M}$.

We next bound the TV deficiency of the reduction $\mathcal{R}$. For notational convenience, we let $f:\real \rightarrow \real$ be the square function, i.e., $f(t) = t^{2}$ for all $t\in \real$ and let
\[
		\overline{M} : = \frac{\mu^{2}}{\varT^{2} k} |r||c|^{\top} + \NORMAL(0,1)^{\otimes n \times n},
\]
where $|r|,|c|\in \real^{n}$ denotes the vectors obtained by taking the entrywise absolute value of $r$ and $c$, respectively.
Note that for all $r,c\in S_{n,k}$, we have 
\[
 	r_{i}^{2}  c_{j}^{2} = \frac{ |r_{i}| \cdot |c_{j}| }{k} \quad \text{for all} \quad i,j \in [n].
\]
Then observe that for all $i,j \in [n]$,
\[
	M_{ij} \sim \NORMAL(\mu r_{i}c_{j},1) \quad \text{and} \quad \overline{M}_{ij} \sim \NORMAL\Big( f\Big( \frac{\mu r_{i}c_{j}}{\varT} \Big),1 \Big).
\]
Consequently, since $\big(\mu r_{i}c_{j} \big)^{2} \leq \mu^{2}/k^{2}$, applying Corollary~\ref{example:monomial} with $\epsilon = \delta/n^{2}$ for any $\delta\in(0,1)$ and using the condition on $\varT$~\eqref{condition-tau-ROS-BC} yield the TV guarantee
\[
		\mathsf{d}_{\mathsf{TV}} \Big( \mathsf{K}(M_{ij}), \overline{M}_{ij} \Big) \leq \frac{\delta}{n^{2}} \quad \text{for all }i,j \in [n].
\]
Putting the pieces together yields
\begin{align*}
	\mathsf{d}_{\mathsf{TV}} \Big( \mathcal{R}(M), \overline{M}  \Big) \leq \sum_{i=1}^{n} \sum_{j=1}^{n} \mathsf{d}_{\mathsf{TV}} \Big( \mathsf{K}(M_{ij}), \overline{M}_{ij}  \Big) \leq \delta. 
\end{align*}
Note that the bound in the display is independent of $r,c \in S_{n,k}$ and $\mu \in \real$. Taking the supremum over $r,c \in S_{n,k}$ and $\mu \in \real$ on both sides yields the desired TV deficiency guarantee. 
The computational complexity of the reduction $\mathcal{R}$ follows from the computational complexity of the scalar reduction $\mathsf{K}$ in Corollary~\ref{example:monomial}, and we use this reduction $n^{2}$ times in total. \qed

\small
\bibliographystyle{abbrvnat}
\bibliography{refs}

\begin{thebibliography}{86}
\providecommand{\natexlab}[1]{#1}
\providecommand{\url}[1]{\texttt{#1}}
\expandafter\ifx\csname urlstyle\endcsname\relax
  \providecommand{\doi}[1]{doi: #1}\else
  \providecommand{\doi}{doi: \begingroup \urlstyle{rm}\Url}\fi

\bibitem[Andoni et~al.(2017)Andoni, Hsu, Shi, and
  Sun]{andoni2017correspondence}
A.~Andoni, D.~Hsu, K.~Shi, and X.~Sun.
\newblock Correspondence retrieval.
\newblock In \emph{Conference on Learning Theory}, pages 105--126. PMLR, 2017.

\bibitem[Arpino and Venkataramanan(2023)]{arpino2023statistical}
G.~Arpino and R.~Venkataramanan.
\newblock {S}tatistical-computational tradeoffs in mixed sparse linear
  regression.
\newblock In \emph{The Thirty Sixth Annual Conference on Learning Theory},
  pages 921--986. PMLR, 2023.

\bibitem[Axelrod et~al.(2024)Axelrod, Garg, Han, Sharan, and
  Valiant]{axelrod2024statistical}
B.~Axelrod, S.~Garg, Y.~Han, V.~Sharan, and G.~Valiant.
\newblock On the statistical complexity of sample amplification.
\newblock \emph{The Annals of Statistics}, 52\penalty0 (6):\penalty0
  2767--2790, 2024.

\bibitem[Bandeira et~al.(2018)Bandeira, Perry, and Wein]{bandeira2018notes}
A.~Bandeira, A.~Perry, and A.~S. Wein.
\newblock Notes on computational-to-statistical gaps: predictions using
  statistical physics.
\newblock \emph{Portugaliae mathematica}, 75\penalty0 (2):\penalty0 159--186,
  2018.

\bibitem[Bangachev et~al.(2025)Bangachev, Bresler, Tiegel, and
  Vaikuntanathan]{bangachev2025near}
K.~Bangachev, G.~Bresler, S.~Tiegel, and V.~Vaikuntanathan.
\newblock Near-optimal time-sparsity trade-offs for solving noisy linear
  equations.
\newblock In \emph{Proceedings of the 57th Annual ACM Symposium on Theory of
  Computing}, pages 1910--1920, 2025.

\bibitem[Barbier et~al.(2019)Barbier, Krzakala, Macris, Miolane, and
  Zdeborov{\'a}]{barbier2019optimal}
J.~Barbier, F.~Krzakala, N.~Macris, L.~Miolane, and L.~Zdeborov{\'a}.
\newblock Optimal errors and phase transitions in high-dimensional generalized
  linear models.
\newblock \emph{Proceedings of the National Academy of Sciences}, 116\penalty0
  (12):\penalty0 5451--5460, 2019.

\bibitem[Ben~Arous et~al.(2019)Ben~Arous, Mei, Montanari, and
  Nica]{arous2019landscape}
G.~Ben~Arous, S.~Mei, A.~Montanari, and M.~Nica.
\newblock The landscape of the spiked tensor model.
\newblock \emph{Communications on Pure and Applied Mathematics}, 72\penalty0
  (11):\penalty0 2282--2330, 2019.

\bibitem[Ben~Arous et~al.(2020)Ben~Arous, Gheissari, and
  Jagannath]{arous2020algorithmic}
G.~Ben~Arous, R.~Gheissari, and A.~Jagannath.
\newblock Algorithmic thresholds for tensor {PCA}.
\newblock \emph{The Annals of Probability}, 48\penalty0 (4):\penalty0
  2052--2087, 2020.

\bibitem[Ben~Arous et~al.(2023)Ben~Arous, Wein, and Zadik]{arous2023free}
G.~Ben~Arous, A.~S. Wein, and I.~Zadik.
\newblock Free energy wells and overlap gap property in sparse {PCA}.
\newblock \emph{Communications on Pure and Applied Mathematics}, 76\penalty0
  (10):\penalty0 2410--2473, 2023.

\bibitem[Berthet and Rigollet(2013{\natexlab{a}})]{berthet2013complexity}
Q.~Berthet and P.~Rigollet.
\newblock Complexity theoretic lower bounds for sparse principal component
  detection.
\newblock In \emph{Conference on learning theory}, pages 1046--1066. PMLR,
  2013{\natexlab{a}}.

\bibitem[Berthet and Rigollet(2013{\natexlab{b}})]{berthet2013sPCA}
Q.~Berthet and P.~Rigollet.
\newblock {Optimal detection of sparse principal components in high dimension}.
\newblock \emph{The Annals of Statistics}, 41\penalty0 (4):\penalty0 1780 --
  1815, 2013{\natexlab{b}}.
\newblock \doi{10.1214/13-AOS1127}.
\newblock URL \url{https://doi.org/10.1214/13-AOS1127}.

\bibitem[Blackwell(1951)]{blackwell1951comparison}
D.~Blackwell.
\newblock Comparison of experiments.
\newblock In \emph{Proceedings of the Second Berkeley Symposium on Mathematical
  Statistics and Probability}, volume~2, pages 93--103. University of
  California Press, 1951.

\bibitem[Blackwell(1953)]{blackwell1953equivalent}
D.~Blackwell.
\newblock Equivalent comparisons of experiments.
\newblock \emph{The Annals of Mathematical Statistics}, pages 265--272, 1953.

\bibitem[Brennan and Bresler(2019)]{brennan2019optimal}
M.~Brennan and G.~Bresler.
\newblock Optimal average-case reductions to sparse {PCA}: From weak
  assumptions to strong hardness.
\newblock In \emph{Conference on Learning Theory}, pages 469--470. PMLR, 2019.

\bibitem[Brennan and Bresler(2020)]{brennan2020reducibility}
M.~Brennan and G.~Bresler.
\newblock Reducibility and statistical-computational gaps from secret leakage.
\newblock In \emph{Conference on Learning Theory}, pages 648--847. PMLR, 2020.

\bibitem[Brennan et~al.(2018)Brennan, Bresler, and
  Huleihel]{brennan2018reducibility}
M.~Brennan, G.~Bresler, and W.~Huleihel.
\newblock Reducibility and computational lower bounds for problems with planted
  sparse structure.
\newblock In \emph{Conference On Learning Theory}, pages 48--166. PMLR, 2018.

\bibitem[Brennan et~al.(2019)Brennan, Bresler, and
  Huleihel]{brennan2019universality}
M.~Brennan, G.~Bresler, and W.~Huleihel.
\newblock Universality of computational lower bounds for submatrix detection.
\newblock In \emph{Conference on Learning Theory}, pages 417--468. PMLR, 2019.

\bibitem[Brennan et~al.(2021)Brennan, Bresler, Hopkins, Li, and
  Schramm]{brennan2021statistical}
M.~Brennan, G.~Bresler, S.~Hopkins, J.~Li, and T.~Schramm.
\newblock Statistical query algorithms and low degree tests are almost
  equivalent.
\newblock In \emph{Conference on Learning Theory}, pages 774--774. PMLR, 2021.

\bibitem[Bresler and Harbuzova(2025)]{bresler2025computational}
G.~Bresler and A.~Harbuzova.
\newblock Computational equivalence of spiked covariance and spiked wigner
  models via gram-schmidt perturbation.
\newblock \emph{arXiv preprint arXiv:2503.02802}, 2025.

\bibitem[Bruna et~al.(2021)Bruna, Regev, Song, and Tang]{bruna2021continuous}
J.~Bruna, O.~Regev, M.~J. Song, and Y.~Tang.
\newblock Continuous {LWE}.
\newblock In \emph{Proceedings of the 53rd Annual ACM SIGACT Symposium on
  Theory of Computing}, pages 694--707, 2021.

\bibitem[Cai et~al.(2016)Cai, Li, and Ma]{cai2016optimal}
T.~Cai, X.~Li, and Z.~Ma.
\newblock Optimal rates of convergence for noisy sparse phase retrieval via
  thresholded {W}irtinger flow.
\newblock \emph{The Annals of Statistics}, 44\penalty0 (82):\penalty0
  2221--2251, 2016.

\bibitem[Celentano et~al.(2020)Celentano, Montanari, and
  Wu]{celentano2020estimation}
M.~Celentano, A.~Montanari, and Y.~Wu.
\newblock The estimation error of general first order methods.
\newblock In \emph{Conference on Learning Theory}, pages 1078--1141. PMLR,
  2020.

\bibitem[Chen(2019)]{chen2019phase}
W.-K. Chen.
\newblock Phase transition in the spiked random tensor with {R}ademacher prior.
\newblock \emph{The Annals of Statistics}, 47\penalty0 (5):\penalty0
  2734--2756, 2019.

\bibitem[Chen and Xu(2016)]{chen2016statistical}
Y.~Chen and J.~Xu.
\newblock Statistical-computational tradeoffs in planted problems and submatrix
  localization with a growing number of clusters and submatrices.
\newblock \emph{Journal of Machine Learning Research}, 17\penalty0
  (27):\penalty0 1--57, 2016.

\bibitem[Damian et~al.(2024)Damian, Pillaud-Vivien, Lee, and
  Bruna]{damian2024computational}
A.~Damian, L.~Pillaud-Vivien, J.~D. Lee, and J.~Bruna.
\newblock Computational-statistical gaps in {G}aussian single-index models.
\newblock \emph{arXiv preprint arXiv:2403.05529}, 2024.

\bibitem[Diakonikolas and Kane(2023)]{diakonikolas2023algorithmic}
I.~Diakonikolas and D.~M. Kane.
\newblock \emph{Algorithmic high-dimensional robust statistics}.
\newblock Cambridge university press, 2023.

\bibitem[Diakonikolas et~al.(2017)Diakonikolas, Kane, and
  Stewart]{diakonikolas2017statistical}
I.~Diakonikolas, D.~M. Kane, and A.~Stewart.
\newblock Statistical query lower bounds for robust estimation of
  high-dimensional {G}aussians and {G}aussian mixtures.
\newblock In \emph{2017 IEEE 58th Annual Symposium on Foundations of Computer
  Science (FOCS)}, pages 73--84. IEEE, 2017.

\bibitem[Ding et~al.(2024)Ding, Kunisky, Wein, and
  Bandeira]{ding2024subexponential}
Y.~Ding, D.~Kunisky, A.~S. Wein, and A.~S. Bandeira.
\newblock Subexponential-time algorithms for sparse {PCA}.
\newblock \emph{Foundations of Computational Mathematics}, 24\penalty0
  (3):\penalty0 865--914, 2024.

\bibitem[Dudeja and Hsu(2021)]{dudeja2021statistical}
R.~Dudeja and D.~Hsu.
\newblock Statistical query lower bounds for tensor pca.
\newblock \emph{Journal of Machine Learning Research}, 22\penalty0
  (83):\penalty0 1--51, 2021.

\bibitem[Dudeja and Hsu(2024)]{dudeja2024statistical}
R.~Dudeja and D.~Hsu.
\newblock Statistical-computational trade-offs in tensor pca and related
  problems via communication complexity.
\newblock \emph{The Annals of Statistics}, 52\penalty0 (1):\penalty0 131--156,
  2024.

\bibitem[Fan et~al.(2018)Fan, Liu, Wang, and Yang]{fan2018curse}
J.~Fan, H.~Liu, Z.~Wang, and Z.~Yang.
\newblock Curse of heterogeneity: {C}omputational barriers in sparse mixture
  models and phase retrieval.
\newblock \emph{arXiv preprint arXiv:1808.06996}, 2018.

\bibitem[Gamarnik et~al.(2021)Gamarnik, Jagannath, and
  Sen]{gamarnik2021overlap}
D.~Gamarnik, A.~Jagannath, and S.~Sen.
\newblock The overlap gap property in principal submatrix recovery.
\newblock \emph{Probability Theory and Related Fields}, 181\penalty0
  (4):\penalty0 757--814, 2021.

\bibitem[Goldstein and Reinert(2005)]{goldstein2005distributional}
L.~Goldstein and G.~Reinert.
\newblock Distributional transformations, orthogonal polynomials, and {S}tein
  characterizations.
\newblock \emph{Journal of Theoretical Probability}, 18\penalty0 (1):\penalty0
  237--260, 2005.

\bibitem[Gupte et~al.(2022)Gupte, Vafa, and
  Vaikuntanathan]{gupte2022continuous}
A.~Gupte, N.~Vafa, and V.~Vaikuntanathan.
\newblock Continuous {LWE} is as hard as {LWE} and applications to learning
  {G}aussian mixtures.
\newblock In \emph{2022 IEEE 63rd Annual Symposium on Foundations of Computer
  Science (FOCS)}, pages 1162--1173. IEEE, 2022.

\bibitem[Gupte et~al.(2024)Gupte, Vafa, and Vaikuntanathan]{gupte2024sparse}
A.~Gupte, N.~Vafa, and V.~Vaikuntanathan.
\newblock Sparse linear regression and lattice problems.
\newblock In \emph{Theory of Cryptography Conference}, pages 276--307.
  Springer, 2024.

\bibitem[Hajek et~al.(2015)Hajek, Wu, and Xu]{hajek2015computational}
B.~Hajek, Y.~Wu, and J.~Xu.
\newblock Computational lower bounds for community detection on random graphs.
\newblock In \emph{Conference on Learning Theory}, pages 899--928. PMLR, 2015.

\bibitem[Hajek et~al.(2018)Hajek, Wu, and Xu]{hajek2018submatrix}
B.~Hajek, Y.~Wu, and J.~Xu.
\newblock Submatrix localization via message passing.
\newblock \emph{Journal of Machine Learning Research}, 18\penalty0
  (186):\penalty0 1--52, 2018.

\bibitem[Han et~al.(2022)Han, Luo, Wang, and Zhang]{han2022exact}
R.~Han, Y.~Luo, M.~Wang, and A.~R. Zhang.
\newblock Exact clustering in tensor block model: Statistical optimality and
  computational limit.
\newblock \emph{Journal of the Royal Statistical Society Series B: Statistical
  Methodology}, 84\penalty0 (5):\penalty0 1666--1698, 2022.

\bibitem[Hansen and Torgersen(1974)]{hansen1974comparison}
O.~H. Hansen and E.~N. Torgersen.
\newblock Comparison of linear normal experiments.
\newblock \emph{The Annals of Statistics}, pages 367--373, 1974.

\bibitem[Hopkins(2018)]{hopkins2018statistical}
S.~Hopkins.
\newblock \emph{Statistical inference and the sum of squares method}.
\newblock PhD thesis, Cornell University, 2018.

\bibitem[Hopkins et~al.(2015)Hopkins, Shi, and Steurer]{hopkins2015tensor}
S.~B. Hopkins, J.~Shi, and D.~Steurer.
\newblock Tensor principal component analysis via sum-of-square proofs.
\newblock In \emph{Conference on Learning Theory}, pages 956--1006. PMLR, 2015.

\bibitem[Hopkins et~al.(2016)Hopkins, Schramm, Shi, and
  Steurer]{hopkins2016fast}
S.~B. Hopkins, T.~Schramm, J.~Shi, and D.~Steurer.
\newblock Fast spectral algorithms from sum-of-squares proofs: tensor
  decomposition and planted sparse vectors.
\newblock In \emph{Proceedings of the forty-eighth annual ACM symposium on
  Theory of Computing}, pages 178--191, 2016.

\bibitem[Hopkins et~al.(2017)Hopkins, Kothari, Potechin, Raghavendra, Schramm,
  and Steurer]{hopkins2017power}
S.~B. Hopkins, P.~K. Kothari, A.~Potechin, P.~Raghavendra, T.~Schramm, and
  D.~Steurer.
\newblock The power of sum-of-squares for detecting hidden structures.
\newblock In \emph{2017 IEEE 58th Annual Symposium on Foundations of Computer
  Science (FOCS)}, pages 720--731. IEEE, 2017.

\bibitem[Jagannath et~al.(2020)Jagannath, Lopatto, and
  Miolane]{jagannath2020statistical}
A.~Jagannath, P.~Lopatto, and L.~Miolane.
\newblock Statistical thresholds for tensor {PCA}.
\newblock \emph{The Annals of Applied Probability}, 30\penalty0 (4):\penalty0
  1910--1933, 2020.

\bibitem[Karlin and Rubin(1956)]{karlin1956theory}
S.~Karlin and H.~Rubin.
\newblock The theory of decision procedures for distributions with monotone
  likelihood ratio.
\newblock \emph{The Annals of Mathematical Statistics}, pages 272--299, 1956.

\bibitem[Kearns(1998)]{kearns1998efficient}
M.~Kearns.
\newblock Efficient noise-tolerant learning from statistical queries.
\newblock \emph{Journal of the ACM (JACM)}, 45\penalty0 (6):\penalty0
  983--1006, 1998.

\bibitem[Kim et~al.(2017)Kim, Bandeira, and Goemans]{kim2017community}
C.~Kim, A.~S. Bandeira, and M.~X. Goemans.
\newblock Community detection in hypergraphs, spiked tensor models, and
  sum-of-squares.
\newblock In \emph{2017 international conference on sampling theory and
  applications (sampta)}, pages 124--128. IEEE, 2017.

\bibitem[Kunisky(2025)]{kunisky2025low}
D.~Kunisky.
\newblock Low coordinate degree algorithms {I}: {U}niversality of computational
  thresholds for hypothesis testing.
\newblock \emph{The Annals of Statistics}, 53\penalty0 (2):\penalty0 774--801,
  2025.

\bibitem[Kunisky et~al.(2019)Kunisky, Wein, and Bandeira]{kunisky2019notes}
D.~Kunisky, A.~S. Wein, and A.~S. Bandeira.
\newblock Notes on computational hardness of hypothesis testing: {P}redictions
  using the low-degree likelihood ratio.
\newblock In \emph{ISAAC Congress (International Society for Analysis, its
  Applications and Computation)}, pages 1--50. Springer, 2019.

\bibitem[Lecu{\'e} and Mendelson(2015)]{lecue2015minimax}
G.~Lecu{\'e} and S.~Mendelson.
\newblock Minimax rate of convergence and the performance of empirical risk
  minimization in phase retrieval.
\newblock \emph{Electron. J. Probab}, 20\penalty0 (57):\penalty0 1--29, 2015.

\bibitem[Lehmann(1988)]{lehmann2011comparing}
E.~L. Lehmann.
\newblock {Comparing Location Experiments}.
\newblock \emph{The Annals of Statistics}, 16\penalty0 (2):\penalty0 521 --
  533, 1988.

\bibitem[Lesieur et~al.(2015)Lesieur, Krzakala, and
  Zdeborov{\'a}]{lesieur2015phase}
T.~Lesieur, F.~Krzakala, and L.~Zdeborov{\'a}.
\newblock Phase transitions in sparse {PCA}.
\newblock In \emph{2015 IEEE International Symposium on Information Theory
  (ISIT)}, pages 1635--1639. IEEE, 2015.

\bibitem[Lesieur et~al.(2017{\natexlab{a}})Lesieur, Krzakala, and
  Zdeborov{\'a}]{lesieur2017constrained}
T.~Lesieur, F.~Krzakala, and L.~Zdeborov{\'a}.
\newblock Constrained low-rank matrix estimation: {P}hase transitions,
  approximate message passing and applications.
\newblock \emph{Journal of Statistical Mechanics: Theory and Experiment},
  2017\penalty0 (7):\penalty0 073403, 2017{\natexlab{a}}.

\bibitem[Lesieur et~al.(2017{\natexlab{b}})Lesieur, Miolane, Lelarge, Krzakala,
  and Zdeborov{\'a}]{lesieur2017statistical}
T.~Lesieur, L.~Miolane, M.~Lelarge, F.~Krzakala, and L.~Zdeborov{\'a}.
\newblock Statistical and computational phase transitions in spiked tensor
  estimation.
\newblock In \emph{2017 ieee international symposium on information theory
  (isit)}, pages 511--515. IEEE, 2017{\natexlab{b}}.

\bibitem[Li and Schramm(2024)]{li2024some}
S.~Li and T.~Schramm.
\newblock Some easy optimization problems have the overlap-gap property.
\newblock \emph{arXiv preprint arXiv:2411.01836}, 2024.

\bibitem[Li and Voroninski(2013)]{li2013sparse}
X.~Li and V.~Voroninski.
\newblock Sparse signal recovery from quadratic measurements via convex
  programming.
\newblock \emph{SIAM Journal on Mathematical Analysis}, 45\penalty0
  (5):\penalty0 3019--3033, 2013.

\bibitem[Lindley(1956)]{lindley1956measure}
D.~V. Lindley.
\newblock On a measure of the information provided by an experiment.
\newblock \emph{The Annals of Mathematical Statistics}, 27\penalty0
  (4):\penalty0 986--1005, 1956.

\bibitem[Lou et~al.(2025)Lou, Bresler, and Pananjady]{lou2025computationally}
M.~Lou, G.~Bresler, and A.~Pananjady.
\newblock Computationally efficient reductions between some statistical models.
\newblock \emph{IEEE Transactions on Information Theory}, 71\penalty0
  (9):\penalty0 7097--7133, 2025.
\newblock \doi{10.1109/TIT.2025.3587354}.

\bibitem[Luo and Zhang(2022)]{luo2022tensor}
Y.~Luo and A.~R. Zhang.
\newblock Tensor clustering with planted structures: {S}tatistical optimality
  and computational limits.
\newblock \emph{The Annals of Statistics}, 50\penalty0 (1):\penalty0 584--613,
  2022.

\bibitem[Ma and Wigderson(2015)]{ma2015sum}
T.~Ma and A.~Wigderson.
\newblock Sum-of-squares lower bounds for sparse {PCA}.
\newblock \emph{Advances in Neural Information Processing Systems}, 28, 2015.

\bibitem[Ma and Wu(2015)]{ma2015computational}
Z.~Ma and Y.~Wu.
\newblock {Computational barriers in minimax submatrix detection}.
\newblock \emph{The Annals of Statistics}, 43\penalty0 (3):\penalty0 1089 --
  1116, 2015.
\newblock \doi{10.1214/14-AOS1300}.

\bibitem[Macris et~al.(2020)Macris, Rush, et~al.]{macris2020all}
N.~Macris, C.~Rush, et~al.
\newblock All-or-nothing statistical and computational phase transitions in
  sparse spiked matrix estimation.
\newblock \emph{Advances in Neural Information Processing Systems},
  33:\penalty0 14915--14926, 2020.

\bibitem[Montanari and Richard(2014)]{montanari2014statistical}
A.~Montanari and E.~Richard.
\newblock A statistical model for tensor {PCA}.
\newblock \emph{Advances in neural information processing systems}, 27, 2014.

\bibitem[Montanari et~al.(2015)Montanari, Reichman, and
  Zeitouni]{montanari2015limitation}
A.~Montanari, D.~Reichman, and O.~Zeitouni.
\newblock On the limitation of spectral methods: {F}rom the {G}aussian hidden
  clique problem to rank-one perturbations of {G}aussian tensors.
\newblock \emph{Advances in Neural Information Processing Systems}, 28, 2015.

\bibitem[Pananjady and Samworth(2022)]{pananjady2022isotonic}
A.~Pananjady and R.~J. Samworth.
\newblock Isotonic regression with unknown permutations: Statistics,
  computation and adaptation.
\newblock \emph{The Annals of Statistics}, 50\penalty0 (1):\penalty0 324--350,
  2022.

\bibitem[Perry et~al.(2020)Perry, Wein, and Bandeira]{perry2020statistical}
A.~Perry, A.~S. Wein, and A.~S. Bandeira.
\newblock Statistical limits of spiked tensor models.
\newblock 2020.

\bibitem[Polyanskiy and Wu(2025)]{polyanskiy2023book}
Y.~Polyanskiy and Y.~Wu.
\newblock \emph{Information theory: From coding to learning}.
\newblock Cambridge university press, 2025.

\bibitem[Quandt and Ramsey(1978)]{quandt1978estimating}
R.~E. Quandt and J.~B. Ramsey.
\newblock Estimating mixtures of normal distributions and switching
  regressions.
\newblock \emph{Journal of the American Statistical Association}, 73\penalty0
  (364):\penalty0 730--738, 1978.

\bibitem[Raghavendra et~al.(2018)Raghavendra, Schramm, and
  Steurer]{raghavendra2018high}
P.~Raghavendra, T.~Schramm, and D.~Steurer.
\newblock High dimensional estimation via sum-of-squares proofs.
\newblock In \emph{Proceedings of the International Congress of Mathematicians:
  Rio de Janeiro 2018}, pages 3389--3423. World Scientific, 2018.

\bibitem[Raginsky(2011)]{raginsky2011shannon}
M.~Raginsky.
\newblock Shannon meets {B}lackwell and {L}e {C}am: {C}hannels, codes, and
  statistical experiments.
\newblock In \emph{2011 IEEE International Symposium on Information Theory
  Proceedings}, pages 1220--1224. IEEE, 2011.

\bibitem[Ros et~al.(2019)Ros, Ben~Arous, Biroli, and Cammarota]{ros2019complex}
V.~Ros, G.~Ben~Arous, G.~Biroli, and C.~Cammarota.
\newblock Complex energy landscapes in spiked-tensor and simple glassy models:
  {R}uggedness, arrangements of local minima, and phase transitions.
\newblock \emph{Physical Review X}, 9\penalty0 (1):\penalty0 011003, 2019.

\bibitem[Schramm and Wein(2022)]{schramm2022computational}
T.~Schramm and A.~S. Wein.
\newblock Computational barriers to estimation from low-degree polynomials.
\newblock \emph{The Annals of Statistics}, 50\penalty0 (3):\penalty0
  1833--1858, 2022.

\bibitem[Shah et~al.(2019)Shah, Balakrishnan, and Wainwright]{shah2019feeling}
N.~B. Shah, S.~Balakrishnan, and M.~J. Wainwright.
\newblock Feeling the {B}ern: Adaptive estimators for {B}ernoulli probabilities
  of pairwise comparisons.
\newblock \emph{IEEE Transactions on Information Theory}, 65\penalty0
  (8):\penalty0 4854--4874, 2019.

\bibitem[Shiryaev and Spokoiny(2000)]{shiryaev2000statistical}
A.~N. Shiryaev and V.~G. Spokoiny.
\newblock \emph{Statistical Experiments And Decisions, Asymptotic Theory},
  volume~8.
\newblock World Scientific, 2000.

\bibitem[St{\"a}dler et~al.(2010)St{\"a}dler, B{\"u}hlmann, and van~de
  Geer]{stadler2010l1}
N.~St{\"a}dler, P.~B{\"u}hlmann, and S.~van~de Geer.
\newblock $\ell_1$-penalization for mixture regression models.
\newblock \emph{Test}, 19:\penalty0 209--256, 2010.

\bibitem[Stone(1961)]{stone1961non}
M.~Stone.
\newblock Non-equivalent comparisons of experiments and their use for
  experiments involving location parameters.
\newblock \emph{The Annals of Mathematical Statistics}, pages 326--332, 1961.

\bibitem[Szego(1939)]{szeg1939orthogonal}
G.~Szego.
\newblock \emph{Orthogonal polynomials}, volume~23.
\newblock American Mathematical Soc., 1939.

\bibitem[Torgersen(1991)]{torgersen1991comparison}
E.~Torgersen.
\newblock \emph{Comparison of statistical experiments}, volume~36.
\newblock Cambridge University Press, 1991.

\bibitem[Wainwright(2014)]{wainwright2014constrained}
M.~J. Wainwright.
\newblock Constrained forms of statistical minimax: Computation, communication
  and privacy.
\newblock In \emph{Proceedings of the International Congress of
  Mathematicians}, pages 13--21, 2014.

\bibitem[Wang et~al.(2017)Wang, Zhang, Giannakis, Ak{\c{c}}akaya, and
  Chen]{wang2017sparse}
G.~Wang, L.~Zhang, G.~B. Giannakis, M.~Ak{\c{c}}akaya, and J.~Chen.
\newblock Sparse phase retrieval via truncated amplitude flow.
\newblock \emph{IEEE Transactions on Signal Processing}, 66\penalty0
  (2):\penalty0 479--491, 2017.

\bibitem[Wang et~al.(2024)Wang, Lou, and Pananjady]{wang2023algorithms}
G.~Wang, M.~Lou, and A.~Pananjady.
\newblock Do algorithms and barriers for sparse principal component analysis
  extend to other structured settings?
\newblock \emph{IEEE Transactions on Signal Processing}, 72:\penalty0
  3187--3200, 2024.

\bibitem[Wang et~al.(2019)Wang, Yang, and Wang]{wang2019statistical}
L.~Wang, Z.~Yang, and Z.~Wang.
\newblock Statistical-computational tradeoff in single index models.
\newblock \emph{Advances in neural information processing systems}, 32, 2019.

\bibitem[Wang and Zhang(2022)]{wang2022simultaneous}
R.~Wang and Z.~Zhang.
\newblock Simultaneous optimal transport.
\newblock \emph{arXiv preprint arXiv:2201.03483}, 2022.

\bibitem[Wein(2025)]{wein2025computational}
A.~S. Wein.
\newblock Computational {C}omplexity of {S}tatistics: {N}ew {I}nsights from
  {L}ow-{D}egree {P}olynomials.
\newblock \emph{arXiv preprint arXiv:2506.10748}, 2025.

\bibitem[Wein et~al.(2019)Wein, El~Alaoui, and Moore]{wein2019kikuchi}
A.~S. Wein, A.~El~Alaoui, and C.~Moore.
\newblock The {K}ikuchi hierarchy and tensor {PCA}.
\newblock In \emph{2019 IEEE 60th Annual Symposium on Foundations of Computer
  Science (FOCS)}, pages 1446--1468. IEEE, 2019.

\bibitem[Zadik et~al.(2022)Zadik, Song, Wein, and Bruna]{zadik2022lattice}
I.~Zadik, M.~J. Song, A.~S. Wein, and J.~Bruna.
\newblock Lattice-based methods surpass sum-of-squares in clustering.
\newblock In \emph{Conference on Learning Theory}, pages 1247--1248. PMLR,
  2022.

\end{thebibliography}
\normalsize

\newpage
\appendix

\begin{center}
\LARGE{\bf Appendix}
\end{center}


\section{Supplementary results}
We state the computational hardness conjectures and prove some auxiliary lemmas and claims in this section.

\subsection{Explicit statements of hardness conjectures}\label{sec:hardness_conjectures}
In this section, we explicitly state the conjectures that are used in Section~\ref{sec:applications}.

Given two positive integers $m$ and $n$, let $\mathcal{G}(m,n,1/2)$ denote the distribution on bipartite graphs $G$ with parts of size $m$ and $n$ wherein each edge between the two parts is included independently with probability $1/2$. Let $k_{m} \in [m]$ and $k_{n} \in [n]$ such that $k_{n}$ divides $n$. Let $E = \{E_1,\dots,E_{k_n}\}$ be a partition of $[n]$ into $k_{n}$ parts, each of size $n/k_{n}$. We define the distribution
\[
    k\text{-}\mathrm{BPC}_{E}(m,n,k_{m},k_{n},1/2)
\]
on bipartite graphs with left side $[m]$ and right side $[n]$ as follows:
\begin{enumerate}
    \item Sample $G \sim \mathcal{G}(m,n,1/2)$.
    \item Choose a subset $S \subseteq [m]$ of size $k_m$ uniformly at random.
    \item Choose a subset $T \subseteq [n]$ of size $k_n$ uniformly at random subject to the constraint
    \[
        |T \cap E_j| = 1 \quad \text{for all } j \in [k_n].
    \]
    \item Plant the complete bipartite graph between $S$ and $T$ in the graph $G$, i.e. set all edges $\{(u,v): u \in S, v \in T\}$ to be present.
\end{enumerate}
The resulting distribution on $G$ is denoted $k\text{-}\mathrm{BPC}_{E}(m,n,k_{m},k_{n},1/2)$. We define the following two hypotheses:
\begin{align}\label{k-BPC}
  \mathcal{H}_{0}: G \sim \mathcal{G}(m,n,1/2) \qquad \text{and} \qquad \mathcal{H}_{1}: G \sim  k\text{-}\mathrm{BPC}_{E}(m,n,k_{m},k_{n},1/2).
\end{align}
We now state the $k$-$\mathrm{BPC}$ conjecture.
\begin{conjecture}[k-BPC]\label{k-BPC-conjecture}
Let $m$ and $n$ be two positive integers and be polynomial in one another. Let $k_{m} \in [m]$ and $k_{n} \in [n]$ such that $k_{n}$ divides $n$.
Let $E$ be any partition of $[n]$ into $k_{n}$ parts, each of size $n/k_{n}$, and $E$ is released. Let $\mathcal{H}_{0}$ and $\mathcal{H}_{1}$ be defined in Eq.~\eqref{k-BPC}. If
\[
  k_{m} = o(\sqrt{m}) \qquad \text{and} \qquad k_{n} = o(\sqrt{n}),
\]
then there cannot exist an algorithm $\mathcal{A}: \{0,1\}^{m \times n} \rightarrow \{0,1\}$ that runs in time $\mathsf{poly}(n)$ and satisfies
\[
    \mathbb{P}_{\mathcal{H}_{0}} \Big\{ \mathcal{A}(G) = 1 \Big\} + \mathbb{P}_{\mathcal{H}_{1}} \Big\{ \mathcal{A}(G) = 0 \Big\} = 1 - \Omega_{n}(1).
\]
\end{conjecture}

Let $n$ and $s$ be non-negative integers with $s \geq 3$. Let $\mathcal{G}^{s}(n,1/2)$ denote the distribution on $s$-uniform hypergraphs $G$ with $n$ vertices and each hyperedge is included independently with probability $1/2$. Let $k\in [n]$ and $k$ divides $n$. Let $E = \{E_1,\dots,E_{k}\}$ be a partition of $[n]$ into $k$ parts, each of size $n/k$. We define the distribution $k\text{-}\mathrm{HPC}_{E}^{s}(n,k,1/2)$ as follows:
\begin{enumerate}
    \item Sample $G \sim \mathcal{G}^{s}(n,1/2)$.
    \item Choose a subset $S \subseteq [n]$ of size $k$ uniformly at random subject to the constraint
    \[
        |S \cap E_j| = 1 \quad \text{for all } j \in [k].
    \]
    \item Plant the complete $s$-uniform hypergraph with vertex set $S$ in the graph $G$, i.e., add all $\binom{k}{s}$ hyperedges within $S$ to $G$.
\end{enumerate}
The resulting distribution on $G$ is denoted $k\text{-}\mathrm{HPC}_{E}^{s}(n,k,1/2)$.
We define two hypotheses:
\begin{align}\label{k-HPC}
  \mathcal{H}_{0}: G \sim \mathcal{G}^{s}(n,1/2) \qquad \text{and} \qquad \mathcal{H}_{1}: G \sim  k\text{-}\mathrm{HPC}_{E}^{s}(n,k,1/2).
\end{align}
We now state the $k$-$\mathrm{HPC}^{s}$ conjecture.
\begin{conjecture}[k-HPC]\label{k-HPC-conjecture}
Let $n$ be a positive integer and let $k \in [n]$ divides $n$.
Let $E$ be any partition of $[n]$ into $k$ parts, each of size $n/k$, and $E$ is released. Let $\mathcal{H}_{0}$ and $\mathcal{H}_{1}$ be defined in Eq.~\eqref{k-HPC}. If $k = o(\sqrt{n})$, then there cannot exist an algorithm $\mathcal{A}: \{0,1\}^{n^{s}} \rightarrow \{0,1\}$ that runs in time $\mathsf{poly}(n)$ and satisfies
\[
    \mathbb{P}_{\mathcal{H}_{0}} \Big\{ \mathcal{A}(G) = 1 \Big\} + \mathbb{P}_{\mathcal{H}_{1}} \Big\{ \mathcal{A}(G) = 0 \Big\} = 1 - \Omega_{n}(1).
\]
\end{conjecture}

\subsection{Inefficiency of the plug-in approach}\label{sec:plug_in_approach}

We next show that a contraction of the signal-to-noise ratio (equivalently, an inflation of the variance) by a factor of order polynomially in $1/\epsilon$ is necessary for the plug-in approach. To compare with the guarantee of Theorem~\ref{thm:gaussian-mean-nonlinear}, we consider the source model $\mathcal{U}$ and target model $\mathcal{V}$ defined as
\[
    \mathcal{U} = (\real, \{ \NORMAL(\theta,1)\}_{\theta \in \Theta} ) \quad \text{and} \quad \mathcal{V} = \big(\real, \big\{ \NORMAL\big(\theta^{2}/\varT^{2},1 \big) \big\}_{\theta \in \Theta} \big) \quad \text{with} \quad \Theta = [-1,1].
\]
Given $X_{\theta} \sim \NORMAL(\theta,1)$, the plug-in approach returns the random variable
\begin{align}\label{plug-in-reduction}
      \widetilde{\mathsf{K}}_{\mathsf{plugin}}(X_{\theta}): = \frac{X_{\theta}^{2}}{\varT^{2}} + G, \quad G \sim \NORMAL(0,1).
\end{align}
In the above, the random variable $G$ is independent of $X_{\theta}$. Operationally, the plug-in approach uses the single sample $X_{\theta}$ as an estimator of the true mean $\theta$ and attempts to approximate the random variable $Y_{\theta} \sim \NORMAL(\theta^{2}/\varT^{2},1)$. 
\begin{lemma} \label{lem:plugin-bad}
For all positive signal-to-noise ratio parameters $\tau^2$, we have
\begin{align*}
     \sup_{\theta \in \Theta} \; \mathsf{d}_{\mathsf{TV}} \big( \widetilde{\mathsf{K}}_{\mathsf{plugin}}(X_{\theta}), Y_{\theta} \big) \geq 
     \frac{ (\sqrt{2}-1) e^{-3/2} }{2\pi} \cdot \min\{1 , \varT^{-2} \}.
\end{align*}
Consequently, it is necessary to take $\varT^{2} \gtrsim \epsilon^{-1}$ to guarantee $\epsilon$-TV deficiency.
\end{lemma}
Note that Lemma~\ref{lem:plugin-bad} illustrates the inefficiency of the plug-in approach. In contrast, Theorem~\ref{thm:gaussian-mean-nonlinear} and Corollary~\ref{example:monomial} show that it suffices to take $\varT^{2} = \mathcal{O}(\mathsf{polylog}(1/\epsilon))$ for the approach we proposed.
\begin{proof}
By definition of TV, we obtain 
\begin{align*}
  \sup_{\theta \in \Theta} \; \mathsf{d}_{\mathsf{TV}} \big( \widetilde{\mathsf{K}}_{\mathsf{plugin}}(X_{\theta}), Y_{\theta} \big) &\geq \sup_{\theta \in \Theta} \; \sup_{t \in \real}\; \Big| \mathbb{P} \big\{ \widetilde{\mathsf{K}}_{\mathsf{plugin}}(X_{\theta}) \leq t \big\} - \mathbb{P} \big\{ Y_{\theta} \leq t \big\} \Big| \\
  &\overset{\1}{\geq} \sup_{\theta \in \Theta} \; \sup_{t \in \real}\; \bigg| \mathbb{P} \bigg\{ \frac{\theta^{2} + 2\theta Z + Z^{2}}{\varT^{2}}  + G \leq t \bigg\} - \mathbb{P} \bigg\{ \frac{\theta^{2}}{\varT^{2}} + G \leq t \bigg\}  \bigg| \\
  & \overset{\2}{\geq} 
  \bigg| \mathbb{P} \bigg\{ \frac{ Z^{2}}{\varT^{2}}  + G \leq 0 \bigg\} - \frac{1}{2}  \bigg|
\end{align*}
where in step $\1$ we use that for $Z,G \overset{\mathsf{i.i.d.}}{\sim} \NORMAL(0,1)$, we have
\[
X_{\theta} \overset{(d)}{=} \theta + Z, \quad Y_{\theta} \overset{(d)}{=} \frac{\theta^{2}}{\varT^{2}} + G, \quad \widetilde{\mathsf{K}}_{\mathsf{plugin}}(X_{\theta}) \overset{(d)}{=} \frac{ \theta^{2} + 2\theta Z +Z^{2}}{\varT^{2}} + G,
\]
and in step $\2$ we evaluate the supremum at $t = \theta^{2}/\varT^{2}$ and $\theta = 0$. Continuing, we have 
\begin{align*}
  \mathbb{P} \bigg\{ \frac{ Z^{2}}{\varT^{2}}  + G \leq 0 \bigg\} &= \EE_{Z \sim \NORMAL(0,1)}\bigg[ \int_{-\infty}^{-Z^{2}/\varT^{2}} \frac{e^{-\frac{x^{2}}{2}}}{\sqrt{2\pi}} \mathrm{d}x \bigg] = \frac{1}{2} - \EE_{Z \sim \NORMAL(0,1)}\bigg[ \int_{0}^{Z^{2}/\varT^{2}} \frac{e^{-\frac{x^{2}}{2}}}{\sqrt{2\pi}} \mathrm{d}x \bigg].
\end{align*}
Putting the two pieces together yields
\begin{align*}
  \sup_{\theta \in \Theta} \; \mathsf{d}_{\mathsf{TV}} \big( \mathsf{K}_{\mathsf{plugin}}(X_{\theta}), Y_{\theta} \big) &\geq \EE_{Z \sim \NORMAL(0,1)}\bigg[ \int_{0}^{Z^{2}/\varT^{2}} \frac{e^{-\frac{x^{2}}{2}}}{\sqrt{2\pi}} \mathrm{d}x \bigg] \\
  &= \int_{-\infty}^{\infty} \frac{e^{-\frac{t^{2}}{2}}}{\sqrt{2\pi}} \int_{0}^{t^{2}/\varT^{2}} \frac{e^{-\frac{x^{2}}{2}}}{\sqrt{2\pi}} \mathrm{d}x \; \mathrm{d}t \\
  &\geq \int_{1}^{\sqrt{2}} \frac{e^{-1}}{\sqrt{2\pi}} \mathrm{d}t \int_{0}^{1/\varT^{2}} \frac{e^{-\frac{x^{2}}{2}}}{\sqrt{2\pi}} \mathrm{d}x  = \frac{ (\sqrt{2}-1) e^{-1} }{\sqrt{2\pi}} \int_{0}^{1/\varT^{2}} \frac{e^{-\frac{x^{2}}{2}}}{\sqrt{2\pi}} \mathrm{d}x.
\end{align*}
Continuing, we have 
\[
    \int_{0}^{1/\varT^{2}} \frac{e^{-\frac{x^{2}}{2}}}{\sqrt{2\pi}} \mathrm{d}x \geq \begin{cases} \frac{e^{-\frac{1}{2}}}{\sqrt{2\pi}}, &\text{if } \varT^{2} \leq 1 \\  \frac{1}{\varT^2} \frac{e^{-\frac{1}{2}}}{\sqrt{2\pi}}, & \text{if } \varT^{2}>1.
     \end{cases}
\]
Combining the two inequalities in the display above yields
\begin{align*}
     \sup_{\theta \in \Theta} \; \mathsf{d}_{\mathsf{TV}} \big( \widetilde{\mathsf{K}}_{\mathsf{plugin}}(X_{\theta}), Y_{\theta} \big) \geq 
     \frac{ (\sqrt{2}-1) e^{-3/2} }{2\pi} \cdot \min\{1 , \varT^{-2} \} \quad \text{for all } \varT \in \real.
\end{align*}
Clearly, this implies that one must take $\varT^{2} \gtrsim \epsilon^{-1}$ to guarantee $\epsilon$-TV deficiency. 
\end{proof}

\section{Proofs of supporting lemmas}

In this section, we provide additional proofs deferred from the main text.

\subsection{Proof of Lemma~\ref{lemma:linear-combination}}\label{sec:pf_lemma_linear_combi}
By definition, we have 
\[
    f^{(n)}(x) = \sum_{k=1}^{K} c_{k} f_{k}^{(n)}(x) \quad \text{for all } n \in \mathbb{N} \text{ and } x\in \real.
\]
Consequently, applying the triangle inequality yields that for all $t \in \real$,
\begin{align*}
    \bigg(\sum_{j=1}^{n} \frac{|f^{(j)}(t)|}{\big(\ClinkT\big)^{j}j!}\bigg)^{n} &\leq  \bigg( \sum_{k=1}^{K} |c_{k}| \cdot \sum_{j=1}^{n} \frac{|f_{k}^{(j)}(t)|}{\big(\ClinkT\big)^{j}j!} \bigg)^{n} \\
    & \leq \bigg(\sum_{\ell=1}^{K}|c_{\ell}|\bigg)^{n} \sum_{k=1}^{K} \bigg(\frac{|c_{k}|}{\sum_{\ell=1}^{K}|c_{\ell}| }\bigg)^{n} \cdot \bigg( \sum_{j=1}^{n} \frac{|f_{k}^{(j)}(t)|}{\big(\ClinkT\big)^{j}j!} \bigg)^{n} \\
    & \leq \sum_{k=1}^{K} \bigg(\frac{|c_{k}|}{\sum_{\ell=1}^{K}|c_{\ell}| }\bigg)^{n} \cdot \bigg( \sum_{j=1}^{n} \frac{|f_{k}^{(j)}(t)|}{\big(\widetilde{C}_{f_{k}} \big)^{j}j!} \bigg)^{n},
\end{align*}
where in the last step we use $\ClinkT = \sum_{\ell=1}^{K} |c_{\ell}| \cdot \max_{k \in [K]} \{ \widetilde{C}_{f_{k}} \}$. By letting $t(x) = (\theta + \sqrt{2}\sigma x)/\varT$ and integrating over $x\in \real$ yields
\begin{align*}
  \int_{\real} \bigg(\sum_{j=1}^{n} \frac{|f^{(j)}(t(x))|}{\big(\ClinkT\big)^{j}j!}\bigg)^{n} \frac{e^{-x^{2}}}{\sqrt{2\pi}} \mathrm{d}x &\leq \sum_{k=1}^{K}\bigg(\frac{|c_{k}|}{\sum_{\ell=1}^{K}|c_{\ell}| }\bigg)^{n} \int_{\real} \bigg( \sum_{j=1}^{n} \frac{|f_{k}^{(j)}(t)|}{\big(\widetilde{C}_{f_{k}} \big)^{j}j!} \bigg)^{n} \frac{e^{-x^{2}}}{\sqrt{2\pi}} \mathrm{d}x \\
  & \overset{\1}{\leq} \sum_{k=1}^{K}\bigg(\frac{|c_{k}|}{\sum_{\ell=1}^{K}|c_{\ell}| }\bigg)^{n} \cdot (C_{f_{k}}n)! \\
  & \leq \Big( \max_{k \in [K]} \{ C_{f_{k}} \} \cdot n \Big)!,
\end{align*}
where in step $\1$ we use the fact that $f_{k}$ satisfies Assumption~\ref{assump:link} with parameters $\widetilde{C}_{f_{k}}$ and $C_{f_{k}}$. This verifies that $f$ satisfies Assumption~\ref{assump:link}. We next bound $\Llink(n)$ and $\LlinkT(n,\varT)$. By definition~\eqref{def-quantity-link}, we obtain
\begin{align*}
   \Llink(n) = \sup_{|x| \leq 1} \sum_{j=1}^{n} \frac{|f^{(j)}(x)|}{\big(\ClinkT \big)^{j} j!} &\leq \sum_{k=1}^{K} |c_{k}| \sup_{|x| \leq 1} \sum_{j=1}^{n} \frac{|f_{k}^{(j)}(x)|}{\big(\ClinkT \big)^{j} j!} \\
    & \overset{\1}{\leq} \sum_{k=1}^{K} \frac{|c_{k}|}{\sum_{\ell = 1}^{K} |c_{\ell}| }  \sup_{|x| \leq 1} \sum_{j=1}^{n} \frac{|f_{k}^{(j)}(x)|}{\big(\widetilde{C}_{f_{k}} \big)^{j} j!} \\
    & = \sum_{k=1}^{K} \frac{|c_{k}|}{\sum_{\ell = 1}^{K} |c_{\ell}| } L_{f_{k}}(n),
\end{align*}
where in step $\1$ we use $\ClinkT = \sum_{\ell=1}^{K} |c_{\ell}| \cdot \max_{k \in [K]} \{ \widetilde{C}_{f_{k}} \}$. Following identical steps yields the desired bound for $\LlinkT(n,\varT)$. 
\qed

\subsection{Proof of Lemma~\ref{lemma:aux-px-qx-nongaussian}}\label{sec:pf_lemma_aux_px_qx_nongaussian}
We prove each part separately.\\
\noindent \underline{Proof of part (i).} Applying Faà di Bruno's formula as in Eq.~\eqref{eq-target-derivative-nongaussian} yields that for all $n \in \mathbb{N}$ and all $x,y\in \real$, we have
\begin{align}\label{eq-target-derivative-nongaussian-1}
  \nabla_{x}^{(n)} v(y;x) =  \frac{1}{\varT} e^{-\psi(\frac{y-x}{\varT})} \cdot n! \cdot \sum_{m_{1},\dots,m_{n}} \frac{(-1)^{m_{1}+\cdots+m_{n}}}{m_{1}!\, m_{2}! \cdots m_{n}!} \prod_{j=1}^{n} \bigg( \frac{ \psi^{(j)}(\frac{y-x}{\varT})}{(-\varT)^{j}j!} \bigg)^{m_{j}},
\end{align}
where the sum is over all $n$-tuples of nonnegative integers $(m_{1},\dots,m_{n})$ satisfying $\sum_{j=1}^{n} jm_{j} = n$. Consequently, using definition~\eqref{truncated-signed-kernel} yields
\begin{align*}
  \mathcal{S}_{N}^{\star}(y\mid x) &= \frac{1}{\varT} e^{-\psi(\frac{y-x}{\varT})} \Bigg( 1 + \sum_{k=1}^{N} \frac{(-1)^{k} \sigma^{2k} (2k)!}{(2k)!!} \sum_{m_{1},\dots,m_{2k}} \frac{(-1)^{m_{1}+\cdots+m_{2k}}}{m_{1}!\,m_{2}! \cdots m_{2k}!}  \prod_{j=1}^{2k} \bigg( \frac{ \psi^{(j)}(\frac{y-x}{\varT})}{(-\varT)^{j}j!} \bigg)^{m_{j}} \Bigg) \\
  & \overset{\1}{\geq} \frac{1}{\varT} e^{-\psi(\frac{y-x}{\varT})} \Bigg( 1 - \sum_{k=1}^{N} \frac{ \sigma^{2k} (2k)!}{(2k)!!} \frac{\lambda^{2k}}{\varT^{2k}} \sum_{m_{1},\dots,m_{2k}} \frac{1}{m_{1}!\,m_{2}! \cdots m_{2k}!}  \prod_{j=1}^{2k} \bigg( \frac{\big| \psi^{(j)}(\frac{y-x}{\varT}) \big|}{\lambda^{j}j!} \bigg)^{m_{j}} \Bigg) \\
  & \overset{\2}{\geq} 
  \frac{1}{\varT} e^{-\psi(\frac{y-x}{\varT})} \Bigg( 1 - \sum_{k=1}^{N} \frac{ \sigma^{2k} (2k)!}{(2k)!!} \frac{\lambda^{2k}}{\varT^{2k}} \sum_{m_{1},\dots,m_{2k}} \bigg(\sum_{j=1}^{2k} \frac{\big| \psi^{(j)}(\frac{y-x}{\varT}) \big|}{\lambda^{j}j!} \bigg)^{m_{1}+\cdots+m_{2k}} \Bigg),
\end{align*}
where in step $\1$ we apply the triangle inequality and use $\sum_{j=1}^{2k} jm_{j} = 2k$ so that
\[
    \frac{\lambda^{2k}}{\varT^{2k}} \prod_{j=1}^{2k} \bigg(\frac{1}{\lambda^{j}}\bigg)^{m_{j}} =   \prod_{j=1}^{2k} \bigg(\frac{1}{\varT^{j}}\bigg)^{m_{j}},
\]
and in step $\2$ we use the multinomial theorem. Continuing, note that we have $(y-x)/\varT \in \Set(\psi,2N,\lambda)$~\eqref{def-set-complement} by assumption, and so we obtain
\[
    \sum_{j=1}^{2k} \frac{\big| \psi^{(j)}(\frac{y-x}{\varT}) \big|}{\lambda^{j}j!} \leq 1 \quad \text{for all } k \in [N].
\]
Consequently,
\begin{align*}
   \mathcal{S}_{N}^{\star}(y\mid x) &\geq \frac{1}{\varT} e^{-\psi(\frac{y-x}{\varT})} \bigg( 1 - \sum_{k=1}^{N} \frac{ \sigma^{2k} (2k)!}{(2k)!!} \frac{\lambda^{2k}}{\varT^{2k}} \sum_{m_{1},\dots,m_{2k}} 1 \bigg) \\
   &\overset{\1}{\geq}  \frac{1}{\varT} e^{-\psi(\frac{y-x}{\varT})} \bigg( 1 - \sum_{k=1}^{N} \frac{ \sigma^{2k} (2k)^{k} \lambda^{2k} e^{2k} }{\varT^{2k}} \bigg) \\
   &\overset{\2}{\geq} \frac{1}{\varT} e^{-\psi(\frac{y-x}{\varT})} \bigg( 1 - \sum_{k=1}^{N} 4^{-k} \bigg) \geq 0.
\end{align*}
Here in step $\1$, we use the fact there are at most $e^{2k}$ terms in the summation and that $(2k)!/(2k)!! \leq (2k)^{k}$. On the other hand, step $\2$ holds under the assumption that $\varT \geq 8e\lambda \sigma N$, which is the condition in the lemma statement. This proves part (i) of Lemma~\ref{lemma:aux-px-qx-nongaussian}. \\

\noindent \underline{Proof of part (ii).} We first verify the first part of Eq.~\eqref{ineq1:integral-derivative-density}. Note that $\rho(t) = g(\psi(t))$, where $g(t) = e^{-t}$. Applying Faà di Bruno's formula yields
\begin{align*}
  \rho^{(n)}(t) = \sum_{m_{1},\dots,m_{n}} \frac{n!}{m_{1}!\,m_{2}! \cdots m_{n}!} (-1)^{m_{1}+\cdots+m_{n}} e^{-\psi(t)} \prod_{j=1}^{n} \bigg(\frac{ \psi^{(j)}(t)}{j!} \bigg)^{m_{j}}.
\end{align*}
Applying the triangle inequality yields
\begin{align}\label{ineq-derivative-rho-up}
  \big| \rho^{(n)}(t) \big| &\leq \sum_{m_{1},\dots,m_{n}} \frac{n!}{m_{1}!\,m_{2}! \cdots m_{n}!}  e^{-\psi(t)} \prod_{j=1}^{n} \bigg(\frac{ | \psi^{(j)}(t) | }{j!} \bigg)^{m_{j}} \nonumber \\
  & \overset{\1}{\leq} n! \cdot \lambda^{n} \sum_{m_{1},\dots,m_{n}} \frac{e^{-\psi(t)}}{m_{1}!\,m_{2}! \cdots m_{n}!}   \prod_{j=1}^{n} \bigg(\frac{ | \psi^{(j)}(t) | }{\lambda^{j}j!} \bigg)^{m_{j}} \nonumber \\
  & \leq n! \cdot \lambda^{n} \sum_{m_{1},\dots,m_{n}} e^{-\psi(t)} \bigg(  \sum_{j=1}^{n}\frac{ | \psi^{(j)}(t) | }{\lambda^{j}j!}  \bigg)^{m_{1}+\cdots+m_{n}},
\end{align}
where in step $\1$ we use $\sum_{j=1}^{n}jm_{j} = n$ so that $\prod_{j=1}^{n} (\frac{1}{\lambda^{j}})^{m_{j}}  = \lambda^{-n}$ and in step $\2$ we use the multinomial theorem. By setting $\lambda = \CdensityT$, we obtain
\begin{align*}
  \int_{\real} \big| \rho^{(n)}(t) \big| \mathrm{d}t &\leq n! \cdot \big(\CdensityT \big)^{n} \sum_{m_{1},\dots,m_{n}} \int_{\real} e^{-\psi(t)} \bigg(  \sum_{j=1}^{n}\frac{ | \psi^{(j)}(t) | }{\big(\CdensityT \big)^{j}j!}  \bigg)^{m_{1}+\cdots+m_{n}} \mathrm{d}t \\
  & \overset{\1}{\leq} n! \cdot \big(\CdensityT \big)^{n} \sum_{m_{1},\dots,m_{n}} (\Cdensity n)! \overset{\2}{\leq} n! \cdot \big(\CdensityT \big)^{n} e^{n} \cdot (\Cdensity n)! <+\infty,
\end{align*}
where in step $\1$ we use Assumption~\ref{assump:log-density} and in step $\2$ we use the fact that there are at most $e^{n}$ terms in the summation. Since $|\rho^{(n)}(t)|$ is continuous and nonnegative, we must have $\lim_{t \uparrow +\infty} \rho^{(n)}(t) = \lim_{t \downarrow -\infty} \rho^{(n)}(t) = 0$. Consequently, we obtain for $n\geq 1$ that
\[
    \int_{\real} \rho^{(n)}(t) \mathrm{d}t = \int_{\real} \mathrm{d} \rho^{(n-1)}(t) = \Big[ \rho^{(n-1)}(t) \Big] \Big|_{-\infty}^{+\infty} = 0.
\]
This proves the first part of Eq.~\eqref{ineq1:integral-derivative-density}.

We now turn to prove the second part of Eq.~\eqref{ineq1:integral-derivative-density}. Using Eq.~\eqref{eq-target-derivative-nongaussian-1} and applying triangle inequality, we obtain for all $n\geq 1$ that
\begin{align}\label{ineq-up-derivative-vabs}
    \big| \nabla_{x}^{(n)} v(y;x) \big| &\leq \frac{1}{\varT} e^{-\psi(\frac{y-x}{\varT})} \cdot n! \cdot \sum_{m_{1},\dots,m_{n}} \frac{1}{m_{1}!\, m_{2}! \cdots m_{n}!} \prod_{j=1}^{n} \bigg( \frac{ |\psi^{(j)}(\frac{y-x}{\varT})|}{\varT^{j}j!} \bigg)^{m_{j}} \nonumber \\
    &\leq \frac{1}{\varT} e^{-\psi(\frac{y-x}{\varT})} \cdot n! \cdot \frac{\lambda^{n}}{\varT^{n}} \sum_{m_{1},\dots,m_{n}}   \bigg( \sum_{j=1}^{n} \frac{ |\psi^{(j)}(\frac{y-x}{\varT})|}{\lambda^{j}j!} \bigg)^{m_{1}+\cdots+m_{n}},
\end{align}
where in the last step we use $(\lambda/\varT)^{n} \prod_{j=1}^{n}(\frac{1}{\lambda^{j}})^{m_{j}} = \prod_{j=1}^{n}(\frac{1}{\varT^{j}})^{m_{j}}$ and then apply the multinomial theorem. Integrating over $y$ yields
\begin{align*}
  &\int_{\frac{y-x}{\varT} \in \Set^{\complement}(\psi,2N,\lambda)} \big| \nabla_{x}^{(n)} v(y;x) \big| \mathrm{d}y \\
  & \leq n! \cdot \frac{\lambda^{n}}{\varT^{n}} \sum_{m_{1},\dots,m_{n}} \int_{\frac{y-x}{\varT} \in \Set^{\complement}(\psi,2N,\lambda)} \frac{1}{\varT} e^{-\psi(\frac{y-x}{\varT})}   \bigg( \sum_{j=1}^{n} \frac{ |\psi^{(j)}(\frac{y-x}{\varT})|}{\lambda^{j}j!} \bigg)^{m_{1}+\cdots+m_{n}} \mathrm{d}y \\
  & \overset{\1}{\leq} n! \cdot \frac{\lambda^{n}}{\varT^{n}} \sum_{m_{1},\dots,m_{n}} \int_{\Set^{\complement}(\psi,2N,\lambda)} e^{-\psi(t)} \bigg(  \sum_{j=1}^{n}\frac{ | \psi^{(j)}(t) | }{\lambda^{j}j!}  \bigg)^{m_{1}+\cdots+m_{n}} \mathrm{d}t \\
  & \leq  n! \cdot \frac{\lambda^{n}}{\varT^{n}} \cdot e^{n} \cdot \err(\psi,2N,\lambda),
\end{align*}
where in step $\1$ we use the change of variables $t = (y-x)/\varT$. This proves Eq.~\eqref{ineq1:integral-derivative-density}.
\qed

\subsection{Proof of Lemma~\ref{lemma:aux-derivative-bound-nonlinear}}\label{sec:pf_lemma_aux_derivative_bound_nonlinear}
Note that by definition~\eqref{target:density} and Eq.~\eqref{shorthand-h-x}, we have $v(y;x) = \phi(g(x))$. Applying Faà di Bruno's formula yields
\begin{align*}
  \nabla_{x}^{(n)} v(y;x) &= \sum_{m_{1}, \dots,m_{n}} \frac{n!}{m_{1}!\,m_{2}! \cdots m_{n}!} \cdot \phi^{(m_{1} + \cdots +m_{n})}(g(x)) \cdot \prod_{j=1}^{n} \bigg( \frac{g^{(j)}(x) }{j!} \bigg)^{m_{j}} \\
  & \overset{\1}{=} n! \cdot \sum_{m_{1}, \dots,m_{n}} \frac{(-1)^{m_{1} + \cdots+m_{n}}}{m_{1}!\,m_{2}! \cdots m_{n}!} \cdot \phi(g(x)) \cdot H_{m_{1}+\cdots+m_{n}}(g(x))  \cdot \prod_{j=1}^{n} \bigg( \frac{f^{(j)}(x/\varT) }{2\varT^{j}j!} \bigg)^{m_{j}} \\
  & \overset{\2}{=} \frac{(\ClinkT)^{n}}{\varT^{n}} \cdot n! \cdot \sum_{m_{1}, \dots,m_{n}} \frac{(-1)^{m_{1} + \cdots+m_{n}}}{m_{1}!\,m_{2}! \cdots m_{n}!} \cdot \phi(g(x)) \cdot H_{m_{1}+\cdots+m_{n}}(g(x))  \cdot \prod_{j=1}^{n} \bigg( \frac{f^{(j)}(x/\varT) }{2(\ClinkT)^{j}j!} \bigg)^{m_{j}},
\end{align*}
where in all steps the summation is over all nonnegative integers $(m_{1},\cdots,m_{n})$ satisfying $\sum_{j=1}^{n} jm_{j} = n$. In step $\1$, we use the definition of Hermite polynomials~\eqref{def:hermite} and in step $\2$, we use $\prod_{j=1}^{n} (\frac{1}{\varT^{j}})^{m_{j}} =  \frac{(\ClinkT)^{n}}{\varT^{n}} \prod_{j=1}^{n} (\frac{1}{(\ClinkT)^{j}})^{m_{j}}$. 
Using the definition of $h_{j}(x)$ in Eq.~\eqref{shorthand-h-x} and combining the pieces then yields Eq.~\eqref{eq:derivative-gaussian-nonlinear}. 

To prove the second claim of the lemma, note that from Ineq.~\eqref{Hermite-boundness}, we have
\begin{align*}
  \big| \phi(g(x)) \cdot H_{m_{1}+\cdots+m_{n}}(g(x)) \big| &\leq \frac{e^{-g(x)^{2}/2}}{\sqrt{2\pi}} \pi^{1/4} \cdot (m_{1} + \cdots + m_{n})! \cdot 2^{m_{1} + \cdots+m_{n}} \\
  & \leq  e^{-g(x)^{2}/2} \cdot 2^{n} \cdot (m_{1} + \cdots + m_{n})!.
\end{align*}
Now applying the triangle inequality to the display before yields
\begin{align*}
   \big| \nabla_{x}^{(n)} v(y;x) \big| &\leq e^{-\frac{g(x)^{2}}{2}} \cdot \frac{(2\ClinkT)^{n}}{\varT^{n}} \cdot n! \sum_{m_{1}, \dots,m_{n}} \frac{(m_{1} + \cdots + m_{n})!}{m_{1}!\,m_{2}! \cdots m_{n}!} \prod_{j=1}^{n} | h_{j}(x) |^{m_{j}} \\
   & \overset{\1}{\leq} e^{-\frac{g(x)^{2}}{2}} \cdot \frac{(2\ClinkT)^{n}}{\varT^{n}} \cdot n! \sum_{m_{1}, \dots,m_{n}} \bigg( \sum_{j=1}^{n} |h_{j}(x)| \bigg)^{m_{1} + \cdots + m_{n}} \\
   & \leq e^{-\frac{g(x)^{2}}{2}} \cdot \frac{(2\ClinkT)^{n}}{\varT^{n}} \cdot n! \cdot e^{n} \cdot \bigg[ \bigg( \sum_{j=1}^{n} |h_{j}(x)| \bigg)^{n} \vee 1 \bigg]
\end{align*}
where in step $\1$ we use the multinomial theorem so that
\[
    \frac{(m_{1} + \cdots + m_{n})!}{m_{1}!\,m_{2}! \cdots m_{n}!} \prod_{j=1}^{n} | h_{j}(x) |^{m_{j}} \leq \bigg( \sum_{j=1}^{n} |h_{j}(x)| \bigg)^{m_{1} + \cdots + m_{n}},
\]
and in the last step we use there are at most $e^{n}$ in the summation and $m_{1}+\cdots+m_{n} \leq n$.
This concludes the proof. \qed

\subsection{Proof of Lemma~\ref{lemma:positivity-nonlinear}}\label{sec:pf_lemma_positivity_nonlinear}
We prove each part separately. Recall the shorthand $h_{j}(x)$, $g(x)$, and $\phi(x)$ from Eq.~\eqref{shorthand-h-x}.
\paragraph{Proof of Lemma~\ref{lemma:positivity-nonlinear}(a).}

From Eq.~\eqref{Hermite-closed-form}, we obtain for all $n \in \PosInt$, $n\geq 1$, and $x\in \real$ that
\[
   |H_{n}(x)| \leq n! \cdot \frac{n}{2} \cdot (1 \vee 2|x|)^{n}.
\]
Consequently, if $|g(x)| = |f(x/\varT) - y|/2 \leq \frac{\sqrt{\varT}}{2}$, then for $\varT \geq 4$ and $1\leq m_{1}+\cdots+m_{n}\leq n$, we have
\begin{align*}
   \big|H_{m_{1}+\cdots+m_{n}}(g(x))\big| &\leq \frac{m_{1}+\cdots+m_{n}}{2} \cdot (m_{1}+\cdots+m_{n})! \cdot \varT^{\frac{m_{1}+\cdots+m_{n}} {2}} \\
   &\leq \frac{n}{2} \cdot (m_{1}+\cdots+m_{n})! \cdot \varT^{\frac{n}{2}}.
\end{align*}
Using the bound in the above display along with Eq.~\eqref{eq:derivative-gaussian-nonlinear}, we obtain for $n\geq 1$ that
\begin{align}\label{ineq:upbound-derivative-abs-nonlinear}
  \big| \nabla_{x}^{(n)}v(y;x) \big| &\leq \phi(g(x)) \cdot \frac{n}{2} \cdot \varT^{\frac{n}{2}} \cdot n! \cdot \frac{(\ClinkT)^{n}}{\varT^{n}} \sum_{m_{1},\dots,m_{n}} \frac{(m_{1}+\cdots+m_{n})!}{m_{1}! \, m_{2}! \cdots m_{n}!} \cdot \prod_{j=1}^{n} \big| h_{j}(x) \big|^{m_{j}} \nonumber \\
  &\leq  \phi(g(x)) \cdot \frac{n}{2} \cdot \varT^{\frac{n}{2}} \cdot n! \cdot \frac{(\ClinkT)^{n}}{\varT^{n}} \cdot e^{n} \cdot \bigg[ \bigg( \sum_{j=1}^{n} |h_{j}(x)| \bigg)^{n} \vee 1 \bigg],
\end{align}
where in the last step we use the multinomial theorem so that
\[
    \frac{(m_{1}+\cdots+m_{n})!}{m_{1}! \, m_{2}! \cdots m_{n}!} \cdot \prod_{j=1}^{n} \big| h_{j}(x) \big|^{m_{j}} \leq \bigg( \sum_{j=1}^{n} |h_{j}(x)| \bigg)^{m_{1} + \cdots +m_{n}} \leq \bigg( \sum_{j=1}^{n} |h_{j}(x)| \bigg)^{n} \vee 1,
\]
and we also use there are at most $e^{n}$ terms in the summation (see Ineq.~\eqref{ineq:size-n-tuples}). Continuing, by using the signed kernel definition~\eqref{truncated-signed-kernel} for $\sigma = 1$, we obtain for $|x| \leq \varT$ that
\begin{align*}
  \mathcal{S}_{N}^{\star}(y \mid x) &= \sum_{k=0}^{N} \frac{(-1)^{k} }{(2k)!!} \nabla_{x}^{(2k)} v(y;x) \\
  & \geq v(y;x) - \sum_{k=1}^{N}  \frac{1}{(2k)!!} \big| \nabla_{x}^{(2k)} v(y;x) \big| \\
  & \overset{\1}{\geq} \phi(g(x)) - \sum_{k=1}^{N}  \frac{1}{(2k)!!} \cdot  \phi(g(x)) \cdot \frac{2k}{2} \cdot (2k)! \cdot \varT^{k} \cdot \bigg( \frac{\ClinkT}{\varT} \bigg)^{2k} e^{2k}  \bigg[ \bigg( \sum_{j=1}^{2k} |h_{j}(x)| \bigg) \vee 1 \bigg]^{2k} \\
  & \overset{\2}{\geq} \phi(g(x)) \cdot \Bigg( 1 - \sum_{k=1}^{N} \frac{ \Big( e^{3}  \big( \ClinkT \big)^{2} \big( \Llink(2k) \vee 1 \big) 2k  \Big)^{k} }{\varT^{k}} \Bigg) \\
  & \overset{\3}{\geq} \phi(g(x)) \cdot \Bigg( 1 - \sum_{k=1}^{N} \frac{1}{2^{k}} \Bigg) \geq 0.
\end{align*}
In step $\1$ above, we use Ineq.~\eqref{ineq:upbound-derivative-abs-nonlinear} for $n = 2k$. In step $\2$, we use the inequalities $k \leq e^{k}$ and $(2k)!/(2k)!! = (2k-1)!! \leq (2k)^{k}$, and also that
\[
   \sum_{j=1}^{2k} |h_{j}(x)| \leq \Llink(2k)  \quad \text{for all } |x| \leq \varT,
\]
which in turn follows from combining Eq.~\eqref{def-quantity-link} and Eq.~\eqref{shorthand-h-x}.
Step $\3$ holds if 
\[
   \varT \geq 2 e^{3} \big( \ClinkT \big)^{2} \big( \Llink(2N) \vee 1 \big) 2N \geq 2e^{3} \big( \ClinkT \big)^{2} \big( \Llink(2k) \vee 1 \big) 2k  \quad \text{for all } k\in[N],  
\]
as assumed in the lemma.
This concludes the proof of part (a).

\paragraph{Proof of Lemma~\ref{lemma:positivity-nonlinear}(b)}
Using Ineq.~\eqref{ineq:derivative-gaussian-nonlinear} and noting that $g(x) = (f(x/\varT)-y)/2$ (Eq.~\eqref{shorthand-h-x}), we obtain
\begin{align*}
    &\int_{f(\frac{x}{\varT}) + \sqrt{\varT}}^{+\infty} \big| \nabla_{x}^{(2k)}v(y;x) \big| \mathrm{d}y + \int_{-\infty}^{f(\frac{x}{\varT}) - \sqrt{\varT}} \big| \nabla_{x}^{(2k)}v(y;x) \big| \mathrm{d}y  \\
    &\leq  (2k)! \cdot e^{2k} \cdot \bigg( \frac{2\ClinkT}{\varT}\bigg)^{2k} \cdot \bigg[ \bigg( \sum_{j=1}^{2k} \big| h_{j}(x) \big| \bigg)^{2k} \vee 1 \bigg] \cdot 2\int_{\sqrt{\varT}}^{+\infty} e^{-t^{2}/8} \mathrm{d}t  \\
    &\overset{\1}{\leq} 2\sqrt{2\pi} \exp(-\varT/8) \cdot (2k)! \cdot \bigg( \frac{2e\ClinkT}{\varT}\bigg)^{2k} \cdot \big( \Llink(2k) \vee 1 \big)^{2k},
\end{align*}
where in step $\1$ we use the fact that $|x| \leq \varT$ as well as the definitions of $h_{j}(x)$~\eqref{shorthand-h-x} and $\Llink(n)$~\eqref{def-quantity-link}. In particular, these yield
\[
    \sum_{j=1}^{2k} \big| h_{j}(x) \big| \leq \Llink(2k) \quad \text{and}\quad  2\int_{\sqrt{\varT}}^{+\infty} e^{-t^{2}/8} \mathrm{d}t \leq 2\sqrt{2\pi}\exp(-\varT/8). 
\]

\paragraph{Proof of Lemma~\ref{lemma:positivity-nonlinear}(c)}
Using the closed-form expression for $\nabla_{x}^{(n)}v(y;x)$ (Eq.~\eqref{eq:derivative-gaussian-nonlinear}) and noting that $g(x) = (f(x/\varT)-y)/2$ (Eq.~\eqref{shorthand-h-x}), we obtain
\begin{align*}
    \int_{\real} \nabla_{x}^{(n)}v(y;x) \mathrm{d}y &= \sum_{m_{1},\dots,m_{n}} \frac{n! \cdot (-1)^{m_{1}+\cdots+m_{n}}}{m_{1}!\,m_{2}! \cdots m_{n}!}   \bigg( \frac{\ClinkT}{\varT} \bigg)^{n}  \prod_{j=1}^{n} \Big(h_{j}(x)/2\Big)^{m_{j}} 2\int_{\real} H_{m_{1}+\cdots+m_{n}}( t  ) \cdot \phi(t) \mathrm{d}t = 0.
\end{align*}
In the last step, we use the facts that $m_{1}+\cdots+m_{n} \geq 1$ and $H_{0}(x) = 1$ along with the orthogonal property of Hermite polynomials. In particular, these imply that
\begin{align*}
    \int_{\real} H_{m_{1}+\cdots+m_{n}}( t  ) \cdot \phi(t) \mathrm{d}t = \frac{1}{\sqrt{2\pi}} \big \langle H_{m_{1}+\cdots+m_{n}}(\cdot), H_{0}(\cdot)  \big \rangle = 0.
\end{align*}
This concludes the proof of the final part of the lemma.
\qed

\subsection{Proof of Lemma~\ref{lemma:reduction-tensor-BB20}}\label{sec:pf_claims_TPCA}
Note that there are two differences between the tensor $\widetilde{T}$ in Eq.~\eqref{tensor-BB20} and the tensor $T$ in Eq.~\eqref{tensor-continuous-spike}: (i) the additive Gaussian noise $\widetilde{\zeta}$ is asymmetric, and (ii) the spike $\widetilde{v}$ is discrete. We now construct two procedures $\mathcal{R}_{1}, \mathcal{R}_{2}: \mathbb{R}^{n^{\otimes s}} \to \mathbb{R}^{n^{\otimes s}}$, each designed to resolve one of these two differences. Our overall reduction is the composition of the two, i.e.,
$
    \mathcal{R}\big( \widetilde{T} \big): = \mathcal{R}_{2} \circ \mathcal{R}_{1}\big( \widetilde{T} \big).
$

\paragraph{Construction of reduction $\mathcal{R}_{1}$ to make instance symmetric.}
We let
\[
    \mathcal{R}_{1}\big( \widetilde{T} \big) = \frac{1}{s!}\sum_{\pi \in \in \mathfrak{S}_{s}} \big( \widetilde{T} \big)^{\pi},
\]
where $\big(\widetilde{T}\big)^{\pi} \in \real^{n^{\otimes s}}$ denotes the permutation of the tensor $\widetilde{T}$ so that
\[
  \big(\widetilde{T}\big)^{\pi}_{i_{1}, \dots,i_{s}} =  \big(\widetilde{T}\big)^{\pi}_{i_{\pi(1)}, \dots,i_{\pi(s)}} 
  \quad \text{for all indices} \quad i_{1},\cdots, i_{s} \in [n].
\] 
Consequently, by the definition of $\widetilde{T}$ in Eq.~\eqref{tensor-BB20}, we obtain
\[
    \mathcal{R}_{1}\big( \widetilde{T} \big) = \beta(\sqrt{n} \widetilde{v} )^{\otimes s} + \zeta, \quad where \quad \widetilde{v} \sim \mathsf{Unif}\big( \big\{\tfrac{-1}{\sqrt{n}}, \tfrac{1}{\sqrt{n}} \big\}^{n}  \big)
\]
and the noise tensor $\zeta$ is now symmetric and satisfies
\[
    \zeta \overset{(d)}{=} \frac{1}{s!} \sum_{\pi \in \in \mathfrak{S}_{s}} G^{\pi},
\]
where $G \in \mathbb{R}^{n^{\otimes s}} $ is an asymmetric tensor with entries that are $\mathsf{i.i.d.}$ $\mathcal{N}(0,1)$, and $G^{\pi} \in \real^{n^{\otimes s}}$ satisfies $G^{\pi}_{i_{1}, \dots,i_{s}} = G_{i_{\pi(1)}, \dots,i_{\pi(s)}}$ for all indices $i_{1},\cdots, i_{s} \in [n]$. 

\paragraph{Construction of reduction $\mathcal{R}_{2}$ to spread the signal spike.}
We next construct the second procedure $\mathcal{R}_{2}$, which transforms the tensor $\mathcal{R}_{1}\big(\widetilde{T}\big)$ into one whose spike $v$ is uniformly distributed on the unit sphere and whose additive Gaussian noise remains independent of $v$ and symmetric Gaussian. For notational convenience, we let 
\[
  \overline{T} := \mathcal{R}_{1}\big(\widetilde{T}\big) = \beta(\sqrt{n} \widetilde{v} )^{\otimes s} + \zeta.
\]
We first review some definitions. For a tensor $T \in \real^{n \otimes s}$, define its vectorization as
\[
    \mathsf{vec}(T) = \sum_{i_{1}=1}^{n} \cdots \sum_{i_{s}=1}^{n} T_{i_{1},\dots, i_{s}} \times e_{i_{1}} \otimes e_{i_{2}} \otimes \cdots \otimes e_{i_{s}},
\]
where $e_{i}$ is $i$-th standard basis vector in $\real^{n}$ and $\otimes$ denotes the Kronecker product. For any permutation $\pi : [s] \to [s]$, define $P_{\pi} \in \mathbb{R}^{n^{s} \times n^{s}}$ by its action on the Kronecker basis:
\[
    P_{\pi} \big( e_{i_{1}} \otimes e_{i_{2}} \otimes \cdots \otimes e_{i_{s}} \big)
    \;=\;
    e_{i_{\pi^{-1}(1)}} \otimes e_{i_{\pi^{-1}(2)}} \otimes \cdots \otimes e_{i_{\pi^{-1}(s)}},
    \quad i_{1},\dots,i_{s} \in [n].
\]
Since $P_{\pi}$ maps the standard basis of $\mathbb{R}^{n^{s}}$ bijectively to itself, it is a permutation matrix.

We now construct the reduction $\mathcal{R}_{2}: \real^{n \otimes s} \rightarrow  \real^{n \otimes s}$: Let $U \in \real^{n \times n}$ be a unitary matrix that is drawn uniformly at random from the set $\{U \in \real^{n\times n}: UU^{\top} = I_{n}\}$. Given $\overline{T} \in \real^{n^{\otimes s}}$, construct $\mathcal{R}_{2}\big(\overline{T}\big)$ such that
\[
    \mathsf{vec}\big(\mathcal{R}_{2}\big(\overline{T}\big)\big) = (U \otimes U \otimes \cdots \otimes U)\; \mathsf{vec}\big( \overline{T} \big).
\]
The tensor $\mathcal{R}_{2}\big(\overline{T}\big)$ by suitably stacking this vector into the appropriate tensor dimension.

\paragraph{The reduction achieves zero-deficiency.}
To verify $\mathcal{R}_{2}\big(\overline{T}\big) \overset{(d)}{=} T$, it is equivalent to verify $\mathsf{vec}\big(\mathcal{R}_{2}\big(\overline{T}\big)\big) \overset{(d)} = \mathsf{vec}(T)$. Note that we have by definition,
\begin{align*}
  \mathsf{vec}\big(\mathcal{R}_{2}\big(\overline{T}\big)\big) &= (U \otimes U \otimes \cdots \otimes U)\; \mathsf{vec}\big( \overline{T} \big) \\
  & = (U \otimes U \otimes \cdots \otimes U)\; \Big( \beta \times (\sqrt{n}\widetilde{v}) \otimes (\sqrt{n}\widetilde{v}) \cdots \otimes (\sqrt{n}\widetilde{v}) + \mathsf{vec}\big(\zeta\big) \Big) \\
  & = \beta \times (\sqrt{n}U \widetilde{v}) \otimes (\sqrt{n}U \widetilde{v}) \otimes \cdots \otimes (\sqrt{n}U \widetilde{v}) + 
  (U \otimes U \otimes \cdots \otimes U)\; \mathsf{vec}\big(\zeta\big).
\end{align*}
We now claim two equalities, whose proofs we defer to the end:
\begin{subequations} \label{eq:two-claims-last-lemma}
\begin{align}\label{claim1-vec-xi}
  &\mathsf{vec}\big(\zeta\big) = \frac{1}{s!} \sum_{\pi \in \mathfrak{S}_{s}} P_{\pi}\; \mathsf{vec} \big(G\big) \quad \text{and}\\
  \label{claim2-vec-xi}
  &(U \otimes U \otimes \cdots \otimes U)\; \mathsf{vec}\big(\zeta\big) = 
  \frac{1}{s!} \sum_{\pi \in \mathfrak{S}_{s}} P_{\pi}\;  (U \otimes U \otimes \cdots \otimes U)\;\mathsf{vec} \big(G\big).
\end{align}
\end{subequations}
Putting the pieces together yields
\begin{align}\label{eq:vec-reduction-T-bar}
    \mathsf{vec}\big(\mathcal{R}_{2}\big(\overline{T}\big)\big) = 
    \beta \times (\sqrt{n}U \widetilde{v}) \otimes ( \sqrt{n}U \widetilde{v}) \otimes \cdots \otimes ( \sqrt{n}U \widetilde{v}) +  \frac{1}{\sqrt{s!}} \sum_{\pi \in \mathfrak{S}_{s}} P_{\pi}\;  (U \otimes U \otimes \cdots \otimes U)\;\mathsf{vec} \big(G\big),
\end{align}
and by definition of $T$~\eqref{tensor-continuous-spike} and Eq.~\eqref{claim1-vec-xi}, we obtain
\begin{align}\label{eq:vec-T}
    \mathsf{vec}\big(T\big) = \beta \times (\sqrt{n}v) \otimes (\sqrt{n}v) \otimes \cdots \otimes (\sqrt{n}v) + \frac{1}{s!} \sum_{\pi \in \mathfrak{S}_{s}} P_{\pi}\; \mathsf{vec} \big(G\big).
\end{align}
We now compare Eq.~\eqref{eq:vec-reduction-T-bar} and Eq.~\eqref{eq:vec-T} to verify that $\mathsf{vec}\big(\mathcal{R}\big(\overline{T}\big)\big) \overset{(d)} = \mathsf{vec}(T)$. First, note that $\widetilde{v} \sim \mathsf{Unif}\big(\big\{\tfrac{-1}{\sqrt{n}}, \tfrac{1}{\sqrt{n}} \big\}^{n}\big)$, $U$ is a unitary matrix that is selected uniformly at random and is independent of $\widetilde{v}$, and $v \sim \mathsf{Unif}\big(\mathbb{S}^{n-1}\big)$. We have $U \widetilde{v}  \overset{(d)}{=} v$. So the first terms in the RHS of Eq.~\eqref{eq:vec-reduction-T-bar} and Eq.~\eqref{eq:vec-T} are equally distributed. Second, note that
$\mathsf{vec} \big(G\big) \sim \NORMAL(0,I_{n^{s}})$ and we have
\begin{align*}
    &\EE \big[  (U \otimes U \otimes \cdots \otimes U)\;\mathsf{vec} \big(G\big) \big] = 0 \quad \text{and} \\
    &\EE \Big[ \big( (U \otimes U \otimes \cdots \otimes U)\;\mathsf{vec} \big(G\big) \big)^{\top} \big( (U \otimes U \otimes \cdots \otimes U)\;\mathsf{vec} \big(G\big) \big)  \Big] \\
    & = \EE \Big[ \big(\mathsf{vec} \big(G\big) \big)^{\top} \big(U^{\top} \otimes U^{\top} \otimes \cdots \otimes U^{\top} \big)  \big(U \otimes U \otimes \cdots \otimes U\big)\;\mathsf{vec} \big(G\big)   \Big] \\
    & = \EE \Big[ \big(\mathsf{vec} \big(G\big) \big)^{\top} \big( (U^{\top} U) \otimes (U^{\top} U) \otimes \cdots \otimes (U^{\top} U) \big)  \;\mathsf{vec} \big(G\big)   \Big] \\
    & = \EE \Big[ \big(\mathsf{vec} \big(G\big) \big)^{\top}  \;\mathsf{vec} \big(G\big)   \Big] = I_{n^{s}}.
\end{align*}
Consequently, we obtain
\[
    (U \otimes U \otimes \cdots \otimes U)\;\mathsf{vec} \big(G\big) \sim \NORMAL(0,I_{n^{s}}).
\]
This proves that the second terms in the RHS of Eq.~\eqref{eq:vec-reduction-T-bar} and Eq.~\eqref{eq:vec-T} are equally distributed. Therefore, we have verified $\mathsf{vec}\big(\mathcal{R}\big(\overline{T}\big)\big) \overset{(d)} = \mathsf{vec}(T)$.
It remains to prove the claimed equalities~\eqref{eq:two-claims-last-lemma}.

\paragraph{Proof of Eq.~\eqref{claim1-vec-xi}}
Note that we have 
\begin{align*}
  \mathsf{vec}\big( \zeta \big) = \frac{1}{s!} \sum_{\pi \in \mathfrak{S}_{s}} \mathsf{vec} \big(G^{\pi}\big).
\end{align*}
By by definition, we have
\begin{align*}
  \mathsf{vec}\big( G^{\pi} \big) &= 
  \sum_{i_{1} =1}^{n} \cdots \sum_{i_{s} =1}^{n} G_{i_{\pi(1)}, \dots, i_{\pi(s)} } \; e_{i_{1}} \otimes e_{i_{1}} \otimes \cdots \otimes e_{i_{s}} \\
  & \overset{\1}{=} 
  \sum_{i_{\pi(1)} =1}^{n} \cdots \sum_{i_{\pi(s)} =1}^{n}  G_{i_{\pi(1)}, \dots, i_{\pi(s)} }\; e_{i_{1}} \otimes e_{i_{2}} \otimes \cdots \otimes e_{i_{s}} \\
  & \overset{\2}{=}  
  \sum_{ j_{1}=1}^{n} \cdots \sum_{ j_{s} =1}^{n}  G_{j_{1}, \dots, j_{s}}\;  e_{j_{\pi^{-1}(1)}} \otimes e_{j_{\pi^{-1}(2)}} \otimes \cdots \otimes e_{j_{\pi^{-1}(s)}} \\
  & \overset{\3}{=}  
  \sum_{ j_{1}=1}^{n} \cdots \sum_{ j_{s} =1}^{n}  G_{j_{1}, \dots, j_{s}}\; P_{\pi} \Big( e_{j_{1}} \otimes e_{j_{2}} \otimes \cdots \otimes e_{j_{s}} \Big) \\
  & = P_{\pi}\; \mathsf{vec} \big( G \big),
\end{align*}
where in step $\1$ we reorder the summation sequence, in step $\2$ we rename the indices by letting
\[
    j_{1} = i_{\pi(1)},\quad  j_{2} = i_{\pi(2)}, \quad \cdots,  \quad j_{s} = i_{\pi(s)}
\]
so that 
\[
   i_{1} = j_{\pi^{-1}(1)}, \quad  i_{2} = j_{\pi^{-1}(2)},  \quad \cdots,  \quad i_{s} = j_{\pi^{-1}(s)},
\]
and in step $\3$ we use the definition of $P_{\pi}$. Summing over permutations completes the proof.

\paragraph{Proof of Eq.~\eqref{claim2-vec-xi}}
We have
\begin{align*}
  (U \otimes U \otimes \cdots \otimes U)\; \mathsf{vec}\big( \zeta \big) = \frac{1}{s!} \sum_{\pi \in \mathfrak{S}_{s}} (U \otimes U \otimes \cdots \otimes U)\;  P_{\pi}\; \mathsf{vec} \big(G\big).
\end{align*}
Note that the operators $(U \otimes U \otimes \cdots \otimes U)$ and $P_{\pi}$ commute, since
\begin{align*}
  (U \otimes U \otimes \cdots \otimes U)\;  P_{\pi}\; (e_{i_{1}} \otimes e_{i_{2}} \otimes \cdots \otimes e_{i_{s}}) &= (U \otimes U \otimes \cdots \otimes U)\; \big(e_{i_{\pi^{-1}(1)}} \otimes e_{i_{\pi^{-1}(2)}} \otimes \cdots \otimes e_{i_{\pi^{-1}(s)}} \big) \\
  &= \big(U e_{i_{\pi^{-1}(1)}}\big) \otimes \big(U e_{i_{\pi^{-1}(2)}}\big) \otimes \cdots \otimes \big(U e_{i_{\pi^{-1}(s)}}\big) \\
  & = P_{\pi} \Big( \big(U e_{i_{1}}\big) \otimes \big(U e_{i_{2}}\big) \otimes \cdots \otimes \big(U e_{i_{s}}\big) \Big) \\
  & = P_{\pi} \; (U \otimes U \otimes \cdots \otimes U) \; (e_{i_{1}} \otimes e_{i_{2}} \otimes \cdots \otimes e_{i_{s}}).
\end{align*}
Consequently, we have 
\begin{align*}
   (U \otimes U \otimes \cdots \otimes U)\; \mathsf{vec}\big( \zeta \big) = \frac{1}{s!} \sum_{\pi \in \mathfrak{S}_{s}} P_{\pi}\; (U \otimes U \otimes \cdots \otimes U)\; \mathsf{vec} \big(G\big),
\end{align*}
as claimed. This concludes the proof. \qed

\end{document}